\documentclass[12pt]{amsart}
\usepackage{latexsym,fancyhdr,amssymb,color,amsmath,amsthm,graphicx,listings,comment}
\usepackage[OT2,T1]{fontenc}
\usepackage[section]{placeins}
\newtheorem*{satz}{Theorem} 
\newtheorem{coro}{Corollary}
\definecolor{red1}{rgb}{1,0.9,0.9} \definecolor{blue1}{rgb}{0.9,0.9,1} \definecolor{green1}{rgb}{0.9,1,0.9} 
\definecolor{yellow1}{rgb}{1,1,0.8} \definecolor{yellow2}{rgb}{1,1,0.8}
\definecolor{color11}{rgb}{1,0.9,0.6}
\def\chapter#1{ \vspace{2mm} \begin{center} \fcolorbox{green1}{green1}{ \parbox{16.2cm}{{\Large {\bf #1}}}} \vspace{2mm} \end{center} }

\def\satz#1{ \vspace{2mm} \begin{center} \fcolorbox{yellow1}{yellow1}{ \parbox{14.0cm}{{\bf Theorem:} #1}} \vspace{2mm} \end{center} }

\newcommand{\NN}{\mathbb{N}} \newcommand{\ZZ}{\mathbb{Z}} \newcommand{\RR}{\mathbb{R}} \newcommand{\CC}{\mathbb{C}}
  \newcommand{\QQ}{\mathbb{Q}} 
\newcommand{\A}{\mathcal{A}}  \newcommand{\R}{\mathcal{R}}

\definecolor{yellow2}{rgb}{1,1,0.9}
\def\tweet#1{ \vspace{2mm} \begin{center} \fcolorbox{yellow2}{yellow2}{ \parbox{11.2cm}{#1}} \vspace{2mm} \end{center} }

\title{Some Fundamental Theorems in Mathematics}
\author{Oliver Knill}
\date{7/22/2018, last update 03/02/2022}
\address{Department of Mathematics \\ Harvard University \\ Cambridge, MA, 02138 }

\makeindex

\begin{document}
\maketitle

\begin{abstract}
An expository hitchhikers guide to some theorems in mathematics. 
\end{abstract} 

Criteria for the current list of 250 theorems are whether the result 
can be formulated elegantly, whether it is beautiful or useful and whether
it could serve as a guide \cite{HitchhikersGuide} without leading to panic.
The order is not a ranking but ordered along a time-line when things were written
down. Since \cite{RotaBeauty} stated
``a mathematical theorem only becomes beautiful if presented as a crown jewel 
within a context" we try sometimes to give some context. 
Of course, any such list of theorems is a matter of personal preferences, taste
and limitations. The number of theorems is arbitrary, the initial obvious goal was 42 
but that number got eventually surpassed as it is hard to stop, once started. 
As a compensation, there are 42 ``tweetable" theorems with included proofs. 
More comments on the choice of the theorems is included in an epilogue.
For literature on general mathematics, see 
\cite{princetonguide,princetonguideapplied,atiyah2000,Garrity,Gowers2002,Stillwell2016,KrantzPark,ConwayCapstone},
for history \cite{Fink1890,Struik1948,Klein1979,BourbakiHistory,Bell,Eves,Kline2,Katz2007,Wussing,Ceruzzi,
Stillwell2010,boyer,guinness,Jackson2012}, for popular, beautiful or elegant things 
\cite{AigZie,Mathbook,Elwes,Dunham,AlsinaNelson,Wells1,Wells2,BehrendsFuenfMinutenMath,EricksonBeautiful,
30SecondMath,MathInMinutes,MacCormick, 17Equations,Henshaw,Elwes,100theorems,CochraneEquations,Crilly50,Cochrane,Neunhauserer,BogdanGrechuk,
MathBook2019}.
For comprehensive overviews in large parts of mathematics, 
\cite{BourbakiElements,DieudonnePanorama,Dieudonne1,BergerPanorama,Simon2017} or predictions on developments 
\cite{MathematicsUnlimited}. For reflections about mathematics in general
\cite{CourantRobbins,RobertoMartinez2018,BellDevelopment,Hersh1997,Livio2009,Byers,RuelleBrain}.
Encyclopedic source examples are
\cite{EDM,Zwillinger,CRCEncyclopedia,CavagnaroHaight,KrantzDictionary,Downing,FrancoiseNaberTsun,
KrantzComprehensive,OxfordDictionary,Tanton}. \\

This is a live document which is in the process of being extended. Thanks so far to 
Johan Commelin, Mikhail Katz, David McCarthy, Kapil Paranjape, 
Jordan Stoyanov, Michael Somos, Ross Rosenwald for some valuable comments or 
corrections. 

\section{Arithmetic}

Let $\NN=\{0,1,2,3, \dots \}$ be the set of {\bf natural numbers}.
A number $p \in \NN, p>1$ is {\bf prime} if $p$ has no factors different from $1$ and $p$. 
With a {\bf prime factorization} $n=p_1 \dots p_n$, we understand the prime factors $p_j$
of $n$ to be ordered as $p_i \leq p_{i+1}$. The {\bf fundamental theorem of arithmetic} is
\index{natural numbers}
\index{primes}
\index{prime factors} 
\index{prime factorization}

\satz{
Every $n \in \NN, n>1$ has a unique prime factorization.
}

Euclid anticipated the result. Carl Friedrich Gauss gave in 1798 the first proof in his monograph
``Disquisitiones Arithmeticae". Within abstract algebra, the result is the statement that the ring of
integers $\ZZ$ is a {\bf unique factorization domain}. For a literature source, see \cite{Kaluzhnin}.
For more general number theory literature, see \cite{Hua1982,Chandrasekharan1968}.
\index{Gauss}
\index{Unique prime factorization}
\index{factorization domain}

\section{Geometry}

Given an {\bf inner product space} $(V,\cdot)$ with {\bf dot product} $v \cdot w$ leading to 
{\bf length} $|v|=\sqrt{v.v}$, three non-zero vectors $v,w,v-w$ define 
a {\bf right angle triangle} if $v$ and $w$ are {\bf perpendicular} meaning that $v \cdot w=0$. 
If $a=|v|,b=|w|,c=|v-w|$ are the lengths of the three vectors, then the {\bf Pythagoras theorem} is
\index{inner product}
\index{Pythagoras theorem}
\index{perpendicular}
\index{right angle triangle}

\satz{
$a^2+b^2=c^2$.
}

Anticipated by Babylonians mathematicians in examples, it appeared independently also in Chinese
mathematics \cite{WasPythagorasChinese} and might have been proven first by Pythagoras \cite{Strathern}
but already early source express uncertainty (see e.g. \cite{KahnPythagoras} p. 32).
The theorem is used in many parts of mathematics like in 
the {\bf Perseval equality} of Fourier theory or that for uncorrelated random variables the 
variance is additive ${\rm Var}[X]+{\rm Var}[Y]={\rm Var}[X+Y]$. In linear algebra it generalizes
to the {\bf Lagrange identity} ${\rm det}(F^T F) = \sum_{|P|=m} \det^2(F_P)$ which holds for all $n \times m$
matrices, where the sum to the right is over all $m \times m$ sub-matrices $P$ of $F$ \cite{HuppertWillems},
a formula which in calculus becomes 
$|\vec{v}|^2 |\vec{w}|^2 - (\vec{v} \cdot \vec{w})^2 = |\vec{v} \wedge \vec{w}|^2$.
See \cite{HiddenHarmonies,PosamentierPythagoras,Maor,Katz2007}.
\index{Perseval equality}
\index{Fourier theory}
\index{Variance}

\section{Calculus}

Let $f$ be a function of one variables which is {\bf continuously differentiable}, meaning that the {\bf limit}
$g(x)= \lim_{h \to 0} [f(x+h)-f(x)]/h$ exists at every point $x$ and defines a continuous function $g$. 
For any such function $f$, we can form the {\bf integral} $\int_a^b f(t) \; dt$ and the {\bf derivative}
$d/dx f(x) = f'(x)$. 
\index{differentiable}
\index{integral}
\index{limit}
\index{derivative}
\index{continuously differentiable}

\satz{
$\int_a^b f'(x) dx = f(b)-f(a), \; \; \; \;   \frac{d}{dx} \int_0^x f(t) dt = f(x)$
}

Newton and Leibniz discovered the result independently, Gregory wrote down the first proof in his
``Geometriae Pars Universalis" of 1668. The result generalizes to higher dimensions in the form
of the {\bf Green-Stokes-Gauss-Ostogradski theorem} $\int_M dF = \int_{\delta M} F$ which holds for
$n$-forms $F$ with {\bf exterior derivative} $dF$ and compact $(n+1)$-manifolds $M$ with boundary 
$\delta M$. \cite{EisenmanJoke} tells the ``tongue in the cheek" proof: as the derivative is 
a limit of {\bf quotient} of {\bf differences}, the anti-derivative must be a limit of 
{\bf sums} of {\bf products}. For history, see \cite{Katz79,EdwardsHistory}.
\index{Newton}
\index{Leibniz}
\index{Stokes}
\index{Gregory}

\section{Algebra}

A {\bf polynomial} is a complex-valued function of the form $f(x) = a_0 + a_1 x + \cdots + a_n x^n$, 
where the entries $a_k$ are in the complex plane $\CC$. The space of all polynomials is denoted by $\CC[x]$. 
The largest non-negative integer $n$ for which $a_n \neq 0$ is called the {\bf degree} of the polynomial. 
Degree $1$ polynomials are {\bf linear}, degree $2$ polynomials are called {\bf quadratic} etc. 
The {\bf fundamental theorem of algebra} is 
\index{polynomial}
\index{degree}
\index{linear}
\index{fundamental theorem of algebra}

\satz{
Every $f \in \CC[x]$ of degree $n$ can be factored into $n$ linear factors.
}

This result was anticipated during the 17th century. The first 
author to assert that any n'th degree polynomial has a root is Peter Roth
in 1600 \cite{Carrera1992}. This was proven first by Carl Friedrich Gauss and finalized 
in 1920 by Alexander Ostrowski who fixed a topological mistake in Gauss proof. 
The theorem assures that the field of complex numbers $\CC$ 
is {\bf algebraically closed}. For history and many proofs see \cite{FineRosenberger}.
\index{Ostrowski}
\index{Gauss}
\index{algebraic closure}
\index{polynomial}

\section{Probability}

Given a sequence $X_k$ of {\bf independent random variables} on a probability space $(\Omega,\mathcal{A},{\rm P})$ 
which all have the same {\bf cumulative distribution functions} 
$F_X(t) = {\rm P}[X \leq t]$. The {\bf normalized random variable} $\overline{X}=$ is
$(X-{\rm E}[X])/\sigma[X]$, where ${\rm E}[X]$ is the {\bf mean} $\int_{\Omega} X(\omega) dP(\omega)$ 
and $\sigma[X]  = {\rm E}[(X-{\rm E}[X])^2]^{1/2}$ is the standard deviation.
A sequence of random variables $Z_n \to Z$ {\bf converges in distribution} to $Z$ if 
$F_{Z_n}(t) \to F_Z(t)$ for all $t$ as $n \to \infty$.
If $Z$ is a {\bf Gaussian random variable}
with zero mean ${\rm E}[Z]=0$ and standard deviation $\sigma[Z]=1$, 
the {\bf central limit theorem} is:
\index{central limit theorem}
\index{mean}
\index{standard deviation}
\index{cumulative distribution function}
\index{independent random variable}
\index{convergence in distribution}
\index{Gaussian random variables}

\satz{
$\overline{(X_1+X_2+ \cdots + X_n)} \to Z$ in distribution.
}

Proven in a special case by Abraham De-Moivre for discrete random variables 
and then by Constantin Carath\'eodory and 
Paul L\'evy, the theorem explains the importance and ubiquity of the 
{\bf Gaussian density function} $e^{-x^2/2}/\sqrt{2\pi}$ 
defining the {\bf normal distribution}. The Gaussian distribution was first
considered by Abraham de Moivre in 1738.  See \cite{Stroock,knillprobability}.
\index{De Moivre}
\index{Caratheodory}
\index{Levy}

\section{Dynamics}

Assume $X$ is a {\bf random variable} on a {\bf probability space} 
$(\Omega,\mathcal{A},{\rm P})$ for which $|X|$ has 
finite mean ${\rm E}[|X|]$. This means $X: \Omega \to \mathbb{R}$ is measurable and 
$\int_\Omega |X(x)| d{\rm P}(x)$ is finite. Let $T$ be an ergodic, measure-preserving 
transformation from $\Omega$ to $\Omega$.
{\bf Measure preserving} means that $P[T^{-1}(A)]=P[A]$ for all {\bf measurable sets} $A \in \mathcal{A}$.  
{\bf Ergodic} means that that $T(A)=A$ implies ${\rm P}[A]=0$ or ${\rm P}[A]=1$ for all $A \in \mathcal{A}$. The
{\bf ergodic theorem} states, that for an ergodic transformation $T$ on has:
\index{Ergodic}
\index{measure preserving}
\index{probability space}
\index{random variable}

\satz{
$[X(x)+X(Tx) + \cdots + X(T^{n-1}(x))]/n \to {\rm E}[X]$ for almost all $x$.
}

This theorem from 1931 is due to George Birkhoff and is called {\bf Birkhoff's
pointwise ergodic theorem}. It assures that ``time averages" are equal to ``space averages". 
A draft of the {\bf von Neumann mean ergodic theorem} which appeared in 1932 by John von Neumann has motivated
Birkhoff, but the mean ergodic version is weaker. See \cite{Zund} for history. 
A special case is the {\bf law of large numbers}, in which case the random variables $x \to X(T^k(x))$ 
are independent with equal distribution (IID). The theorem belongs to ergodic theory 
\cite{Halmos,CFS,Sinai}.
\index{Birkhoff}
\index{Birkhoff theorem}
\index{law of large numbers}
\index{time average}
\index{space average}

\section{Set theory}

A {\bf bijection} is a map from a set $X$ to a set $Y$ which is {\bf injective}: $f(x)=f(y) \Rightarrow x=y$ and 
{\bf surjective}: for every $y \in Y$, there exists $x \in X$ with $f(x)=y$. 
Two sets $X,Y$ have the {\bf same cardinality}, if there exists a bijection from $X$ to $Y$. 
Given a set $X$, the {\bf power set} $2^X$ is the set of all subsets of $X$, including the {\bf empty
set} and $X$ itself. If $X$ has $n$ elements, the power set has $2^n$ elements. Cantor's theorem is 
\index{cardinality}
\index{bijective}
\index{surjective}
\index{injective}
\index{power set}
\index{empty set}

\satz{For any set $X$, the sets $X$ and $2^X$ have different cardinality. }

The result is due to Cantor. Taking for $X$ the natural numbers, then every $Y \in 2^X$ 
defines a real number $\phi(Y)=\sum_{y \in Y} 2^{-y} \in [0,1]$. As $Y$ and $[0,1]$ have
the same cardinality (as {\bf double counting pair cases} 
like $0.39999999 \dots = 0.400000 \dots$ form a 
countable set), the interval $[0,1]$ is uncountable. 
There are different types of infinities leading to {\bf countable infinite sets}
and {\bf uncountable infinite sets}. 
In order to compare sets, the {\bf Schr\"oder-Bernstein theorem} is important. If there exist
injective functions $f:X \to Y$ and $g: Y \to X$, then there exists also a bijection $X \to Y$. 
This result was used by Cantor already. For literature, see \cite{HalmosNaive}.
\index{Cantor}
\index{countable}
\index{uncountable}
\index{Schroeder-Bernstein}

\section{Statistics}

A {\bf probability space} $(\Omega,\mathcal{A},{\rm P})$ consists of a set $\Omega$,
a {\bf $\sigma$-algebra} $\mathcal{A}$ and a {\bf probability measure} ${\rm P}$. A 
$\sigma$-algebra is a collection of subset of $\Omega$ which contains the empty set
and which is closed under the operations of taking {\bf complements}, {\bf countable unions} and
{\bf countable intersections}. The function ${\rm P}$ on $\mathcal{A}$ takes values in the
interval $[0,1]$, satisfies ${\rm P}[\Omega]=1$ and 
${\rm P}[ \bigcup_{A \in S} A ] = \sum_{A \in S} {\rm P}[A]$ for any finite or countable set 
$S \subset \mathcal{A}$ of pairwise disjoint sets. 
The elements in $\mathcal{A}$ are called {\bf events}. Given two events $A,B$ where
$B$ satisfies ${\rm P}[B]>0$, one can define the {\bf conditional probability} 
${\rm P}[A|B]={\rm P}[A \cap B]/{\rm P}[B]$. {\bf Bayes theorem} states:
\index{probability measure}
\index{$\sigma$-algebra}
\index{probability space}
\index{conditional probability}
\index{Bayes theorem}

\satz{${\rm P}[A|B] = {\rm P}[B|A] {\rm P}[A]/{\rm P}[B]$}

The setup stated the {\bf Kolmogorov axioms} by Andrey Kolmogorov who wrote in 1933 the ``Grundbegriffe der
Wahrscheinlichkeitsrechnung" \cite{Kolmogorov} based on measure theory built by Emile Borel and Henry Lebesgue. 
For history, see \cite{ShaferVovk}, who report that 
``Kolmogorov sat down to write the Grundbegriffe, in a rented cottage on the Klyaz'ma River in November 1932".
Bayes theorem is rather a fantastically clever definition and not really a theorem. 
There is almost nothing to prove as multiplying with ${\rm P}[B]$ gives ${\rm P}[A \cap B]$ on both sides. 
It essentially restates that $A \cap B = B \cap A$, the Abelian property of the product in the ring
$\mathcal{A}$. More general is the statement that if $A_1, \dots, A_n$ is a disjoint 
set of events whose union is $\Omega$, then 
${\rm P}[A_i|B] = {\rm P}[B|A_i] {\rm P}[A_i]/(\sum_j {\rm P}[B|A_j] {\rm P}[A_j]$.
Bayes theorem was first noticed in 1763 by Thomas Bayes. It is by some considered 
to the theory of probability what the Pythagoras theorem is to geometry.
``Monty Hall" type stories \cite{Rosenhouse} illustrate that conditional expectation is not always intuitive. 

\section{Graph theory}

A {\bf finite simple graph} $G=(V,E)$ is a finite collection $V$ of {\bf vertices} connected by a 
finite collection $E$ of {\bf edges}, which are un-ordered pairs $(a,b)$ with $a,b \in V$. 
{\bf Simple} means that no {\bf self-loops} nor {\bf multiple connections} are present in the graph.
The {\bf vertex degree} $d(x)$ of $x \in V$ is the number of edges containing $x$. 
\index{graph}
\index{finite simple graph}
\index{vertices}
\index{edges}
\index{vertex degree}
\index{self-loops}
\index{multiple connections}

\satz{
$\sum_{x \in V} d(x)/2 = |E|$.
}

This formula is also called the {\bf Euler handshake formula} because every edge in a graph 
contributes exactly two handshakes. It can be seen as a {\bf Gauss-Bonnet formula} for the {\bf valuation}
$G \to v_1(G)$ counting the number of edges in $G$. A {\bf valuation} $\phi$ is a function defined on 
{\bf sub-graphs} with the property that $\phi(A \cup B) = \phi(A) + \phi(B) - \phi(A \cap B)$. 
Examples of valuations are the number $v_k(G)$ of {\bf complete sub-graphs} of dimension $k$ of $G$. 
An other example is the {\bf Euler characteristic} $\chi(G)=v_0(G)-v_1(G)+v_2(G)-v_3(G)+ \cdots + (-1)^d v_d(G)$.
If we write $d_k(x)=v_k(S(x))$, where $S(x)$ is the unit sphere of $x$, then 
$\sum_{x \in V} d_k(x)/(k+1) = v_k(G)$ is the {\bf generalized handshake formula}, the Gauss-Bonnet
result for $v_k$. The Euler characteristic then satisfies $\sum_{x \in V} K(x) = \chi(G)$, 
where $K(x) = \sum_{k=0}^{\infty} (-1)^k v_k(S(x))/(k+1)$. 
This is the {\bf discrete Gauss-Bonnet result}. The handshake result is the special case for the 
valuation $v_1(G)$ counting the number of edges of a graph $G$. 
It was found by Euler and has by some called the {\bf fundamental theorem of graph theory}.
For more about graph theory, \cite{BM,Merris,BR,handbookgraph} about Euler: \cite{Fellmann}. 
\index{Euler handshake}
\index{valuation}
\index{Gauss-Bonnet}
\index{Euler characteristic}
\index{subgraph}
\index{Generalized handshake}

\section{Polyhedra}

A {\bf finite simple graph} $G=(V,E)$ is given by a finite vertex set $V$ and edge set $E$.
A subset $W$ of $V$ {\bf generates} the sub-graph $(W,\{ \{a,b\} \in E \; | \; a,b \in W \})$. 
The {\bf unit sphere} of $v \in V$ is the sub graph generated by
$S(x) = \{ y \in V \; | \; \{x,v\} \in E \}$. 
The {\bf empty graph} $0=(\emptyset,\emptyset)$ is called the {\bf $(-1)$-sphere}.
The $1$-point graph $1=(\{1\},\emptyset)=K_1$ is the smallest contractible graph.
Inductively, a graph $G$ is {\bf contractible}, if it is either $1$ or if
there exists $x \in V$ such that both $G-x$ and $S(x)$ are contractible. 
Inductively, a graph $G$ is a {\bf $d$-sphere}, if it is either $0$ or
if every $S(x)$ is a $(d-1)$-sphere and if there exists a vertex $x$ such that $G-x$ is contractible.
Let $v_k$ denote the number of complete sub-graphs $K_{k+1}$
of $G$. The vector $(v_0,v_1, \dots)$ is the {\bf $f$-vector} of $G$ and
$\chi(G)=v_0-v_1+v_2- \dots $ is the {\bf Euler characteristic} of $G$. The generalized {\bf Euler gem}
formula due to Schl\"afli is:
\index{Euler characteristic}
\index{sphere}
\index{Euler gem formula}
\index{Euler polyhedron formula}
\index{contractible}
\index{unit sphere}

\satz{For $d=2$, $\chi(G)=v-e+f=2$. For $d$-spheres, $\chi(G) = 1+(-1)^d$.}

Convex Polytopes were studied already in ancient Greece. The Euler characteristic
relations were discovered in dimension $2$ by Descartes \cite{Aczel} and interpreted topologically
by Euler who proved the case $d=2$. This is written as $v-e+f=2$, where $v=v_0,e=v_1,f=v_2$. The
two-dimensional case can be stated for {\bf planar graphs}, where one has a clear notion of what
the two dimensional cells are and can use the topology of the ambient sphere in which the graph 
is embedded. Historically there had been confusions \cite{cromwell,Richeson} about the definitions.
It was Ludwig Schl\"afli \cite{Schlafli} who covered the higher dimensional case. The 
above set-up is a modern reformulation of his set-up, due essentially to Alexander Evako. 
Multiple refutations \cite{lakatos} can be blamed to ambiguous definitions. 
Polytopes are often defined through convexity \cite{gruenbaum,Ziegler} and there is not much consensus 
on a general definition \cite{Gruenbaum2003}, which was the reason in this entry to formulate
Schl\"afli's theorem in a rather restrictive case (where all cells are simplices),
but where we have a simple combinatorial definition of what a ``sphere" is.  See also 
\cite{symmetries}. 

\section{Topology}

The {\bf Zorn lemma} assures that the Cartesian product of a non-empty family of non-empty 
sets is non-empty. The {\bf Zorn lemma} is equivalent to the {\bf axiom of choice} C in the 
Zermelo-Frenkel {\bf ZFC axiom system} and also equivalent to the {\bf Tychonov theorem} in topology.
Let $X=\prod_{i \in I} X_i$ denote the {\bf product} of topological spaces. The {\bf product topology}
is the {\bf weakest topology} on $X$ which renders all {\bf projection functions} 
$\pi_i: X \to X_i$ continuous. Here is Tychonov;s theorem
\index{Zorn lemma}
\index{axiom of choice}
\index{ZFC axiom system}
\index{Tychonov theorem}
\index{product topology}
\index{projection function}

\satz{
If all $X_i$ are compact, then  $\prod_{i \in I} X_i$ is compact. 
}

{\bf Zorn's lemma} is due to Kazimierz Kuratowski in 1922 and Max August Zorn in 1935. 
Andrey Nikolayevich Tykhonov proved his theorem
in 1930. One application of the Zorn lemma is the {\bf Hahn-Banach theorem} in functional analysis,
the existence of {\bf spanning trees} in infinite graphs or the fact that commutative rings 
with units have {\bf maximal ideals}. For literature, see \cite{Jech2008}. 
\index{Kuratowski}
\index{Hahn-Banach}
\index{Tychonov theorem}
\index{Axiom of choice}
\index{Zermele-Frenkel}
\index{product topology}

\section{Algebraic geometry}

The {\bf algebraic set} $V(J)$ of an {\bf ideal} $J$ in the commutative ring $R=k[x_1,\dots, x_n]$ over an 
{\bf algebraically closed} {\bf field} $k$ defines the ideal $I(V(J))$ containing all polynomials 
that vanish on $V(J)$. 
The {\bf radical} $\sqrt{J}$ of an ideal $J$ is the set of polynomials in $R$ such that $r^n \in J$ for
some positive $n$. [An {\bf ideal} $J$ in a ring $R$ is a subgroup of the additive group of $R$ 
such that $r x \in I$ for all $r \in R$ and all $x \in I$. It defines the {\bf quotient ring} $R/I$ 
and is so the kernel of a ring homomorphism from $R$ to $R/I$. 
The algebraic set $V(J)=\{ x \in k^n \; | \; f(x)=0, \forall f \in J\}$ 
of an ideal $J$ in the polynomial ring $R$ is the set of common roots of all these functions $f$.
The algebraic sets are the closed sets in the {\bf Zariski topology} of $R$. 
The ring $R/I(V)$ is the {\bf coordinate ring}
of the algebraic set $V$.]  The {\bf Hilbert Nullstellensatz} is
\index{Algebraic set}
\index{radical}
\index{ideal}
\index{commutative ring}
\index{Nullstellensatz}

\satz{$I(V(J)) = \sqrt{J}$.}

The theorem is due to Hilbert from 1893 \cite{Hilbert1893} (page 320).  Of course, Hilbert did not yet
use the language of ideals but in terms of having ``ganze rationalen homogene Funktionen" of several variables.
A simple example is when $J=\langle p \rangle = \langle x^2-2xy+y^2 \rangle$ is the ideal 
$J$ generated by $p$ in $\RR[x,y]$; then $V(J) = \{ x=y \}$ and $I(V(J))$ is the ideal generated by $x-y$. 
For literature, see \cite{hartshorne,ZariskiCommutativeAlgebra}. 
\index{Hilbert}
\index{Nullstellenstatz}

\section{Cryptology}

An integer $p>1$ is {\bf prime} if $1$ and $p$ are the only factors of $p$. 
The number $k \; {\rm mod} \; p$ is the {\bf reminder} when dividing $k$ by $p$. 
For example $18 {\rm mod} 7 = 4$. 
{\bf Fermat's little theorem} is 
\index{reminder}
\index{prime}
\index{Fermat's little theorem}

\satz{$a^p=a \; {\rm mod} \; p$ for every prime $p$ and every integer $a$.}

The theorem was found by Pierre de Fermat in 1640. A first proof appeared in 1683 by Leibniz. 
Euler in 1736 published the first proof. The result is used in the {\bf Diffie-Hellman key exchange}, 
where a large public prime $p$ and 
a public base value $a$ are taken. Ana chooses a number $x$ and publishes $X=a^x {\rm mod} p$ and 
Bob picks $y$ publishing $Y=a^y {\rm mod} p$. Their secret key is $K=X^y = Y^x$. 
An adversary Eve who only knows $a,p, X$ and $Y$ can from this not get $K$ due to the difficulty
of the {\bf discrete log  problem}.
More generally, for possibly composite numbers $n$, the theorem extends to the fact that $a^{\phi(n)}=1$ modulo $p$, 
where the {\bf Euler's totient function} $\phi(n)$ counts the number of positive integers less than $n$ which are 
{\bf coprime} to $n$. The generalized Fermat theorem is the key for RSA {\bf crypto systems}: 
in order for Ana and Bob to communicate.
Bob publishes the product $n=pq$ of two large primes as well as some base integer $a$. Neither Ana nor any 
third party Eve do know the factorization. Ana communicates a message $x$ to Bob by 
sending $X=a^x {\rm mod} n$ using {\bf modular exponentiation}. 
Bob, who knows $p,q$, can find $y$ such that $xy=1 \; {\rm mod} \; \phi(n)$. This
is because of Fermat $a^{(p-1)(q-1)} = a \; {\rm mod} \; n$. 
Now, he can compute $x=y^{-1} {\rm mod} \; \phi(n)$. Not even Ana herself could recover $x$ from $X$.

\index{coprime}
\index{crypto system}
\index{Diffie-Hellman system}
\index{RSA crypto system}
\index{discrete log problem}
\index{Euler totient function}

\section{Spectral theorem}

A bounded linear operator $A$ on a {\bf Hilbert space} is called {\bf normal} if $AA^* = A^*A$, 
where $A^*=\overline{A}^T$ is the {\bf adjoint} and $A^T$ is the {\bf transpose} and $\overline{A}$ 
is the {\bf complex conjugate}. Examples of normal operators are {\bf self-adjoint} operators 
(meaning $A=A^*$) or {\bf unitary operators} (meaning $A A^*=1$). 
\index{Hilbert space}
\index{self-adjoint operator}
\index{bounded linear operator}
\index{normal operator}
\index{unitary operator}
\index{adjoint}
\index{complex conjugate}

\satz{$A$ is normal if and only if $A$ is unitarily diagonalizable. }

In finite dimensions, any unitary $U$ diagonalizing $A$ using $B=U^* A U$ contains an {\bf orthonormal eigenbasis} 
of $A$ as column vectors. The theorem is due to Hilbert. In the self-adjoint case, all the eigenvalues are real
and in the unitary case, all eigenvalues are on the unit circle. 
The result allows a {\bf functional calculus} for normal operators: for any continuous
function $f$ and any bounded linear operator $A$, one can define $f(A) = U f(B) U^*$, 
if $B=U^* A U$. See \cite{ConwayFunctionalAnalysis}.
\index{Diagonalizable}
\index{functional calculus}
\index{eigenbasis}
\index{orthonormal eigenbasis}
\index{eigenvalues}

\section{Number systems}

A {\bf monoid} is a set $X$ equipped with an {\bf associative operation} $*$ and an 
{\bf identity element} $1$ satisfying $1*x=x$ for all $x \in X$. {\bf Associativity} 
means $x*(y*z)=(x*y)*z$ for all $x,y,z \in X$. The monoid structure belongs to a collection of
mathematical structures {\bf magmas} $\supset$ {\bf semigroups} $\supset$ {\bf monoids} $\supset$ {\bf groups}. 
A monoid is {\bf commutative}, if $x*y=y*x$ for all $x,y \in X$. 
A {\bf group} is a monoid in which every element $x$ has an {\bf inverse} $y$ satisfying $x*y=y*x=1$. 
\index{monoid}
\index{group}
\index{inverse}
\index{associativity}
\index{commutative}
\index{identity element}

\satz{
Every commutative monoid can be extended to a group.
}

The general result is due to Alexander Grothendieck from around 1957. A more precise statement is that
there is a group containing a homomorphic image of the monoid. It is for {\bf cancellative monoids}
(the statement $a*x=b*x$ implies $a=b$ in the monoid) that the monoid is also contained isomorphically inside 
the group. In general, like for a {\bf zero monoid} with 3 or more elements defined by $x*y=1$ for all $x,y$
which is not cancellative, such a collapse already appears. The group is called the {\bf Grothendieck group completion}
of the monoid. For example, the additive monoid of natural numbers can be extended to the group of integers, the
multiplicative monoid of non-zero integers can be extended to the group of rational numbers. 
The construction of the group is used in {\bf K-theory} \cite{Atiyah1967,Karoubi}
For insight about the philosophy of Grothendieck's mathematics, see \cite{risingsea}. 
\index{Grothendick group completion}
\index{Grothendieck}
\index{group completion}
\index{K-Theory}
\index{monoid}
\index{cancellative monoid}

\section{Combinatorics}

Let $|X|$ denote the {\bf cardinality} of a finite set $X$. This means that
$|X|$ is the number of elements in $X$. A function $f$ from a set $X$ to a set $Y$
is called {\bf injective} if $f(x)=f(y)$ implies $x=y$. The {\bf pigeonhole principle} tells:
\index{cardinality}
\index{pigeonhole principle}

\satz{
If $|X|>|Y|$ then no function $X \to Y$ can be injective. 
}

This implies that if we place $n$ items into $m$ boxes and $n>m$, then one box
must contain more than one item. The principle is believed to be formalized first by Peter Dirichlet.
Despite its simplicity, the principle has many applications, like proving that something exists. 
An example is the statement that there are two trees in New York City streets which have the same 
number of leaves. The reason is that the U.S. Forest services states 592'130 trees in the year 2006
and that a mature, healthy tree has about 200'000 leaves. 
One can also use it for less trivial statements like that in a cocktail party there are at least two
with the same number of friends present at the party. A mathematical application is the
{\bf Chinese remainder Theorem} stating that that there exists a solution to $a_i x = b_i \; {\rm mod} \; m_i$
all disjoint pairs $m_i, m_j$ and all pairs $a_i,m_i$ are relatively prime \cite{DingPeiSalomaa,Martzloff}.
The principle generalizes to infinite set if $|X|$ is the cardinality. It implies then for example
that there is no injective function from the real numbers to the integers.  For literature,
see for example \cite{Brualdi2004}, which states also a stronger version which for example allows
to show that any sequence of real $n^2+1$ real numbers contains either an increasing subsequence of length $n+1$
or a decreasing subsequence of length $n+1$. 
\index{Dirichlet}
\index{Chinese Remainder Theorem}

\section{Complex analysis}

Assume $f$ is an {\bf analytic function} in an {\bf open domain} $G$ of the {\bf complex plane} $\CC$. 
Such a function is also called {\bf holomorphic} in $G$. 
Holomorphic means that if $f(x+iy) = u(x+iy) + i v(x+iy)$, then the {\bf Cauchy-Riemann} differential 
equations $u_x = v_y, u_y = -v_x$ hold in $G$. Assume $z$ is in $G$ and assume $C \subset G$ is a {\bf circle} $a+r e^{i \theta}$ 
centered at $z$ which is bounding a disc $D=\{ z \in \CC \; | \; |z-a| < r \} \subset G$. 
\index{analytic function}
\index{open domain}
\index{Cauchy-Riemann differential equation}
\index{circle}
\index{complex plane}
\index{holomorphic function}

\satz{
For analytic $f$ and circle $C \subset G$, one has $f(a) = \frac{1}{2\pi i} \int_C \frac{f(z) dz}{(z-a)}$.
}

This {\bf Cauchy integral formula} of Cauchy is used for other results and estimates.
It implies for example the {\bf Cauchy integral theorem} assuring that $\int_C f(z) dz=0$ for any simple closed 
curve $C$ in $G$ bounding a simply connected region $D \subset G$. {\bf Morera's theorem} 
assures that for any domain $G$, if $\int_C f(z) \; dz=0$ for all simple closed smooth curves 
$C$ in $G$, then $f$ is holomorphic in $G$. An other generalization is {\bf residue calculus}: 
For a simply connected region $G$ and a function $f$ which is analytic except in a finite set $A$ 
of points. If $C$ is piecewise smooth continuous closed curve not intersecting $A$, then 
$\int_C f(z) \; dz = 2\pi i \sum_{a \in A} I(C,a) {\rm Res}(f,a)$, where $I(C,a)$ is the 
{\bf winding number} of $C$ with respect to $a$ and ${\rm Res}(f,a)$ is the {\bf residue} of $f$ at $a$ which 
is in the case of poles given by $\lim_{z \to a} (z-a) f(z)$. 
See \cite{CauchyCoursDAnalyse,AhlforsComplexAnalysis,Conway1978}.
\index{Cauchy integral formula}
\index{Cauchy integral theorem}
\index{Morera's theorem}
\index{Residue calculus}
\index{residue}
\index{winding number}

\section{Linear algebra}

If $A$ is a $m \times n$ {\bf matrix} with {\bf image} ${\rm ran}(A)$ and {\bf kernel} ${\rm ker}(A)$. 
If $V$ is a linear subspace of $\RR^m$, then $V^\perp$ denotes the {\bf orthogonal complement} 
of $V$ in $\RR^m$, the
linear space of vectors perpendicular to all $x \in V$. 
\index{matrix}
\index{image}
\index{kernel}
\index{orthogonal complement}
\index{perpendicular}

\satz{
${\rm dim}({\rm ker} A)+{\rm dim}({\rm ran} A) = n,
 {\rm dim}(({\rm ran} A)^\perp) = {\rm dim}({\rm ker} A^T)$.
}

The result is used in {\bf data fitting} for example when understanding the {\bf least square solution}
$x=(A^T A)^{-1} A^T b$ of a {\bf system of linear equations} $A x = b$. 
It assures that $A^T A$ is invertible
if $A$ has a trivial kernel. The result is a bit stronger than the {\bf rank-nullity theorem}
${\rm dim}({\rm ran}(A)) + {\rm dim}({\rm ker}(A)) = n$ alone 
and implies that for finite $m \times n$ matrices the 
{\bf index} ${\rm dim}({\rm ker} A) - {\rm dim}({\rm ker} A^*)$ is always $n-m$, 
which is the value for the $0$ matrix. 
For literature, see \cite{Strang1993}. The result has an abstract generalization in the form 
of the group isomorphism theorem for a group homomorphism $f$ stating that $G/{\rm ker}(f)$
is isomorphic to $f(G)$. It can also be described using the 
{\bf singular value decomposition} $A=U D V^T$.
The number $r={\rm ran} A$ has as a basis the first $r$ columns of $U$.
The number $n-r={\rm ker} A$ has as a basis the last $n-r$ columns of $V$. 
The number ${\rm ran} A^T$ has as a basis the first $r$ columns of $V$.
The number ${\rm ker} A^T$ has as a basis the last $m-r$ columns of $U$.
\index{rank-nullity}
\index{index}
\index{least square solution}
\index{data fitting}
\index{singular value decomposition}

\section{Differential equations}

A {\bf differential equation} $\frac{d}{dt} x=f(x)$ and $x(0)=x_0$ in a 
{\bf Banach space} $(X,|| \cdot ||)$ (a normed, complete vector space) 
defines an {\bf initial value problem}: we look for a solution $x(t)$ 
satisfying the equation and given initial condition $x(0)=x_0$ and $t \in (-a,a)$ for some $a>0$. 
A function $f$ from $\RR$ to $X$ is called {\bf Lipschitz}, if there exists a constant $C$ 
such that for all $x,y \in X$ the inequality $||f(x)-f(y)|| \leq C |x-y|$ holds. 
\index{differential equation}
\index{Banach space}
\index{initial value theorem}
\index{Lipschitz}

\satz{
If $f$ is Lipschitz, a unique solution of $x'=f(x),x(0)=x_0$ exists. 
}

This result is due to Picard and Lindel\"of from 1894. Replacing the Lipschitz condition 
with continuity still gives an {\bf existence theorem} which is due to Giuseppe Peano in 1886, 
but uniqueness can fail like for $x'=\sqrt{x},x(0)=0$ with solutions $x=0$ and $x(t)=t^2/4$. 
The example $x'(t)=x^2(t),x(0)=1$ with solution $1/(1-t)$ shows that we can not have solutions
for all $t$. The proof is a simple application of the Banach fixed point theorem. 
For literature, see \cite{CoddingtonLevinson}.
\index{Picard}
\index{Lidelof}
\index{Peano}

\section{Logic}

An {\bf axiom system} $A$ is a collection of formal statements assumed to be true.
We assume it to contain the basic {\bf Peano axioms} of arithmetic. (One only needs
first order Peano arithmetic PA, for the first incompletness theorem one can even do with 
the weaker Robinson arithmetic.)
An axiom system is {\bf complete}, if every true statement can be proven within the system.  
The system is {\bf consistent} if one can not prove $1=0$ within the system. 
It is {\bf provably consistent} if one can prove a theorem 
{\bf "The axiom system $A$ is consistent."} within the system. 
It is important that the axiom system is strong enough to contain the Peano arithmetic
as there are interesting and widely studied theories that happen to be complete, 
such as the theory of {\bf real closed fields}.
\index{axiom system}
\index{Peano axioms}
\index{formal rules}
\index{consistent}
\index{complete}
\index{inconsistent}
\index{provably consistent}

\satz{
An axiom system is neither complete nor provably consistent.
}

The result is due to Kurt G\"odel who proved it in 1931. In this thesis, G\"odel had proven a
completeness theorem of first order predicate logic. 
The incompleteness theorems of 1931 destroyed the dream of {\bf Hilbert's program} which 
aimed for a complete and consistent {\bf axiom system} for mathematics. A commonly assumed
axiom system is the {\bf Zermelo-Frenkel axiom system} together with the axiom of choice ZFC.
Other examples are Quine's {\bf new foundations} NF or Lawvere's 
{\bf elementary theory of the category of sets} ETCS.
For a modern view on Hilbert's program, see \cite{Tapp}. For G\"odel's theorem \cite{Franzen,NagelNewman}. 
Hardly any other theorem had so much impact outside of mathematics. 
\index{G\"odel}
\index{axiom system}
\index{Hilbert's program}
\index{new foundations}
\index{elementary theory of the category of sets}

\section{Representation theory}

For a {\bf finite group} or {\bf compact topological group} $G$, one can look at 
{\bf representations}, group homomorphisms from $G$ to the automorphisms of a {\bf vector space} $V$. 
A representation of $G$ is {\bf irreducible} if the only $G$-invariant
subspaces of $V$ are $0$ or $V$. The {\bf direct sum} of of two representations $\phi,\psi$ is defined
as $\phi \oplus \psi(g)(v \oplus w) = \phi(g)(v) \oplus \phi(g)(w)$. A representation is {\bf semi simple} 
if it is a unique direct sum of irreducible finite-dimensional representations:
\index{finite group}
\index{compact group}
\index{representation}
\index{vector space}
\index{direct sum of representations}
\index{irreducible representation}
\index{semi simple representation}

\satz{
Representations of compact topological groups are semi simple.
}

For representation theory, see \cite{FultonHarris}. Pioneers in representation theory were
Ferdinand Georg Frobenius, Herman Weyl, and \'Elie Cartan.
Examples of compact groups are {\bf finite group}, or {\bf compact Lie groups} (a smooth manifold
which is also a group for which the multiplications and inverse operations are smooth) like
the {\bf torus group} $T^n$, the orthogonal groups $O(n)$ of all orthogonal 
$n \times n$ matrices or the {\bf unitary groups} $U(n)$ of all unitary $n \times n$ matrices
or the group ${\rm Sp}(n)$ of all {\bf symplectic} $n \times n$ matrices. Examples of groups that
are not Lie groups are the groups $Z_p$ of {\bf $p$-adic integers}, which are examples of 
{\bf pro-finite groups}.
\index{finite group} 
\index{symplectic}
\index{p-adic integers}
\index{unitary group}
\index{orthogonal group}
\index{profinite group}

\section{Lie theory}

Given a {\bf topological group} $G$, a {\bf Borel measure} $\mu$ on $G$ is called 
{\bf left invariant} if $\mu(g A) = \mu(A)$ 
for every $g \in G$ and every measurable set $A \subset G$. A left-invariant measure on $G$ is also called
a {\bf Haar measure}. A topological space is called {\bf locally compact}, if every point has 
a compact neighborhood. 
\index{topological group}
\index{Haar measure}
\index{Borel measure}
\index{left invariant}
\index{locally compact}
\index{character}

\satz{
A locally compact group has a unique Haar measure.
}

Alfr\'ed Haar showed the existence in 1933 and John von Neumann proved that it is unique.
In the compact case, the measure is finite, leading to an inner product
and so to {\bf unitary representations}. 
Locally compact {\bf Abelian} groups $G$ can be understood by their {\bf characters}, 
continuous group homomorphisms from $G$ to the circle group $\mathbb{T}=\mathbb{R}/\mathbb{Z}$. The set of 
characters defines a new locally compact group $\hat{G}$, the {\bf dual} of $G$. The multiplication 
is the pointwise multiplication, the inverse is the complex conjugate and the topology is the one of {\bf uniform
convergence} on compact sets. If $G$ is compact, then $\hat{G}$ is discrete, and if $G$ is discrete, then $\hat{G}$
is compact. In order to prove {\bf Pontryagin duality} $\hat\hat{G}=G$, one needs a generalized
{\bf Fourier transform} $\hat{f}(\chi) = \int_G f(x) \overline{\chi(x)} d\mu(x)$ which uses the Haar measure.
The {\bf inverse Fourier transform} gives back $f$ using the {\bf dual Haar measure}. The Haar measure is
also used to define the {\bf convolution} $f \star g(x) = \int_G f(x-y) g(y) d\mu(y)$ rendering $L^1(G)$ a 
{\bf Banach algebra}. The Fourier transform then produces a homomorphism from $L^1(G)$ to $C_0(\hat{G})$ or a unitary 
transformation from $L^2(G)$ to $L^2(\hat{G})$. For literature, see \cite{Chevalley,Warner}. 
\index{unitary representation}
\index{Haar measure} 
\index{Banach algebra}
\index{Pontryagin duality}
\index{character}
\index{convolution}
\index{Fourier transform}

\section{Computability}

The class of {\bf general recursive functions} is the smallest class of functions which
allows {\bf projection}, {\bf iteration}, {\bf composition} and {\bf minimization}. The class of {\bf Turing 
computable functions} are the functions which can be implemented by a {\bf Turing machine} 
possessing finitely many states. Turing introduced this in 1936 \cite{AnnotatedTuring}.
\index{recursive function}
\index{general recursive function}
\index{Turing computable}
\index{Turing machine}

\satz{
The generally recursive class is the Turing computable class.
}

Kurt G\"odel and Jacques Herbrand defined the class of general recursive functions around 1933. They were motivated
by work of Alonzo Church who then created {\bf $\lambda$ calculus} later in 1936.
Alan Turing developed the idea of a {\bf Turing machine} which allows to replace Herbrand-G\"odel recursion and
$\lambda$ calculus. The {\bf Church thesis} or {\bf Church-Turing thesis} 
states that everything we can compute is generally recursive. As ``whatever we can compute" is not 
formally defined, this always will remain a thesis unless some more effective computation concept would emerge.  
\index{Church}
\index{Church thesis}
\index{Church-Turing thesis}

\section{Category theory}

Given an element $A$ in a {\bf category} $C$, let $h^A$ denote the {\bf functor} which 
assigns to a set $X$ the set ${\rm Hom}(A,X)$ of all {\bf morphisms} from $A$ to $X$. 
Given a {\bf functor} $F$ from $C$ to the category $S={\rm Set}$, let $N(G,F)$ be the
set of {\bf natural transformations} from $G=h^A$ to $F$. 
(A {\bf natural transformation} between two functors $G$ and $F$ from $C$ to $S$ 
assigns to every object $x$ in $C$ a morphism $\eta_x: G(x) \to F(x)$ such that for every 
morphism $f: x \to y$ in $C$ we have $\eta_y \circ G(f) = F(f) \circ \eta_x$.) 
The {\bf functor category} defined by $C$ and $S$ has as objects the functors $F$ and
as morphisms the natural transformations. The {\bf Yoneda lemma} is 
\index{category}
\index{functor}
\index{set category}
\index{Yoneda lemma} 
\index{natural transformation}
\index{morphism}

\satz{$N(h^A,F)$ can be identified with $F(A)$.}

Category theory was introduced in 1945 by Samuel Eilenberg and Sounders Mac Lane. 
The lemma above is due to Nobuo Yoneda from 1954. It allows to see
a category embedded in a {\bf functor category} which is a {\bf topos} and serves as a 
sort of completion. One can identify 
a set $S$ for example with ${\rm Hom}(1,S)$. An other example is {\bf Cayley's theorem} 
stating that the category of groups can be completely understood by looking at the 
group of permutations of $G$. For category theory, see \cite{McLarty,MacLaneCategory}. 
For history, \cite{Kroemer}. 
\index{pre-sheave}
\index{Cayley theorem}

\section{Perturbation theory}

A function $f$ of several variables is called {\bf smooth} if one can take {\bf first partial derivatives}
like $\partial_x,\partial_y$ and second partial derivatives like $\partial_x \partial_y  f(x,y) = f_{xy}(x,y)$ and still 
have continuous functions. Assume $f(x,y)$ is a {\bf smooth function} of two Euclidean variables $x,y \in \RR^n$.
If $f(a,0)=0$, we say $a$ is a {\bf root} of $x \to f(x,y)$. If $f_y(x_0,y)$ is invertible, the root is 
called {\bf non-degenerate}. If there is a solution $f(g(y),y)=0$ such that $g(0)=a$ and $g$ 
is continuous, the root $a$ has a {\bf local continuation} and say that it {\bf persists} 
under perturbation.
\index{smooth function}
\index{non-degenerate}
\index{local continuation}
\index{root}


\satz{A non-degenerate root persists under perturbation.}

This is the {\bf implicit function theorem}. There are concrete and fast algorithms 
to compute the continuation.  An example is the {\bf Newton method} which iterates
$T(x) = x -  f(x,y)/f_x(x,y)$ to find the roots of $x \to f(x,y)$ for fixed $y$. 
The importance of the implicit function theorem is both theoretical as well as applied. 
The result assures that one can makes statements about a complicated 
theory near some model, which is understood. There are related situations, like if we
want to continue a solution of $F(x,y) = (f(x,y),g(x,y))=(0,0)$ giving {\bf equilibrium
points} of the {\bf vector field} $F$. Then the Newton
step $T(x,y) = (x,y) - dF^{-1}(x,y) \cdot F(x,y)$ method allows a continuation if $dF(x,y)$
is invertible. This means that small deformations of $F$ do not lead to changes of the nature
of the equilibrium points. When equilibrium points change, the system exhibits {\bf bifurcations}.
This in particular applies to $F(x,y) = \nabla f(x,y)$, where equilibrium points
are {\bf critical points}. The derivative $dF$ of $F$ is then the {\bf Hessian}. 
\cite{KrantzParksImplicitFunction} call it one of the most important and oldest
pradigms in modern mathematics for which the germ of the idea was already formed 
in the writings of Isaac Newton and Gottfried Leibniz but only riped under
Augustin-Louis Cauchy to the theorem we know today.
\index{implicit function theorem}
\index{Newton method}
\index{Vector field}
\index{bifurcation}
\index{equilibrium point} 

\section{Counting}

A {\bf simplicial complex} $X$ is a finite set of non-empty sets that is closed under the operation
of taking finite non-empty subsets. The {\bf Euler characteristic} $\chi$ of a simplicial complex 
$G$ is defined as $\chi(X)=\sum_{x \in X} (-1)^{{\rm dim}(x)}$, where the {\bf dimension} 
${\rm dim}(x)$ of a set $x$ is its cardinality $|x|$ minus $1$. 
\index{dimension}
\index{cardinality}
\index{simplicial complex}
\index{Euler characteristic}

\satz{
$\chi(X \times Y) = \chi(X) \chi(Y)$.
}

For {\bf zero-dimensional simplicial complexes} $G$, (meaning that all sets in $G$ have cardinality $1$), 
we get the {\bf rule of product}: if you have $m$ ways to do 
one thing and $n$ ways to do an other, then there are $m n$ 
ways to do both. This {\bf fundamental counting principle} is used in 
probability theory for example. 
The {\bf Cartesian product} $X \times Y$ of two complexes is defined as the set-theoretical product of the 
two finite sets. It is not a simplicial complex any more in general but has the same Euler characteristic 
than its Barycentric refinement $(X \times Y)_1$, which is a simplicial complex. 
The maximal dimension of $A \times B$ is ${\rm dim}(A) + {\rm dim}(B)$ and $p_X(t) = \sum_{k=0}^n v_k(X) t^k$ 
is the generating function of $v_k(X)$, then $p_{X \times Y}(t) = p_X(t) p_Y(t)$ implying the counting principle as $p_X(-1) = \chi(X)$. 
The function $p_X(t)$ is called the {\bf Euler polynomial} of $X$. 
The importance of Euler characteristic as a {\bf counting tool}
lies in the fact that only $\chi(X)=p_X(-1)$ is invariant under 
{\bf Barycentric subdivision} $\chi(X)=X_1$, where $X_1$ is the complex which consists of
the vertices of all complete subgraphs of the graph in which the sets of $X$ are the vertices
and where two are connected if one is contained in the other. The concept of Euler characteristic goes so over
to continuum spaces like {\bf manifolds} where the product property holds too. 
See for example \cite{alexandroff}. 
\index{fundamental counting principle}
\index{rule of product}
\index{zero dimensional complexes}
\index{cardinality}
\index{Cartesian product} 
\index{Euler polynomial}
\index{Barycentric subdivision}

\section{Metric spaces}

A continuous map $T:X \to X$, where $(X,d)$ is a {\bf complete} non-empty {\bf metric space} is called
a {\bf contraction} if there exists a real number $0<\lambda<1$ such that $d(T(x),T(y)) \leq \lambda d(x,y)$
for all $x,y \in X$. The space is called {\bf complete} if every {\bf Cauchy sequence} 
in $X$ has a {\bf limit}.  (A sequence $x_n$ in $X$ is called {\bf Cauchy} if for all $\epsilon>0$, there exists
$n>0$ such that for all $i,j>n$, one has $d(x_i,x_j)<\epsilon$.)
\index{complete metric space}
\index{metric space}
\index{limit} 
\index{contraction}
\index{Cauchy sequence}

\satz{
A contraction has a unique fixed point in $X$.
}

This result is the {\bf Banach fixed point theorem} proven by Stefan Banach from 1922. 
The example case $T(x)=(1-x^2)/2$ on $X=\QQ \cap [0.3,0.6]$ having contraction rate 
$\lambda=0.6$ and $T(X)=\QQ \cap [0.32,0.455] \subset X$ shows that completeness is 
necessary. The unique fixed point of $T$ in $X$ is $\sqrt{2}-1=0.414...$ which is not in 
$\QQ$ because $\sqrt{2}=p/q$ would imply $2q^2=p^2$, which is not possible for integers 
as the left hand side has an odd number of prime factors $2$ while the right hand side 
has an even number of prime factors.  
\index{Banach fixed point theorem}
See \cite{PalaysBanach}
\index{completeness} 
\index{Banach fixed point theorem}


\section{Dirichlet series}

The {\bf abscissa of simple convergence} of a {\bf Dirichlet series}
$\zeta(s) = \sum_{n=1}^{\infty} a_n e^{-\lambda_n s}$ is
$\sigma_0 = \inf \{ a \in \R \; | \; \zeta(z)$ converges for all ${\rm Re}(z)>a \; \}$.
For $\lambda_n=n$ we have the {\bf Taylor series} $f(z) = \sum_{n=1}^{\infty} a_n z^n$ with $z=e^{-s}$. For
$\lambda_n=\log(n)$ we have the {\bf standard Dirichlet series} $\sum_{n=1}^{\infty} a_n/n^s$. For example, for $a_n=z^n$,  
one gets the {\bf poly-logarithm} ${\rm Li}_s(z)=\sum_{n=1}^{\infty} z^n/n^s$ and especially ${\rm Li}_s(1)=\zeta(s)$, the 
{\bf Riemann zeta function} or the {\bf Lerch transcendent} $\Phi(z,s,a) = \sum_{n=1}^{\infty} z^n/(n+a)^s$. 
Define $S(n) = \sum_{k=1}^n a_k$. The {\bf Cahen's formula} applies if the series $S(n)$ does not converge.
\index{abscissa of convergence}
\index{polylogarithm}
\index{Lerch transcendent}
\index{Taylor series}
\index{Dirichlet series}
\index{Cahen formula}
\index{standard Dirichlet series}

\satz{$\sigma_0 = \limsup_{n \to \infty} \frac{\log|S(n)|}{\lambda_n} \; . $}

There is a similar formula for the {\bf abscissa of absolute convergence} of $\zeta$ which is defined as
$\sigma_a = \inf \{ a \in \R \; | \; \zeta(z)$ converges absolutely for all ${\rm Re}(z)>a \; \}$. The
result is $\sigma_a = \limsup_{n \to \infty} \frac{\log(\overline{S}(n))}{\lambda_n}$,
For example, for the {\bf Dirichlet eta function} $\zeta(s) = \sum_{n=1}^{\infty} (-1)^{n-1}/n^s$ has the abscissa of
convergence $\sigma_0 = 0$ and the absolute abscissa of convergence $\sigma_a = 1$. The series
$\zeta(s) = \sum_{n=1}^{\infty} e^{i n^{\alpha}}/n^s$ has $\sigma_a=1$ and $\sigma_0 = 1-\alpha$.
If $a_n$ is multiplicative $a_{n+m}=a_n a_m$ for relatively prime $n,m$, then 
$\sum_{n=1}^{\infty} a_n/n^s = \prod_p (1+a_p/p^s+a_{p^2}/p^{2s} + \cdots )$ 
generalizes the {\bf Euler golden key formula} $\sum_n 1/n^s=\prod_p (1-1/p^s)^{-1}$. 
See \cite{Hardydivergent,HardyRiesz}.
\index{convergence}
\index{Euler golden key}
\index{Dirichlet Eta function}
\index{abscissa of convergence}
\index{abscissa of absolute convergence}

\section{Trigonometry}

Mathematicians had a long and painful struggle with the concept of {\bf limit}. One of the first to ponder
the question was Zeno of Elea around 450 BC \cite{Mazur2007}. Archimedes of Syracuse made some progress around 250 BC. 
Since Augustin-Louis Cauchy \cite{CauchyCoursDAnalyse} 
one uses the notion of limits. See also \cite{Grabiner1983}. Today,
one defines the {\bf limit} $\lim_{x \to a} f(x) = b$ {\bf to exist}
if and only if for all $\epsilon>0$, there exists a $\delta>0$ such that if $|x-a|<\delta$, then $|f(x)-b|<\epsilon$. 
A place where limits appear are when computing {\bf derivatives} $g'(0) = \lim_{x \to 0} [g(x)-g(0)]/x$.
In the case $g(x)=\sin(x)$, one has to understand the limit of 
the function $f(x)=\sin(x)/x$ which is the {\bf sinc} function. A prototype result
is the {\bf fundamental theorem of trigonometry} (called as such in some calculus texts like \cite{Bretscher2006}).
\index{limit}
\index{sinc function}
\index{trigonometry}

\satz{
$\lim_{x \to 0} \sin(x)/x = 1$. 
}

It appears strange to give weight to such a special result but it
explains the difficulty of limit and the {\bf l'H\^opital rule} of 1694,
which was formulated in a book of Bernoulli commissioned to H\^opital: the limit can
be obtained by differentiating both the denominator and nominator and taking the limit of
the quotients. The result allows to derive (using trigonometric identities) that in general
$\sin'(x)=\cos(x)$ and $\cos'(x)=-\sin(x)$. One single limit is the gateway.
It is important also culturally because it embraces thousands of years of struggle. 
It was Archimedes, who used the theorem when computing the {\bf circumference of the circle formula} $2\pi r$ 
using {\bf exhaustion} using regular polygons from the inside and outside. Comparing the 
lengths of the approximations essentially battled that fundamental theorem of trigonometry. 
The identity is therefore the epicenter around the development 
of {\bf trigonometry}, {\bf differentiation} and {\bf integration}.
\index{trigonometric functions}
\index{L'Hopital rule}
\index{fundamental theorem of trigonometry}
\index{exhaustion}
\index{circumference of circle}

\section{Logarithms}

The {\bf natural logarithm} is the inverse of the {\bf exponential function} $\exp(x)$
establishing so a {\bf group homomorphism} from the additive group $(\mathbb{R},+)$ to the
multiplicative group $(\mathbb{R}^+,*)$. We have:
\index{logarithm}
\index{natural logarithm}

\satz{$\log(u v) = \log(u) + \log(v)$.}

This follows from $\exp(x+y) = \exp(x) \exp(y)$ and $\log(\exp(x))=\exp(\log(x))=x$ by
plugging in $x=\log(u), y = \log(v)$.
The logarithms were independently discovered by Jost B\"urgi around 1600 and John Napier in 1614
\cite{StaudacherBuergi}. The {\bf logarithm} with base $b>0$ is denoted by $\log_b$. It
is the inverse of $x \to b^x = e^{x \log(b)}$. The concept of logarithm has been extended
in various ways: in any {\bf group} $G$, one can define the {\bf discrete logarithm} $\log_b(a)$ to base $b$
as an {\bf integer} $k$ such that $b^k = a$ (if it exists). For complex numbers the {\bf complex logarithm}
$\log(z)$ as any solution $w$ of $e^w = z$. It is {\bf multi-valued}
as $\log(|z|) + i {\rm arg}(z) + 2\pi i k$ all solve this with some integer $k$,
where ${\rm arg}(z) \in (-\pi,\pi)$. The identity $\log(u v) = \log(u) + \log(v)$ is now only true
up to $2\pi k i$. Logarithms can also be defined for matrices. Any matrix $B$ solving $\exp(B)=A$ is called a {\bf logarithm}
of $A$. For $A$ close to the identity $I$, can define $\log(A) = (A-I) - (A-I)^2/2 + (A-I)^3/3 - ... $,
which is a {\bf Mercator series}. For {\bf normal invertible matrices}, one can define logarithms using
the {\bf functional calculus} by diagonalization.  On a {\bf Riemannian manifold} $M$, one also has an
exponential map: it is a diffeomorphim from a small ball $B_r(0)$ in the {\bf tangent space}
$x \in M$ to $M$. The map $v \to \exp_x(v)$ is obtained by defining $\exp_x(0)=x$ and by taking for
$v \neq 0$ a {\bf geodesic} with initial direction $v/|v|$ and running it for time $|v|$. The
logarithm $\log_x$ is now defined on a {\bf geodesic ball} of radius $r$ and defines an element in the
tangent space. In the case of a Lie group $M=G$, where the points are matrices, each tangent space is
its {\bf Lie algebra}.
\index{discrete log}
\index{principal logarithm}
\index{argument}
\index{complex logarithm}
\index{diagonalization}
\index{exponential map}
\index{Lie group}

\section{Geometric probability}

A subset $K$ of $\RR^n$ is called {\bf compact} if it is {\bf closed} and {\bf bounded}. By {\bf Bolzano-Weierstrass}
this is equivalent to the fact that every infinite sequence $x_n$ in $K$ has a {\bf subsequence} which converges. 
A subset $K$ of $\RR^n$ is called {\bf convex}, if for any two given points $x,y \in K$, the interval 
$\{ x+t(y-x), t \in [0,1]  \}$ is a subset of $K$. 
Let $G$ be the set of all {\bf compact convex subsets} of $\RR^n$. An {\bf invariant valuation}
$X$ is a function $X: G \to \RR$ satisfying $X(A \cup B) + X(A \cap B)= X(A)+X(B)$,
which is continuous in the {\bf Hausdorff metric} 
$d(K,L)= {\rm max}( \sup_{x \in K} \inf_{y \in L} d(x,y) + 
\sup_{y \in K} \inf_{x \in L} d(x,y))$ and invariant under {\bf rigid motion}
generated by rotations, reflections and translations in the linear space $\RR^n$. 
\index{Hausdorff metric}
\index{valuation}
\index{invariant valuation}
\index{rigid motion}

\satz{The space of valuations is $(n+1)$-dimensional.}

The theorem is due to Hugo Hadwiger from 1937. 
The coefficients $a_j(G)$ of the polynomial
${\bf Vol}(G+tB) = \sum_{j=0}^n a_j t^j$ are a basis, where $B$ is the 
{\bf unit ball} $B=\{ |x| \leq 1 \}$.  See \cite{KlainRota}.
\index{Hadwiger}
\index{unit ball}

\section{Partial differential equations}

A {\bf quasilinear partial differential equation} is a differential equation
of the form $u_t(x,t)= F(x,t,u) \cdot \nabla_x u(x,t) + f(x,t,u)$ 
with analytic initial condition $u(x,0)=u_0(x)$ and an {\bf analytic} {\bf vector 
field} $F$. It defines a {\bf quasi-linear Cauchy problem}. 
\index{partial differential equation}
\index{quasi-linear Cauchy problem}
\index{vector field}

\satz{A quasi-linear Cauchy problem has a unique analytic solution.}

This is the {\bf Cauchy-Kovalevskaya theorem}. 
It was initiated by Augustin-Louis Cauchy in 1842
and proven in 1875 by Sophie Kowalevskaya. Analyticity is important, 
smoothness alone is not enough.  
If $F$ is analytic in each variable, one can look at equations like
the Cauchy problem $u_t = F(t,x,u,u_x,u_{xx})$. Examples are
partial differential equations like the heat equation $u_t=u_{xx}$ or the
wave equation $u_{tt} = u_{xx}$. Given an initial condition 
$u(0,x) = u_0(x)$ one then deals with an ordinary differential equation in a function
space. One can then try to approach the Cauchy-Kovalevskaya problem by Picard-Lindel\"of.
The problem is that the Lipschitz condition fails because
the corresponding operators are unbounded. Even Cauchy-Peano (which does not ask for 
uniqueness) fails. And this even in an analytic setting. \cite{Mouhot} gives
the example $u_t = u_{xx}$ with initial condition $u(0,x) = 1/(1+x^2)$ for which 
the entire series solving the problem has a zero radius of convergence in $x$ for any $t>0$.
Texts like \cite{TaylorPDE,Mouhot} give full versions of the Cauchy-Kovalevskaya theorem 
for real-analytic Cauchy initial data on a real analytic hypersurface satisfying a 
non-characteristic condition for the partial differential equation. 
For a shorter introduction to partial differential equations, see \cite{Arnold2004}.
\index{Kowalevskaya}
\index{Cauchy} 
\index{Cauchy-Kovalevskaya theorem}

\section{Game theory}

If $S=(S_1,\dots, S_n)$ are $n$ {\bf players} and $f=(f_1,\dots, f_n)$ is
a {\bf payoff function} defined on a {\bf strategy profile} $x=(x_1, \dots, x_n)$. 
A point $x^*$ is called an {\bf equilibrium} if $f_i(x^*)$ is {\bf maximal}
with respect to changes of $x_i$ alone in the profile $x$ for every player $i$. 
\index{player}
\index{payoff function}
\index{equilibrium} 
\index{maximal equilibrium} 
\index{strategy profile}

\satz{There is an equilibrium for any game with mixed strategy}

The equilibrium is called a {\bf Nash equilibrium}. It tells us what we would see in 
a world if everybody is doing their best, given what everybody else is doing. 
John Forbes Nash used in 1950 the {\bf Brouwer fixed point theorem} and later in 1951 the 
{\bf Kakutani fixed point theorem} to prove it. The Brouwer fixed point theorem itself is generalized
by the {\bf Lefschetz fixed point theorem} which equates the super trace of the induced
map on cohomology with the sum of the indices of the fixed points.  About John Nash and some history
of game theory, see \cite{SiegfriedNash}: game theory started maybe with Adam Smith's ``the 
Wealth of Nations" published in 1776, Ernst Zermelo in 1913 (Zermelo's theorem), \'Emile Borel in the 1920s 
and John von Neumann in 1928 pioneered mathematical game theory. Together with Oskar Morgenstern, John
von Neumann merged game theory with economics in 1944. Nash published his thesis in a paper
of 1951. For the mathematics of games, see \cite{WebbGameTheory}.
\index{Nash equilibrim} 
\index{Brouwer fixed point theorem} 
\index{Kakutani fixed point theorem}
\index{Kakutani fixed point theorem}
\index{Lefschetz fixed point theorem}
\index{Kakutani fixed point theorem}

\section{Measure theory}

A topological space with open sets $\mathcal{O}$ defines
the {\bf Borel $\sigma$-algebra}, the smallest $\sigma$ algebra
which contains $\mathcal{O}$. For the metric space $(\mathbb{R},d)$
with $d(x,y)=|x-y|$, already the intervals generate the Borel $\sigma$
algebra $\mathcal{A}$. A {\bf Borel measure} is a  measure defined on a
Borel $\sigma$-algebra. Every {\bf Borel measure} $\mu$ on the real line $\RR$ can be 
decomposed uniquely into an {\bf absolutely continuous} part 
$\mu_{ac}$, a {\bf singular continuous} part $\mu_{sc}$ and
a {\bf pure point} part $\mu_{pp}$:
\index{Borel measure}
\index{absolutely continuous measure}
\index{singular continuous measure}
\index{pure point measure}

\satz{ $\mu = \mu_{ac} + \mu_{sc} + \mu_{pp}$. }  

This is called the {\bf Lebesgue decomposition theorem}. 
It uses the {\bf Radon-Nikodym theorem}. The decomposition 
theorem implies the decomposition theorem of the {\bf spectrum}
of a linear operator. See \cite{Simon2017} (like page 259). 
Lebesgue's theorem was published in 1904. A generalization 
due to Johann Radon and Otto Nikodym was done in 1913. 
\index{Lebesgue measure}
\index{Radon-Nikodym theorem}
\index{spectrum} 

\section{Geometric number theory}

If $\Gamma$ is a {\bf lattice} in $\RR^n$, denote with $\RR^n/\Gamma$ the 
{\bf fundamental region} and by $|\Gamma|$ its {\bf volume}. A set $K$ is {\bf convex} 
if $x,y \in K$ implies $x+t(x-y) \in K$ for all $0 \leq t \leq 1$. A set 
$K$ is {\bf centrally symmetric} if $x \in K$ implies $-x \in K$. A region is {\bf Minkowski}
if it is convex and centrally symmetric. Let $|K|$ denote the volume of $K$.
\index{lattice}
\index{fundamental region} 
\index{volume}
\index{convex}
\index{Minkowski}
\index{centrally symmetric}

\satz{If $K$ is Minkowski and $|K| >2^n |\Gamma|$ then $K \cap \Gamma \neq \emptyset$. }

The theorem is due to Hermann Minkowski in 1896. It lead to a field called {\bf geometry
of numbers}. \cite{OldsLaxDavidoff}. It has many applications in number theory
and {\bf Diophantine analysis} \cite{BurgerNumberJungle,Hua1982}
\index{Diophantine analysis}
\index{Number theory}

\section{Fredholm}

An {\bf integral kernel} $K(x,y) \in L^2([a,b]^2)$ 
defines an {\bf integral operator} $A$ defined by $A f(x) = \int_a^b K(x,y) f(y) \; dy$ with
adjoint $T^* f(x) = \int_a^b \overline{K(y,x)} f(y) \; dy$. The $L^2$ assumption makes the 
function $K(x,y)$ what one calls a {\bf Hilbert-Schmidt kernel}.
Fredholm showed that the {\bf Fredholm equation} $A^*f=(T^*-\overline{\lambda}) f = g$ has a solution $f$ if and 
only if $f$ is perpendicular to the kernel of $A=T-\lambda$. 
This identity ${\rm ker}(A)^{\perp} = {\rm im}(A^*)$ is in finite dimensions
part of the {\bf fundamental theorem of linear algebra}. The {\bf Fredholm alternative} reformulates
this in a more catchy way as an {\bf alternative}: 
\index{integral operator}

\satz{Either $\exists f \neq 0$ with $Af=0$ or for all $g$, $\exists f$ with $Af=g$. }

In the second case, the solution depends continuously on $g$. The alternative 
can be put more generally by stating that if $A$
is a {\bf compact operator} on a Hilbert space and $\lambda$ is not an eigenvalue of $A$, 
then the {\bf resolvent} $(A-\lambda)^{-1}$ is bounded. A 
bounded operator $A$ on a Hilbert space $H$ is called {\bf compact}
if the image of the unit ball is relatively compact (has a compact closure).  
The Fredholm alternative is part of {\bf Fredholm theory}. It was developed by 
Ivar Fredholm in 1903. 
\index{Hilbert Schmidt kernel}
\index{Fredholm alternative}
\index{Compact operator}
\index{Fredholm theory}

\section{Prime distribution}

The {\bf Dirichlet theorem} about
the primes along an arithmetic progression tells that if $a$ and $b$ are {\bf relatively prime}
meaning that there largest common divisor is $1$, then
there are infinitely many primes of the form $p=a \; {\rm mod} \; b$. 
The Green-Tao theorem strengthens this. We say that a set $A$ contains {\bf arbitrary long
arithmetic progressions} if for every $k$ there exists an {\bf arithmetic progression}
$\{ a + b j, j=1,\cdots, k\}$ within $A$. 
\index{Dirichlet prime number theorem}
\index{arithmetic progression} 
\index{Green Tao Theorem}

\satz{The set of primes contains arbitrary long arithmetic progressions.}

The {\bf Dirichlet prime number theorem} was found in 1837. 
The {\bf Green-Tao theorem} was done in 2004 and appeared in 2008 \cite{GreenTao}. 
It uses {\bf Szemer\'edi's theorem} \cite{FurstenbergRecurrence} which
shows that any set $A$ of positive upper density $\limsup_{n \to \infty} |A \cap \{ 1 \cdots n\}|/n$ 
has arbitrary long arithmetic progressions. So, any subset $A$ of the primes $P$ for which the 
{\bf relative density}  $\limsup_{n \to \infty} |A \cap \{1 \cdots n\}|/|P \cap \{ 1 \cdots n\}|$ is positive
has arbitrary long arithmetic progressions.  For non-linear sequences of numbers the problems are
wide open. The {\bf Landau problem} of the infinitude of primes of the form $x^2+1$ illustrates this. 
The Green-Tao theorem gives hope to tackle the {\bf Erd\"os conjecture on arithmetic progressions}
telling that a sequence $\{ x_n \}$ of integers satisfying $\sum_n x_n=\infty$ contains arbitrary long
arithmetic progressions.
\index{Landau problem}
\index{relative density}
\index{Szemer\'edi theorem}
\index{Erd\"os conjecture on arithmetic progressions}

\section{Riemannian geometry}

A {\bf Riemannian manifold} is a smooth finite dimensional manifold $M$ equipped with 
a {\bf smooth}, {\bf symmetric}, {\bf positive definite tensor} $g$ defining on each 
{\bf tangent space} $T_xM$ an {\bf inner product} $(u,v)_x = (g(x) u,v) = \sum_{i,j} g_{ij}(x) u^i v^j$.
Let $\Omega$ be the space of {\bf smooth vector fields} $X$ on $M$. A vector field $X$ acts on smooth 
functions $f$ as directional derivative $X f = \delta_X f$. Given two vector fields $X,Y$ on $M$,
one has at each point $x \in M$ a number $g(X,Y) = (g(x) X(x),Y(x))$ so that $g(X,Y)$ is a smooth function
on $M$. A {\bf connection} is a bilinear map $(X,Y) \to \nabla_X Y$ from $\Omega \times \Omega$
to $\Omega$ satisfying the differentiation rules $\nabla_{fX} Y = f \nabla_X Y$ and
{\bf Leibniz rule} $\nabla_{X} (fY) = df(X) Y + f \nabla_X Y$. It is {\bf compatible with the metric}
if the {\bf Lie derivative} satisfies $\delta_X g(Y,Z) = g(\Gamma_X Y,Z) + g(Y,\Gamma_X Z)$.
It is {\bf torsion-free} if $\nabla_X Y - \nabla_Y X = [X,Y]$ is the {\bf Lie bracket} on $\Omega$.
\index{Riemannian manifold}
\index{symmetric tensor}
\index{positive definite tensor}
\index{inner product}
\index{connection}
\index{Leibniz rule} 
\index{torsion free}
\index{smooth vector field}

\satz{There is exactly one torsion-free connection compatible with $g$.}

This is the {\bf fundamental theorem of Riemannian geometry.} 
The connection is called the {\bf Levi-Civita connection}, named after Tullio Levi-Civita.
One proof goes by establishing the {\bf Koszul formula} which determines $\nabla_X Y$
explicitely
$$ 2g(\nabla_X Y,Z) =   X g(Y,Z) +   Y g(X,Z)  -   Z g(X,Y) 
                    - g(X, [Y,Z])- g(Y, [X,Z]) + g(Z, [X,Y]) \; . $$
See for example \cite{DFN,AMR,Spivak1999,Cycon}.

\index{fundamental theorem of Riemannian geometry}
\index{Levi Civita}
\index{Koszul formula}
\index{Levi-Civita connection}

\section{Symplectic geometry}

A {\bf symplectic manifold} $(M,\omega)$ is a smooth $2n$-manifold $M$ equipped with a non-degenerate
closed $2$-form $\omega$. The later is called a {\bf symplectic form}. 
As a 2-form, it satisfies $\omega(x,y)= -\omega(y,x)$. 
{\bf Non-degenerate} means $\omega(u,v)=0$ for all $v$ implies $u=0$. 
The {\bf standard symplectic form} is $\omega_0 = \sum_{i<j} dx_i \wedge dx_j$. 
\index{symplectic manifold}
\index{symplectic form}
\index{non-degenerate 2-form}
\index{standard symplectic form} 

\satz{Every symplectic form is locally diffeomorphic to $\omega_0$.}

This theorem is due to Jean Gaston Darboux from 1882. Modern proofs use {\bf Moser's trick}
from 1965 (i.e. \cite{HoferZehnder1994}). The Darboux theorem assures that locally, 
two symplectic manifolds of the same dimension are symplectic equivalent. 
It also implies that {\bf symplectic matrices} $A$ ($2n \times 2n$ matrices satisfying $A^T J A = J$ 
with skew symmetric $J=\left[\begin{array}{cc} 0 & I \\  -I & 0 \end{array} \right]$)
have {\bf determinant} $1$ which is not obvious as applying the determinant to $A^T J A = J$ 
only establishes ${\rm det}(A)^2=1$.
In contrast, for {\bf Riemannian manifolds}, one can not trivialize the Riemannian
metric in a neighborhood one can only render it the standard metric at the point itself. 
\index{Moser trick}
\index{Darboux}
\index{determinant} 
\index{symplectic matrix} 

\section{Differential topology}

Given a {\bf smooth function} $f$ on a {\bf differentiable manifold} $M$. Let $df$ denote
the {\bf gradient} of $f$. A point $x$ is called a {\bf critical point}, if $df(x)=0$. 
We assume $f$ has only finitely many {\bf critical points} and that
all of them are {\bf non-degenerate}. The later means that the {\bf Hessian}
$d^2 f(x)$ is invertible at $x$. One calls such functions {\bf Morse functions}. 
The {\bf Morse index} of a critical point $x$ is the number of negative eigenvalues of $d^2f$. 
The {\bf Morse inequalities} relate the number $c_k(f,K)$ of critical points of index $k$
of $f$ with the {\bf Betti numbers} $b_k(M)$, defined as the nullity of the
 {\bf Hodge star operator} $d d^* + d^* d$ restricted to $k$-forms 
$\Omega_k$, where $d_k: \Omega_k  \to \Omega_{k+1}$ 
is the {\bf exterior derivative}.
\index{exterior derivative}
\index{Morse index}
\index{critical point} 
\index{non-degenerate critical point}
\index{manifold}
\index{differentiable manifold}
\index{Hodge operator} 
\index{Betti number}
\index{Morse inequality}
\index{Morse function}

\satz{$c_k -c_{k-1} + \cdots + (-1)^k c_0 \geq b_k-b_{k-1} + \cdots + (-1)^k b_0$.}

These are the {\bf Morse inequalities} due to Marston Morse from 1934. 
It implies in particular the {\bf weak Morse inequalities} $b_k \leq c_k$. 
Modern proofs use {\bf Witten deformation} \cite{Cycon} of the exterior derivative $d$. 
\index{strong Morse inequality} 
\index{weak Morse inequalities}
\index{Witten deformation}

\section{Non-commutative geometry}

A {\bf spectral triple} $(A,H,D)$ is given by a {\bf Hilbert space} $H$, a 
{\bf $C^*$-algebra} $A$ of operators on $H$ and a densely defined {\bf self-adjoint} operator $D$
satisfying $||[D,a]|| < \infty$ for all $a \in A$ the operator $e^{-tD^2}$ is 
{\bf trace class}. The operator $D$ is called a {\bf Dirac operator}. The set-up generalizes 
Riemannian geometry because of the following result dealing with the 
{\bf exterior derivative} $d$ on a Riemannian manifold $(M,g)$, where $A=C(M)$ is the 
$C^*$-algebra of continuous functions and $D=d+d^*$ is the Dirac operator, defining
a spectral triple for $(M,g)$. Let $\delta$ denote the {\bf geodesic distance} in $(M,g)$: 
\index{Spectral triple}
\index{Hilbert space}
\index{C star algebra}
\index{exterior derivative}
\index{Dirac operator}

\satz{ $\delta(x,y) = {\rm sup}_{f \in A, ||[D,f]|| \leq 1} |f(x)-f(y)|$. } 

This formula of Alain Connes tells that the spectral triple determines the geodesic distance in $(M,g)$
and so the metric $g$. It justifies to look at spectral triples as 
non-commutative generalizations of Riemannian geometry. See \cite{Connes}.
\index{Connes formula} 
\index{non-commutative geometry} 
\index{geodesic distance} 


\section{Polytopes}

A {\bf convex polytop} $P$ in dimension $n$ is the {\bf convex hull} of finitely many points in $R^n$. 
One assumes all vertices to be {\bf extreme points}, points which do not lie in an open line segment
of $P$. The {\bf boundary} of $P$ is formed by $(n-1)$ dimensional boundary facets. The notion
of {\bf Platonic solid} is recursive. 
A convex polytope is {\bf Platonic}, if all its facets are Platonic $(n-1)$-dimensional polytopes and
vertex figures. Let $p=(p_2,p_3,p_4, \dots )$ encode the number of Platonic solids meaning 
that $p_d$ is the number of Platonic polytops in dimension $d$. 
\index{Platonic polytop}

\satz{There are 5 platonic solids and $p=(\infty,5,6,3,3,3,\dots)$}

In dimension 2, there are infinitely many. They are the {\bf regular polygons}. The list of Platonic
solids is ``octahedron", ``dodecahedron", ``icosahedron", ``tetrahedron" and ``cube" 
has been known by the Greeks already. Ludwig Schl\"afli first classified the higher dimensional case.
There are six in dimension 4: they are the ``5 cell", the ``8 cell" ({\bf tesseract}), the ``16 cell", the ``24 cell", 
the ``120 cell" and the ``600 cell". There are only three regular polytopes 
in dimension 5 and higher, where only the analog of the 
tetrahedron, cube and octahedron exist.  For literature, see \cite{gruenbaum,Ziegler,Richeson}.
\index{regular polygon}
\index{octahedron}
\index{dodecahedron}
\index{icosahedron}
\index{tetrahedron}
\index{cube}
\index{tesseract}

\section{Descriptive set theory} 

A {\bf metric space} $(X,d)$ is a set with a {\bf metric} $d$ (a function $X \times X \to [0,\infty)$
satisfying {\bf symmetry} $d(x,y)=d(y,x)$, the {\bf triangle inequality} $d(x,y) + d(y,z) \geq d(x,z)$, and
$d(x,y)=0 \leftrightarrow x=y$.)
A metric space $(X,d)$ is {\bf complete} if every {\bf Cauchy sequence} converges
in $X$. A metric space is of {\bf second Baire category} if the intersection of a 
countable set of open dense sets is dense. The {\bf Baire Category theorem} tells 
\index{Second Baire category}
\index{Complete metric space}
\index{Generic}
\index{Baire category theorem}

\satz{Complete metric spaces are of second Baire category.}

One calls the intersection $A$ of a countable set of open dense sets $A$ in $X$ also a
{\bf generic set} or {\bf residual set}. The complement of a generic 
set is also called a {\bf meager set} or {\bf negligible} or a set of {\bf first category}. 
Such a set is the union of countably many nowhere dense sets. Like measure theory, 
Baire category theory can be used to get existence results. There can be surprises: a generic
continuous function is not differentiable for example. For descriptive set theory, 
see \cite{Kechris}. The frame work for classical descriptive set theory often are
{\bf Polish spaces}, which are separable complete metric spaces. See \cite{Bredon}.
\index{first Baire category}
\index{meager set}
\index{meagre set} 
\index{Polish space}
\index{residual set}

\section{Calculus of variations}

Let $X$ be the vector space of {\bf smooth}, {\bf compactly supported} functions $h$ 
on an interval $(a,b)$. The {\bf fundamental lemma of calculus of variations}
tells
\index{compactly supported}
\index{smooth}
\index{fundamental lemma of calculus of variations}

\satz{ $\int_a^b f(x) g(x) dx = 0$ for all $g \in X$, then $f=0$. }

The result is due to Joseph-Louis Lagrange. 
One can restate this as the fact that if $f=0$ {\bf weakly} then $f$ is actually zero. 
It implies that if $\int_a^b f(x) g'(x) \; dx=0$ for all $g \in X$, then 
$f$ is constant. This is nice as $f$ is not assumed to be differentiable. The
result is used to prove that extrema to a {\bf variational problem} $I(x) = \int_a^b L(t,x,x') \; dt$
are weak solutions of the {\bf Euler Lagrange equations} $L_x = d/dt L_{x'}$. 
See \cite{GH,MoserVariations}.
\index{Euler-Lagrange equations}
\index{variational problem}
\index{weak solutions}

\section{Integrable systems}

Given a {\bf Hamilton differential equation} $x' = J \nabla H(x)$ on a compact 
{\bf symplectic $2n$-manifold} $(M,\omega)$. The {\bf almost complex structure} $J: T^*M \to TM$ 
is tied to $\omega$ using a Riemannian metric $g$ by $\omega(v,w)=\langle v,J g \rangle$.
A function $F:M \to \mathbb{R}$ is called an {\bf first integral} if $d/dt F(x(t)) = 0$ for all $t$. An example
is the {\bf Hamiltonian function} $H$ itself. A set of integrals $F_1, \dots, F_k$ {\bf Poisson commutes}
if $\{ F_j,F_k \} = J \nabla F_j \cdot \nabla F_k =0$ for all $k,j$. They are {\bf linearly independent}, 
if at every point the vectors $\nabla F_j$ are linearly independent in the sense of linear algebra. 
A system is {\bf Liouville integrable} if there
are $d$ linearly independent, Poisson commuting integrals. 
The following theorem due to Liouville and Arnold
characterizes the {\bf level surfaces} $\{ F=c \} = \{ F_1=c_1, \dots F_d=c_d \}$:
\index{Liouville integrable} 
\index{first integral}
\index{Poisson commute}
\index{Almost complex structure}
\index{symplectic manifold}

\satz{For a Liouville integrable system, level surfaces $F=c$ are tori.}

An example how to get integrals is to write the system as an {\bf isospectral deformation} of an operator $L$.
This is called a {\bf Lax system}. 
Such a differential equation has the form $L' = [B,L]$, where $B=B(L)$ is skew symmetric. An example is the 
{\bf periodic Toda system} $\dot{a}_{n}=a_{n}(b_{n+1}-b_{n})$, $\dot{b}_{n}=2(a_{n}^{2}-a_{n-1}^{2})$,
where $(Lu)_n=a_n u_{n+1}+a_{n-1} u_{n-1} + b_n u_n$ and $(Bu)_n=a_n u_{n+1}-a_{n-1} u_{n-1}$. 
An other example is the motion of a {\bf rigid body} in $n$ dimensions if the center of mass
is fixed. See \cite{Arnold1980}. 
\index{isospectral deformation}
\index{Toda system}
\index{Lax system}
\index{rigid body}
\index{center of mass}

\section{Harmonic analysis}

On the vector space $X$ of continuously differentiable $2\pi$ periodic, complex-
valued functions, define the {\bf inner product} $(f,g) = (2\pi)^{-1} \int f(x) \overline{g}(x)  \; dx$. 
The {\bf Fourier coefficients} of $f$ are $\hat{f}_n = (f,e_n)$, where $\{ e_n(x) = e^{i n x}\}_{n \in \mathbb{Z}}$ 
is the {\bf Fourier basis}. The {\bf Fourier series} of $f$ is the sum 
$\sum_{n \in \mathbb{Z}} \hat{f}_n e^{i n x}$. 
\index{inner product}
\index{Fourier basis}

\satz{ The Fourier series of $f \in X$ converges point-wise to $f$. }

Already Fourier claimed this always to be true in his ``Th\'eorie Analytique de la Chaleur".
After many fallacious proofs, Dirichlet gave the first proof of convergence \cite{KoernerFourier}.
The case is subtle because there are continuous functions for which the convergence fails
at some points. Lip\'ot F\'ejer was able to show that for a continuous function $f$, 
the coefficients $\hat{f}_n$ nevertheless determine the function using 
{\bf C\'esaro convergence}. See \cite{Katznelson}. 
\index{Fourier coefficients}
\index{Fourier series}
\index{C\'esaro convergence}
\index{F\'ejer kernel}

\section{Jensen inequality}

If $V$ is a {\bf vector space}, a set $X$ is called {\bf convex} if for
all points $a,b \in X$, the {\bf line segment} 
$\{ t b+ (1-t) a \; | \; t \in [0,1] \}$ is contained in $X$.
A real-valued function $\phi:X \to \mathbb{R}$  is called {\bf convex}
if $\phi(tb + (1-t) a) \leq t \phi(b) +(1-t) \phi(a)$ for all $a,b \in X$ 
and all $t \in [0,1]$. 
Let now $(\Omega,\mathcal{A},{\rm P})$ be a {\bf probability space},
and $f \in L^1(\Omega,{\rm P})$ an {\rm integrable function}. We write
${\rm E}[f] = \int_{\omega} f(x) \; dP(x)$ for the {\bf expectation}
of $f$. For any convex $\phi: \mathbb{R} \to \mathbb{R}$ and
$f \in L^1(\Omega,P)$, we have the {\bf Jensen inequality}
\index{integrable function}
\index{convex}
\index{probability space}
\index{line segment}

\satz{$\phi({\rm E}[f])  \leq {\rm E}[ \phi(f) ]$. }

For $\phi(x) = \exp(x)$ and a finite probability space
$\Omega=\{1,2, \dots, n\}$ with $f(k) = x_k=\exp(y_k)$
and ${\rm P}[\{x\}] = 1/n$, this gives the 
{\bf arithmetic mean- geometric mean inequality}
$(x_1 \cdot x_2 \cdots x_n)^{1/n} \leq (x_1 + x_2 + \cdots + x_n)/n$.
The case $\phi(x)=e^x$ is useful in general as it 
leads to the inequality $e^{{\rm E}[f]} \leq {\rm E}[e^f]$ if $e^f \in L^1$.
For $f \in L^2(\omega,P)$ one gets $({\rm E}[f])^2  \leq {\rm E}[f^2]$ which reflects 
the fact that ${\rm E}[f^2]-({\rm E}[f])^2={\rm E}[(f-{\rm E}[f])^2]={\rm Var}[f] \geq 0$
where ${\rm Var}[f]$ is the {\bf variance} of $f$. 
\index{arithmetic mean}
\index{geometric mean}
\index{variance}
\index{line segment}
\index{convex set}
\index{convex function}
\index{Jensen inequality}

\section{Jordan curve theorem}

A {\bf closed curve} in the image of a continuous map $r: \mathbb{T} \to \mathbb{R}^2$.
It is called {\bf simple}, if this map $r$ is injective. One then calls the map an {\bf embedding}
and the image a {\bf topological 1-sphere}, meaning that it is homeomorphic to the standard circle
$x^2+y^2=1$ in $\mathbb{R}^2$. The image is then called a {\bf Jordan curve}.
The {\bf Jordan curve theorem} deals with such simple closed curves $S$ in the two-dimensional plane.
\index{Jordan curve}
\index{embedding}
\index{simple}

\satz{ A simple closed curve divides the plane into two regions. }

The Jordan curve theorem is due to Camille Jordan. His proof \cite{Jordan} was objected at
first \cite{Kline42} but rehabilitated in \cite{HalesJordanProof}.
The theorem can be strengthened, a {\bf theorem of Schoenflies} tells that each of the two regions
is homeomorphic to the disk $\{ (x,y) \in \mathbb{R}^2 \; | \; x^2+y^2 <1 \}$. In the smooth
case, it is even possible to extend the map to a diffeomorphism in the plane.
In higher dimensions, one knows that an embedding of the $(d-1)$ dimensional sphere in a $\mathbb{R}^d$
divides space into two regions. This is the {\bf Jordan-Brouwer} separation theorem.
It is no more true in general that the two parts are homeomorphic to
$\{ x \in \mathbb{R}^d \; | \;|x| <1 \}$: a counter example is the {\bf Alexander horned sphere}
which is a topological $2$-sphere but where the unbounded component is not simply connected and so not
homeomorphic to the complement of a unit ball. See \cite{Bredon}. 
\index{Schoenflies theorem}
\index{Jordan-Brouwer separation theorem}

\section{Chinese remainder theorem}

Given integers $a,b$, a {\bf linear modular equation} or {\bf congruence}
$a x + b = 0 \; {\rm \; mod \; m}$ asks to find an integer $x$ such
that $a x + b$ is divisible by $m$. This linear equation can always be solved if $a$ and
$m$ are coprime. The {\bf Chinese remainder theorem} deals with the {\bf system of linear modular equations}
$x=b_1 \; {\rm mod} \; m_1, x=b_2 \; {\rm mod} \; m_2, \dots, x=b_n \; {\rm mod} \; m_n$, where $m_k$
are the {\bf moduli}. More generally, for an integer $n \times n$ matrix $A$ 
we call $A x = b {\rm mod} \; m$ a {\bf Chinese remainder theorem system} or shortly {\bf CRT system}
if the $m_j$ are pairwise relatively prime and in each row there is a matrix element $A_{ij}$ relatively
prime to $m_i$. 
\index{CRT system}
\index{moduli}

\satz{Every Chinese remainder theorem system has a solution.}

The classical single variable case case is when $A_{i1}=1$ and $A_{ij}=0$ for $j>1$. 
Let $M=m_1 \cdots m_2 \cdots m_n$ be the product. In this one-dimensional case,
the result implies that $x {\rm mod} \; M$
$\to (x \; {\rm mod} \;$ $m_1, \dots, (x \; {\rm mod} \; m_n)$ is a ring isomorphism.
Define $M_i=M/m_i$.  An explicit algorithm is to finding numbers $y_i,z_i$ with
$y_i M_i + z_i m_i = 1$ (finding $y,z$ solving $a y + b z = 1$ for coprime $a,b$ is computed
using the {\bf Euclidean algorithm}), then finding $x=b_1 m_1 y_1 + \cdots + b_n m_n y_n$.
\cite{DingPeiSalomaa,Martzloff}. The multi-variable version appeared in 2005 
\cite{Brandeis,KnillChinese} and can be found also in \cite{Sury2015}.
\index{moduli}
\index{coprime}
\index{Euclidean algorithm}

\section{B\'ezout's theorem} 

A polynomial is {\bf homogeneous} if the total degree of all its {\bf monomials} is the same. 
A {\bf homogeneous polynomial} $f$ in $n+1$ variables of degree $d \geq 1$ defines 
a {\bf projective hypersurface} $f=0$. 
Given $n$ projective irreducible hypersurfaces $f_k=c_k$ of degree $d_k$ in a {\bf projective space} $\mathbb{P}^n$
we can look at the solution set $\{ f = c \} = \{f_1=c_1, \cdots, f_k=c_k \}$ of a system of nonlinear equations. 
The {\bf B\'ezout's bound} is $d=d_1 \cdots d_k$ the product of the degrees. 
{\bf B\'ezout's theorem} allows to count the number of solutions
of the system, where the number of solutions is counted with multiplicity. 
\index{hypersurface}
\index{B\'ezout's bound}
\index{projective space}
\index{monomial}
\index{homogeneous polynomial}
\index{hypersurface}

\satz{The set $\{f=c\}$ is either infinite or has $d$ elements.}

B\'ezout's theorem was stated in the ``Principia" of Newton in 1687 but was proven fist only in 1779 by \'Etienne B\'ezout.
If the hypersurfaces are all {\bf irreducible} and in ``general position", then there are exactly $d$ solutions and
each has multiplicity $1$. This can be used also for affine surfaces. If $y^2-x^3-3 x - 5=0$ is an {\bf elliptic
curve} for example, then $y^2 z -x^3-3xz^2-5z^3=$ is a projective hypersurface, its {\bf projective completion}.
B\'ezout's theorem implies part the fundamental theorem of algebra as for $n=1$, when we have only
one homogeneous equation we have $d$ roots to a polynomial of degree $d$. The theorem implies for example that the 
intersection of two {\bf conic sections} have in general $2$ intersection points.
The example $x^2-yz=0,x^2+z^2-yz=0$ has only the solution $x=z=0,y=1$ but with multiplicity $2$. 
As non-linear systems of equations appear frequently in {\bf computer algebra} this theorem gives a lower 
bound on the computational complexity for solving such problems. 
\index{general position}
\index{conic section}
\index{projective completion}
\index{elliptic curve}
\index{computer algebra}

\section{Group theory}

A {\bf finite group} $(G,*,1)$ is a finite set containing a {\bf unit} $1 \in G$
and a binary operation $*: G \times G \to G$ satisfying the {\bf associativity property} $(x * y) * z  = x * (y * z)$
and such that for every $x$, there exists a unique $y=x^{-1}$ such that $x*y=y*x=1$. The
{\bf order} $n$ of the group is the number of elements in the group. An element $x \in G$ generates
a {\bf subgroup} formed by $1,x,x^2=x*x, \dots$. This is the {\bf cyclic subgroup} $C(x)$ generated by $x$. 
{\bf Lagrange's theorem} tells
\index{group}
\index{subgroup}
\index{Lagrange theorem}
\index{associativity}
\index{finite group}
\index{cyclic subgroup}

\satz{ $|C(x)|$ is a factor of $|G|$}

The origins of group theory go back to Joseph Louis Lagrange, Paulo Ruffini and \'Evariste Galois. 
The concept of abstract group appeared first in the work of Arthur Cayley. 
Given a subgroup $H$ of $G$, the {\bf left cosets} of $H$ are the equivalence classes of the 
equivalence relation $x \sim y$ if there exists $z \in H$ with $x = z*y$. The equivalence classes
$G/N$ partition $G$. The number $[G:N]$ of elements in $G/H$ is called the {\bf index} of $H$ in $G$. 
It follows that $|G| = |H| [G:H]$ and more generally that if $K$ is a subgroup of $H$ and $H$ is a subgroup of $G$
then $[G:K] = [G:H] [H:K]$. 
The group $N$ generated by $x$ is a called a {\bf normal group} $N \triangleleft G$ if for 
all $a \in N$ and all $x$ in $G$ the element $x * a * x^{-1}$ is in $N$. This can be rewritten
as $H*x = x*H$. If $N$ is a normal group, then $G/H$ is again a group, the {\bf quotient group}.
For example, if $f: G \to G'$ is a group homomorphism, then the kernel of $f$ is a normal subgroup and
$|G| = |{\rm ker}(f)| |{\rm im}(f)|$ because of the {\bf first group isomorphism theorem}.
\index{partition}
\index{index}
\index{quotient group}
\index{equivalence class}

\section{Primes}

A {\bf rational prime} (or simply ``prime") is an integer larger than $1$ which is only divisible by $1$ or itself.
{\bf The Wilson theorem}
allows to define a prime as a number $n$ for which $(n-1)!+1$ is divisible by $n$. Euclid already knew that there
are infinitely many primes (if there were finitely many $p_1, \dots, p_n$, the new number $p_1 p_2 \cdots p_n+1$
would have a prime factor different from the given set). It also follows from the {\bf divergence} of the 
{\bf harmonic series} $\zeta(1) = \sum_{n=1}^{\infty} 1/n = 1+1/2+1/3 + \cdots$ and the
{\bf Euler golden key} or 
{\bf Euler product} $\zeta(s) = \sum_{n=1}^{\infty} 1/n^2 = \sum_{p \; {\rm prime}} (1-1/p^s)^{-1}$ for 
the {\bf Riemann zeta function} $\zeta(s)$ that there are infinitely many primes as otherwise, the product
to the right would be finite. 
\index{Riemann zeta function}
\index{Euler product}
\index{Euler golden key}
\index{Wilson theorem}
\index{harmonic series}
\index{divergence}
\index{prime}

Let $\pi(x)$ be the {\bf prime-counting function} which gives the number of primes
smaller or equal to $x$. Given two functions $f(x),g(x)$ from the integers to the integers,
we say $f \sim g$, if $\lim_{x \to \infty} f(x)/g(x)=1$. The {\bf prime number theorem} tells

\satz{ $\pi(x) \sim x/\log(x)$. } 

The result was investigated experimentally first by Anton Ferkel and Jurij Vega, Adrien-Marie Legendre
first conjectured in 1797 a law of this form. Carl Friedrich Gauss wrote in 1849 that he 
experimented independently around 1792 with such a law. The theorem was proven in 1896 by 
Jacques Hadamard and Charles de la Vall\'ee Poussin. Proofs without complex analysis were put forward
by Atle Selberg and Paul Erd\"os in 1949. A simple analytic proof was given by Donald Newman in 
\cite{Newman1980}.  The prime number theorem also assures that there are
infinitely many primes but it makes the statement {\bf quantitative} in that it gives an idea
how fast the number of primes grow asymptotically. Under the assumption of the Riemann hypothesis, 
Lowell Schoenfeld proved $|\pi(x)-{\rm li}(x)| < \sqrt{x} \log(x)/(8 \pi)$, 
where ${\rm li}(x)= \int_0^x dt/\log(t)$ is the {\bf logarithmic integral}.  
\index{Prime number theorem}
\index{prime counting function}
\index{natural logarithm}
\index{logarithm}
\index{Gauss}
\index{Hadamard}
\index{de la Vallee-Poussin}

\section{Cellular automata}

A finite set $A$ called {\bf alphabet} and an integer $d \geq 1$ defines the compact topological
space $\Omega=A^{\mathbb{Z}^d}$ of all infinite d-dimensional configurations.
The topology is the product topology which is compact by the Tychonov theorem.
The translation maps $T_i(x)_n=x_{n+e_i}$ are homeomorphisms of $\Omega$
called {\bf shifts}. A closed $T$ invariant subset $X \subset \Omega$ defines a {\bf subshift}    
$(X,T)$. An automorphism $T$ of $\Omega$ which commutes with the translations $T_i$ is called a   
{\bf cellular automaton}, abbreviated $CA$. An example of a cellular automaton is a map
$T x_n = \phi(x_{n+u_1}, \dots x_{n+u_k})$ where $U = \{ u_1, \dots u_k \} \subset \mathbb{Z}^d$ is a
fixed finite set. It is called an {\bf local automaton} because it is defined by a finite rule
so that the status of the cell $n$ at the next step depends only on the status of the ``neighboring cells" 
$\{ n+u  \; | \; u \in U \}$. The following result is the {\bf Curtis-Hedlund-Lyndon theorem}:
\index{Bernoulli shift}
\index{alphabet}  
\index{cellular automaton}
\index{local automaton}
\index{Curtis-Hedlund-Lyndon}

\satz{Every cellular automaton is a local automaton.}

Cellular automata were introduced by John von Neumann and mathematically
in 1969 by Hedlund \cite{Hed69}. The result appears there. Hedlund saw cellular automata also as
maps on subshifts. One can so look at cellular automata on subclasses of subshifts. For example,
one can restrict the cellular automata map $T$ on almost periodic configurations, which are subsets $X$ of $\Omega$
on which $(X,T_1,\dots T_j)$ has only invariant measures $\mu$ for which the Koopman operators 
$U_i f = f(T_i)$ on $L^2(X,\mu)$ have pure point spectrum. 
A particularly well studied case is $d=1$ and $A=\{0,1\}$, if $U=\{ -1,0,1\}$, where the automaton is called an 
{\bf elementary cellular automaton}. The {\bf Wolfram numbering} labels the $2^8$ possible 
elementary automata with a number between 1 and 255.
The {\bf game of life} of Conway is a case for $d=2$ and $A=\{ -1,0,1\} \times \{-1,0,1 \}$.
For literature on cellular automata see \cite{Wolfram86} or as part of complex systems \cite{Wolfram2002}
or evolutionary dynamics \cite{Nowak}. For topological dynamics, see \cite{DGS}.

\index{Wolfram numbering}
\index{Hedlund-Curtis-Lyndon}
\index{Game of life}

\section{Topos theory}

A {\bf category} has {\bf objects} as {\bf nodes} and {\bf morphisms} as {\bf arrows}
going from one object to an other object. There can be multiple connections and self-loops so that one can
visualize a category as a {\bf quiver}. Every object has the identity arrow $1_A$. A {\bf topos} $X$ is a {\bf Cartesian closed}
category $C$ in which {\bf finite limits} exists and which has a {\bf sub-object classifier} $\Omega$ allowing
to identify sub-objects with morphisms from $X$ to $\Omega$.
{\bf Cartesian closed} means that one can define for any pair of objects $A,B$ in $C$ the {\bf product} $A \times B$
and an {\bf equalizer} representing solutions $f=g$ to arrows $f:A \to B, G:A \to B$ as well as an 
{\bf exponential} $B^A$ representing all arrows from $A$ to $B$. An example is the topos of sets.
An example of a sub-object classifier is $\Omega=\{0,1\}$ encoding ``true or false". 
\index{category}
\index{subobject classifier}
\index{terminal object}
\index{Cartesian closed}
\index{exponential of categories}
\index{equalizer}

The {\bf slice category} $E/X$ of a category $E$ with an object $X$ in $E$ is a category, 
where the objects are the arrows from $E \to X$. An $E/X$ arrow between objects 
$f: A \to X$ and $g: B \to X$ is a map $s: A \to B$ which produces a commutative 
triangle in $E$. The composition is pasting triangles together.
The {\bf fundamental theorem of topos theory} is:

\satz{The slice category $E/X$ of a topos $E$ is a topos.}

For example, if $E$ is the topos of sets, then the slice category is the category of {\bf pointed sets}:
the objects are then sets together with a function selecting a point as a ``base point".
A morphism $f:A \to B$ defines a functor $E/B \to E/A$ which preserves exponentials
and the {\bf subobject classifier} $\Omega$.
Topos theory was motivated by geometry (Grothendieck), physics (Lawvere), topology (Tierney)
and algebra (Kan). It can be seen as a generalization and even a replacement of set theory:
the Lawvere's {\bf elementary theory of the category of sets} ETCS is seen as part of ZFC
which are less likely to be inconsistent \cite{Leinster2012}. For a short introduction \cite{Illusie},
for textbooks \cite{McLarty,Caramello}, for history of topos theory in particular, see 
\cite{McLarty1990}.

\section{Transcendentals}

A {\bf root} of an equation $f(x)=0$ with integer polynomial 
$f(x)=a_n x^n + a_{n-1} x^{n-1} + \cdots + a_0$ with $n \geq 0$ and $a_j \in \mathbb{Z}$ 
is called an {\bf algebraic number}. The set $A$ of {\bf algebraic numbers} is sub-field of the field
$\mathbb{R}$ of {\bf real numbers}. The field $A$ is the {\bf algebraic closure} of the rational numbers $\mathbb{Q}$. 
It is of number theoretic interest as it contains all {\bf algebraic number fields}, finite degree
field extensions of $\mathbb{Q}$. The complement $\RR \setminus A$ is the set of {\bf transcendental
numbers}. Transcendental numbers are necessarily irrational because every rational number $x=p/q$
is algebraic, solving $q x-p=0$. Because the set of algebraic numbers is countable and the  
real numbers are not, most numbers are transcendental. The group of all automorphisms of $A$ which fix 
$\mathbb{Q}$ is called the {\bf absolute Galois group} of $\mathbb{Q}$. 

\satz{$\pi$ and $e$ are transcendental}

\index{algebraic numbers}
\index{transcendental numbers}
\index{integer polynomial}
\index{algebraic closure}
\index{algebraic number field}
\index{absolute Galois group}

This result is due to Ferdinand von Lindemann. He proved that $e^x$ is transcendental for every 
non-zero algebraic number $x$. This immediately implies $e$ is transcendental. 
Now, if $\pi$ were algebraic, then $\pi i$ would be algebraic and 
$e^{i \pi} = -1$ would be transcendental. But $-1$ is rational. 
Lindemann's result was extended in 1885 by Karl Weierstrass to the statement telling that if 
$x_1, \dots x_n$ are linearly independent algebraic numbers, then $e^{x_1}, \dots e^{x_n}$ 
are algebraically independent.  The transcendental property of $\pi$ also proves that 
$\pi$ is irrational. This is easier to prove directly. See \cite{Hua1982}.

\section{Recurrence}

A {\bf homeomorphism} $T: X \to X$ of a compact topological space $X$ defines a {\bf topological
dynamical system} $(X,T)$. We write $T^j(x) = T(T( \dots T(x)))$ to indicate that the map $T$ is applied $j$ times. 
For any $d>0$, we get from this a set $(T_1,T_2, \dots, T_d)$ of commuting
homeomorphisms on $X$, where $T_j(x) = T^j x$.  A point $x \in X$ is called {\bf multiple recurrent} 
for $T$ if for every $d>0$, there exists a sequence $n_1<n_2<n_3< \cdots $ of integers $n_k \in \mathbb{N}$ for which 
$T_j^{n_k}x \to x$ for $k \to \infty$ and all $j=1, \dots, d$. F\"urstenberg's {\bf multiple recurrence theorem} 
states:
\index{homeomorphism}
\index{topological dynamical system}
\index{multiple recurrent} 

\satz{Every topological dynamical system is multiple recurrent.}

It is known even that the set of multiple recurrent points are Baire generic. Hillel 
F\"urstenberg proved this result in 1975. 
There is a parallel theorem for {\bf measure preserving systems}: an automorphism $T$ of a
probability space $(\Omega,\mathcal{A},{\rm P})$ is called {\bf multiple recurrent} if
there exists $A \in \mathcal{A}$ and an integer $n$ such that 
${\rm P}[ A \cap T_1(A) \cap \cdots \cap T_d(A)] >0$. 
This generalizes the {\bf Poincar\'e recurrence theorem}, which is the case $d=1$. 
Recurrence theorems are related to the {\bf Szemer\'edi theorem} telling that a
subset $A$ of $\mathbb{N}$ of positive {\bf upper density} contains arithmetic 
progressions of arbitrary finite length. See \cite{FurstenbergRecurrence}.

\index{Poincar\'e recurrence}
\index{Szemer\'edi theorem}
\index{F\"urstenberg theorem}
\index{Upper density}
\index{arithmetic progressions}

\section{Solvability}

A basic task in mathematics is to solve {\bf polynomial equations}
$p(x)=a_n x^n + a_{n-1} x^{n-1} + \cdots + a_1 x + a_0=0$
with complex coefficients $a_k$ using explicit formulas involving {\bf roots}. 
One calls this finding an {\bf explicit algebraic solution}. 
The linear case $a x + b=0$ with $x=-b/a$, the quadratic case $a x^2 + b x + c =0$ with 
$x=(-b \pm \sqrt{b^2-4ac})/(2a)$ were known since antiquity. 
The cubic $x^3+ax^2+bx+C=0$ was solved by Niccolo Tartaglia and 
Cerolamo Cardano: a first substitution $x=X-a/3$ produces the 
{\bf depressed cubic} $X^3 + p X + q$ (first solved by Scipione dal Ferro). The substitution
$X=u-p/(3u)$ then produces a quadratic equation for $u^3$. 
Lodovico Ferrari solved finally the quartic by reducing it to the cubic. 
It was Paolo Ruffini, Niels Abel and \'Evariste Galois who realized
that there are no algebraic solution formulas any more for polynomials
of degree $n \geq 5$. 
\index{Depressed cubic}
\index{quartic equation}
\index{cubic equation}

\satz{Explicit algebraic solutions to $p(x)=0$ exist if and only if $n \leq 4$. }

The quadratic case was settled over a longer period in independent 
developments in Babylonian, Egyptian, Chinese and Indian mathematics. 
The cubic and quartic discoveries were dramatic culminating with Cardano's book
of 1545, marking the beginning of modern algebra. After centuries of failures
of solving the quintic, Paolo Ruffini published the first proof in 1799, a proof
which had a gap but who paved the way for Niels Hendrik Abel and \'Evariste Galois. 
For further discoveries see \cite{MazurImagining,LivioSymmetry,DuelAtDawn}. 

\section{Galois theory}

If $F$ is sub-field of $E$, then $E$ is a vector space over $F$.
The dimension of this vector space is called the {\bf degree} $[E:F]$ of the 
{\bf field extension} $E/F$. The field extension is called {\bf finite}
if $[E:F]$ is finite. A field extension is called {\bf transcendental}
if there exists an element in $E$ which is not a root of an integral
polynomial $f$ with coefficients in $F$. Otherwise, the extension is called
{\bf algebraic}. In the later case, there exists a unique monic polynomial
$f$ which is irreducible over $F$ and the field extension is finite. 
An algebraic field extension $E/F$ is called {\bf normal}
if every irreducible polynomial over $K$ with at least
one root in $E$ {\bf splits} over $F$ into linear factors.  
An algebraic field extension $E/F$ is called {\bf separable} if the 
associated irreducible polynomial $f$ is separable, meaning that $f'$ is not
zero. This means, that $F$ has zero characteristic or that $f$ is not of 
the form $\sum_k a_k x^{p k}$ if $F$ has characteristic $p$. 
A field extension is called {\bf Galois} if it normal and separable. 
Let ${\bf Fields}(E/F)$ be the set of subfields of $E/F$
and ${\bf Groups}(E/F))$ the set of subgroups of the automorphism group ${\rm Aut}(E/F)$. 
The {\bf Fundamental theorem of Galois theory} assures:
\index{field extension}
\index{Galois extension}
\index{normal extension}
\index{monique polynomial}
\index{algebraic extension}
\index{transcendental extension}

\satz{${\bf Fields}(E/F) \overset{bijective}{\leftrightarrow} {\bf Groups}(E/F)$ if $E/F$ is Galois.}

The {\bf intermediate fields} of $E/F$ are so described by groups.
It implies the {\bf Abel-Ruffini theorem} about the non-solvability of the quintic by radicals.
The fundamental theorem demonstrates that solvable extensions correspond to solvable groups.
The {\bf symmetry groups} of permutations of $5$ or more elements are no more solvable. 
See \cite{StewartTall}.

\section{Metric spaces}

A {\bf topological space} $(X,\mathcal{O})$ is given by a set $X$
and a finite collection $\mathcal{O}$ of subsets of $X$ with the
property that the {\bf empty set} $\emptyset$ and $\Omega$ both
belong to $\mathcal{O}$ and that $\mathcal{O}$ is closed under 
arbitrary unions and finite intersections. The sets in $\mathcal{O}$ are called
{\bf open sets}. {\bf Metric spaces} $(X,d)$ are special topological spaces. In that
case, $\mathcal{O}$ consists of all sets $U$ such that for every $x \in U$
there exists $r>0$ such that the {\bf open ball} 
$B_r(x) = \{ y \in X \; | \; d(x,y)<r \}$ is contained in $U$. 
Two topological spaces $(X,\mathcal{O})$, $(Y,\mathcal{Q})$ are {\bf homeomorphic}
if there exists a bijection $f:X \to Y$, such that $f$ and $f^{-1}$ are both
continuous. A function $f:X \to Y$ is {\bf continuous} if $f^{-1}(A) \in \mathcal{O}$ for all $A \in Q$. 
When is a topological space homeomorphic to a metric space? The
{\bf Urysohn metrization theorem} gives an answer: we need the {\bf regular Hausdorff property}
meaning that a closed set $K$ and a point $x$ can be separated by disjoint neighborhoods
$K \subset U, y \in V$. We also need the space to be {\bf second countable} meaning that
there is a countable topological base (a {\bf topological base} in $\mathcal{O}$ is a subset $\mathcal{B} \subset \mathcal{O}$
such that every $U \in \mathcal{O}$ can be written as a union of elements in $\mathcal{B}$.)
\index{topological spaces}
\index{continuous map}
\index{open set}
\index{homeomorphism}
\index{second countable}
\index{regular Hausdorff}
\index{topological base}
\index{Urysohn metrization}

\satz{A second countable regular Hausdorff space is metrizable.}

The result was proven by Pavel Urysohn in 1925 with ``regular" replaced by ``normal" and by 
Andrey Tychonov in 1926. It follows that a compact Hausdorff space is metrizable if and only if it is
second countable. For literature, see \cite{Bredon}.

\section{Fixed point}

Given a continuous {\bf transformation} $T:X \to X$ of a compact topological space $X$, one can look for
the {\bf fixed point set} ${\rm Fix}_T(X) = \{ x \; | \; T(x)=x \}$. This is useful for finding {\bf periodic points}
as fixed points of $T^n=T\circ T \circ T \cdots \circ T$ are periodic points of period $n$. 
If $X$ has a finite {\bf cohomology} like if $X$ is a compact $d$-manifold with boundary, one
can look at the {\bf linear map} $T_p$ induced on the cohomology groups $H^p(X)$. The {\bf super
trace} $\chi_T(X)=\sum_{p=0}^d (-1)^p {\rm tr}(T_p)$ is called the {\bf Lefschetz number} of $T$ on $X$.
If $T$ is the identity, this is the {\bf Euler characteristic}.
Let ${\rm ind}_T(x)$ be the {\bf Brouwer degree} of the map $T$ induced on a small $(d-1)$-sphere $S$ 
centered at $x$.
This is the {\bf trace} of the linear map $T_{d-1}$ induced from $T$ on the cohomology group $H^{d-1}(S)$ 
which is an integer. If $T$ is differentiable and $dT(x)$ is invertible, the Brouwer degree is 
${\rm ind}_T(x) = {\rm sign}({\rm det}(dT))$.
Let ${\rm Fix}_T(X)$ denote the set of fixed points of $T$. The {\bf Lefschetz-Hopf fixed point theorem} is 
\index{transformation}
\index{fixed points}
\index{cohomology} 
\index{Brouwer degree}
\index{Lefschetz-Hopf fixed point theorem}

\satz{If ${\rm Fix}_T(X)$ is finite, then $\chi_T(X) = \sum_{x \in {\rm Fix}_T(X)} {\rm ind}_T(x)$. }

A special case is the {\bf Brouwer fixed point theorem}:
if $X$ is a compact convex subset of Euclidean space. In that case $\chi_T(X)=1$
and the theorem assures the existence of a fixed point. In particular, if $T: D \to D$ is a continuous
map from the disc $D=\{ x^2+y^2 \leq 1 \}$ onto itself, then $T$ has a fixed point.
This {\bf Brouwer fixed point theorem} was proved in 1910 by Jacques Hadamard and 
Luitzen Egbertus Jan Brouwer. The {\bf Schauder fixed point theorem} from 1930 generalizes the result
to convex compact subsets of Banach spaces. The Lefschetz-Hopf fixed point theorem was given in 1926.
For literature, see \cite{Dold,Border}.
\index{Schauder fixed point theorem}
\index{Brouwer fixed point theorem}
\index{Brouwer fixed point theorem}

\section{Quadratic reciprocity}

Given a prime $p$, a number $a$ is called a {\bf quadratic residue}
if there exists a number $x$ such that $x^2$ has remainder $a$ modulo $p$.
In other words, quadratic residues are the squares in the field $\mathbb{Z}_p$.
The {\bf Legendre symbol} $(a|p)$ is defined by be $0$ if $a$ is $0$ or a multiple of $p$
and $1$ if $a$ is a non-zero residue of $p$ and $-1$ if it is not.
While the integer $0$ is sometimes considered to be a quadratic residue we don't
include it as it is a special case. Also, in the multiplicative group $\mathbb{Z}_p^*$
without zero, there is a symmetry: there are the same number of quadratic
residues and non-residues. This is made more precise in the {\bf law of quadratic reciprocity}
\index{law of quadratic reciprocity}
\index{Legendre symbol}
\index{quadratic residue}
\index{quadratic non-residue}
\index{quadratic reciprocity}

\satz{ For any two odd primes $(p|q) (q|p) = (-1)^{\frac{p-1}{2} \frac{q-1}{2}}$.}

This means that $(p|q) = -(q|p)$ if and only if both $p$ and $q$ have remainder $3$
modulo $4$. The odd primes with of the form $4k+3$ are also prime in the Gaussian integers.
To remember the law, one can think of them as ``Fermions" and quadratic reciprocity 
tells they Fermions are anti-commuting. The odd primes of the form $4k+1$ factor by the 
{\bf 4-square theorem} in the Gaussian plane to $p=(a+ib) (a-ib)$ and are as a product of 
two Gaussian primes and are therefore Bosons. One can remember the rule because Bosons commute both other 
particles so that if either $p$ or $q$ or both are ``Bosonic", then 
$(p|q)=(q|p)$. The law of quadratic reciprocity was first conjectured by Euler and
Legendre and published by Carl Friedrich Gauss in his Disquisitiones Arithmeticae of 
1801. (Gauss found the first proof in 1796).  \cite{HardyWright,Hua1982}.
\index{quadratic residue}
\index{four square theorem}

\section{Quadratic map}

Every quadratic map $z \to f(z)=z^2+b z + d$ in the complex plane is conjugated by a linear transformation
to one of the quadratic family maps $T_c(z)=z^2+c$. 
The {\bf Mandelbrot set} $M = \{ c \in \mathbb{C}$, $T_c^n(0)$ stays bounded $\}$
is also called the {\bf connectedness locus} of the quadratic family because for $c \in M$, the
{\bf Julia set} $J_c = \{ z \in \mathbb{C}; T^n(z)$ stays bounded $\}$
is connected and for $c \notin M$, the Julia set $J_c$ is a {\bf Cantor set}. The fundamental
theorem for quadratic dynamical systems is:
\index{Julia set}
\index{Mandelbrot set}
\index{quadratic family}
\index{connectedness locus}
\index{Cantor set}

\satz{The Mandelbrot set is connected.}

Mandelbrot first thought after doing experiments and picturing the set using a computer and
printing it out that it was disconnected. The theorem is due to Adrien Duady and
John Hubbard in 1982. One can also look at the connectedness locus for $T(z)=z^d+c$, which
leads to {\bf Multibrot sets} or the map $z \to \overline{z}+c$, which leads to the
{\bf tricorn} or {\bf mandelbar} which is not path connected.
One does not know whether the Mandelbrot set $M$ is locally connected, nor whether it
is path connected. See \cite{MilnorNotes,Carlson,Beardon}
\label{tricorn set}
\label{mandelbar set}
\index{Multibrot set}
\index{locally connected}

\section{Differential equations}

Let us say that a differential equation $x'(t) = F(x(t))$ is {\bf integrable} if
a trajectory $x(t)$ either converges to infinity, or to an {\bf equilibrium point}
or to a {\bf limit cycle} or to a {\bf limiting torus}, where it is a 
periodic or almost periodic trajectory. We assume that $F$ has global solutions
meaning that a unique solution $x(t), t \geq 0$ solving $x'=F(x)$ exists for all times
The {\bf Poincar\'e-Bendixon} theorem is:
\index{Poincar\'e Bendixon}
\index{limit cycle} 

\satz{Any differential equation in the plane is integrable.}

This changes in dimensions $3$ and higher. The {\bf Lorenz attractor}
or the {\bf R\"ossler attractor} are examples of {\bf strange attractors},
limit sets on which the dynamics can have positive topological entropy
and is therefore no more integrable. The theorem also does not hold any more
if $\mathbb{R}^2$ is replaced by the 2-dimensional torus $\mathbb{T}^2$ because
there can be recurrent non-periodic orbits and even weak mixing situations can
occur generically in smooth situations. The proof of the Poincar\'e-Bendixon 
theorem relies on the {\bf Jordan curve theorem} which states that a simple closed
curve has an interior and exterior in $\mathbb{R}^2$. 
\cite{CoddingtonLevinson,KH}. 
\index{Roessler attractor}
\index{strange attractor}
\index{Lorenz attractor}

\section{Approximation theory}

A function $f$ on a closed interval $I=[a,b]$ is called {\bf continuous} if for every $\epsilon>0$
there exists a $\delta>0$ such that if $|x-y|<\delta$ then $|f(x)-f(y)|<\epsilon$. 
In the space $X=C(I)$ of all continuous functions, one can define a {\bf distance}
$d(f,g) = {\rm max}_{x \in I} |f(x)-g(x)|$. A subset $Y$ of $X$ is called {\bf dense} if for
every $\epsilon>0$ and every $x \in X$, there exists $y \in Y$ with $d(x,y)<\epsilon$. 
Let $P$ denote the class of {\bf polynomials} in $X$. The {\bf Weierstrass approximation theorem} 
tells that 
\index{continuous function}
\index{polynomial}
\index{Weierstrass approximation theorem}

\satz{Polynomials $P$ are dense in continuous functions $C(I)$. }

The Weierstrass theorem has been proven in 1885 by Karl Weierstrass. A constructive proof suggested 
by Sergey Bernstein in 1912 uses
{\bf Bernstein polynomials} $f_n(x) = \sum_{k=0}^n f(k/n) B_{k,n}(x)$ 
with $B_{k,n}(x) = B(n,k) x^k (1-x)^{n-k}$, where $B(n,k)$ denote the Binomial coefficients.
The result has been generalized to compact Hausdorff spaces $X$ and more general subalgebras of $C(X)$.
The {\bf Stone-Weierstrass approximation theorem} was proven by Marshall Stone in 1937 and simplified
in 1948 by Stone. In the complex, there is {\bf Runge's theorem} from 1885 approximating functions 
holomomorphic on a bounded region $G$ with rational functions uniformly on a compact subset $K$ of $G$ and 
{\bf Mergelyan's theorem} from 1951 allowing approximation uniformly on a compact subset with polynomials
if the region $G$ is simply connected. In {\bf numerical analysis} one has the task to approximate a given function space
by functions from a simpler class. Examples are approximations of smooth functions
by polynomials, trigonometric polynomials. There is also the {\bf interpolation problem}
of approximating a given data set with polynomials or piecewise polynomials like {\bf splines}
or {\bf B\'ezier curves}. See \cite{ToddConstructiveTheory,NatansonConstructiveTheory}.
\index{splines}
\index{Binomial coefficients}
\index{Bernstein polynomials}
\index{Stone-Weierstrass theorem}
\index{B\'ezier curves}
\index{Mergelyan theorem}
\index{Runge theorem}
\index{Numerical analysis}

\section{Diophantine approximation} 

An {\bf algebraic number} is a root of a polynomial $p(x)= a_n x^n + a_{n-1} x^{n-1} + 
\cdots + a_1 x + a_0$ with {\bf integer coefficients} $a_k$. 
A real number $x$ is called {\bf Diophantine} if there exists $\epsilon>0$
and a positive constant $C$ such that the {\bf Diophantine condition} $|x-p/q| > C/q^{2+\epsilon}$ 
is satisfied for all $p$, and all $q>0$. {\bf Thue-Siegel-Roth theorem} tells:
\index{Algebraic number}
\index{Thue-Siegel-Roth Theorem}
\index{Diophantine number}
\index{Diophantine condition}

\satz{Any irrational algebraic number is Diophantine.}

The {\bf Hurwitz's theorem} from 1891 assures that there are infinitely  many 
$p,q$ with $|x-p/q| < C/q^2$ for $C=1/\sqrt{5}$. This shows that the Tue-Siegel-Roth
Theorem can not be extended to $\epsilon=0$. The {\bf Hurwitz constant} $C$ is optimal. 
For any $C<1/\sqrt{5}$ one can with the {\bf golden ratio} 
$x=(1+\sqrt{5})/2$ have only finitely many $p,q$ with $|x-p/q| < C/q^2$.  
The set of {\bf Diophantine numbers} has full Lebesgue measure. 
A slightly larger set is the {\bf Brjuno set} of all numbers for which the continued
fraction {\bf convergent} $p_n/q_n$ satisfies $\sum_n \log(q_{n+1})/q_n < \infty$. A Brjuno
rotation number assures the {\bf Siegel linearization theorem} still can be 
proven. For quadratic polynomials, Jean-Christophe Yoccoz showed that linearizability 
implies the rotation number must be a Brjuno number. \cite{Carlson,Her79}
\index{Diophantine number}
\index{Siegel theorem}
\index{Brjuno number}
\index{Golden ratio}

\section{Almost periodicity}

If $\mu$ is a {\bf probability measure} of compact support on $\mathbb{R}$,
then $\hat{\mu}_n = \int e^{i n x} \; d\mu(x)$  are the {\bf Fourier 
coefficients} of $\mu$. The {\bf Riemann-Lebesgue lemma} tells that 
if $\mu$ is absolutely continuous, then $\hat{\mu}_n$ goes
to zero. The pure point part can be detected with the
following {\bf Wiener theorem}: 
\index{Riemann Lebesgue theorem}
\index{Wiener theorem}
\index{Fourier coefficients}

\satz{ $\lim_{n \to \infty} \frac{1}{n} \sum_{k=1}^n |\hat{\mu}_k|^2 = \sum_{x \in \mathbb{T}} |\mu(\{x\})|^2$. }

This looks a bit like the {\bf Poisson summation formula} 
$\sum_n f(n) = \sum_n \hat{f}(n)$, where $\hat{f}$ is the Fourier transform of $f$. 
[The later follows from $\sum_n e^{2\pi i k x} = \sum_n \delta(x-n)$, where 
$\delta(x)$ is a Dirac delta function. The Poisson formula holds if $f$ is uniformly continuous
and if both $f$ and $\hat{f}$ satisfy the growth condition $|f(x)| \leq C/|1+|x||^{1+\epsilon}$. ] 
More generally, one can read off the {\bf Hausdorff dimension} from decay rates of 
the Fourier coefficients. See \cite{Katznelson,Sternberg2019}.

\index{Poisson summation} 
\index{Dirac delta function}
\index{Hausdorff dimension}

\section{Shadowing}

Let $T$ be a {\bf diffeomorphism} on a smooth {\bf Riemannian manifold} $M$ with geodesic metric $d$.
A {\bf $T$-invariant set} is called {\bf hyperbolic} if for each $x \in K$, the 
tangent space $T_xM$ splits into a {\bf stable and unstable bundle}
$E^+_x \oplus E^-_x$ such that for some $0<\lambda<1$ and constant $C$, one has
$dT E_x^\pm = E_{Tx}^\pm$ and $|d T^{\pm n} v| \leq C \lambda^{n}$
for $v \in E^\pm$ and $n \geq 0$. An {\bf $\epsilon$-orbit} is a sequence
$x_n$ of points in $M$ such that $x_{n+1} \in B_{\epsilon}(T(x_n))$, where $B_{\epsilon}$
is the geodesic ball of radius $\epsilon$. 
Two sequences $x_n,y_n \in M$ are called {\bf $\delta$-close} 
if $d(y_n,x_n) \leq \delta$ for all $n$. 
We say that a set $K$ has the {\bf shadowing property}, if there exists an open neighborhood $U$ 
of $K$ such that for all $\delta>0$ there exists $\epsilon>0$ such that every $\epsilon$-pseudo 
orbit of $T$ in $U$ is $\delta$-close to true orbit of $T$. 
\index{Diffeomorphism}
\index{hyperbolic set}
\index{pseudo orbit}
\index{close orbits}
\index{shadowing property}

\satz{ Every hyperbolic set has the shadowing property. }

This is only interesting for infinite $K$ as if $K$ is a finite periodic hyperbolic orbit,
then the orbit itself is the orbit. It is interesting however for a hyperbolic invariant
set like a {\bf Smale horse shoe} or in the {\bf Anosov case}, which is the 
situation when the entire manifold is hyperbolic.
See \cite{KH}.
\index{horse shoe}
\index{Anosov} 
\index{Smale horse shoe}

\section{Partition function}

Let $p(n)$ denote the number of ways we can write $n$ as a sum of positive integers
without distinguishing the order. For example, $p(4)=5$ because $4=1+3=2+2=1+1+2=1+1+1+1$
can be written in 4 different ways as a sum of positive integers.
Euler used its {\bf generating function}
which is $\sum_{n=0}^\infty p(n) x^n = \prod_{k=1}^{\infty} (1-x^k)^{-1}$.
The reciprocal function $(1-x)(1-x^2)+(1-x^3) \cdots $ is called the {\bf Euler function}
and generates the {\bf generalized Pentagonal number theorem}
$\sum_{k \in \mathbb{Z}} (-1)^k x^{k(3k-1)/2} = 
1-x-x^2+x^5-x^7-x^{12} - x^{15} \cdots $ leading to the
recursion $p(n) = p(n-1) + p(n-2) - p(n-5)-p(n-7) + p(n-12) + p(n-15)  \cdots $.
The {\bf Jacobi triple product} identity is

\satz{ $\prod_{n=1}^{\infty} (1-x^{2m}) (1-x^{2m-1} y^2) (1-x^{2m-1} y^{-2}) = 
       \sum_{n=-\infty}^{\infty} x^{n^2} y^{2n}$. }

The formula was found in 1829 by Jacobi. For $x=z \sqrt{z}$ and
$y^2=-\sqrt{z}$ the identity reduces to the {\bf pentagonal number theorem} of Euler.
See \cite{AndrewsPartitions}.
\index{Pentagonal number theorem} 
\index{Generating function} 
\index{Jacobi triple product}         

\section{Burnside lemma}

If $G$ is a finite group acting on a finite set $X$, let $X/G$ denote the number of
disjoint {\bf orbits} and $X^g = \{ x \in X \; | \; g.x = x, \forall g \in G \}$ 
the {\bf fixed point set} of elements which are fixed by $g$. The number $|X/G|$ 
of orbits and the {\bf group order} $|G|$ and the size of the {\bf fixed point 
sets} are related by the {\bf Burnside lemma}: 
\index{Burnside lemma}
\index{Group}
\index{orbit}
\index{fixed point set}

\satz{$|X/G| = \frac{1}{|G|} \sum_{g \in G} |X^g|$}

The result was first proven by Frobenius in 1887. Burnside popularized it
in 1897 \cite{Burnside}. 

\section{Taylor series}

A complex-valued function $f$ which is {\bf analytic} in a disc $D=D_r(a)=\{ |x-a|<r\}$ can be
written as a series involving the $n$'th derivatives $f^{(n)}(a)$ of $f$ at $a$.
If $f$ is real valued on the real axes, the function is called {\bf real analytic}
in $(x-a,x+a)$. In several dimensions we can use multi-index notation 
$a=(a_1, \dots, a_d)$, $n=(n_1, \dots, n_d)$, $x=(x_1, \dots, x_d)$ and $x^n = x_1^{n_1} \cdots x_d^{n_d}$
and $f^{(n)}(x) = \partial_{x_1}^{n_1} \cdots \partial_{x_d}^{n_d}$ and use a {\bf polydisc}
$D=D_r(a) = \{ |x_1-a_1|<r_1, \dots |x_d-a_d|<r_d \}$.  The {\bf Taylor series formula} is:
\index{real analytic}
\index{Taylor series}

\satz{ For analytic $f$ in $D$, $f(x) = \sum_{n=0}^{\infty} \frac{f^{(n)}(a)}{n!} (x-a)^n$. }

Here, $T_r(a) = \{ |x_i-a_1|=r_1 \dots |x_d-a_d|=r_d \}$ is the boundary torus.
For example, for $f(x)=\exp(x)$, where $f^{(n)}(0)=1$, one has $f(x)=\sum_{n=0}^{\infty} x^n/n!$. 
Using the {\bf differential operator} $D f(x) = f'(x)$, one 
can see $f(x+t) = \sum_{n=0}^{\infty} \frac{f^{(n)}(x)}{n!} t^n = e^{D t} f(x)$ as a solution of
the {\bf transport equation} $f_t = D f$. 
One can also represent $f$ as a {\bf Cauchy formula} for polydiscs
$1/(2\pi i)^d \int_{|T_r(a)|} f(z)/(z-a)^d dz$ integrating along the boundary torus. 
Finite Taylor series hold in the case if $f$ is $m+1$ times differentiable. In that case one has
a finite series $S(x) = \sum_{n=0}^{m} \frac{f^{(n)}(a)}{n!} (x-a)^n$ such that  the {\bf Lagrange rest term}
is $f(x)-S(x)= R(x) = f^{m+1}(\xi) (x-a)^{m+1}/((m+1)!)$, where $\xi$ is between $x$ and $a$.
This generalizes the {\bf mean value theorem} in the case $m=0$, where $f$ is only differentiable.
The remainder term can also be written as $\int_a^x f^{(m+1)}(s) (x-a)^m/m! \; ds$. 
Brook Taylor did state but not justify the formula in 1715.
In 1742 Colin Maclaurin uses the modern form. \cite{Krantz2001}.
\index{Mean value theorem}
\index{Lagrange rest term} 
\index{polydisk}
\index{transport equation}
\index{Differential operator}

\section{Isoperimetric inequality}

Given a smooth surface $S$ in $\mathbb{R}^n$ homeomorphic to a sphere and bounding a region
$B$. Assume that the {\bf surface area} $|S|$ is fixed. How large can the {\bf  volume}
$|B|$ of $B$ become? If $B$ is the unit ball $B_1$ with
volume $|B_1|$ the answer is given by the {\bf isoperimetric inequality}:
\index{isoperimetric inequality}
\index{volume}

\satz{ $n^n |B|^{n-1} \leq |S|^n/|B_1|$.}

If $B=B_1$, this gives $n |B| \leq |S|$, which is an equality as
then the {\bf volume of the ball} $|B|=\pi^{n/2}/\Gamma(n/2+1)$
and the {\bf surface area of the sphere} $|S|=n \pi^{n/2}/\Gamma(n/2+1)$
which Archimedes first got in the case $n=3$, where $|S|=4\pi$ and $|B|=4\pi/3$. 
The classical {\bf isoperimetric problem} is $n=2$, where we are in
the plane $\mathbb{R}^2$. The inequality tells then $4 |B| \leq |S|^2/\pi$
which means $4\pi {\rm Area} \leq {\rm Length}^2$. The ball $B_1$ with area
$1$ maximizes the functional. For $n=3$, with usual Euclidean space $\mathbb{R}^3$, the inequality
tells $|B|^2 \leq (4\pi)^3/(27 \cdot 4\pi/3)$  
which is $|B| \leq 4 \pi/3$. The first proof in the case $n=2$ was attempted
by Jakob Steiner in 1838 using the {\bf Steiner symmetrization} process which is a 
refinement of the {\bf Archimedes-Cavalieri principle}.
In 1902 a proof by Hurwitz was given using Fourier series. The result has been
extended to geometric measure theory \cite{Federer}.
One can also look at the discrete problem to maximize the area defined by a polygon:
if $\{ (x_i,y_i), i =0,\dots n-1\}$ are the points of the polygon, then the area is
given by Green's formula as $A=\sum_{i=0}^{n-1} x_i y_{i+1}-x_{i+1}y_i$  and the length is
$L=\sum_{i=0}^{n-1} (x_i-x_{i+1})^2 + (y_i-y_{i+1})^2$ with
$(x_n,y_n)$ identified with $(x_0,y_0)$. The {\bf Lagrange equations} for $A$ under the
constraint $L=1$ together with a fix of $(x_0,y_0)$ and $(x_1=1/n,0)$ produces two maxima
which are both {\bf regular polygons}. A generalization to $n$-dimensional 
Riemannian manifolds is given by the L\'evi-Gromov isoperimetric inequality.
  
\index{volume of ball}
\index{surface area of sphere}
\index{isoperimetric inequality}
\index{area of polygon}
\index{length of polygon}
\index{Lagrange equations} 
\index{Cavalieri principle}

\section{Riemann Roch}

A {\bf Riemann surface} is a one-dimensional complex manifold. It is a two-dimensional real analytic
manifold but it has also a {\bf complex structure} forcing it to be orientable for example. 
Let $G$ be a compact connected {\bf Riemann surface} of Euler characteristic 
$\chi(G)=1-g$, where $g=b_1(G)$ is the {\bf genus}, the number of handles of $G$ (and $1=b_0(G)$
indicates that we have only one connected component). 
A {\bf divisor} $D=\sum_i a_i z_i$ on $G$ is an element of the
free Abelian group on the points of the surface. These are finite formal sums of points $z_i$ in $G$,
where $a_i \in \mathbb{Z}$ is the multiplicity of the point $z_i$. The {\bf degree} of the divisor
is defined as ${\rm deg}(D)=\sum_i a_i$. Let us write $\chi(D)={\rm deg}(D) + \chi(G)= {\rm deg}(D)+1-g$ and 
call this the {\bf Euler characteristic} of the divisor $D$ as one can see a divisor as a geometric
object by itself generalizing the complex manifold $X$ (which is the case $D=0$). 
A {\bf meromorphic function} $f$ on $G$ defines the {\bf principal divisor}
$(f)=\sum_i a_i z_i - \sum_j b_j w_j$, where $a_i$ are the multiplicities of the {\bf roots}
$z_i$ of $f$ and $b_j$ the multiplicities of the {\bf poles} $w_j$ of $f$. The principal divisor
of a global meromorphic 1-form $dz$ which is called the {\bf canonical divisor} $K$.
Let $l(D)$ be the dimension of the linear space of meromorphic functions $f$ on $G$ for which
$(f) + D \geq 0$. (The notation $\geq 0$ means that all coefficients are non-negative. 
One calls such a divisor {\bf effective}). The {\bf Riemann-Roch} theorem is
\index{divisor}
\index{meromorphic function}
\index{principal divisor}
\index{canonical divisor}
\index{degree of a divisor}
\index{genus}
\index{Riemann surface}
\index{effective divisor}

\satz{$l(D) - l(K-D) = \chi(D)$}

The idea of a Riemann surfaces was defined by Bernhard Riemann.
Riemann-Roch was proven for Riemann surfaces by Bernhard Riemann in 1857 and Gustav Roch in 1865.
It is possible to see this as a {\bf Euler-Poincar\'e type relation} by identifying the left
hand side as a signed cohomological Euler characteristic and the right hand side as 
a combinatorial Euler characteristic. There are various generalizations, to 
arithmetic geometry or to higher dimensions. See \cite{GriffithsHarris,Schenk}.
\index{Riemann-Roch theorem}

\section{Optimal transport}

Given two probability spaces $(X,P),(Y,Q)$
and a continuous {\bf cost function} $c: X \times Y \to [0,\infty]$,
the {\bf optimal transport problem} or
{\bf Monge-Kantorovich minimization problem} is to find the minimum
of $\int_X c(x,T(x)) \; dP(x)$ among all {\bf coupling transformations}
$T: X \to Y$ which have the property that it transports the measure $P$
to the measure $Q$. More generally, one looks at a measure $\pi$ on
$X \times Y$ such that the projection of $\pi$ onto $X$ it is $P$ and
the projection of $\pi$ onto $Y$ is $Q$. The function to optimize is then
$I(\pi) = \int_{X \times Y} c(x,y) \; d\pi(x,y)$. One of the fundamental results
is that optimal transport exists. The technical assumption is that if 
the two probability spaces $X,Y$ are {\bf Polish} (=separable complete metric spaces) and that the cost function
$c$ is continuous.
\index{Polish space}
\index{Monge-Kantorovich}
\index{coupling transformation}
\index{optimal transport}
\index{cost function}

\satz{For continuous cost functions $c$, there exists a minimum of $I$.}

In the simple set-up of probability spaces, this just follows from the 
compactness (the Alaoglu theorem for balls in the weak star topology of a Banach 
space) of the set of probability measures: any sequence $\pi_n$
of probability measures on $X \times Y$ has a convergent subsequence. Since
$I$ is continuous, picking a sequence $\pi_n$ with $I(\pi_n)$ decreasing produces 
to a minimum. The problem was formalized in 1781 by Gaspard Monge and worked on by
Leonid Kantorovich. Hirisho Tanaka in the 1970ies produced connections
with partial differential equations like the Bolzmann equation. There are
also connections to {\bf weak KAM theory} in the form of Aubry-Mather theory.
The above existence result is true under substantial less regularity.
The question of uniqueness or the existence of a Monge coupling given in the
form of a transformation $T$ is subtle \cite{VillaniTransport}.
\index{Bolzmann equation}
\index{Optimal transport problem}
\index{Aubry-Mather theory}
\index{KAM theory}
\index{Alaoglu theorem}

\section{Structure from motion}

Given $m$ hyper planes in $\mathbf{R}^d$ serving as retinas or photographic plates
for {\bf affine cameras}
and $n$ points in $\mathbf{R}^d$. The {\bf affine structure from motion} problem is to
understand under which conditions it is possible to recover both the points and
planes when knowing the orthogonal projections onto the planes. It is a model
problem for the task to reconstruct both the scene as well as the camera positions
if the scene has $n$ points and $m$ camera pictures were taken. 
Ullman's theorem is a prototype result with $n=3$ different cameras and $m=3$ points which 
are not collinear. Other setups are {\bf perspective cameras} or {\bf omni-directional cameras}.
The {\bf Ullman} map $F$ is a nonlinear map from $R^{d \cdot 2} \times SO_d^2$ to $(R^{3d-3})^2$ which
is a map between equal dimensional spaces if $d=2$ and $d=3$. The group $SO_d$ is the rotation
group in $\mathbb{R^d}$ describing the possible ways in which the affine camera can be positioned. 
Affine cameras capture the same picture when translated so that the planes can all go through the
origin. 
In the case $d=2$, we get a map from $R^4 \times SO_2^2$ to $R^6$ and
in the case $d=3$, $F$ maps $\mathbf{R}^6 \times SO_3^2$ into $\mathbf{R}^{12}$.
\index{affine camera}
\index{Ullman's theorem}
\index{Orthogonal projection}
\index{structure from motion}

\satz{The structure from motion map is locally invertible.}

In the case $d=2$, there is a reflection ambiguity. In dimension $d=3$, the number of ambiguities
is typically $64$. Ullman's theorem appeared in 1979 in \cite{Ullman}. Ullman states the theorem
for d=3 with 4 points as adding a four point cuts the number of ambiguities from 64 to 2.
See \cite{KnillHerranB} both in dimension d=2 and d=3 the Jacobean $dF$ of the Ullman map is
seen to be invertible and the inverse of $F$ is given explicitly. For structure from
motion problems in computer vision in general, see \cite{faugeras96,hartley,Trucco}. In            
applications one takes $n$ and $m$ large and reconstructs both the points as well as 
the camera parameters using {\bf statistical data fitting}.
\index{reflection ambiguity}
\index{data fitting}

\section{Poisson equation}

What functions $u$ solve the {\bf Poisson equation} $-\Delta u = f$, a 
partial differential equation? The right hand side can be
written down for $f \in L^1$ as $K_f(x) = \int_{{\mathbb R}^n} G(x,y) f(y) \; dy + h$,
where $h$ is {\bf harmonic}. If $f=0$, then the Poisson equation is the {\bf Laplace equation}.
The function $G(x,y)$ is the {\bf Green's function}, an {\bf integral kernel}. It satisfies
$-\Delta G(x,y) = \delta(y-x)$, where $\delta$ is the {\bf Dirac delta function}, a distribution.
It is given by $G(x,y) = -\log|x-y|/(2\pi)$ for $n=2$ or $G(x,y) =|x-y|^{-1}/(4\pi)$ for $n=3$.
In {\bf elliptic regularity theory}, one replaces the Laplacian
$-\Delta$ with an {\bf elliptic} second order {\bf differential operator}
$L=A(x) \cdot D \cdot D + b(x) \cdot D + V(x)$
where $D=\nabla$ is the gradient and $A$ is a positive definite matrix,
$b(x)$ is a vector field and $c$ is a scalar field.
\index{Differential operator}
\index{elliptic regularity}
\index{Green's function}
\index{Poisson equation}
\index{Laplace equation}

\satz{ For $f \in L^p$ and $p>n$, then $K_f$ is differentiable. }

The result is much more general and can be extended.
If $f$ is in $C^k$ and has compact support for example, then $K_f$ is in $C^{k+1}$.
An example of the more general set up is the {\bf Schr\"odinger operator}
$L=-\Delta + V(x)-E$. The solution to $L u = 0$, solves then an eigenvalue problem.
As one looks for solutions in $L^2$, the solution only exists if $E$ is an {\bf eigenvalue}
of $L$. The Euclidean space $\mathbb{R}^n$ can be replaced by a bounded domain $\Omega$ of
$\mathbb{R}^n$ where one can look at boundary conditions like of Dirichlet or von Neumann type.
Or one can look at the situation on a general Riemannian manifold $M$ with or without
boundary. On a Hilbert space, one has then {\bf Fredholm theory}. The equation
$u=\int G(x,y) f(y) dy$ is called a {\bf Fredholm integral equation} and
${\rm det}(1 - s G) = \exp(-\sum_n s^n {\rm tr}(G^n)/n!)$ the {\bf Fredholm determinant}
leading to the {\bf zeta function} $1/{\rm det}(1 - s G)$. See \cite{ReedSimon,LiebLoss}. 
\index{Eigenvalue}
\index{Schroedinger equation} 
\index{Fredholm theory}
\index{Schroedinger operator}

\section{Four square theorem}

{\bf Waring's problem} asked whether there exists for every $k$ an integer $g(k)$
such that every positive integer can be written as a sum of $g(k)$ powers
$x_1^k  + \cdots + x_{g(k)}^k$. Obviously $g(1)=1$. 
David Hilbert proved in 1909, that $g(k)$ is finite. This is the {\bf Hilbert-Waring theorem}.
The following {\bf theorem of Lagrange} tells that $g(2)=4$:
\index{Waring problem}
\index{Hilbert-Waring theorem}
\index{Theorem of Lagrange}

\satz{Every positive integer is a sum of four squares}.

The result needs only to be verified for prime numbers as $N(a,b,c,d)=a^2+b^2+c^2+d^2$
is a norm for {\bf quaternions} $q=(a,b,c,d)$ which has the property $N(p q) = N(p) N(q)$.
This property can be seen also as a {\bf Cauchy-Binet formula}, when writing quaternions
as complex $2 \times 2$ matrices.
The four-square theorem had been conjectured already by Diophantus, but was proven first
by Lagrange in 1770. The case $g(3)=9$ was done by Wieferich in 1912. 
It is conjectured that $g(k) = 2^k + [(3/2)^k] - 2$, where $[x]$ is the integral part
of a real number. See \cite{dicksonI,dicksonII,Hua1982}.
\index{four square theorem}

\section{Knots}

A {\bf knot} is a closed curve in $\mathbb{R}^3$, an embedding of the circle in three
dimensional Euclidean space. One also draws knots in the $3$-sphere $S^3$. 
As the {\bf knot complement} $S^3-K$ of a knot $K$ characterizes the knot up to mirror reflection, the theory of knots is
part of {\bf 3-manifold theory}. The {\bf HOMFLYPT} polynomial $P$ of a knot or {\bf link} $K$
is defined recursively using {\bf skein relations} $l P(L_+) + l^{-1} P(L^-) + m P(L_0) = 0$.
Let $K \# L$ denote the {\bf knot sum} which is a {\bf connected sum}.
Oriented knots form with this operation a commutative monoid with {\bf unknot} 
as unit. It features a unique prime factorization.
The {\bf unknot} has $P(K)=1$, the {\bf unlink} has $P(K)=0$. The {\bf trefoil knot} 
has $P(K)=2l^2-l^4+l^2 m^2$.    
\index{knot}
\index{unknot}
\index{link}
\index{knot sum}
\index{HOMFLYPT}
\index{Alexander polynomial}
\index{Jones polynonial}

\satz{ $P(K \# L) = P(K) P(L)$. }

The {\bf Alexander polynomial} was discovered in 1928 and initiated classical knot theory.
John Conway showed in the 60ies how to compute the Alexander polynomial using a recursive 
{\bf skein relations} (skein comes from French escaigne=hank of yarn). 
The Alexander polynomial allows to compute an invariant for knots by looking at the projection.
The Jones polynomial found by Vaughan Jones came in 1984. This is generalized by the
HOMFLYPT polynomial named after Jim Hoste, Adrian Ocneanu, Kenneth Millett, 
Peter J. Freyd and W.B.R. Lickorish from 1985 and J. Przytycki and P. Traczyk from 1987.  
See \cite{AdamsKnot}. 
Further invariants are {\bf Vassiliev invariants} of 1990 and {\bf Kontsevich invariants} of 1993.

\section{Hamiltonian dynamics}

Given a probability space $(M,\mathcal{A},m)$ and a smooth Lie manifold $N$
with potential function $V:N \to \mathbb{R}$, the              
{\bf Vlasov Hamiltonian differential equations} on all maps $X=(f,g): M \to T^*N$
is $f'=g,g'=\int_N \nabla V(f(x)-f(y)) \; dm(y)$.
Starting with $X_0=Id$, we get a flow $X_t$ and by push forward an evolution
$P^t=X_t^* m$ of probability measures on $N$.      
The Vlasov intro-differential equations on measures in $T^*N$ are
$\dot{P}^t(x,y) + y \cdot \nabla_x P^t(x,y) - W(x) \cdot \nabla_y P^t(x,y) = 0$
with $W(x)=\int_{M} \nabla_x V(x-x') P^t(x',y') ) \; dy' dx'$.                     
Note that while $X_t$ is an infinite dimensional {\bf ordinary differential equations}
evolving maps $M \to T^*N$, the path $P^t$ is an {\bf integro differential equation} 
describing the evolution of measures on $T^*N$.
\index{Vlasov system}
\index{Vlasov dynamics}

\satz{If $X_t$ solves the Vlasov Hamiltonian, then $P^t=X_t^* m$ solves Vlasov.}

This is a result which goes back to James Clerk Maxwell. Vlasov dynamics was 
introduced in 1938 by Anatoly Vlasov. An existence result was proven by W. Brown and Klaus Hepp in 1977.
The maps $X_t$ will stay perfectly smooth if smooth initially. However, even if $P^0$ is smooth, the
measure $P^t$ in general rather quickly develops singularities so that the partial differential 
equation has only {\bf weak solutions}. The analysis of $P$ directly would involve complicated function spaces.
The {\bf fundamental theorem of Vlasov dynamics} therefore plays the role of 
the {\bf method of characteristics} in this field. If $M$ is a finite probability space, then the 
Vlasov Hamiltonian system is the {\bf Hamiltonian $n$-body problem} on $N$. An other
example is $M=T^*N$ and where $m$ is an initial phase space measure. Now $X_t$ is a one parameter
family of diffeomorphisms $X_t: M \to T^*N$ pushing forward $m$ to a measure $P^t$ on the cotangent bundle. 
If $M$ is a circle then $X^0$ defines a closed curve on $T^*N$. In particular, 
if $\gamma(t)$ is a curve in $N$ and $X^0(t) = (\gamma(t),0)$, we have a continuum 
of particles initially at rest which evolve by interacting with a force $\nabla V$.
About interacting particle dynamics, see \cite{Spohn}.
\index{Weak solution}
\index{n-body problem}
\index{fundamental theorem of Vlasov dynamics}
\index{method of characteristics}

\section{Hypercomplexity}

A {\bf hypercomplex algebra} is a finite dimensional algebra over $\mathbb{R}$ which is {\bf unital}
and distributive. The classification of hypercomplex algebras (up to isomorphism) 
of two-dimensional hypercomplex algebras over the reals are the {\bf complex numbers} $x+iy$
with $i^2=-1$, the {\bf split complex numbers} $x+jy$ with $j^2=-1$ and the {\bf dual numbers}
(the exterior algebra) $x+\epsilon y$ with $\epsilon^2=0$.
A {\bf division algebra} over a field $F$ is an algebra over $F$ in which division is possible.
{\bf Wedderburn's little theorem} tells that a finite division algebra must be a finite field.
Only $\mathbb{C}$ is the only two dimensional {\bf division algebra} over $\mathbb{R}$.
The following theorem of Frobenius classifies the class $\mathcal{X}$ of
finite dimensional associative division algebras over $\mathbb{R}$:
\index{division algebra}
\index{exterior algebra}
\index{split algebra}
\index{dual numbers}
\index{Wedderburn little theorem}
\index{hypercomplex algebra}

\satz{ $\mathcal{X}$ consists of the algebras $\mathbb{R},\mathbb{C}$ and $\mathbb{H}$.}

Hypercomplex numbers like {\bf quaternions}, {\bf tessarines} or {\bf octonions} extend the algebra of
complex numbers. Cataloging them started with Benjamin Peirce 1872 "Linear associative algebra".
{\bf Dual numbers} were introduced in 1873 by William Clifford.
The {\bf Cayley-Dickson constructions} generates iteratively algebras of twice the dimensions:
like the complex numbers from the reals, the quaternions from the complex numbers or the octonions
from the quaternions (for octonions associativity is lost). 
The next step leads to {\bf sedenions} but the later are not even an alternative algebra any more.
The Hurwitz and Frobenius theorems limit the number in the case of real normed division algebras.
Ferdinand George Frobenius classified in 1877 the finite-dimensional associative division algebras.
Adolf Hurwitz proved in 1923 (posthumously) that unital finite dimensional real algebra endowed with 
a positive-definite quadratic form (a {\bf real normed division algebra} must be 
$\mathbb{R},\mathbb{C},\mathbb{H}$ or $\mathbb{O}$). 
These four are the only {\bf Euclidean Hurwitz algebras}.
In 1907, Joseph Wedderburn classified simple algebras 
(simple meaning that there are no non-trivial two-sided ideals and so that $ab=0$ implies $a=0$ or $b=0$).
In 1958 J. Frank Adams showed topologically
that $\mathbb{R},\mathbb{C},\mathbb{H},\mathbb{O}$ are the only finite dimensional
real division algebras. In general, division algebras have dimension $1,2,4$ or $8$ as Michel Kervaire 
and Raoul Bott and John Milnor have shown in 1958 by relating the problem 
to the {\bf parallelizability of spheres}. 
The problem of classification of division algebras over a field $F$
led Richard Brauer to the {\bf Brauer group} $BR(F)$, which Jean Pierre Serre identified it with
{\bf Galois cohomology} $H^2(K,K^*)$, where $K^*$ is the multiplicative group of $K$ seen as an algebraic group.
Each Brauer equivalence class among central simple algebras ({\bf Brauer algebras}) contains a unique
division algebra by the Artin-Wedderburn theorem. Examples: the Brauer group of an algebraically closed field
or finite field is trivial, the Brauer group of $\mathbb{R}$ is $Z_2$. Brauer groups were later defined
for commutative rings by Maurice Auslander and Oscar Goldman and by Alexander Grothendieck in 1968 for schemes.
Ofer Gabber extended the Serre result to schemes with ample line bundles. The finiteness of the Brauer
group of a proper integral scheme is open.  See \cite{Baez2002,FarbDennis}.
\index{hypercomplex numbers}
\index{tessarines}
\index{octonions}
\index{sedenions}
\index{Dual numbers} 
\index{Brauer group}

\section{Approximation}

The {\bf Kolmogorov-Arnold superposition theorem} shows that {\bf continuous functions} 
$C(\mathbb{R}^n)$ of several variables can be written as a composition of continuous 
functions of two variables:
\index{continuous functions}
\index{composition of functions}
\index{superposition theorem}

\satz{Every $f \in C(\mathbb{R}^n)$ composition of continuous functions in $C(\mathbb{R}^2)$. }

More precisely, it is now known since 1962 that there exist functions $f_{k,l}$ and a function $g$
in $C(\mathbb{R})$ such that 
$f(x_1, \dots, x_n) = \sum_{k=0}^{2n} g(f_{k,1}(x_1) + \cdots + f_{k,n} x_n)$.
As one can write finite sums using functions of two variables like $h(x,y)=x+y$  or $h(x+y,z)=x+y+z$
two variables suffice. The above form was given by by George Lorentz in 1962.
Andrei Kolmogorov reduced the problem in 1956 to functions of three variables. 
Vladimir Arnold showed then (as a student at Moscow State university) in 1957 that one can do with two variables. 
The problem came from a more specific problem in algebra, the problem of
finding roots of a polynomial $p(x)=x^n+a_1 x^{n-1} + \cdots a_n$ using radicals
and arithmetic operations in the coefficients is not possible in general for $n \geq 5$. 
Erland Samuel Bring shows in 1786 that a quintic can be reduced to $x^5+ax+1$. 
In 1836 William Rowan Hamilton showed that the sextic can be reduced to $x^6+ax^2+bx+1$
to $x^7+ax^3+bx^2+cx+1$ and the degree 8 to a 4 parameter problem $x^8+ax^4+bx^3+cx^2+dx+1$. 
Hilbert conjectured that one can not do better. They are the {\bf Hilbert's 13th problem}, the 
{\bf sextic conjecture} and {\bf octic conjecture}. In 1957, Arnold and Kolmogorov
showed that no topological obstructions exist to reduce the number of variables. Important progress
was done in 1975 by Richard Brauer.  
Some history is given in \cite{FarbWolfson}:
\index{octic conjecture}
\index{sextic conjecture}

\section{Determinants}

The {\bf determinant} of a $n \times n$ matrix $A$ is defined as the sum
$\sum_{\pi} (-1)^{{\rm sign}(\pi)} A_{1\pi(1)} \cdots A_{n \pi(n)}$, where
the sum is over all $n!$ permutations $\pi$ of $\{1,\dots, n\}$ and ${\rm sign}(\pi)$
is the {\bf signature} of the permutation $\pi$. The determinant functional satisfies the
{\bf product formula} ${\rm det}(A B) = {\rm det}(A) {\rm det}(B)$.
As the determinant is the constant coefficient of the {\bf characteristic polynomial}
$p_A(x) = {\rm det}(A-x 1)  = p_0 (-x)^n
+p_1 (-x)^{n-1} + \cdots + p_k (-x)^{n-k} + \cdots + p_n$
of $A$, one can get the coefficients of the product $F^T G$ of two $n \times m$ 
matrices $F,G$ as follows:

\satz{$p_k = \sum_{|P|=k} \det(F_P) \det(G_P)$.}

The right hand side is a sum over all minors of length $k$ including the empty one $|P|=0$,
where $\det(F_P) \det(G_P)=1$. This implies $\det(1+F^T G) = \sum_P \det(F_P) \det(G_P)$
and so $ \det(1+F^T F) = \sum_P \det^2(F_P)$. The classical Cauchy-Binet theorem
is the special case $k=m$, where $\det(F^T G) = \sum_P \det(F_P) {\rm det}(G_P)$ is a sum
over all $m \times m$ patterns if $n \geq m$. It has as even more special case the
Pythagorean consequence $\det(A^T A) = \sum_P \det(A_P^2)$. The determinant product 
formula is the even more special case when $n=m$. 
\cite{Shafarevich,CauchyBinetKnill,HoffmanWu}.

\section{Triangles}

A {\bf triangle} $T$ on a two-dimensional surface $S$ is defined by three points $A,B,C$ joined by three
geodesic paths. (It is assumed that the three geodesic paths have no self-intersections nor
other intersections besides $A,B,C$ so that $T$ is a topological disk with a piecewise
geodesic boundary).  If $\alpha,\beta,\gamma$ are the 
{\bf inner angles} of a {\bf triangle} $T$ located on a surface
with {\bf curvature} $K$, there is the Gauss-Bonnet formula 
$\int_S K(x) dA(x) = \chi(S)$, where $dA$ denotes the {\bf area element} 
on the surface. This implies a relation between the integral of the curvature
over the triangle and the angles: 
\index{area element}
\index{Guss-Bonnet formula} 

\satz{ $\alpha+\beta+\gamma = \int_T K \; dA + \pi$ }

This can be seen as a special Gauss-Bonnet result for {\bf Riemannian manifolds with boundary}
as it is equivalent to $\int_T K \; dA +\alpha'+\beta +\gamma'=2\pi$
with {\bf complementary angles}
$\alpha'=\pi-\alpha, \beta'=\pi-\beta,\gamma' = \pi - \gamma$. One can think of
the vertex contributions as {\bf boundary curvatures} (generalized function). 
In the case of {\bf constant curvature} $K$, 
the formula becomes $\alpha + \beta + \gamma = K A + \pi$, where $A$ is the {\bf area of the triangle}.
Since antiquity, one knows the flat case $K=0$, where $\pi=\alpha+\beta+\gamma$ taught in elementary school.
On the {\bf unit sphere} this is $\alpha+ \beta + \gamma = A + \pi$, result of
Albert Girard which was predated by Thomas Harriot.
In the {\bf Poincar\'e disk model} $K=-1$, this is $\alpha+ \beta + \gamma = -A + \pi$ which is usually stated that
the area of a triangle in the disk is $\pi - \alpha-\beta-\gamma$. This was proven by Johann Heinrich Lambert.
See \cite{Brummelen} for spherical geometry and 
\cite{AndersonHyperbolic} for hyperbolic geometry, which are both part of 
{\bf non-Euclidean geometry} and now part of {\bf Riemannian geometry}.  \cite{BergerGostiaux,Jost}
\index{Klein model}
\index{sphere} 
\index{complementary angles}
\index{area triangle}
\index{spherical geometry}
\index{hyperbolic geometry}
\index{non-Euclidean geometry} 
\index{Riemannian geometry}
\index{Poincar\'e disk}

\section{KAM}

An {\bf area preserving map} $T(x,y) = (2x-y + cf(x),x)$ has an orbit $(x_{n+1},x_n)$
on $\mathbb{T}^2=(\mathbb{R}/\mathbb{Z})^2$ which satisfies the recursion $x_{n+1}-2x_n+x_{n-1} = cf(x_n)$. 
The $1$-periodic function $f$ is assumed to be real-analytic, non-constant satisfying $\int_0^1 f(x)\;dx=0$. 
In the case $f(x)=\sin(2\pi x)$, one has the {\bf Standard map}. 
When looking for invariant curves $(q(t+\alpha),q(t))$ with smooth $q$, we seek
a solution of the nonlinear equation $F(q) = q(t+\alpha)-2 q(t)+ q(t-\alpha) - cf(q(t))=0$.
For $c=0$, there is the solution $q(t)=t$. The {\bf linearization}
$dF(q)(u) = L u = u(t+\alpha)-2 u(t)+u(t-\alpha) -c f'(q(t)) u(t)$ is a bounded linear operator 
on $L^2(\mathbb{T})$ but not invertible for $c=0$ so that the {\bf implicit function theorem} does 
not apply. The map $L u = u(t+\alpha) -2 u(t) + u(t-\alpha)$ becomes after a Fourier transform
the diagonal matrix $\hat{L} \hat{u}_n = [2 \cos(n \alpha) - 2] \hat{u}_n$ which
has the inverse diagonal entries $[2 \cos(n \alpha) - n]^{-1}$ leading to {\bf small
divisors}. A real number $\alpha$ is called {\bf Diophantine} if there exists a
constant $C$ such that for all integers $p,q$ with $q \neq 0$, we have
$|\alpha - p/q| \geq C/q^2$. {\bf KAM theory} assures that the solution $q(t)=t$ persists
and remains smooth if $c$ is small. With {\bf solution} the theorem means a {\bf smooth solution}. 
For real analytic $F$, it can be real analytic. The following result is a special case of the
{\bf twist map theorem}.
\index{Standard map}
\index{small divisors}
\index{Diophantine number}
\index{implicit function theorem}
\index{KAM theory}

\satz{For Diophantine $\alpha$, there is a solution of $F(q)=0$ for small $|c|$.}

The KAM theorem was predated by the {\bf Poincar\'e-Siegel theorem} in complex dynamics
which assured that if $f$ is analytic near $z=0$ and $f'(0)=\lambda = \exp(2 \pi i \alpha)$
with Diophantine $\alpha$, then there exists $u(z) = z+q(z)$ such that
$f(u(z)) = u(\lambda z)$ holds in a small disk $0$: there is an analytic solution
$q$ to the {\bf Schr\"oder equation} $\lambda z + g(z+q(z)) = q(\lambda z)$. 
\index{Poincar\'e-Sigel theorem}
\index{Siegel linearization theorem}
\index{Schr\"oder equation}
The question about the existence of invariant curves is important as it determines
the {\bf stability}. The twist map theorem result follows also from a {\bf strong implicit function theorem}
initiated by John Nash and J\"urgen Moser. 
For larger $c$, or non-Diophantine $\alpha$, the solution $q$ still exists but
it is no more continuous. This is {\bf Aubry-Mather theory}.
For $c \neq 0$, the operator $\hat{L}$ is an almost periodic 
{\bf Toeplitz matrix} on $l^2(\mathbb{Z})$ which is a special kind of 
{\bf discrete Schr\"odinger operator}. The decay rate of the off diagonals depends on the
smoothness of $f$. Getting control of the inverse can be technical \cite{Bourgain2005}.
Even in the {\bf Standard map} case $f(x)=\sin(x)$, the composition $f(q(t))$
is no more a trigonometric polynomial so that $\hat{L}$ appearing here is not a {\bf Jacobi matrix
in a strip}. The first breakthrough of the theorem in a frame work of Hamiltonian differential
equations was done in 1954 by Andrey Kolmogorov. J\"urgen Moser proved the discrete twist map version and
Vladimir Arnold in 1963 proved the theorem for Hamiltonian systems. The above stated result generalizes 
to higher dimensions where one looks for {\bf invariant tori} called {\bf KAM tori}. 
one needs some non-degeneracy conditions
See \cite{Carlson,MoserStableRandom,MoserVariations}. For the story on KAM, see \cite{Dumas2014}.
\index{KAM}
\index{Aubry-Mather theory}
\index{Schroeder equation}
\index{Schroedinger operator}
\index{Jacobi matrix}
\index{Standard map}
\index{Toeplitz matrix}
\index{stability}
\index{strong implicit function theorem}

\section{Continued Fraction}

Given a positive {\bf square free} integer $d$, the {\bf Diophantine equation} $x^2-d y^2=1$
is called {\bf Pell's equation}. Solving it means to find a nontrivial unit in the ring 
$\mathbb{Z}[\sqrt{d}]$ because $(x+y \sqrt{d})(x-y \sqrt{d})=1$. 
The trivial solutions are $x=\pm 1, y=0$.  
Solving the equation is therefore part of the {\bf Dirichlet unit problem} from algebraic number theory.
Let $[a_0;a_1,\dots]$ denote the {\bf continued fraction expansion} of $x=\sqrt{d}$. 
This means $a_0=[x]$ is the integer part and $[1/(x-a_0)]=a_1$ etc. 
If $x=[a_0; a_1,\dots, a_n+b_n]$, then $a_{n+1}=[1/b_n]$.  Let $p_n/q_n = [a_0; a_1,a_2, \dots, a_n]$ denote
the $n$'th {\bf convergent} to the regular continued fraction of $\sqrt{d}$. A solution $(x_1,y_1)$ 
which minimizes $x$ is called the {\bf fundamental solution}. 
The theorem tells that it is of the form $(p_n,q_n)$:
\index{Pell's equation}
\index{continued faction expansion}
\index{Dirichlet unit problem}

\satz{Any solution to the Pell's equation is a convergent $p_n/q_n$.} 

One can find more solutions recursively because the ring of units in                                      
$\mathbb{Z}[\sqrt{d}]$ is $\mathbb{Z}_2 \times C_n$ for some cyclic group $C_n$. 
The other solutions $(x_k,y_k)$ can be obtained from $x_k + \sqrt{d} y_k = (x_1+\sqrt{d} y_1)^k$.
One of the first instances, where the equation appeared is in the {\bf Archimedes cattle problem}
which is $x^2-410286423278424y^2 = 1$.
The equation is named after John Pell, who has nothing to do with the equation. It was Euler who
attributed the solution by mistake to Pell. It was first found by William Brouncker. 
The approach through continued fractions started with Euler and Lagrange.
See \cite{Riesel,BressoudPrimality,Lenstra2008}.

\section{Gauss-Bonnet-Chern}

Let $(M,g)$ be a {\bf Riemannian manifold} of dimension $d$ with {\bf volume element $d\mu$}.
If $R^{ij}_{kl}$ is {\bf Riemann curvature tensor} with respect to the metric $g$, define
the constant $C=((4\pi)^{d/2} (-2)^{d/2} (d/2)!)^{-1}$ and 
the {\bf curvature} $K(x) = C \sum_{\sigma,\pi} {\rm sign}(\sigma) {\rm sign}(\pi) 
R^{\sigma(1) \sigma(2)}_{\pi(1) \pi(2)} \cdots R^{\sigma(d-1) \sigma(d)}_{\pi(d-1) \pi(d)}$,
where the sum is over all permutations $\pi,\sigma$ of $\{1, \dots, d\}$.
It can be interpreted as a {\bf Pfaffian}. In odd dimensions, the curvature is zero.
Denote by $\chi(M)$ the {\bf Euler characteristic} of $M$.
\index{Pfaffian}
\index{Riemann curvature tensor}
\index{curvature}

\satz{   $\int_M K(x) \; d\mu(x)  = 2\pi \chi(M)$. }

The case $d=2$ was solved by Carl Friedrich Gauss and by Pierre Ossian Bonnet in 1848.
Gauss knew the theorem but never published it. In the case $d=2$, the curvature $K$
is the {\bf Gaussian curvature } which is the product of the {\bf principal curvatures} 
$\kappa_1,\kappa_2$ at a point. For a sphere of radius $R$ for example, 
the Gauss curvature is $1/R^2$ and $\chi(M)=2$. 
The {\bf volume form} is then the usual {\bf area element} normalized so that 
$\int_M 1 \; d\mu(x)=1$.
Allendoerfer-Weil in 1943 gave the first proof, based on previous work of
Allendoerfer, Fenchel and Weil. Chern finally, in 1944 proved the theorem
independent of an embedding.
\cite{Cycon} features a proof of Vijay Kumar Patodi. A more classical approach is in
in \cite{TuDifferentialGeometry}.
\index{Gaussian curvature}
\index{principal curvature}
\index{Gauss-Bonnet-Chern}

\section{Atiyah-Singer}

Assume $M$ is a compact orientable finite dimensional {\bf manifold} of dimension $n$ and 
assume $D$ is an {\bf elliptic} {\bf differential operator} $D:E \to F$ between two smooth 
{\bf vector bundles} $E,F$ over $M$. Using multi-index notation 
$D^k=\partial_{x_1}^{k_1} \cdots \partial_{x_n}^{k_n}$,
a {\bf differential operator} $\sum_{k} a_k(x) D^k x$ is called {\bf elliptic} if for all $x$,      
its {\bf symbol} the polynomial $\sigma(D)(y)=\sum_{|k|=n} a_k(x) y^k$ is not zero for nonzero $y$. 
{\bf Elliptic regularity} assures that both the kernel of $D$ and the kernel of the
{\bf adjoint} $D^*: F \to E$ are both finite dimensional. The {\bf analytical index} of $D$ is
defined as $\chi(D) = {\rm dim}({\rm ker}(D)) - {\rm dim}( {\rm ker}(D^*))$. We think of it as
the Euler characteristic of $D$. The {\bf topological index} of $D$ is
defined as the integral of the $n$-form $K_D=(-1)^n {\rm ch}(\sigma(D)) \cdot {\rm td}(TM)$, over $M$. 
This $n$-form is the cup product $\cdot$ of the {\bf Chern character} ${\rm ch}(\sigma(D))$ and the
{\bf Todd class} of the complexified tangent bundle $TM$ of $M$. We think about $K_D$ as a {\bf curvature}.
Integration is done over the {\bf fundamental class} $[M]$ of $M$ which is the natural {\bf volume form} on $M$.          
The Chern character and the Todd classes are both mixed rational cohomology classes.
On a complex vector bundle $E$ they are both given by 
concrete power series of {\bf Chern classes} $c_k(E)$ like 
${\rm ch}(E)=e^{a_1(E)} + \cdots + e^{a_n(E)}$ and
${\rm td}(E)=a_1 (1+e^{-a_1})^{-1}  \cdots a_n (1+e^{-a_n})^{-1}$ with $a_i=c_1(L_i)$ if       
$E=L_1 \oplus \cdots \oplus L_n$ is a direct sum of {\bf line bundles}.
\index{differential operator}
\index{symbol}
\index{elliptic regularity}
\index{differential operator}
\index{Todd class}
\index{Chern character}
\index{line bundle}
\index{Chern classes} 
\index{Kronecker pairing}
\index{fundamental class}
\index{Analytical index}
\index{Topological index}

\satz{The analytic index and topological indices agree: $\chi(D)=\int_M K_D$.}

In the case when $D=d+d^*$ from the vector bundle of even forms $E$ to the vector bundle
of odd forms $F$, then $K_D$ is the Gauss-Bonnet curvature and $\chi(D)=\chi(M)$.
Israil Gelfand conjectured around 1960 that the analytical index should have a 
topological description. The Atiyah-Singer index theorem has been proven in 1963 by 
Michael Atiyah and Isadore Singer. The result generalizes the Gauss-Bonnet-Chern
and Riemann-Roch-Hirzebruch theorem. According to \cite{Rognes2004}, 
``the theorem is valuable, because it connects 
analysis and topology in a beautiful and insightful way". See \cite{PalaisSeminar}. 

\section{Complex multiplication}

A {\bf n'th root of unity} is a solution to the equation $z^n=1$ in the complex plane $\mathbb{C}$.
It is called {\bf primitive} if it is not a solution to $z^k=1$ for some $1 \leq k<n$.
A {\bf cyclotomic field} is a number field $\mathbb{Q}(\zeta_n)$ which is obtained
by adjoining a complex {\bf primitive root of unity} $\zeta_n$ to $\mathbb{Q}$.
Every cyclotomic field is an Abelian field extension of the field of rational numbers $\mathbb{Q}$.
The {\bf Kronecker-Weber} theorem reverses this. It is also called
the main theorem of {\bf class field theory over $\mathbb{Q}$}
\index{Kronecker-Weber theorem}
\index{cyclotomic field}
\index{root of unity}  
\index{primitive root of unity}

\satz{Every Abelian extension $L/\mathbb{Q}$ is a subfield of a cyclotomic field.}

Abelian field extensions of $\mathbb{Q}$ are also called {\bf class fields}.
It follows that any {\bf algebraic number field} $K/Q$ with Abelian {\bf Galois group}
has a {\bf conductor}, the smallest $n$ such that $K$ lies in the field generated by $n$'th
roots of unity. Extending this theorem to other base number fields is {\bf Kronecker's Jugendtraum}
or {\bf Hilbert's twelfth problem}. The theory of {\bf complex multiplication} does the
generalization for {\bf imaginary quadratic fields}.
The theorem was stated by Leopold Kronecker in 1853 and proven by Heinrich Martin
Weber in 1886. A generalization to {\bf local fields} was done by Jonathan Lubin and
John Tate in 1965 and 1966. (A {\bf local field} is a locally compact topological field with respect
to some non-discrete topology. The list of local fields is $\mathbb{R},\mathbb{C}$, field extensions of the
{\bf p-adic numbers} $\mathbb{Q}_p$, or formal Laurent series $F_q((t))$ over a finite field $F_q$.)
The study of {\bf cyclotomic fields} came from elementary geometric problems like the construction of a regular
$n$-gon with {\bf ruler and compass}. Gauss constructed a regular 17-gon and showed that a {\bf regular $n$-gon}
can be constructed if and only if $n$ is a {\bf Fermat prime} $F_n=2^{2^n}+1$ 
(the known ones are $3,7,17,257,65537$ and a problem of Eisenstein of 1844 asks whether 
there are infinitely many). 
Further interest came in the context of {\bf Fermat's last theorem} because $x^n+y^n=z^n$ can be written as 
$x^n+y^n=(x+y) (x+\zeta y) \cdots (x+\zeta^{n-1} y)$, where $\zeta$ is an $n$'th root of unity
for $n>2$. 
\index{conductor}
\index{class fields}
\index{local field} 
\index{regular n-gon}
\index{Fermat prime}
\index{imaginary quadratic field}
\index{Hilbert's 12th problem}
\index{complex multiplication}
\index{Fermat's last theorem}
\index{Cyclotomic field} 
\index{p-adic field}
\index{17-gon}

\section{Choquet theory}

Let $K$ be a {\bf compact} and {\bf convex} set in a Banach space $X$. A point $x \in K$
is called {\bf extreme} if $x$ is not in an open interval $(a,b)$ with $a,b \in K$.
Let $E$ be the set of extreme points in $K$.
The {\bf Krein-Milman theorem}, proven in 1940 by Mark Krein and David Milman,
assures that $K$ is the convex hull of $E$.
 Given a probability measure $\mu$ on $E$, it defines the point $x=\int y d\mu(y)$.
We say that $x$ is the {\bf Barycenter} of $\mu$. The {\bf Choquet theorem} is
\index{Barycenter}
\index{extreme point}
\index{convex}
\index{Banach space}

\satz{Every point in $K$ is a Barycenter of its extreme points.}

This result of Choquet implies the Krein-Milman theorem. It generalizes to {\bf locally compact}
{\bf topological spaces}. The measure $\mu$ is not unique in general. It is in finite
dimensions if $K$ is a simplex. But in general, as shown by Heinz Bauer in 1961,
for an extreme point $x \in K$ the measure $\mu_x$ is unique.
It has been proven by {\bf Gustave Choquet} in 1956 and
was generalized by Erret Bishop and Karl de Leeuw in 1959. 
\cite{Phelps}
\index{locally compact topological space}

\section{Helly's theorem}

Given a family $\mathcal{K}=\{K_1, \dots K_n\}$ of {\bf convex} sets 
$K_1,K_2, \dots, K_n$ in the {\bf Euclidean space} $\mathbb{R}^d$ and 
assume that $n>d$. Let $\mathcal{K}_{m}$ denote the set of subsets of
$\mathcal{K}$ which have exactly $m$ elements. We say that $\mathcal{K}_m$ has the
{\bf intersection property} if every of its elements has a non-empty 
common intersection. The {\bf theorem of Helly} assures that 
\index{Theorem of Helly}

\satz{$\mathcal{K}_n$ has the intersection property if $\mathcal{K}_{d+1}$ has. } 

The theorem was proven in 1913 by Eduard Helly. It generalizes to an infinite
collection of compact, convex subsets. This theorem led Johann Radon to prove
in 1921 the {\bf Radon theorem} which states that any set of $d+2$ points in 
$\mathbb{R}^d$ can be partitioned into two disjoint subsets whose convex hull 
intersect. A nice application of Radon's theorem is the {\bf Borsuk-Ulam theorem}
which states that a continuous function $f$ from the $d$-dimensional sphere
$S^n$ to $\mathbb{R}^d$ must some pair of {\bf antipodal points} to the same point:
$f(x)=f(-x)$ has a solution. For example, if $d=2$, this implies that
on earth, there are at every moment two antipodal points on the Earth's surface
for which the temperature and the pressure are the same. The {\bf Borsuk-Ulam}
theorem appears first have been stated in work of 
Lazar Lyusternik and Lev Shnirelman in 1930, and 
proven by Karol Borsuk in 1933 who attributed it to Stanislav Ulam. 
\index{Radon theorem} 
\index{antipodal point}
\index{Borsuk-Ulam theorem} 

\section{Weak Mixing}

An {\bf automorphism} $T$ of a probability space $(X,\mathcal{A},m)$ is a measure preserving
invertible measurable transformation from $X$ to $X$. It is called {\bf ergodic}
if $T(A)=A$ implies $m(A)=0$ or $m(A)=1$. It is called {\bf mixing} if
$m(T^n(A) \cap B) \to m(A) \cdot m(B)$ for $n \to \infty$ for all $A,B$.
It is called {\bf weakly mixing} if $n^{-1} \sum_{k=0}^{n-1} |m(T^k(A) \cap B) - m(A) \cdot m(B)| \to 0$
for all $A,B \in \mathcal{A}$ and $n \to \infty$. This is equivalent to the fact that the unitary operator $Uf=f(T)$ on $L^2(X)$ has
no point spectrum when restricted to the orthogonal complement of the constant functions. A topological transformation 
(a continuous map on a locally compact topological space) with a
weakly mixing invariant measure is {\bf not integrable} as for integrability, one wants every invariant
measure to lead to an operator $U$ with pure point spectrum and conjugating it so to a group translation.
Let $\mathcal{G}$ be the complete topological group of automorphisms of
$(X,\mathcal{A},m)$ with the weak topology: $T_j$ converges to $T$ {\bf weakly}, if
$m(T_j(A) \Delta T(A)) \rightarrow 0$ for all $A \in \mathcal{A}$;
this topology is metrizable and completeness is defined with respect to an
equivalent metric.
\index{mixing}
\index{weakly mixing}
\index{ergodic}
\index{metrizable}
\index{Weak convergence}
\index{spectral integrability}
\index{integrable}

\satz{A generic $T$ is weakly mixing and so ergodic.}

Anatol Katok and Anatolii Mikhailovich Stepin in 1967 \cite{KaSt67} proved that 
purely singular continuous spectrum of $U$ is generic. A new proof was given by \cite{ChNa90} and a short proof in
using {\bf Rokhlin's lemma}, Halmos conjugacy lemma and a Simon's {\bf ``wonderland theorem"}
establishes both genericity of weak mixing and genericity of singular spectrum.
On the topological side, a generic volume preserving homeomorphism of a manifold
has purely singular continuous spectrum which
strengthens Oxtoby-Ulam's theorem \cite{Oxtoby} about generic ergodicity. \cite{KaSt70,Halmos}
The Wonderland theorem of Simon \cite{Sim95} also allowed to prove that a generic invariant measure
of a shift is singular continuous \cite{Kni97} or that zero-dimensional singular continuous
spectrum is generic for open sets of flows on the torus allowing also to show that open sets
of Hamiltonian systems contain generic subset with both quasi-periodic as well as weakly 
mixing invariant tori \cite{KnillTori}.
\index{invariant tori}
\index{Wonderland theorem}
\index{Rokhlin lemma}
\index{Quasiperiodic}

\section{Universality}

The space $X$ of {\bf unimodular maps} is the set of 
twice continuously differentiable even maps $f: [-1,1] \to [-1,1]$
satisfying $f(0)=1$ $f''(x)<0$ and $\lambda=f(1)<0$. The 
{\bf Feigenbaum-Cvitanovi\'c functional equation} (FCE) is $g=Tg$
with $T(g)(x) = \frac{1}{\lambda} g(g(\lambda x))$. The map $T$ is a
{\bf renormalization map.} 
\index{Unimodular map}
\index{Feigenbaum-Civtanovic functional equation}

\satz{There exists an analytic hyperbolic fixed point of $T$.}

The first proof was given by Oscar Lanford III in 1982 (computer
assisted). See \cite{Lanford82,Lanford84}. 
That proof also established that the fixed point is 
hyperbolic with a one-dimensional unstable manifold and positive
expanding eigenvalue. This explains
some {\bf universal features} of unimodular maps found 
experimentally in 1978 by Mitchell Feigenbaum and which is now
called {\bf Feigenbaum universality}. 
The result has been ported to area preserving maps \cite{EKW}.
\index{Feigenbaum universality}
\index{Universality}
\index{computer assisted proof}

\section{Compactness}

Let $X$ be a compact metric space $(X,d)$. The Banach space $C(X)$ of real-valued continuous
functions is equipped with the supremum norm. A closed subset $F \subset C(X)$ is called
{\bf uniformly bounded} if for every $x$ the supremum of all values $f(x)$ with $f \in F$
is bounded. The set $F$ is called {\bf equicontinuous} if for every $x$ and every $\epsilon>0$
there exists $\delta>0$ such that if $d(x,y) < \delta$, then 
$|f(x)-f(y)|<\epsilon$ for all $f \in F$.
A set $F$ is called {\bf precompact} if its closure is compact. 
The {\bf Arzel\`a-Ascoli theorem} is:
\index{Arzela-Ascoli}
\index{equicontinuous}
\index{supremum norm}
\index{uniformly bounded}

\satz{Equicontinuous uniformly bounded sets in $C(X)$ are precompact.}          

The result also holds on {\bf Hausdorff spaces} and not only metric spaces. 
In the complex, there is a variant called {\bf Montel's theorem} which is the
fundamental normality test for holomorphic functions: 
an uniformly bounded family of holomorphic functions on a complex domain $G$
is {\bf normal} meaning that its closure is compact with respect to the 
{\bf compact-open topology}. The compact-open topology in $C(X,Y)$ is the topology
defined by the {\bf sub-base} of all continuous maps $f_{K,U}: f:K \to U$, where $K$ runs
over all compact subsets of $X$ and $U$ runs over all open subsets of $Y$. 
\index{holomorphic}
\index{Hausdorff space}
\index{Montel theorem}
\index{normal family}
\index{compact-open topology}

\section{Geodesic}

The {\bf geodesic distance} $d(x,y)$ between two points $x,y$ on a
{\bf Riemannian manifold} $(M,g)$ is defined as the length
of the shortest geodesic $\gamma$ connecting $x$ with $y$. This
renders the manifold a metric space $(M,d)$. We assume it is
{\bf locally compact}, meaning that every point $x \in M$ has a compact
neighborhood. A metric space is called
{\bf complete} if every {\bf Cauchy sequence} in $M$ has a convergent subsequence.
(A sequence $x_k$ is called a Cauchy sequence if for every $\epsilon>0$, there
exists $n$ such that for all $i,j>n$ one has $d(x_i,x_j)<\epsilon$.)
The local existence of differential equations assures
that the geodesic equations exist for small enough time. This can be
restated that the {\bf exponential map} $v \in T_xM \to M$ assigning to a
point $v \neq 0$ in the tangent space $T_xM$ the solution $\gamma(t)$
with initial velocity $v/|v|$ and $t \leq |v|$, and $\gamma(0)=x$.
A Riemannian manifold $M$ is called {\bf geodesically complete} if the
exponential map can be extended to the entire tangent space $T_xM$ for
every $x \in M$. This means that geodesics can be continued for all times.
The Hopf-Rinow theorem assures:
\index{geodesic distance}
\index{Hopf-Rinov theorem}
\index{geodesic complete}
\index{complete}

\satz{Completeness and geodesic completeness are equivalent.}

The theorem was named after Heinz Hopf and his student Willi Rinow
who published it in 1931. See \cite{HopfRinow,DoCarmo}. 

\section{Crystallography}

A {\bf wall paper group} is a discrete subgroup of the {\bf Euclidean symmetry group} $E_2$ 
of the plane. Wall paper groups classify two-dimensional patterns according to their symmetry. 
In the plane $\mathbb{R}^2$, the underlying group is the group $E_2$ of {\bf Euclidean plane symmetries} 
which contain {\bf translations} {\bf rotations} or {\bf reflections} or {\bf glide reflections}. This
group is the group of rigid motions. It is a three dimensional {\bf Lie group} which according
to Klein's {\bf Erlangen program} characterizes {\bf Euclidean geometry}. Every element in $E_2$ can
be given as a pair $(A,b)$, where $A$ is an orthogonal matrix and $b$ is a vector. 
A subgroup $G$ of $E_2$ is called {\bf discrete} if there is a positive minimal distance between two
elements of the group. This implies the {\bf crystallographic restriction theorem} assuring 
that only rotations of order $2,3,4$ or $6$ can appear. 
This means only rotations by $180,120,90$ or $60$ degrees can occur in a Wall paper group. 
\index{wall paper group}
\index{crystallogrphy}
\index{translation}
\index{rotation}
\index{reflection} 

\satz{There are 17 wallpaper groups}

The first proof was given by Evgraf Fedorov in 1891 and then by George Polya in 1924. 
in three dimensions there are 230 {\bf space groups} and 219 types if {\bf chiral copies}
are identified. In space there are 65 space groups which preserve the orientation.
See \cite{OwenOrnament,BruenbaumShephard,symmetries}.
\index{space group}
\index{chiral copies}

\section{Quadratic forms}

A symmetric square matrix $Q$ of size $n \times n$ with integer entries
defines a {\bf integer quadratic form} $Q(x) = \sum_{i,j=1}^n Q_{ij} x_i x_j$.
It is called {\bf positive} if $Q(x)>0$ whenever $x \neq 0$.
A positive integral quadratic form is called {\bf universal} if its range is
$\mathbb{N}$. For example, by the {\bf Lagrange four square theorem}, the form
$Q(x_1,x_2,x_3,x_4)=x_1^2+x_2^2+x_3^2+x_4^2$ is universal.
The {\bf Conway-Schneeberger fifteen theorem} tells
\index{integer quadratic form}
\index{Conway-Schneeberger fifteen theorem}
\index{fifteen theorem}
\index{Lagrange four square theorem}
\index{universal}

\satz{$Q$ is universal if it has $\{1,\dots 15\}$ in the range.}

The interest in quadratic forms started in the 17'th century especially about
numbers which can be represented as sums $x^2+y^2$. Lagrange, in 1770 proved the
four square theorem. In 1916, Ramajujan listed all diagonal quaternary forms which
are universal. The 15 theorem was proven in 1993 by John Conway and William Schneeberger
(a student of Conway's in a graduate course given in 1993).
There is an analogue theorem for {\bf integral positive quadratic forms}, these
are defined by positive definite matrices $Q$ which take only integer values.
The binary quadratic form $x^2+x y + y^2$ for example is integral but not
an integer quadratic form because the corresponding matrix $Q$ has fractions $1/2$.
In 2005, Manjul Bhargava and Jonathan Hanke proved the 290 theorem, assuring that
an integral positive quadratic form is universal if it contains $\{1, \dots, 290\}$
in its range.  \cite{Conway15}.

\index{universal}
\index{quadratic form}
\index{integral quadratic form}
\index{square matrix}
\index{positive definition}
\index{Lagrange four square theorem}
\index{fifteen theorem}
\index{15 theorem}
\index{290 theorem}

\section{Sphere packing}

A {\bf sphere packing} in $\mathbb{R}^d$ is an arrangement of non-overlapping 
unit spheres in the $d$-dimensional Euclidean space $\mathbb{R}^d$ with volume measure $\mu$.
It is known since \cite{Gro63} that packings with maximal densities exist.
Denote by $B_r(x)$ the ball of radius $r$ centered at $x \in \mathbb{R}^d$.
If $X$ is the set of centers of the sphere and $P=\bigcup_{x \in X} B_1(x)$ is the
union of the unit balls centered at points in $X$, then the {\bf density} of the
packing is defined as
$\Delta_d = \limsup \int_{B_r(0)} P \; d\mu/\int_{B_r(0)} 1 \; d\mu$.
The sphere packing problem is now solved in 5 different cases:
\index{sphere packing}
\index{density of sphere packing}
\index{optimal sphere packing}

\satz{Optimal sphere packings are known for $d=1,2,3,8,24$.}

The one-dimensional case $\Delta_1=1$ is trivial. The case $\Delta_2 = \pi/\sqrt{12}$
was known since Axel Thue in 1910 but proven only by L\'asl\'o Fejes To\'oth in 1943.
The case $d=3$ was called the {\bf Kepler conjecture} as Johannes Kepler conjectured
$\Delta_3 = \pi/\sqrt{18}$. It was settled by Thomas Hales in 1998 using computer assistance.
A complete formal proof appeared in 2015.
The case $d=8$ was settled by Maryna Viazovska who proved in 2017 \cite{Viazovska} that $\Delta_8 = \pi^4/384$ and 
also established uniqueness. The densest packing in the case $d=8$ is the $E_8$ lattice.
The proof is based on linear programming bounds developed by Henry Cohn and
Noam Elkies in 2003. Later with other collaborators, she also covered the case $d=24$.
The densest packing in dimension $24$ is the {\bf Leech lattice}.
For sphere packing see \cite{CS,CoSl95}.

\index{Leech lattice}
\index{$E_8$ lattice}

\section{Sturm theorem}

Given a square free {\bf real-valued polynomial} $p$ let $p_k$ denote the {\bf Sturm chain},
$p_0=p$, $p_1=p'$, $p_2 = p_0 \; {\rm mod} \; p_1$, $p_3 = p_1 \; {\rm mod} \; p_2$ etc.
Let $\sigma(x)$ be the number of {\bf sign changes} ignoring zeros in the sequence
$p_0(x),p_1(x), \dots ,p_m(x)$.

\satz{The number of distinct roots of $p$ in $(a,b]$ is $\sigma(b)-\sigma(a)$.}

Sturm proved the theorem in 1829.
He found his theorem on sequences while studying solutions of differential
equations {\bf Sturm-Liouville theory} and credits Fourier for inspiration.
See \cite{RahmanSchmeisser}.
\index{Sturm-Liouville theory}
\index{sign changes}
\index{Sturm chain}

\section{Smith Normal form}

A integer $m \times n$ matrix $A$ is said to be expressible in {\bf Smith normal form} if there
exists an invertible $m \times m$ matrix $S$ and an invertible $n \times n$ matrix $T$
so that $S M T$ is a diagonal matrix ${\rm Diag}(\alpha_1, \dots,\alpha_r,0,0,0)$ with
$\alpha_i|\alpha_{i+1}$. The integers $\alpha_i$ are called {\bf elementary divisors}. They can be
written as $\alpha_i=d_i(A)/d_{i-1}(A)$, where $d_0(A)=1$ and $d_k(A)$ is the greatest common divisor
of all $k \times k$ minors of $A$. The Smith normal form is called {\bf unique}
if the elementary divisors $\alpha_i$ are determined up to a sign.

\satz{Any integer  matrix has a unique Smith normal form.}

The result was proven by Henry John Stephen Smith in 1861.
The result holds more generally in a {\bf principal ideal domain}, which is an
{\bf integral domain} (a ring $R$ in which $ab=0$ implies $a=0$ or $b=0$) in which
every {\bf ideal} (an additive subgroup $I$ of the ring such that $ab \in I$ if $a \in I$ and $b \in R$)
is generated by a single element.
\index{Smith normal form}
\index{principal ideal domain}

\section{Spectral perturbation}

A complex valued matrix $A$ is {\bf self-adjoint} = Hermitian if $A^*=A$, where
$A^*_{ij} = \overline{A}_{ji}$.  The spectral theorem assures that $A$ has real eigenvalues
Given two selfadjoint complex $n \times n$ matrices $A,B$ with
eigenvalues $\alpha_1 \leq \alpha_2 \leq \dots \leq \alpha_n$ and
$\beta_1 \leq \beta_2 \leq \dots \leq \beta_n$, one has the Lidskii-Last theorem:
\index{Self adjoint}
\index{Harmitian}
\index{Lidskii-Last}

\satz{ $\sum_{j=1}^n |\alpha_j - \beta_j| \leq \sum_{i,j=1}^{n} |A-B|_{ij}$.}

The result has been deduced by Yoram Last (around 1993) from {\bf Lidskii's inequality}
found in 1950 by Victor Lidskii $\sum_j |\alpha_j-\beta_j| \leq \sum_j |\gamma_j|$
where $\gamma_j$ are the eigenvalues of $C=B-A$
(see \cite{SimonTrace} page 14). The original Lidskii inequality also holds for $p \geq 1$:
$\sum_j |\alpha_j-\beta_j|^p \leq \sum_j |\gamma_j|^p$.
Last's spin on it allows to estimate the $l^1$ spectral distance of two
self-adjoint matrices using the $l^1$ distance of the matrices.
This is handy as we often know the matrices $A,B$ explicitly rather than the
eigenvalues $\gamma_j$ of $A-B$.

\section{Radon transform}

In order to solve the {\bf tomography problem} like 
{\bf magnetic resonance imaging} (MRI) of finding the density 
function $g(x,y,z)$ of a three dimensional body, one looks at a {\bf slice}
$f(x,y)=g(x,y,c)$, where $z=c$ is kept constant and measures the {\bf Radon transform} 
$R(f)(p,\theta)=\int_{\{x \cos(\theta) + y \sin(\theta)=p\}} f(x,y) \; ds$.
This quantity is the {\bf absorption rate} due to {\bf nuclear magnetic resonance} 
along the line $L$ of polar angle $\alpha$ in distance $p$ from the center.
Reconstructing $f(x,y)=g(x,y,c)$ for different $c$ allows to recover the 
{\bf tissue density} $g$ and so to ``see inside the body".
\index{tissue density}
\index{nuclear magnetic resonance}
\index{MRI}
\index{Magnetic resonance imaging}
\index{Radon transform}

\satz{The Radon transform can be diagonlized and so pseudo inverted.}

We only need that the Fourier series 
$f(r,\phi)=\sum_n f_n(r) e^{i n \phi}$ converges uniformly for all $r>0$
and that $f_n(r)$ has a Taylor series. The expansion
$f(r,\phi)=\sum_{n \in \ZZ} \sum_{k=1}^{\infty} f_{n,k} \psi_{n,k}$
with $ \psi_{n,k}(r,\phi)=r^{-k} e^{i n \phi}$ is an eigenfunction expansion
with eigenvalues $\lambda_{n,k} = 2 \int_{0}^{\pi/2} \cos(nx) \cos(x)^{(k-1)} \; dx
=\frac{\pi}{2^{k-1} \cdot k} \cdot \frac{\Gamma(k+1)}
                  {\Gamma(\frac{k+n+1}{2}) \Gamma(\frac{k-n+1}{2})}$. 
The {\bf inverse problem} is subtle due to the existence of a {\bf kernel}
spanned by $\{\psi_{n,k} \; | \; (n+k) \; {\rm odd} \; , |n|>k \}$. One calls it
an {\bf ill posed problem} in the sense of Hadamard. 
The Radon transform was first studied by Johann Radon in 1917 \cite{Helgason}.
\index{inverse problem}
\index{Fourier series}
\index{Radon transform}
\index{ill posed problem} 

\section{Linear programming}

Given two vectors $c \in \mathbb{R}^m$ and $b \in \mathbb{R}^n$, and a $n \times m$ matrix $A$,
a {\bf linear program} is the variational problem on $\mathbb{R}^m$ to maximize $f(x)=c \cdot x$ subject to
the linear constraints $Ax \leq b$ and $x \geq 0$. The dual problem is
to minimize $b \cdot y$ subject to to $A^T y \geq c ,y \geq 0$. The {\bf maximum principle}
for linear programming is tells that the solution is on the boundary of the {\bf convex polytop}
formed by the {\bf feasable region} defined by the constraints.
\index{convex polytop}
\index{linear program}
\index{maximum principle}
\index{dual linear programming problem}

\satz{Local optima of linear programs are global and on the boundary}.

Since the solutions are located on the vertices of the polytope defined by the constraints
the {\bf simplex algorithm} for solving linear programs works: start at a vertex
of the polytop, then move along the edges along the gradient until the optimum
is reached. If $A=[2,3]$ and $x=[x_1,x_2]$
and $b=6$ and $c=[3,5]$ we have $n=1,m=2$. The problem is to maximize $f(x_1,x_2)=3x_1+5x_2$
on the triangular region $2 x_1 + 3 x_2 \leq 6, x_1 \geq 0, x_2 \geq 0$.  Start at $(0,0)$,
the best improvement is to go to $(0,2)$ which is already the maximum.
Linear programming is used to solve practical problems in operations research.
The simplex algorithm was formulated by George Dantzig in 1947.
It solves random problems nicely but there are expensive cases in general and
it is possible that cycles occur. One of the open problems of Steven Smale asks for a strongly
polynomial time algorithm deciding whether a solution of a linear programming problem exists.
\index{simplex algorithm}
\cite{MurtyLinearProgramming}

\section{Random Matrices}

A {\bf random matrix} $A$ is given by an $n \times n$ array of independent, identically distributed random variables
$A_{ij}$ of zero mean and standard deviation $1$. The eigenvalues $\lambda_j$ of $A/\sqrt{n}$ define a 
discrete measure $\mu_n = \sum_j \delta_{\lambda_j}$ called {\bf spectral measure} of $A$. The {\bf circular law}
on the complex plane $\mathbb{C}$ is the probability measure $\mu_0 = 1_{D}/\pi$, where $D = \{ |z| \leq 1 \}$ 
is the unit disk. A sequence $\nu_n$ of probability measures converges {\bf weakly} or {\bf in law} to $\nu$ if
for every continuous and bounded function $f: \mathbb{C} \to \mathbb{C}$ one has $\int f(z) \; d\nu_n(z) \to \int f(z) \; d\nu(z)$.
The {\bf circular law} is:
\index{random matrix}
\index{spectral measure}

\satz{Almost surely, the spectral measures converge $\mu_n \to \mu_0$. }  

One can think of $A_n$ as a sequence of larger and larger matrix valued random variables. The circular law tells that the
eigenvalues fill out the unit disk in the complex plane uniformly when taking larger and larger matrices.
It is a kind of central limit theorem. An older version due to Eugene Wigner from 1955 is the {\bf semi-circular law} telling that 
in the self-adjoint case, the now real measures $\mu_n$ converge to a distribution with density $\sqrt{4-x^2}/(2\pi)$ on $[-2,2]$.
The circular law was stated first by Jean Ginibre in 1965 and Vyacheslav Girko 1984. 
It was proven first by Z.D. Bai in 1997.  Various authors have generalized it and removed 
more and more moment conditions. The latest condition was removed by 
Terence Tao and Van Vu in 2010, proving so the above ``fundamental theorem of random matrix theory". See \cite{TaoRandomMarix}.
\index{Wigner semi circle law}
\index{Circular law}
\index{Girko law}
\index{random matrix}
\index{eigenvalues}

\section{Diffeomorphisms}

Let $M$ be a compact Riemannian surface and $T: M \to M$ a $C^2$-diffeomorphism. A Borel
probability measure $\mu$ on $M$ is $T$-invariant if $\mu(T(A))=\mu(A)$ for all
$A \in \mathcal{A}$. It is called {\bf ergodic} if $T(A)=A$ implies $\mu(A)=1$ or
$\mu(A)=0$. The {\bf Hausdorff dimension} ${\rm dim}(\mu)$ of a measure $\mu$ is
defined as the Hausdorff dimension of the smallest Borel set $A$ of full measure
$\mu(A)=1$. The {\bf entropy} $h_{\mu}(T)$ is the {\bf Kolmogorov-Sinai entropy}
of the measure-preserving dynamical system $(X,T,\mu)$. 
For an ergodic surface diffeomorphism, the 
{\bf Lyapunov exponents} $\lambda_1,\lambda_2$ of $(X,T,\mu)$ are the
logarithms of the eigenvalues of $A=\lim_{n \to \infty} [(dT^n(x))^* dT^n(x)]^{1/(2n)}$,
which is a limiting Oseledec matrix and constant $\mu$ almost everywhere due to ergodicity.
Let $\lambda(T,\mu)$ denote the Harmonic mean of $\lambda_1,-\lambda_2$. 
The {\bf entropy-dimension-Lyapunov theorem} tells that for every $T$-invariant ergodic
probability measure $\mu$ of $T$, one has:
\index{Hausdorff dimension of measure}
\index{Lyapunov exponent}
\index{Entropy} 

\satz{$h_{\mu} = {\rm dim}(\mu) \lambda/2$.} 

This formula has become famous because it relates ``entropy", ``fractals" and ``chaos", which
are all ``rock star" notions also outside of mathematics. The theorem
implies in the case of Lebesgue measure preserving 
symplectic transformation, where ${\rm dim}(\mu)=2$ and $\lambda_1=-\lambda_2$ that
``entropy = Lyaponov exponent" which is a {\bf formula of Pesin} given by 
$h_{\mu}(T) = \lambda(T,\mu)$. A similar result holds for {\bf circle     
diffeomorphims} or smooth interval maps, where $h_{\mu}(T) = {\rm dim}(\mu) \lambda(T,\mu)$.
The notion of Hausdorff dimension was introduced by Felix Hausdoff in 1918. Entropy 
was defined in 1958 by Nicolai Kolmogorov and in general by Yakov Sinai in 1959, 
Lyapunov exponents were introduced with the work of Valery Oseledec in 1965.  
The above theorem is due to Lai-Sang Young who proved it in 1982. 
Francois Ledrapier and Lai-Sang Young proved in 1985 that in arbitrary dimensions,      
$h_{\mu} = \sum_j \lambda_j \gamma_j$, where $\gamma_j$ are dimensions of    
$\mu$ in the direction of the Oseledec spaces $E_j$. This is called the 
{\bf Ledrappier-Young formula}. It implies the {\bf Margulis-Ruelle inequality}
$h_{\mu}(T)  \leq \sum_j \lambda_j^+(T)$, where $\lambda_j^+={\rm max}(\lambda_j,0)$
and $\lambda_j(T)$ are the Lyapunov exponents.
In the case of a smooth $T$-invariant measure $\mu$ or more generally, for SRB measures,
there is an equality
$h_{\mu}(T)  = \sum_j \lambda_j^+(T)$ which is called the {\bf Pesin formula}.  
See \cite{KH,RuelleEckmann}.
\index{Pesin formula}
\index{Ledrappier-Young formula}
\index{Margulis-Ruelle inequality}

\section{Linearization}

If $F: M \to M$ is a globally Lipschitz continuous function
on a finite dimensional vector space $M$, then
the differential equation $x'=F(x)$ has a global solution 
$x(t)=f^t(x(0))$ (a local by {\bf Picard-Lindel\"of 's existence theorem} and global by 
the {\bf Gr\"onwall inequality}). An {\bf equilibrium point} of
the system is a point $x_0$ for which $F(x_0)=0$.
This means that $x_0$ is a fixed point of a differentiable
mapping $f=f^1$, the {\bf time-1-map}. We say that $f$ is {\bf linearizable} 
near $x_0$ if there exists a homeomorphism $\phi$ from a neighborhood 
$U$ of $x_0$ to a neighborhood $V$ of $x_0$ such that $\phi \circ f \circ \phi^{-1} = df$. 
The {\bf Sternberg-Grobman-Hartman linearization theorem} is
\index{Sternberg linarization}
\index{Grobman-Hartman linearization}
\index{Gr\"onwall inequality}
\index{Time-1-map}

\satz{If $f$ is hyperbolic, then $f$ is linearizable near $x_0$.}

The theorem was proven by D.M. Grobman  in 1959 Philip Hartman in 1960 and 
by Shlomo Sternberg in 1958.
This implies the existence of {\bf stable and unstable manifolds} passing through $x_0$.
One can show more and this is due to Sternberg who wrote a series of papers starting
1957 \cite{Sternberg1957}:
if $A=df(x_0)$ satisfies {\bf no resonance condition}
meaning that no relation $\lambda_0 = \lambda_1 \cdots \lambda_j$ exists between eigenvalues
of $A$, then a {\bf linearization to order $n$}
is a $C^{n}$ map $\phi(x) = x + g(x)$, with $g(0)=g'(0)=0$
such that $\phi \circ f \circ \phi^{-1}(x) = A x + o(|x|^n)$ near $x_0$. 
We say then that $f$ can be {\bf $n$-linearized} near $x_0$. The generalized result
tells that non-resonance fixed points of $C^n$ maps are
$n$-linearizable near a fixed point. See \cite{Lanford87}.
\index{n-linearization}
\index{resonance condition}

\section{Fractals}

An {\bf iterated function system} is a finite set of contractions $\{ f_i \}_{i=1}^n$ on a complete
metric space $(X,d)$. The corresponding {\bf Huntchingson operator} $H(A) = \sum_{i} f_i(A)$ is 
then a contraction on the {\bf Hausdorff metric} of sets and has a unique fixed point called the 
{\bf attractor} $S$ of the iterated function system. 
The definition of {\bf Hausdorff dimension} is as follows: define
$h_{\delta}^s(A) = \inf_{U \in \mathcal{U}} \sum_i |U_i|^s$, where $\mathcal{U}$ is a $\delta$-cover
of $A$. And $h^s(A)=\lim_{\delta \to 0} H_{\delta}^s(A)$. The {\bf Hausdorff dimension}
${\rm dim}_H(S)$ finally is the value $s$, where $h^s(S)$ jumps from $\infty$ to $0$.
If the contractions are maps with contraction factors $0<\lambda_j<1$
then the Hausdorff dimension of the attractor $S$ can be estimated with the 
the {\bf similarity dimension} of the contraction vector $(\lambda_1, \dots, \lambda_n)$:
this number is defined as the solution $s$ of the equation $\sum_{i=1}^n \lambda_i^{-s} = 1$.
\index{similarity dimension}
\index{iterated function system} 
\index{fractal}

\satz{ ${\rm dim}_{{\rm hausdorff}}(S) \leq {\rm dim}_{{\rm similarity}}(S)$. }

There is an equality if $f_i$ are all affine contractions like $f_i(x)=A_i \lambda x+\beta_i$ with 
the same contraction factor and $A_i$ are orthogonal and $\beta_i$ are vectors (a situation which generates
a large class of popular fractals). For equality one also has to assume that there is an open non-empty set $G$
such that $G_i=f_i(G)$ are disjoint. In the case $\lambda_j=\lambda$ are all the same then $n \lambda^{-{\rm dim}} = 1$
which implies ${\rm dim}(S) = -\log(n)/\log(\lambda)$. For the {\bf Smith-Cantor set} $S$, where
$f_1(x)=x/3+2/3, f_2(x)=x/3$ and $G=(0,1)$. One gets with $n=2$ and $\lambda=1/3$ the dimension
${\rm dim}(S)=\log(2)/\log(3)$. For the {\bf Menger carpet} with $n=8$ affine maps $f_{ij}(x,y)=(x/3+i/3,y/3+j/3)$
with $0 \leq i \leq 2, 0 \leq j \leq 2, (i,j) \neq (1,1)$, the dimension is $\log(8)/\log(3)$. The {\bf Menger sponge}
is the analogue object with $n=20$ affine contractions in $\mathbb{R}^3$ and has dimension $\log(20)/\log(3)$. 
For the {\bf Koch curve} on the interval, where $n=4$ affine contractions of contraction factor $1/3$
exist, the dimension is $\log(4)/\log(3)$. These are all {\bf fractals}, sets with Hausdorff dimension different
from an integer. The modern formulation of iterated function systems is due 
to John E. Hutchingson from 1981.  Michael Barnsley used the concept for a 
{\bf fractal compression algorithms}, which uses the idea that storing the rules for an iterated function system
is much cheaper than the actual attractor. 
Iterated function systems appear in complex dynamics in the case when the {\bf Julia set} is
completely disconnected, they have appeared earlier also in work of Georges de Rham 1957. 
See \cite{Mandelbrot,Falconer}.
\index{Cantor set}
\index{Huntchingson operator}
\index{Hausdorff dimension}
\index{Koch curve}
\index{Menger carpet}
\index{attractor of iterated function system}

\section{Strong law of small numbers}

Like the Bayes theorem or the Pigeon hole principle which both are too simple to qualify as
``theorems" but still are of utmost importance, the 
``Strong law of small numbers" is not really a theorem but a 
{\bf fundamental mathematical principle}. It is more fundamental than a specific theorem as it
applies throughout mathematics. It is for example important in Ramsey theory:
The statement is put in different ways like 
"There aren't enough small numbers to meet the many demands made of them". 
\cite{stronglawofsmallnumbers} puts it in the following catchy way: 
\index{Strong law of small numbers}

\satz{You can't tell by looking.}

The point was made by Richard Guy in \cite{stronglawofsmallnumbers} who states 
two ``corollaries": 
{\bf ``superficial similarities spawn spurious statements"} and
{\bf ``early exceptions eclipse eventual essentials"}. 
The statement is backed up with countless
many examples (a list of 35 are given in \cite{stronglawofsmallnumbers}). 
Famous are Fermat's claim that all {\bf Fermat primes} $2^{2^n}+1$ are prime or 
the claim that the number $\pi_3(n)$ of primes of the form $4k+3$ in $\{1, \dots, n\}$ is larger than 
$\pi_1(n)$ of primes of the form $4k+1$ so that the $4k+3$ primes win the {\bf prime race}.
Hardy and Littlewood showed however $\pi_3(n)-\pi_1(n)$ changes sign infinitely often. The prime number theorem
extended to arithmetic progressions shows $\pi_1(n) \sim n/(2\log(n))$ and $\pi_3(n) \sim n/(2 \log(n))$
but the density of numbers with $\pi_3(n)>\pi_1(n)$ is larger than $1/2$. This is the {\bf Chebyshev
bias}. Experiments then suggested the density to be $1$ but also this is false: the density of numbers 
for which $\pi_3(n)>\pi_1(n)$ is smaller than $1$. 
The principle is important in a branch of combinatorics called {\bf Ramsey theory}. But it not only
applies in discrete mathematics. There are many examples, where one can not tell by looking. 
When looking at the boundary of the Mandelbrot set for example, 
one would tell that it is a fractal with Hausdorff dimension between $1$ and $2$. In reality
the Hausdorff dimension is $2$ by a result of Mitsuhiro Shishikura. 
Mandelbrot himself thought first ``by looking" that the 
Mandelbrot set $M$ is disconnected. Douady and Hubbard proved $M$ to be connected. 

\section{Ramsey Theory}

Let $G$ be the complete graph with $n$ vertices. 
An {\bf edge labeling} with $r$ colors is an assignment of $r$ numbers 
to the {\bf edges} of $G$. A complete sub-graph of $G$ is called a {\bf clique}. 
If it is has $s$ vertices, it is denoted by $K_s$. 
A graph $G$ is called {\bf monochromatic} if all edges in $G$ have 
the same color. (We use in here {\bf coloring} as a short for {\bf edge labeling}
and not in the sense of chromatology where an edge coloring assumes that intersecting 
edges have different colors.) Ramsey's theorem is:
\index{monochromatic}
\index{edges}
\index{coloring}
\index{edge labeling}

\satz{For large $n$, every $r$-colored $K_n$ contains a monochromatic $K_s$.}

So, there exist {\bf Ramsey numbers} $R(r,s)$ such that for $n \geq R(r,s)$, the
edge coloring of one of the $s$-cliques can occur. A famous case is the identity
$R(3,3)=6$. Take $n=6$ people. It defines the complete graph $G$. 
If two of them are friends, color the edge blue, otherwise red. This
{\bf friendship graph} therefore is a $r=2$ coloring of $G$. There are 78 possible
colorings. In each of them, there is a triangle of friends or a triangle of strangers. 
In a group of 6 people, there are either a clique with 3 friends or a clique of 
3 complete strangers. The theorem was proven by Frank Ramsey in 1930. Paul Erdoes asked
to give explicit estimated $R(s)$ which is the least integer
$n$ such that any graph on $n$ vertices contains either a {\bf clique} of size $s$ (a set where all 
are connected to each other) or an independent set of size $s$ (a set where none are
connected to each other). Graham for example asks whether the limit $R(n)^{1/n}$ exists. 
Ramsey theory also deals other sets: {\bf van der Waerden's theorem} 
from 1927 for example tells that if the positive integers $\mathbb{N}$ are colored with $r$ colors, then for every $k$,
there exists an $N$ called $W(r,k)$ such that the finite set $\{1 \dots, N\}$ has
an arithmetic progression with the same color. For example, $W(2,3)=9$. Also here, it is an
open problem to find a formula for $W(r,k)$ or even give good upper bounds. \cite{GrahamRamsey}
\cite{GrahamRudimentsRamsey}
\index{Ramsey theory} 
\index{Friendship problem}
\index{Van der Waerden's theorem}

\section{Poincar\'e Duality}

For a differentiable {\bf Riemannian $n$-manifold} $(M,g)$ there is an
{\bf exterior derivative} $d=d_p$ which maps $p$-forms
$\Lambda^p$ to $(p+1)$-forms $\Lambda^{p+1}$. For $p=0$, the derivative is
called the {\bf gradient}, for $p=1$, the derivative is called the {\bf curl} and for
$p=d-1$, the derivative is the adjoint of {\bf divergence}.
The Riemannian metric defines an inner product $\langle f,h \rangle$
on $\Lambda^p$ allowing so to see $\Lambda^p$ as part of a Hilbert space and to
define the adjoint $d^*$ of $d$. It is a linear map from $\Lambda^{p+1}$ to 
$\Lambda^p$. The exterior derivative defines so the self-adjoint 
{\bf Dirac operator} $D=d+d^*$ and the
{\bf Hodge Laplacian} $L=D^2 = d d^* + d^*d$ which now leaves each $\Lambda^p$ 
invariant. {\bf Hodge theory} assures that ${\rm dim}( {\rm ker} (L| \Lambda^p)) =b_p = {\rm dim}(H^p(M))$,
where $H^p(M)$ are the $p$'th {\bf cohomology group}, the kernel of $d_p$
modulo the image of $d_{p-1}$. {\bf Poincar\'e duality} is:
\index{Poincare duality}
\index{Cohomology}
\index{Dirac operator}

\satz{If $M$ is orientable n-manifold, then $b_k(M)=b_{n-k}(M)$.}

The {\bf Hodge dual} of $f \in \Lambda^p$ is defined as the unique $*g \in 
\Lambda^{n-p}$ satisfying  $\langle f,*g \rangle = \langle f \wedge g,\omega \rangle$
where $\omega$ is the volume form. One has $d^* f = (-1)^{d+dp+1} *d* f$
and $L * f = * Lf$. This implies that $*$ is a unitary map from
${\rm ker}(L| \Lambda^p)$ to
${\rm ker}(L| \Lambda^{d-p})$ proving so the duality theorem. 
For $n=4k$, one has $*^2=1$, allowing to define the
{\bf Hirzebruch signature} $\sigma:={\rm dim}\{ u | Lu=0, *u=u \}
- {\rm dim}(u | Lu=0, *u=-u \}$.
The Poinar\'e duality theorem was first stated by Henri Poincar\'e 
in 1895. It took until the 1930ies to clean out the notions and make it precise. 
The Hodge approach establishing an explicit isomorphism between harmonic
$p$ and $n-p$ forms appears for example in
\cite{Cycon}. 
\index{Hodge dual}

\section{Rokhlin-Kakutani approximation}

Let $T$ be an automorphism of a probability space $(\Omega,\mathcal{A},\mu)$.
This means $\mu(A)=\mu(T(A))$ for all $A \in \mathcal{A}$. The system $T$
is called {\bf aperiodic}, if the set of {\bf periodic points }
$P=\{ x \in \Omega \; | \;  \exists n >0, T^n x= x \}$
has measure $\mu(P)=0$. A set $B \in \mathcal{A}$ which has the property that
$B,T(B), \dots, T^{n-1}(B)$ are disjoint is called a {\bf Rokhlin tower}. If the
measure of the tower is $\mu(B \cup \cdots \cup T^{n-1}(B)) = n \mu(B) = 1-\epsilon$, we call
it an $(1-\epsilon)$-Rokhlin tower. We say $T$ can be {\bf approximated arbitrary well}
by Rokhlin towers, if for all $\epsilon >0$, there is an $(1-\epsilon)$ Rokhlin tower.
\index{Rokholin tower}
\index{aperiodic}
\index{periodic points}

\satz{An aperiodic $T$ can be approximated well by Rokhlin towers. }

The result was proven by Vladimir Abramovich Rokhlin in his thesis 1947 and independently by 
Shizuo Kakutani in 1943. The lemma can be used to build {\bf Kakutani skyscrapers}, which are
nice partitions associated to a transformation. 
This lemma allows to approximate an aperiodic transformation $T$ by a periodic
transformations $T_n$.  Just change $T$ on $T^{n-1}(B)$ so that $T_n^n(x)=x$ for all $x$.
The theorem has been generalized by Donald Ornstein and Benjamin Weiss 
to higher dimensions like $\mathbb{Z}^d$ actions of 
measure preserving transformations where the periodicity assumption is replaced by the
assumption that the action is {\bf free}: for any $n \neq 0$, the set $T^n(x)=x$ has zero 
measure. See \cite{CFS,Friedman,Halmos}.
\index{Kakutani skyscraper}
\index{Free action}

\section{Lax approximation}

On the group $\mathcal{X}$ of all measurable, invertible transformations on
the $d$-dimensional {\bf torus} $X=\mathbb{T}^d$ which preserve the Lebesgue volume measure,
one has the metric
$$ \delta(T,S) = |  \delta(T(x),S(x)) |_{\infty} \; , $$
where ${\rm \delta}$ is the geodesic distance on the flat torus and
where $|\cdot|_{\infty}$ is the $L^{\infty}$ supremum norm. Lets call $(\mathbb{T}^d,T,\mu)$
a {\bf toral dynamical system} if $T$ is a {\bf homeomorphism}, a continuous transformation
with continuous inverse. A {\bf cube exchange transformation} on $\mathbb{T}^d$ is a 
periodic, piecewise affine measure-preserving transformation $T$ which permutes rigidly all the cubes
$\prod_{i=1}^d [k_i/n,(k_i+1)/n]$, where $k_i \in \{0, \dots, n-1 \}$.
Every point in $\mathbb{T}^d$ is $T$ periodic. A cube exchange transformation is determined by
a permutation of the set $\{1, \dots, n \}^d$. If it is
cyclic, the exchange transformation is called {\bf cyclic}.
A theorem of Lax \cite{Lax71} states that every toral dynamical system
can approximated in the metric $\delta$ by cube exchange transformations.
The approximations can even be cyclic \cite{AlPr93}.
\index{toral dynamical system}
\index{cube exchange transformation}
\index{cyclic cube exchange transformations}

\satz{Toral systems can be approximated by cyclic cube exchanges}

The result is due to Peter Lax \cite{Lax71}. The proof of this result uses Hall's marriage theorem
in graph theory (for a 'book proof' of the later theorem, see \cite{AigZie}).
Periodic approximations of symplectic maps work surprisingly
well for relatively small $n$ (see \cite{Ran74}). On the Pesin region
this can be explained in part by the shadowing property \cite{KH}. The
approximation by cyclic transformations make long time
stability questions look different \cite{Hal87}.
\index{Hall Mariage}

\section{Sobolev embedding}

All functions are defined on $\mathbb{R}^n$, integrated
$\int$ over $\mathbb{R}^n$ and assumed to be
{\bf locally integrable} meaning that
for every compact set $K$ the {\bf Lebesgue integral}
$\int_K |f| \; dx$ is finite. For functions in $C_c^{\infty}$
which serve as {\bf test functions}, {\bf partial derivatives}
$\partial_i = \partial/\partial_{x_i}$ and more general
{\bf differential operators}
$D^{k}= \partial_{x_1}^{k_1} \cdots \partial_{x_n}^{k_n}$
can be applied. A function $g$ is a {\bf weak partial derivative}
of $f$ if $\int f \partial_i \phi dx = -\int g \phi dx$ for
all test functions $\phi$. For $p \in [1,\infty)$, the $L^p$ space
is $\{ f \; | \; \int |f|^p dx<\infty \}$.
The {\bf Sobolev space} $W^{k,p}$ is the set of functions for which
all $k$'th weak derivatives are in $L^p$. So $W^{0,p}=L^p$.
The {\bf H\"older space} $C^{r,\alpha}$
with $r \in \mathbb{N},\alpha \in (0,1]$ is
defined as the set of functions for which all $r$'th derivatives
are $\alpha$-H\"older continuous. It is a Banach space with
norm ${\rm max}_{|k| \leq r} ||D^kf||_{\infty} +
      {\rm max}_{|k| = r} ||D^k f||_{\alpha}$, where
$||f||_{\infty}$ is the {\bf supremum norm} and $||f||_{\alpha}$
is the {\bf H\"older coefficient} $\sup_{x \neq y} |f(x)-f(y)|/|x-y|^{\alpha}$.
The {\bf Sobolev embedding theorem} is
\index{Sobolev embedding}
\index{differential operator}
\index{Holder continuity}
\index{test functions}
\index{partial derivatives}

\satz{If $n<p$ and $l=r+\alpha<k-n/p$, one has $W^{k,p} \subset C^{r,\alpha}$. }

(\cite{Simon2017} states this as Theorem 6.3.6) gives some history: 
{\bf generalized functions} appeared first in the
work of Oliver Heaviside in the form of ``operational calculus. Paul Dirac used the
formalism in quantum mechanics. In the 1930s, 
Kurt Otto Friedrichs, Salomon Bocher and Sergei Sobolev define
weak solutions of PDE's. Schwartz used the $C_c^{\infty}$ functions, smooth functions
of compact support. 
This means that the existence of $k$ weak derivatives implies the existence
of actual derivatives. For $p=2$, the spaces $W^k$ are Hilbert spaces and the theory
a bit simpler due to the availability of Fourier theory, where tempered distributions flourished.
In that case, one can define for any real $s>0$ the Hilbert space
$H^s$ as the subset of all $f \in S'$ for which $(1+|\xi|^2)^{s/2} \hat{f}(\xi)$
is in $L^2$. The Schwartz test functions $S$ consists of all $C^{\infty}$ functions
having bounded semi norms
$||\phi||_k=\max_{|\alpha|+|\beta| \leq k} || x^{\beta} D^{\alpha} \phi||_{\infty}<\infty$
where $\alpha,\beta \in \mathbb{N}^n$. 
Since $S$ is larger than the set of smooth functions of compact support, the dual
space $S'$ is smaller. They are {\bf tempered distributions}.
Sobolev emedding theorems like above allow to show that weak solutions of PDE's are smooth:
for example, if the Poisson problem $\Delta f = V f$ with smooth 
$V$ is solved by a distribution $f$, then $f$ is smooth. 
\cite{Brezis2011,Simon2017}
\index{generalized functions}
\index{tempered distributions}

\section{Whitney embedding}

A smooth {\bf $n$-manifold} $M$ is a metric space equipped
with a cover $U_j=\phi_j^{-1}(B)$ with $B=\{ x \in \mathbb{R}^n \; | \; |x|^2<1 \})$
or $U_j=\phi_j^{-1}(H)$ with $H=\{ x \in \mathbb{R}^n \; | \; |x|^2<1, x_0 \geq 0 \})$ with
$\delta H=\{ x \in H \; | \; x_0 = 0 \}$ such that
the homeomorphisms $\phi_j: U_j \to B$ or $\phi_j: U_j \to H$ lead to smooth transition maps
$\phi_{kj} =\phi_j \phi_k^{-1}$ from $\phi_k(U_j \cap U_k)$ to $\phi_j(U_j \cap U_k)$
which have the property that all restrictions of $\phi_{kj}$ from       
$\delta \phi_k(U_j \cap U_k)$ to $\delta \phi_j(U_j \cap U_k)$ are smooth too. The
{\bf boundary} $\delta M$ of $M$ now naturally is a smooth $(n-1)$ manifold, the atlas
being given by the sets $V_j = \phi_j(\delta H)$ for the indices $j$ which map $\phi_j: U_j \to H$. 
Two manifolds $M,N$ are {\bf diffeomorphic} if there is a refinement $\{ U_j,\phi_j \}$ of the atlas in $M$
and a refinement $\{V_j,\psi_j \}$ of the atlas in $N$ such that $\phi_j(U_j) = \psi_j(V_j)$.
A manifold $M$ can be {\bf smoothly embedded} in 
$\mathbb{R}^k$ if there is a smooth injective map $f$ from $M$ to $\mathbb{R}^k$ such that 
the image $f(M)$ is diffeomorphic to $M$. 
\index{manifold}
\index{boundary}
\index{embedding}
\index{diffeomorphic}
\index{cover}

\satz{Any $n$-manifold $M$ can be smoothly embedded in $\mathbb{R}^{2n}$.}

The theorem has been proven by Hassler Whitney in 1926 who also was the first to give a precise
definition of manifold in 1936. The standard assumption is that $M$ is second countable Hausdorff 
but as every smooth finite dimensional manifold can be upgraded to be Riemannian, the simpler 
metric assumption is no restriction of generality. 
The modern point of view is to see $M$ as a {\bf scheme} over Euclidean $n$-space,
more precisely as a {\bf ringed space}, that is locally the spectrum of the commutative ring 
$C^{\infty}(B)$ or $C^{\infty}(H)$. The set of manifolds is a {\bf category} in which the smooth maps
$M \to N$ are the {\bf morphisms}. The cover $U_j$ defines an {\bf atlas} and the transition maps $\phi_j$ allow to
port notions like smoothness from Euclidean space to $M$. The maps $\phi_j^{-1}: B \to M$ or $\phi_j^{-1}: H \to M$ 
parametrize the sets $U_j$. \cite{WhitneyCollected}.
\index{parametrization}
\index{scheme}
\index{commutative ring}

\section{Artificial intelligence}

Like {\bf meta mathematics} or {\bf reverse mathematics}, the field of
{\bf artificial intelligence} (AI) is a part of mathematics which also reflects on 
subject itself. It is related of {\bf data science} (algorithms for data mining,
and statistics)  {\bf computation theory}
(like complexity theory) {\bf language theory} and especially
{\bf grammar} and {\bf evolutionary dynamics},
{\bf optimization problems} (like solving optimal transport
or extremal problems) {\bf solving inverse problems} (like developing
algorithms for computer vision or optical character or speech
recognition), {\bf cognitive science} as well as {\bf pedagogy} in education
(human or machine learning and human motivation). There is no apparent ``fundamental theorem"
of AI, (except maybe for Marvin Minsky's
{\it "The most efficient way to solve a problem is to already know how to solve it."} \cite{minsky},
which is a surprisingly deep and insightful statement as modern AI agents like 
{\bf Alexa}, {\bf Siri}, {\bf Google Home}, {\bf IBM Watson} or {\bf Cortana} demonstrate;
they compute little, they just know or look up - or annoy you to look it up yourself...). But there is a {\bf theorem of Lebowski on 
machine super intelligence} which taps into the rather uncharted territory of {\bf machine motivation} 
\index{Theorem of Lebowski}
\index{data science}
\index{data mining}
\index{complexity theory}
\index{inverse problems}
\index{cognitive science}
\index{optimization problem}
\index{OCR}
\index{computer vision}

\satz{No AI will bother after hacking its own reward function.}

The picture \cite{JasonKottke} is that once the AI has figured out the philosophy of the
``Dude" in the Cohen brothers movie Lebowski, also repeated mischiefs does not bother it
and it ``goes bowling". Objections are brushed away with ``Well, this is your, like, opinion, man".
Two examples of human super intelligent units who have succeeded to hack their 
own reward function are Alexander Grothendieck or Grigori Perelman.
The Lebowski theorem is due to Joscha Bach \cite{JoschaBach},
who stated this {\bf theorem of super intelligence} in a tongue-in-cheek tweet.
From a mathematical point of view, the smartest way to ``solve" an optimal transport
problem is to change the utility function. On a more serious level, the 
smartest way to ``solve" the continuum hypothesis is to change the axiom system. 
This might look like a cheat, but on a meta level, more creativity is possible. 
Precursor's of the Lebowski theme is Stanislav Lem's notion of a {\bf mimicretin}
\cite{FuturolgicalCongress}, a computer that plays stupid in order, once and for all, to be 
left in peace or the machine in \cite{HitchhikersGuide} who develops
humor and enjoys fooling humans with the answer to the ultimate question: ``42". This 
document entry is the analogue to the ultimate question: 
``What is the fundamental theorem of AI"? 
\index{mimiocretin}
\index{Lebowski theorem}
\index{ultimate question}

\section{Stokes Theorem}

On a smooth orientable $n$-dimensional manifold $M$, one has
$\Lambda^p$, the vector bundle of smooth {\bf differential $p$-forms}.
As any $p$-form $F$ induces an {\bf induced volume form} on a $p$-dimensional {\bf sub-manifold} $G$
defining so an {\bf integral} $\int_G F$. The {\bf exterior derivative} $d: \Lambda^p \to \Lambda^{p+1}$
satisfies $d^2=0$ and defines an {\bf elliptic complex}. There is a natural {\bf Hodge duality}
isomorphism given called ``Hodge star" $*: \Lambda^p \to \Lambda^{n-p}$. 
Given a $p$-form $F \in \Lambda^p$ and a $(p+1)$-dimensional compact oriented
sub-manifold $G$ of $M$ with boundary $\delta G$ compatible with the orientation
of $G$, we have {\bf Stokes theorem}:`
\index{Stokes theorem}
\index{differential form}
\index{submanifold}

\satz{ $\langle G, dF \rangle = \int_G dF = \int_{\delta G} F = \langle \delta G, F \rangle$.  }

The theorem states that the exterior derivative $d$ is dual to the boundary operator $\delta$.
If $G$ is a connected $1$-manifold with boundary, it is a curve
with boundary $\delta G=\{A,B\}$. A $1$-form
can be integrated over the curve $G$ by choosing the on $G$ induced volume form $r'(t) dt$ given by a
{\bf curve parametrization} $[a,b] \to G$ and integrate
$\int_a^b F(r(t)) \cdot r'(t) dt$, which is the {\bf line integral}.
Stokes theorem is then the {\bf fundamental theorem of line integrals}. Take a $0$-form
$f$ which is a {\bf scalar function} the derivative $df$ is the gradient $F=\nabla f$. Then
$\int_a^b \nabla f(r(t)) \cdot r'(t) \; dt=f(B)-f(A)$. If $G$ is a two
dimensional surface with boundary $\delta G$ and $F$ is a $1$-form, then
the $2$-form $dF$ is the {\bf curl} of $F$. If $G$ is given as a
{\bf surface parametrization} $r(u,v)$, one can apply $dF$ on the pair of tangent
vectors $r_u,r_v$ and integrate this $dF(r_u,r_v)$ over the surface $G$ to get $\int_G dF$.
The {\bf Kelvin-Stokes theorem} tells that this is the same than the line integral $\int_{\delta G} F$.
In the case of $M=\mathbb{R}^3$, where $F=Pdx+Qdy+Rdz$ can be identified with a
vector field $F=[P,Q,R]$ and $dF = \nabla \times F$ and integration of a $2$-form $H$
over a parametrized manifold $G$ is
$\int \int_R H(r(u,v))(r_u,r_v) = \int \int_R H(r(u,v) \cdot r_u \times r_v du dv$
we get the {\bf classical Kelvin-Stokes theorem.} If $F$ is a $2$-form, then $dF$ is a
$3$-form which can be integrated over a $3$-manifold $G$. As $d: \Lambda^2 \to \Lambda^3$
can via Hodge duality naturally be paired with $d_0^*: \Lambda^1 \to \Lambda^0$, which is the
{\bf divergence}, the {\bf divergence theorem} $\int \int \int_G {\rm div}(F) \; dx dy dz 
= \int\int_{\delta G} F \cdot dS$ relates a triple integral with a flux integral.
Historical milestones start with the development of the {\bf fundamental theorem of calculus}
(1666 Isaac Newton, 1668 James Gregory, Isaac Barrow 1670 and Gottfried Leibniz 1693);
the first rigorous proof was done by Cauchy in 1823 (the first textbook appearance in 1876 by Paul du Bois-Reymond).
See \cite{BressoudPortland}. In 1762, Joseph-Louis Lagrange and in 1813 Karl-Friedrich Gauss look
at special cases of divergence theorem, Mikhail Ostogradsky in 1826 and George Green in 1828
cover the general case. Green's theorem in two dimensions was first stated
by Augustin-Louis Cauchy in 1846 and Bernhard Riemann in 1851. Stokes theorem first appeared
in 1854 as an exam question but the theorem has appeared already in a letter of William Thomson
to Lord Kelvin in 1850, hence also the name {\bf Kelvin-Stokes theorem}.
Vito Volterra in 1889 and Henri Poincar\'e in 1899 generalized the theorems to higher dimensions.
Differential forms were introduced in 1899 by \'Elie Cartan. The $d$ notation for exterior
derivative was introduced in 1902 by Theodore de Donder. The ultimate formulation above is
from Cartan 1945. We followed Katz \cite{HistoryStokes} who noticed that only in 1959,
this version has started to appear in textbooks.

\section{Moments}

The {\bf Hausdorff moment problem} asks for necessary and sufficient conditions for a sequence $\mu_n$ to
be realizable as a moment sequence $\int_0^1 x^n \; d\mu(x)$ for a Borel probability measure on $[0,1]$.
One can study the problem also in higher dimensions: for a multi-index $n=(n_1,\dots, n_d)$ denote
by $\mu_n = \int x_1^{n_1} \dots x_d^{n_d} \; d\mu(x)$ the {\bf $n$'th moment} of a
{\bf signed Borel measure} $\mu$ on the unit cube $I^d=[0,1]^d \subset \RR^d$.
We say $\mu_n$ is a {\bf moment configuration} if there exists a measure $\mu$
which has $\mu_n$ as moments. If $e_i$ denotes the standard basis in $\ZZ^d$,
define the {\bf partial difference} $(\Delta_i a)_n = a_{n-e_i}-a_n$ and $\Delta^k=\prod_i \Delta_i^{k_i}$.
We write $\frac{k}{n}=\prod_{i=1}^n \frac{k_i}{n_i}$ and $\left( \begin{array}{c} n \\ k \\ \end{array} \right)  
= \prod_{i=1}^d \left( \begin{array}{c} n_i \\ k_i \end{array} \right)$ and
$\sum_{k=0}^n = \sum_{k_1=0}^{n_1}  \dots \sum_{k_d=0}^{n_d}$.
We say moments $\mu_n$ are {\bf Hausdorff bounded} if there exists a constant $C$ such that
$\sum_{k=0}^n |\left( \begin{array}{c} n \\ k \\ \end{array} \right) (\Delta^{k} \mu)_n| \leq C$ for all
$n \in \mathbb{N}^d$. The {\bf theorem of Hausdorff-Hildebrandt-Schoenberg} is
\index{Borel measure}
\index{partial differences}
\index{Theorem of Hausdorff-Hildebrandt-Schoenberg}

\satz{Hausdorff bounded moments $\mu_n$ are generated by a measure $\mu$.}

The above result is due to Theophil Henry Hildebrandt and Isaac Jacob Schoenberg from 1933.
\cite{HiSc33}. Moments also allow to compare measures: a measure $\mu$ is called
{\bf uniformly absolutely continuous} with respect to $\nu$
if there exists $f \in L^{\infty}(\nu)$ such that $\mu = f \nu$.
A positive probability measure $\mu$ is uniformly absolutely continuous with respect to
a second probability measure $\nu$ if and only if there exists a constant $C$ such that
$(\Delta^k \mu)_n \leq C \cdot (\Delta^k \nu)_n$ for all $k,n \in \NN^d$.
In particular it gives a generalization of a result of Felix Hausdorff from 1921
\cite{Hau21} assuring that $\mu$ is positive if and only if
$(\Delta^k \mu)_n \geq 0$ for all $k,n \in \NN^d$.
An other special case is that $\mu$ is uniformly absolutely continuous with respect
to Lebesgue measure $\nu$ on $I^d$ if and only if
$|\Delta^k \mu_n| \leq  \left( \begin{array}{c} n \\ k \\ \end{array} \right) (n+1)^d$
for all $k$ and $n$. Moments play an important role in statistics,
when looking at {\bf moment generating functions} $\sum_n \mu_n t^n$  of random variables $X$, where
$\mu_n = {\rm E}[X^n]$ as well as in {\bf multivariate statistics}, when looking
at random vectors $(X_1, \dots, X_d)$, where $\mu_n =  {\rm E}[X_1^{n_1} \cdots X_d^{n_d}]$
are {\bf multivariate moments}.  See \cite{knillprobability,Schmuedgen2017}
\index{multivariate moments}
\index{Moment methods}
\index{moment generating function}

\section{Martingales}

A sequence of random variables $X_1,X_2, \dots$ on a probability space 
$(\Omega,\mathcal{A},{\rm P})$ is called a {\bf discrete time stochastic process}.
We assume the $X_k$ to be in $L^2$ meaning that the expectation ${\rm E}[X_k^2] < \infty$ for all $k$.
Given a sub-$\sigma$ algebra $\mathcal{B}$ of $\mathcal{A}$, the {\bf conditional expectation} ${\rm E}[X|\mathcal{B}]$ 
is the projection of $L^2(\Omega,\mathcal{A},P)$ to $L^2(\Omega,\mathcal{B},P)$. Extreme
cases are ${\rm E}[X|\mathcal{A}]=X$ and ${\rm E}[X|\{\emptyset,\Omega\}] = {\rm E}[X]$.
A finite set $Y_1, \dots, Y_n$ of random variables generates a sub- $\sigma$-algebra $\mathcal{B}$ of 
$\mathcal{A}$, the smallest $\sigma$-algebra for which all 
$Y_j$ are still measurable. Write ${\rm E}[X|Y_1,\cdots ,Y_n]={\rm E}[X|\mathcal{B}]$,
where $\mathcal{B}$ is the $\sigma$-algebra generated by $Y_1, \cdots Y_n$. 
A discrete time stochastic process is called a {\bf martingale} if ${\rm E}[X_{n+1} | X_1,\cdots,X_n] = {\rm E}[X_n]$
for all $n$. If the equal sign is replaced with $\leq$ then the process is called a {\bf super-martingale}, if $\geq$ it is a 
{\bf sub-martingale}. The {\bf random walk} $X_n = \sum_{k=1}^n Y_k$ defined by a sequence of independent $L^2$ random variables 
$Y_k$ is an example of a martingale because independence implies ${\rm E}[X_{n+1} | X_1,\cdots,X_n] = {\rm E}[X_{n+1}]$
which is ${\rm E}[X_n]$ by the identical distribution assumption. 
If $X$ and $M$ are two discrete time stochastic processes, define the {\bf martingale transform} (=discrete Ito integral)
$X \cdot M$ as the process $(X \cdot M)_n  = \sum_{k=1}^{n} X_k(M_k-M_{k-1})$. 
If the process $X$ is {\bf bounded} meaning that there exists a constant $C$ such that ${\rm E}[|X_k|] \leq C$ for all $k$, then
if $M$ is a martingale, also $X \cdot M$ is a martingale. The {\bf Doob martingale convergence theorem} is
\index{Doob martingale convergence}
\index{bounded stochastic process}
\index{stochastic process}
\index{discrete time stochastic process}

\satz{For a bounded super martingale $X$, then $X_n$ converges in $L^1$. }

The convergence theorem can be used to prove the {\bf optimal stopping time theorem} which tells that the expected value
of a {\bf stopping time} is the initial expected value. In finance it is known as the {\bf fundamental theorem of 
asset pricing}. If $\tau$ is a stopping time adapted to a martingale $X_k$, it defines the random variable $X_{\tau}$
and ${\rm E}[X_\tau]={\rm E}[X_0]$. For a super-martingale one has $\geq$ and for a sub-martingale $\leq$. 
The proof is obtained by defining the {\bf stopped process} 
$X^{\tau}_n = X_0 + \sum_{k=0}^{{\rm min}(\tau,n)-1} (X_{k+1}-X_k)$ which is a martingale transform and so a martingale.
The martingale convergence theorem gives a limiting random variable $X_{\tau}$ and because ${\rm E}[X_n^{\tau}]={\rm E}[X_0]$
for all $n$,  ${\rm E}[X_\tau]={\rm E}[X_0]$. This is rephrased as ``you can not beat the system" \cite{WilliamsMartingales}.
A trivial implication is that one can not for example design a strategy allowing to win in a fair game by designing a ``clever
stopping time" like betting on ``red" in roulette if 6 times ``black" in a row has occurred. Or to follow the strategy
to stop the game, if one has a first positive total win, which one can always do by doubling the bet in case of 
losing a game. Martingales were introduced by Paul L\'evy in 1934, the name ``martingale" (referring to the just mentioned 
doubling betting strategy) was added in a 1939 probability book of Jean Ville. 
The theory was developed by Joseph Leo Doob in his book of 1953. 
\cite{DoobStochastic}. See \cite{WilliamsMartingales}. 

\index{stopping time}
\index{asset pricing}
\index{bounded martingale}
\index{optimal stopping time}
\index{conditional expectation}
\index{martingale}
\index{discrete Ito integral}
\index{martingale transform}

\section{Theorema Egregium}

A Riemannian metric on a two-dimensional manifold $S$ defines the
quadratic form $I = E du^2+2Fdudv+G dv^2$ called {\bf first fundamental form}
on the surface. If $r(u,v)$ is a parameterization of $S$, then
$E=r_u \cdot r_u, F=r_u \cdot r_v$ and $G=r_v \cdot r_v$.
The {\bf second fundamental form} of $S$ is 
$II=L du^2 + 2 M du dv + N dv^2$,
where $L=r_{uu} \cdot n, M = r_{uv} \cdot n, N = r_{vv} \cdot n$, 
written using the normal vector $n=(r_u \times r_v)/|r_u \times r_v|$. 
The {\bf Gaussian curvature}
$K={\rm det}(II)/{\rm det}(I) = (LN-M^2)/(EG-F^2)$. 
depends on the embedding $r: R \to S$ in space $\mathbb{R}^3$, but
it actually only depends on the intrinsic metric, 
the first fundamental form. This is the {\bf Theorema egregium} of Gauss:
\index{First fundamental form}
\index{Theorema egregium}
\index{second fundamental form}

\satz{ The Gaussian curvature only depends on the Riemannian metric.}

Gauss himself already gave explicit formulas, but a 
formula of {\bf Brioschi} gives the curvature $K$ explicitly as a ratio of
determinants involving $E,F,G$ as well as and first and second 
derivatives of them. 
In the case when the surface is given as a graph $z=f(x,y)$, one can give 
$K= D/(1+|\nabla f|^2)^2$, where $D=(f_{xx} f_{yy} - f_{xy}^2)$ is 
the {\bf discriminant} and $(1+|\nabla f|^2)^2= {\rm det}(II)$. 
If the surface is rotated in space so that $(u,v)$ is a
critical point for $f$, then the {\bf discriminant} $D$ is equal to the curvature.
One can see the independence of the embedding also 
from the {\bf Puiseux formula} 
$K = 3 (|S_0(r)|-S(r))/(\pi r^3)$, where $|S_0(r)|=2\pi r$ is the 
circumference of the circle $S_0(r)$ in the flat case and 
$|S(r)|$ is the circumference of the {\bf geodesic circle} 
of radius $r$ on $S$. The theorem Egregium also follows from
Gauss-Bonnet as the later allows to write
the curvature in terms of the angle sum of a geodesic infinitesimal triangle 
with the angle sum $\pi$ of a flat triangle. As the 
angle sums are entirely defined intrinsically, the curvature is intrinsic. 
The ``Theorema Egregium" was found by Karl-Friedrich Gauss in 1827 and 
published in 1828 in ``Disquisitiones generales circa superficies curvas". 
It is not an accident, that Gauss was occupied with concrete 
geodesic triangulation problems too.  
\index{Puiseux formula}
\index{Brioschi formula}
\index{discriminant}

\section{Entropy}

Given a random variable $X$ on a probability space $(\Omega,\mathcal{A},{\rm P})$
which is {\bf finite and discrete} in the sense that it takes only finitely many values,
the {\bf entropy} is defined as $S(X)  = -\sum_{x} p_x \log(p_x)$,
where $p_x = {\rm P}[ X=x ]$. To compare, for a random variable $X$ with cumulative distribution function
$F(x) = {\rm P}[X \leq x]$ having a continuous derivative $F'=f$, the entropy is defined
as $S(X) = - \int f(x) \log(f(x)) \; dx$, allowing
the value $-\infty$ if the integral does not converge. (We always read $p \log(p)=0$ if $p=0$.)
In the continuous case, one also calls this the {\bf differential entropy}. 
Two discrete random variables $X,Y$ are called
{\bf independent} if one can realize them on a product probability space $\Omega = A \times B$
so that $X(a,b) = X(a)$ and $Y(a,b) = Y(b)$ for some functions
$X:A \to \mathbb{R}, Y:B \to \mathbb{R}$. Independence implies that the random variables
are uncorrelated, ${\rm E}[X Y ] = {\rm E}[X] {\rm E}[Y]$ and that the {\bf entropy adds up}
$S(X Y) = S(X) + S(Y)$.  We can write $S(X) = {\rm E}[ \log(W(x)) ]$, where
$W$ is the ``Wahrscheinlichkeit" random variable assigning to
$\omega \in \Omega$ the value $W(\omega)=1/p_x$ if $X(\omega)=x$.
Let us say, a functional on discrete random variables is {\bf additive} if it is of
the form $H(X) = \sum_x f(p_x)$ for some continuous function $f$ for which $f(t)/t$ is monotone. 
We say it is {\bf multiplicative} if $H(X Y) = H(X) + H(Y)$ for independent
random variables. The functional is {\bf normalized} if $H(X)=\log(4)$ if $X$ is a random variable
taking two values $\{0,1\}$ with probability $p_0=p_1=1/2$. Shannon's theorem is:
\index{entropy}
\index{Wahrscheinlichkeit}
\index{independence}

\satz{Any normalized, additive and multiplicative $H$ is entropy $S$. }

The word ``entropy" was introduced by Rudolf Clausius in 1850 \cite{RovelliTime}.
Ludwig Bolzmann saw the importance of $\frac{d}{dt} S \geq 0$ in the context of heat
and wrote in 1872 $S=k_B \log(W)$, where $W(x)=1/p_x$ is the inverse ``Wahrscheinlichkeit" 
that a state has the value $x$. His equation is understood as the expectation
$S=k_B {\rm E}[\log(W)] = \sum_x p_x \log(W(x))$ which is the {\bf Shannon entropy}, 
introduced in 1948 by Claude Shannon in the context of information theory.
(Shannon characterized functionals $H$ with the property
that if $H$ is continuous in $p$, then for random variables $H_n$ with $p_x(H_n)=1/n$,
one has $H(X_n)/n \leq H(X_m)/m$ if $n \leq m$ and if $X,Y$ are two
random variables so that the finite $\sigma$-algebras $\mathcal{A}$ defined by $X$
is a sub-$\sigma$-algebra $\mathcal{B}$ defined by $Y$,
then $H(Y) = H(X) + \sum_x p_x H(Y_x)$, where $Y_x(\omega)=Y(\omega)$
for $\omega \in \{ X =x \}$. One can show that these Shannon conditions are 
equivalent to the combination of being additive and multiplicative.
In statistical thermodynamics, where $p_x$ is
the probability of a {\bf micro-state}, then $k_B S$ is also called the {\bf Gibbs entropy},
where $k_B$ is the {\bf Boltzmann constant}. For general
random variables $X$ on $(\Omega,\mathcal{A},{\rm P})$ and a finite $\sigma$-sub-algebra $\mathcal{B}$, 
Gibbs looked in 1902 at {\bf course grained entropy}, which is the entropy of the conditional
expectation $Y={\rm E}[X|\mathcal{B}|$, which is now a random variable $Y$ 
taking only finitely many values so that entropy is defined. 
See \cite{Shannon48}.
\index{Boltzmann constant}
\index{Shannon entropy}

\section{Mountain Pass}

Let $H$ be a {\bf Hilbert space}, and let $f$ be a twice Fr\'echet differentiable
function from $H$ to $\mathbb{R}$. The {\bf Fr\'echet derivative} $A=f'$ at a point $x \in H$
is a linear operator $A$ satisfying $f(x+h)-f(x)- A h = o(h)$ for all $h \to 0$.
A point $x \in H$ is called a {\bf critical point} of $f$ if $f'(x)=0$.
The functional satisfies the {\bf Palais-Smale condition}, if every sequence $x_k$ in $H$
for which $\{ f(x_k) \}$ is bounded and $f'(x_k) \to 0$, has a convergent subsequence in 
the closure of $\{ x_k \}_{k \in \mathbb{N}}$. 
A pair of points $a,b \in H$ defines a {\bf mountain pass}, if there exist $\epsilon>0$
and $r>0$ such that $f(x) \geq f(a)+\epsilon$ on $S_r(a) = \{ x \in H \; | \; ||x-a||=r \}$,
$f$ is not constant on $S_r(a)$ and $f(b) \leq f(a)$. 
A critical point is called a {\bf saddle} if it is neither a maximum nor a minimum of $f$.
\index{Palais-Smale condition}
\index{critical point}
\index{saddle point}
\index{mountain pass} 
\index{Fr\'echet derivative}

\satz{If a Palais-Smale $f$ has a mountain pass, it features a saddle.}

The idea is to look at all continuous paths $\gamma$ from $a$ to $b$
parametrized by $t \in [0,1]$.
For each path $\gamma$, the value $c_{\gamma} = f(\gamma(t))$ has to be maximal
for some time $t \in [0,1]$. The infimum over all these critical values $c_{\gamma}$
is a critical value of $f$. The mountain pass condition leads to a
``mountain ridge" and the critical point is a ``mountain pass", hence the
name. The example $(2\exp(-x^2-y^2)-1)(x^2+y^2)$ with $a=(0,0),b=(1,0)$ shows
that the non-constant condition is necessary for a saddle point on 
$S_{r}(a)$ with $r=1/2$. The reason for sticking with a Hilbert space is that
it is easier to realize the compactness condition due to weak star compactness
of the unit ball. But it is possible to weaken the conditions and work with 
a Banach manifolds $X$ continuous G\^ateaux derivatives: $f': X \to X^*$ if
$X$ has the strong and $X^*$ the weak-$*$ topology.
It is difficult to pinpoint historically the first use of the mountain 
pass principle as it must have been known intuitively since antiquity.
The crucial Palais-Smale {\bf compactness condition} which makes the theorem work 
in infinite dimensions appeared in 1964. \cite{AubinEkeland} calls it condition (C),
a notion which already appeared in the original paper \cite{PalaisSmale}.

\index{mountain pass}
\index{compactness condition}

\section{Exponential sums}

Given a smooth function $f: \mathbb{R} \to \mathbb{R}$ which maps integers to integers,
one can look at {\bf exponential sums} $\sum_{x=a}^b \exp(i \pi f(x))$
An example is the {\bf Gaussian sum} $\sum_{x=0}^{n-1} \exp(i \alpha x^2)$. 
There are lots of interesting relations and estimates. One of the
magical formulas is the {\bf Landsberg-Schaar relations} for the 
finite sums $S(q,p) = \frac{1}{\sqrt{p}} \sum_{x=0}^{p-1} \exp(i \pi x^2 q/p)$.
\index{Landsberg-Schaar relation}
\index{Gauss sums}

\satz{If $p,q$ are positive and odd integers, then $S(2q,p) = e^{i \pi/4}  S(-p,2q) $.}

One has $S(1,p) = (1/\sqrt{p}) \sum_{x=0}^{p-1} \exp(i x^2/p) = 1$ for all positive integers $p$
and $S(2,p) = (e^{i \pi/4}/\sqrt{p}) \sum_{x=0}^{p-1} \exp(2 i x^2/p) = 1$ if $p=4k+1$
and $i$ if $p=4k-1$. The method of exponential sums has been 
expanded especially by Vinogradov's papers  \cite{Vinogradov} and used for 
number theory like for quadratic reciprocity \cite{MurtyPacelli}.
The topic is of interest also outside of number theory. Like in dynamical systems
theory as F\"urstenberg has demonstrated. An ergodic theorist would look 
at the dynamical system $T(x,y) = (x+2y+1,y+1)$ on the 2-torus
$\mathbb{T}^2=\mathbb{R}^2/(\pi \mathbb{Z})^2$ and define $g_{\alpha} (x,y)=\exp(i \pi x \alpha)$.
Since the orbit of this toral map is $T^n(1,1) = (n^2,n)$, the exponential sum 
can be written as a {\bf Birkhoff sum} $\sum_{k=0}^{p-1} g_{q/p}(T^k(1,1))$ which 
is a particular orbit of a dynamical system. Results
as those mentioned above show that the random walk grows like $\sqrt{p}$, similarly as in a 
random setting. Now, since the dynamical system is minimal, the growth rate should not
depend on the initial point and $\pi q/p$ should be replaceable by any irrational $\alpha$
and no more be linked to the length of the orbit. The problem is then to study 
the growth rate of the {\bf stochastic process} $S^t(x,y) = \sum_{k=0}^{t-1} g(T^k(x,y))$ 
(= sequence of random variables) for any 
continuous $g$ with zero expectation which by Fourier boils down to look at exponential sums. 
Of course $S^t(x,y)/t \to 0$ by Birkhoff's ergodic theorem, but as in the law of iterated logarithm
one is interested in precise growth rates. This can be subtle. 
Already in the simpler case of an integrable $T(x)=x+\alpha$ on the 
$1$-torus, there is Denjoy-Koskma theory which shows that the growth rate depends
on Diophantine properties of $\pi \alpha$. Unlike for irrational rotations, the F\"urstenberg type 
skew systems $T$ leading to the theta functions are not integrable: it is not conjugated to a group
translation (there is some randomness, even-so weak as Kolmogorov-Sinai entropy is zero).
The dichotomy between structure and randomness and especially
the similarities between dynamical and number theoretical set-ups has been discussed in
\cite{TaoStructureRandomness}. 

\section{Sphere theorem}

A compact {\bf Riemannian manifold} $M$ is said to have
{\bf positive curvature}, if all {\bf sectional curvatures}
are positive. The {\bf sectional curvature} at a point $x \in M$
in the direction of the 2-dimensional plane $\Sigma \subset T_xM$
is defined as the Gaussian curvature of the surface
$\exp_x(\Sigma) \subset M$ at the point. In terms of the
{\bf Riemannian curvature tensor} $R: T_xM^4 \to \mathbb{R}$
and an orthonormal basis $\{u,v\}$ spanning $\Sigma$,
this is $R(u,v,u,v)$. The curvature is called
{\bf quarter pinched}, if it the sectional curvature is
in the interval $(1,4]$ at all points $x \in M$. In particular, a quarter
pinched manifold is a manifold with positive curvature.
We say here, a compact Riemannian manifold {\bf is a sphere} if it is homeomorphic
to a sphere. The {\bf sphere theorem} is:
\index{sphere theorem}
\index{sectional curvature}
\index{quarter pinched}

\satz{A simply-connected quarter pinched manifold is a sphere}

The theorem was proven by Marcel Berger and Wilhelm Klingenberg
in 1960. That a pinching condition would imply a manifold
to be a sphere had been conjectured already by Heinz Hopf.
Hopf himself proved in 1926 that constant sectional curvature
implies that $M$ is even isometric to a sphere.
Harry Rauch, after visiting Hopf in Z\"urich in the 1940's
proved that a 3/4-pinched simply connected manifold is a sphere.
In 2007, Simon Brendle and Richard Schoen proved that
the theorem even holds if the statement {\bf $M$ is a d-sphere}
(meaning that $M$ is diffeomorphic to the Euclidean d-sphere $\{ |x|^2=1 \} \subset \mathbb{R}^{d+1}$).
This is the {\bf differentiable sphere theorem}. Since John Milnor had given
in 1956 examples of spheres which are homeomorphic but not diffeomorphic
to the standard sphere (so called {\bf exotic spheres}, spheres which carry
a smooth maximal atlas different from the standard one), the 
differentiable sphere theorem is a substantial improvement
on the topological sphere theorem. It needed completely 
new techniques, especially the {\bf Ricci flow} $\dot{g} = - 2 {\rm Ric}(g)$ 
of Richard Hamilton which is a weakly parabolic 
partial differential equation deforming
the metric $g$ and uses the {\bf Ricci curvature} ${\rm Ric}$ of $g$. 
See \cite{BergerPanorama,BrendleRicci}.
\index{differentiable sphere theorem}
\index{exotic sphere}
\index{Ricci flow}

\section{Word problem}

The {\bf word problem} in a {\bf finitely presented group} $G=(g|r)$
with {\bf generators} $g$ and {\bf relations} $r$ is the problem to decide,
whether a given set of two words $v,w$ represent the same group element in $G$
or not. The word problem is not solvable in general. There are concrete finitely
presented groups in which it is not. 
The following theorem of Boone and Higman relates the solvability to
algebra. A group is {\bf simple} if its only {\bf normal subgroup} is
either the trivial group or then the group itself.
\index{simple group}
\index{normal subgroup}
\index{finitely presented group}

\satz{Finitely presented simple groups have a solvable word problem.}

More generally, if $G \subset H \subset K$ where $H$ is simple and $K$
is finitely presented, then $G$ has a solvable word problem.
Max Dehn proposed the word problem in 1911.
Pyotr Novikov in 1955 proved that the word problem is undecidable for
finitely presented groups.
William W. Boone and Graham Higman proved the theorem in 1974 \cite{BooneHigman}.
Higman would in the same year also find an example of an infinite finitely
presented simple group.
The non-solvability of the word problem implies the non-solvability of the
homeomorphism problem for $n$-manifolds with $n \geq 4$. See \cite{BooneCannonito}.

\section{Finite simple groups}

A {\bf finite group} $(G,*,1)$ is a finite set $G$ with an operation $*: G \times G \to G$ and
{\bf $1$ element}, such that the operation is {\bf associative} $(a*b)*c=a*(b*c)$, for all $a,b,c$,
such that $a*1=1*a=a$ for every $a$ and such that every $a$ has an inverse $a^{-1}$ satisfying
$a*a^{-1}=1$. A group $G$ is {\bf simple} if the only {\bf normal subgroups} of $G$
are the {\bf trivial group} $\{1\}$ or the group itself.
A subgroup $H$ of $G$ is called {\bf normal} if $gH=Hg$ for all $g$.
Simple groups play the role of the primes in the set of integers.
A {\bf theorem of Jordan-H\"older} is that the composition series of $G$ (with simple groups
as quotients) is unique up to permutations and isomorphisms. 
The {\bf classification theorem of finite simple groups} is
\index{simple group}
\index{normal group}
\index{classification of finite simple groups}

\satz{Every finite simple group is cyclic, alternating, Lie or sporadic.}

There are 18 so called {\bf regular families} of finite simple groups made of {\bf cyclic}, 
{\bf alternating} and 16 {\bf Lie type} groups.
Then there are 26 so called {\bf sporadic groups}, in which 20 are {\bf happy groups} as they
are subgroups or sub-quotients of the {\bf monster} and 6 are {\bf pariahs}, outcasts which
are not under the spell of the monster. 
The classification was a huge collaborative effort with more than 100 authors, covering 500 journal
articles. According to Daniel Gorenstein, the classification was completed in 1981
and fixes were applied until 2004. (Michael Aschbacher and Stephen Smith resolved the last
problems which lasted several years leading to a full
proof of 1300 pages.) A second generation cleaned-out proof written with more details
is under way and currently has 5000 pages. Some history is given in \cite{SolomonFiniteGroups}.
\index{monster}
\index{pariah}
\index{happy groups}

\section{God number}

Given a finite finitely presented group $G=(g|r)$ like
for example the Rubik group. It defines the {\bf Cayley graph} $\Gamma$
in which the group elements are the nodes and
where two nodes $a,b$ are connected if there is a generator
$x$ in in $g$ such that $x a = b$. The {\bf diameter}
of a graph is the largest geodesic distance between two
nodes in $\Gamma$. It is also called the {\bf God number} of the puzzle.
The {\bf Rubik cube} is an example of a finitely
presented group. The original $3 \times 3 \times 3$ cube allows
to permute the 26 boundary cubes using the 18 possible rotations
of the 6 faces as generators. From the $X=8! 12! 3^8 2^{12}$ possible
ways to physically build the cube, only $|G|=X/12= 43252003274489856000$
are present in the Rubik group $G$. Some of the positions ``quarks" \cite{Golomb}
can not be realized but combinations of them ``mesons" or ``baryons" can.

\satz{The God number of the Rubik cube is 20.}

This means that from any position, one could, in principle
solve the puzzle in 20 moves. Note that one has to specify clearly
the generators of the group as this defines the Cayley
graph and so a metric on the group. The lower bound $18$ had already been known in 1980
because a counting of all the possible moves with $17$ steps produced less elements.
The lower bound 20 came in 1995 when Michael Reid proved that the {\bf super-flip position}
(where the edges are all flipped but corners are correct) needs 20 moves.
In July 2010, using about 35 CPU years, a team around Tomas Rokicki established
that the God number is 20. They partitioned the possible group positions into roughly 2 billion sets of 20
billions positions each. Using symmetry they reduced it to 55 million positions,
then found solutions for any of the positions in these sets. \cite{GodNumber}
It appears silly to put a God number computation as a fundamental theorem, but the
status of the Rubik cube is enormous as it has been one of the most popular puzzles for decades and
is a {\bf prototype} for many other similar puzzles, the choice can be defended.
\footnote{I presented the God number problem in the 80ies as an undergraduate in a logic seminar
of Ernst Specker and the choice of topic had been objected to by Specker himself as a too
``narrow problem". But the Rubik cube and its group properties have 
``cult status". The object was one of the triggers for me to study math. }
One can ask to compute the God number of any finitely presented finite group.
Interesting in general is the complexity of evaluating that functional.
The simplest nontrivial {\bf Rubik cuboid} is the $2 \times 2 \times 1$ one. It has
$6$ positions and $2$ generators $a,b$. The finitely presented group is $\{ a,b | a^2=b^2=(ab)^3=1 \}$
which is the {\bf dihedral group} $D_3$. Its group elements are
$G=\{ 1,a=babab,ab=baba,aba=bab,abab=ba,ababa=b \}$. The group is isomorphic to
the {\bf symmetry group of the equilateral triangle}, generated by the two reflections $a,b$ 
at two altitude lines. The God number of that group is 3 because
the Cayley graph $\Gamma$ is the cyclic graph $C_6$. The puzzle solver has
here ``no other choice than solving the puzzle", because one is forced to make non-trivial move in each step. 
See \cite{Joyner} or \cite{Baumslag} for general combinatorial group theory
and \cite{CubedRubik} for a recent auto biography of Erno Rubik. 
\index{God number}
\index{Rubik cube}
\index{symemtry groyup}
\index{Rubik cuboid}

\section{Sard Theorem}

Let $f: M \to N$ be a smooth map between smooth manifolds $M,N$
of dimension ${\rm dim}(M)=m$ and ${\rm dim}(N)=n$. A
point $x \in M$ is called a {\bf critical point} of $f$,
if the Jacobian $n \times m$ matrix $df(x)$ has rank both
smaller than $m$ and $n$. If $C$ is the set of critical points, then
$f(C) \subset N$ is called the {\bf critical set} of $f$. The
{\bf volume measure} on $N$ is a choice of a volume form, obtained
for example after introducing a Riemannian metric.
{\bf Sard}'s theorem is
\index{Sard theorem}
\index{critical set}
\index{critical points} 

\satz{The critical set of $f:M \to N$ has zero volume measure in $N$.}

The theorem applied to smooth map $f: M \to \mathbb{R}$
tells that for almost all $c$, the set $f^{-1}(c)$ is a smooth hypersurface
of $M$ or then empty. The later can happen if $f$ is constant.
We assumed $C^{\infty}$ but one can relax the smoothness assumption of $f$.
If $n \geq m$, then $f$ needs only to be continuously differentiable.
If $n < m$, then $f$ needs to be in $C^{m-n+1}$.
The case when $N$ is one-dimensional has been covered by Antony Morse
(who is unrelated to Marston Morse) in 1939 and by Arthur Sard in general
in 1942. A bit confusing is that Marston Morse (not Antony) covered the case $m=1,2,3$
and Sard in the case $m=4,5,6$ in unpublished papers before as mentioned in
a footnote to \cite{Sard42}. Sard also notes already that examples of
Hassler Whitney show that the smoothness condition can not be relaxed.
Sard formulated the results for $M=\mathbb{R}^m$ and $N=\mathbb{R}^n$
(by the way with the same choice $f: M \to N$ as done here and not as
in many other places). The manifold case appears for example in
\cite{Sternberg1964}.

\section{Elliptic curves}

An {\bf elliptic curve} is a plane algebraic curve defined
by the points satisfying the {\bf Weierstrass equation} $y^2=x^3+ax+b=f(x)$.
One assumes the curve to be {\bf non-singular}, meaning that
the {\bf discriminant} $\Delta = -16(4a^3+27b^2)$
is not zero. This assures that there are no cusps nor
multiple roots for the simple reason that 
the explicit solution formulas for roots of 
$f(x)=0$ involves only square roots of $\Delta$.
A curve is an {\bf Abelian variety}, if it carries
an Abelian algebraic group structure, meaning that the 
addition of a point defines a morphism of the variety. 
\index{elliptic curve}
\index{non-singular curve}
\index{Weierstrass equation}

\satz{Elliptic curves are Abelian varieties.}

The theorem seems first have been realized by Henri Poincar\'e
in 1901. Weierstrass before had used the Weierstrass $\mathcal{P}$
function earlier in the case of elliptic curves over the complex plane.
To define the group multiplication, one uses the {\bf chord-tangent
construction}: first add point $O$ called the {\bf point at infinity}
which serves as the {\bf zero} in the group. Then define
$-P$ as the point obtained by reflecting at the $x$-axes.
The {\bf group multiplication} between two different points
$P,Q$ on the curve is defined to be $-R$ if $R$ is the
point of intersection of the line through
$P,Q$ with the curve. If $P=Q$, then $R$ is defined to be the intersection
of the tangent with the curve. If there is no intersection,
that is if $P=Q$ is an inflection point, then one defines $P+P=-P$.  Finally,
define $P+O=O+P=P$ and $P+(-P)=0$. This recipe can be explicitly given in coordinates
allowing to define the multiplication in any field of characteristic
different from $2$ or $3$. The group structure on elliptic
curves over finite fields provides a rich source of
{\bf finite Abelian groups} which can be used for
cryptological purposes, the so called {\bf elliptic curve cryptograph} ECC.
Any procedure, like public key, Diffie-Hellman or factorization attacks
on integers can be done using groups given by elliptic curves.
\cite{FultonAlgebraicCurves}.
\index{Elliptic curve cryptography}
\index{chord tangent construction}

\section{Billiards}

Billiards are the geodesic flow on a smooth compact $n$-manifold $M$ with boundary. The
dynamics is extended through the boundary by applying the law of reflection.
While the flow of the geodesic $X^t$ is Hamiltonian on the unit tangent bundle $SM$,
the billiard flow is only piecewise smooth and also the return map to the boundary
is not continuous in general but it is a map preserving a natural volume so that one can
look at ergodic theory. Already difficult are flat 2-manifolds $M$ homeomorphic
to a disc having convex boundary homeomorphic to a circle.
For smooth convex tables this leads to a return map $T$ on the annulus
$X=\mathbb{T} \times [-1,1]$ which is $C^{r-1}$ smooth if the boundary is $C^{r}$ \cite{Dou82}.
It defines a {\bf monotone twist map}: in the sense that it preserves the boundary,
is area and orientation preserving and satisfies the {\bf twist condition} that $y \to T(x,y)$
is strictly monotone. A {\bf Bunimovich stadium} is the 2-manifold with boundary 
obtained by taking the convex hull of two discs of equal radius in $\mathbb{R}$ 
with different center. The billiard map is called {\bf chaotic}, if it is ergodic 
and the {\bf Kolmogorov-Sinai entropy} 
is positive. By Pesin theory, this metric entropy is the {\bf Lyapunov exponent} which is
the exponential growth rate of the Jacobian $dT^n$ (and constant almost everywhere due to 
ergodicity). There are coordinates in the tangent bundle of the annulus $X$ in which $dT$ is
the composition of a horizontal shear with strength $L(x,y)$, where $L$ is the trajectory
length before the impact with a vertical shear with strength $-2 \kappa/\sin(\theta)$ 
where $\kappa(x)$ is the curvature of the curve at the impact $x$ and $y=\cos(\theta)$, with 
{\bf impact angle} $\theta \in [0,\pi]$ between the tangent and the trajectory. 
\index{billiards}
\index{Lyapunov exponent}
\index{entropy}

\satz{The Bunimovich stadium billiard is chaotic.}

Jacques Hadmard in 1898 and Emile Artin in 1924 already 
looked at the geodesic flow on a surface of constant negative curvature. 
Yakov Sinai constructed in 1970 the first chaotic billiards, the Lorentz gas or 
Sinai billiard. An example, where Sinai's result applies is the hypocycloid $x^{1/3} + y^{1/3}=1$.  
The Bernoulli property was established by Giovanni Gallavotti and Donald Ornstein in 1974. 
In 1973, Vladimir Lazutkin proved that a generic smooth convex two-dimensional billiard 
can not be ergodic due to the presence of KAM whisper galleries using Moser's twist map
theorem. These galleries are absent in the presence of flat points (by a theorem of John Mather) 
or points, where the curvature is unbounded (by a theorem of Andrea Hubacher \cite{Hubacher87}). 
Leonid Bunimovich \cite{Bun79} constructed in 1979 the first convex chaotic billiard.
No smooth convex billiard table with positive Kolmogorov-Sinai entropy
is known. A candidate is the real analytic $x^4+y^4=1$. Various generalizations
have been considered like in \cite{Woj86}.
A detailed proof that the Bunimovich stadium is measure theoretically
conjugated to a Bernoulli system (the shift on a product space) 
is surprisingly difficult: one has
to show positive Lyapunov exponents on a set of positive measure. Applying
Pesin theory with singularities (Katok-Strelcyn theory \cite{KatokStrelcyn}) 
gives a Markov process.  One needs then to establish ergodicity using a method of Eberhard Hopf 
of 1936 which requires to understand stable and unstable manifolds 
\cite{ChernovMarkarian}.  
See \cite{Tab95, KozlovTreshchev, MoserVariations, Gole, KH, ChernovMarkarian} for sources on billiards.
\index{Hopf method}
\index{Whisper galleries}

\section{Uniformization}

A {\bf Riemann surface} is a one-dimensional {\bf complex manifold}. 
This means is is a connected two-dimensional real manifold so that the transition 
functions of the atlas are holomorphic mappings of the complex plane. 
It is {\bf simply connected} if its fundamental group is trivial (equivalently, its genus
$b_1$ is zero). Two Riemann surfaces are {\bf conformally equivalent} or simply
{\bf equivalent} if they are equivalent as complex manifolds, that is if there is a bijective morphism $f$
between them. A map $f:S \to S'$ is holmorphic if for every choice of coordinates $\phi:S \to \mathbb{C}$
and $\psi':S' \to \mathbb{C}$, the maps $\phi' \circ f \circ \phi^{-1}$ are holomorphic.
The {\bf curvature} is the Gaussian curvature of the surface. The
{\bf uniformization theorem} is:
\index{Riemann surface}
\index{complex manifold}
\index{conformally equivalent}
\index{uniformization theorem}

\satz{A Riemann surface is equivalent to one with constant curvature.}

This is a ``geometrization statement" and means that 
the universal cover of every Riemann surface is conformally equivalent to 
either a {\bf Riemann sphere} (positive curvature), 
a  {\bf complex plane} (zero curvature) or a {\bf unit disk} (negative curvature).
It implies that any region $G \subset \mathbb{C}$ whose complement contains two
or more points has a universal cover which is the disk. It especially
implies the {\bf Riemann mapping theorem} assuring that any region $U$ homeomorphic to
a disk is conformally equivalent to the unit disk (see \cite{Carlson}). 
For a detailed treatment of compact Riemann surfaces, see \cite{GirondoGonzalezDiez}.
\index{Riemann mapping theorem}
It also follows that all {\bf Riemann surfaces} (without restriction of genus) can be obtained as
quotients of these three spaces: for the sphere one does not have to take any quotient,
the genus 1 surfaces = {\bf elliptic curves} can be obtained as quotients of the complex
plane and any genus $g>1$ surface can be obtained as quotients of the unit disk.
Since every closed $2$-dimensional orientable surface is characterized by their
genus $g$, the uniformization theorem implies that any such surface admits a metric
of constant curvature. Teichm\"uller theory parametrizes the possible metrics,
and there are $3g-3$ dimensional parameters for $g \geq 2$, whereas for $g=0$
there is one and for $g=1$ a moduli space $\mathbb{H}/SL_2(\mathbb{Z})$. 
In higher dimensions, closest to the uniformization theorem is
the {\bf Killing-Hopf theorem} telling that every connected 
complete Riemannian manifold of {\bf constant sectional curvature} 
and dimension $n$ is isometric to the quotient of a sphere $\mathbb{S}^n$, 
Euclidean space $\mathbb{R}^n$ or Hyperbolic $n$-space $\mathbb{H}^n$ restating
that constant curvature geometry is either elliptic, parabolic=Euclidean
or yyperbolic geometry. Complex analysis has rich applications in complex dynamics 
\cite{Beardon,MilnorNotes,Carlson} and relates to much more geometry \cite{McMullen2014}.

\section{Control Theory}

A {\bf Kalman filter} is an optional estimates algorithm of a
linear dynamic system from a series of possibly noisy measurements.
The idea is similar as in a {\bf dynamic Bayesian network} or
{\bf hidden Markov model}. The filter applies both to {\bf differential equations}
$\dot{x}(t) = A x(t) + B u(t) + G z(t)$ as well as {\bf discrete dynamical system}
$x(t+1)=A x(t) + B u(t) + G z(t)$, where $u(t)$ is
{\bf external input} and $z(t)$ {\bf input noise} given by independent identically distributed
usually Gaussian {\bf random variables}. Kalman calls this a {\bf Wiener problem}.
One does not see the {\bf state} $x(t)$ of the system but some {\bf output}
$y(t) = C x(t) + D u(t)$. The filter then ``filters out" or ``learns" the best estimate $x^*(t)$
from the observed data $y(t)$.  The linear space $X$ is defined as
the vector space spanned by the already observed vectors. The optimal solution is given
by a sophisticated dynamical data fitting.
\index{discrete dynamical system}
\index{Wiener problem}

\satz{The optimal estimate $x^*$ is the projection of $y$ onto $X$.}

This formulation is the informal 1-sentence description which can be found already in 
Kalman's article. Kalman then gives 
explicit formulas which generate from the {\bf stochastic difference
equation} a concrete {\bf deterministic linear system}. 
For a modern exposition, see \cite{MarchthalerDinglerKalman}.
The {\bf Kalman filter} is named after Rudolf Kalman who wrote
\cite{Kalman} in 1960. Kalman's paper is one of the most
cited papers in applied mathematics. The ideas were used both in the
Apollo and Space Shuttle program.
Similar ideas have been introduced in statistics by the Danish
astronomer Thorvald Thiele and the radar theoretician Peter Swerling.
There are also nonlinear version of the Kalman filter which is
used in nonlinear state estimation like navigation systems and GPS.
The nonlinear version uses a multi-variate Taylor series expansion
to linearise about a working point. See \cite{EubankKalman,MarchthalerDinglerKalman}.
\index{Kalman filter}
\index{stochastic difference equation}

\section{Zariski main theorem}

A {\bf variety} is called {\bf normal} if it can be covered by open affine varieties
whose rings of functions are normal. A commutative ring is called {\bf normal} if
it has no non-zero nilpotent elements and is integrally closed in its complete ring
of fractions. For a curve, a one-dimensional variety, normality is equivalent to
being non-singular but in higher dimensions, a normal variety still can have singularities.
The normal complex variety is called {\bf unibranch at a point $x \in X$}
if there are arbitrary small neighborhoods $U$ of $x$ such that the set of non-singular points of $U$
is connected. {\bf Zariski's main theorem} can be stated as:

\satz{Any closed point of a normal complex variety is unibranch.}

Oscar Zariski proved the theorem in 1943. To cite \cite{Mumford2013}, {\it ``it was the final result
in a foundational analysis of  birational maps between varieties. The 'main Theorem' asserts in a strong sense
that the normalization (the integral closure) of a variety X is the maximal variety $X'$
birational over X, such that the fibres of the map $X' \to X$ are finite. A generalization
of  this fact became Alexandre Grothendieck's concept of  the 'Stein factorization' of  a map.}
The result has been generalized to schemes $X$, which is called {\bf unibranch} at
a point $x$ if the local ring at $x$ is unibranch.
A generalization is the {\bf Zariski connectedness theorem} from 1957:
if $f:X \to Y$ is a birational projective morphism between Noetherian integral schemes,
then the inverse image of every normal point of $Y$ is connected. Put more colloquially,
the fibres of a birational morphism from a projective variety $X$ to a normal variety $Y$
are connected. It implies that a birational morphism $f:X \to Y$ of algebraic varieties $X,Y$
is an open embedding into a neighbourhood of a normal point
$y$ if $f^{-1}(y)$ is a finite set. Especially, a birational morphism between
normal varieties which is bijective near points is an isomorphism.  \cite{hartshorne,Mumford2013}

\section{Poincar\'e's last theorem}

A {\bf homeomorphism} $T$ of an annulus $X=\mathbb{T} \times [0,1]$ is called
{\bf measure preserving} if it preserves the Lebesgue (area) measure
and preserves the orientation of $X$. 
As a homeomorphism it induces also homeomorphisms on each of the two boundary 
circles. It is called {\bf twist homeomorphism}, if it rotates the boundaries 
in different directions. 
\index{twist homeomorphism}
\index{measure preservation}
\index{area preservation}

\satz{A twist map on an annulus has at least two fixed points.}

This is called the {\bf Poincar\'e-Birkhoff theorem} or Poincar\'e's last theorem.
It was stated by Henri Poincar\'e in 1912 in the context of the {\bf three
body problem}. Poincar\'e already gave an index argument for
the existence of one fixed point gives a second. 
The existence of the first was proven by George Birkhoff in 1913 and
in 1925, where Birkhoff added the precise argument for the existence of the second. 
The twist condition is necessary because the rotation of the 
annulus $(r,\theta) \to (r,\theta+1)$ has no fixed point. 
Also area-preservation is necessary as the example $(r,\theta) \to (r(2-r),\theta+2r-1)$ 
shows.  \cite{Bir25,BrNe75}
\index{three body problem}
\index{twist map}
\index{Poincar\'e's last theorem}

\section{Geometrization}

A {\bf closed manifold} $M$ is a smooth compact manifold without boundary.
A closed manifold is {\bf simply connected} if it is connected and the fundamental group is trivial
meaning that every closed loop in $M$ can be pulled together to a point
within $M$: (if $r: \mathbb{T} \to M$ is a
parametrization of a closed path in $M$, then there exists a continuous
map $R: \mathbb{T} \times [0,1] \to M$ such that $R(0,t)=r(t)$ and $R(1,t)=r(0)$.)
We say that {\bf M is 3-sphere} if $M$ is homeomorphic to the $3$-dimensional unit sphere
$\{ (x_1,x_2,x_3,x_4) \in \mathbb{R}^4 \; | \; x_1^2+x_2^2+x_3^2+x_3^2=1 \}$.
\index{Poincar\'e conjecture}
\index{closed manifold}
\index{3-sphere}
\index{unit sphere}
\index{simply connected}

\satz{A closed simply connected 3-manifold is a 3-sphere.}

Henri Poincar\'e conjectured this in 1904. It remained the {\bf Poincar\'e conjecture}
until its proof by Grigori Perelman in 2006  \cite{Morgan2007}.
In higher dimensions, the statement was known
as the {\bf generalized Poincar\'e conjecture}, the case $n>4$ had been proven by
Stephen Smale in 1961 and the case $n=4$ by Michael Freedman in 1982.
A {\bf $d$-homotopy sphere} is a closed d-manifold that is homotopic to a $d$-sphere.
(A manifold $M$ is {\bf homotopic} to a manifold $N$ if there exists a continuous
map $f: M \to N$ and a continuous map $g: N \to M$ such that the composition
$g \circ f: M \to M$ is homotopic to the identity map on $M$ (meaning that
there exists a continuous map $F:M \times [0,1] \to M$ such that $F(x,0)=g(f(x))$ and
$F(x,1) = x$ and the map $f \circ g: N \to N$ is homotopic to the identity
on $N$.) The Poincar\'e conjecture itself, the case $d=3$,
was proven using a theory built by {\bf Richard Hamilton} who suggested to use the
{\bf Ricci flow} to solve the conjecture and more generally the
{\bf geometrization conjecture} of William Thurston: every closed 3-manifold can
be decomposed into {\bf prime manifolds} which are of 8 types, the so called
{\bf Thurston geometries} 
$S^3,E^3,H^3,S^2 \times R,H^2 \times R,\tilde{SL}(2,R),{\rm Nil},{\rm Solv}$.
If the statement {\bf M is a sphere} is replaced by 
{\bf $M$ is diffeomorphic to a sphere}, one has the 
{\bf smooth Poincar\'e conjecture}. Perelman's proof verifies this also
in dimension $d=3$. The smooth Poincar\'e conjecture is false in dimension
$d \geq 7$ as $d$-spheres then can admit non-standard smooth structures, so
called {\bf exotic spheres} constructed first by John Milnor. 
For $d=5$ it is true following result of Dennis Barden from 1964. 
It is also true for $d=6$. For $d=4$, the smooth Poincar\'e conjecture is open, and called 
``the last man standing among all great problems of classical geometric topology" \cite{FGMW}.
See \cite{MorganTian} for details on Perelman's proof.
\index{Generalized Poincare conjecture}
\index{Smooth Poincare conjecture}
\index{prime manifolds}
\index{homotopica}
\index{Thurston geometry}
\index{exotic sphere}

\section{Steinitz theorem}

A non-empty finite simple connected graph $G$ is called {\bf planar} if it can be
embedded in the plane $\mathbb{R}^2$ without self crossings. The abstract
edges of the graph are then realized as actual curves in the plane connecting
two vertices which are realized as actual points in the plane. 
The embedding of $G$ in the plane subdivides the plane now into a finite collection $F$ of
{\bf simply connected regions} called {\bf faces}. (In the two dimensional plane, a region 
is simply connected if it is homeomorphic to a disc.)
Let $v=|V|$ is the
number of vertices, $e=|E|$ the  number of edges and $f=|F|$ is the
number of faces. A planar graph is called {\bf polyhedral} if it can
be realized as a {\bf convex polyhedron}, a convex hull of finitely
many points in $\mathbb{R}^3$. A graph is called {\bf $3$-connected}, if it
remains connected also after removing one or two of its vertices. A
connected, planar 3-connected graph is also called a {\bf 3-polyhedral graph}.
The {\bf Polyhedral formula of Euler} combined with {\bf Steinitz's theorem} means:
\index{planar}
\index{polyhedral formula}
\index{faces}
\index{3-connected}
\index{Steinitz theorem}

\satz{$G$ planar $\Rightarrow$ $v-e+f=2$. Planar 3-connected $\Leftrightarrow$ polyhedral.}

The Euler polyhedron formula has first been noticed in examples
by Ren\'e Descartes \cite{Aczel} and written down in a secret notebook.
It was realized by Euler in 1750 that the formula works for general
planar graphs. Euler already gave an induction proof (also in 1752) but the first complete proof appears
have been given first by Legendre in 1794. The Steinitz theorem was proven by Ernst
Steinitz in 1922, even so he obtained the result already in 1916.
In general, a planar graph always defines a finite generalized CW complex
in which the faces are the 2-cells, the edges are the 1-cells and the
vertices are the 0-cells. The embedding in the plane defines 
then a {\bf geometric realization} of this combinatorial structure as a topological 2-sphere (as
the 2-sphere is the compactification of the plane).
The structure is not required to be achievable in the form of a convex polyhedron. And it is in general not:
take a {\bf tree graph} for example, a connected graph without triangles and without
closed loops. It is planar but it is not even $2$-connected. 
The number of vertices $v$ and the number of edges $e$ satisfy $v-e=1$. After embedding the
tree in the plane, we have exactly one face, so that $f=1$. The Euler polyhedron formula
$v-e+f=2$ is verified, but the graph is far from polyhedral.
Even in the extreme case, where $G$ is a one-point graph, the Euler formula holds: in that case there are
$v=1$ vertices, $e=0$ edges and $f=1$ faces (given by the complement of the point in the plane) 
so that still $v-e+f=2$ holds The 3-connectedness assures that the realization can be done 
using convex polyhedra. It is then even possible to have force the vertices of the polyhedron to be 
on the integer lattice points  \cite{gruenbaum,Ziegler}. In \cite{gruenbaum}, it is stated that 
the Steinitz theorem is ``the most important and deepest known result for 3-polytopes".
\index{geometric realization}
\index{Steinitz theorem}
\index{tree}
\index{secret notebook}

\section{Hilbert-Einstein action}

Let $(M,g)$ be a smooth $4$-dimensional {\bf Lorentzian manifold} which is {\bf asymptotically flat}. 
(A simplification is that the Riemannian curvature tensor $R$ is flat outside a compact subset of $M$
but this is a bit restrictive as the Schwarzschild solution below indicates.)
A Lorentzian manifold is a $4$-dimensional pseudo Riemannian
manifold of signature $(1,3)$ which in the flat case is
$dx^2+dy^2+dz^2-dt^2$. The technical condition of asymptotic flatness should imply
that the volume form $d\mu$ then has
the property that the {\bf scalar curvature} $R$ is in $L^1(M,d\mu)$ (which is the case if the non-flat
part is compact.) One can now look
at the variational problem to find extrema of the functional $g \to \int_M R d\mu$.
More generally, one can add a {\bf Lagrangian} $L$
one consider the {\bf Hilbert-Einstein functional} $\int_M  R/\kappa + L d\mu$,
where $\kappa = 8\pi G/c^4$ is the {\bf Einstein constant}.
Let $R_{ij}$ be the Ricci tensor, a symmetric tensor,
and $T_{ij}$ the energy-momentum tensor. The {\bf Einstein field equations} are
\index{Ricci tensor}
\index{Scalar curvature}
\index{Einstein constant}
\index{energy momentum tensor}

\satz{$G_{ij} =  R_{ij} - g_{ij} R/2 = \kappa T_{ij}$.}

These are the Euler-Lagrange equations of an infinite-dimensional
extremization problem. The variational problem was proposed by
David Hilbert in 1915. Einstein published in the same year
the {\bf general theory of relativity}. In the case of a {\bf vacuum}: $T=0$,
solutions $g$ of the Einstein equations define {\bf Einstein manifolds} $(M,g)$. An example of a
solution to the vacuum Einstein equations different from the flat space solution is
the {\bf Schwarzschild solution}, which was found also in 1915 and published in 1916.
It is the metric given in spherical coordinates as
$-(1-r/\rho) c^2 dt^2+(1-r/\rho)^{-1} d\rho^2 + \rho^2 d\phi^2+\rho^2 \sin^2\phi d\theta^2$,
where $r$ is the {\bf Schwarzschild radius}, $\rho$ the distance to the singularity, $\theta,\phi$
are the standard {\bf Euler angles} ({\bf longitude} and {\bf colatitude}) in calculus. The metric solves the
Einstein equations for $\rho>r$.
The flat metric $-c^2 dt^2+d\rho^2+\rho^2 d\theta^2+\rho^2 \sin^2\theta d\phi^2$ describes
the vacuum and the Schwarzschild solution describes the gravitational field near a
{\bf massive body}. Intuitively, the metric tensor $g$ is determined by $g(v,v)$,
and the Ricci tensor by $R(v,v)$ which is 3 times the average
sectional curvature over all planes passing through a plane
through $v$. The scalar curvature is 6 times the average over
all sectional curvatures passing through a point.
See \cite{Misner,Ciu95}.
\index{Hilbert action}
\index{Hilbert-Einstein equations}
\index{General theory of relativity}

\section{Hall stable marriage}

Let $X$ be a finite set and $\mathcal{A}$ a family of finite
subsets $A$ of $X$. A {\bf transversal} of $\mathcal{A}$ is an injective function
$f: \mathcal{A} \to X$ such that $f(A) \in A$ for all $A \in \mathcal{A}$.
The set $\mathcal{A}$ satisfies the {\bf marriage condition} if for
every finite subset $\mathcal{B}$ of $\mathcal{A}$, one has
$|\mathcal{B}| \leq |\bigcup_{A \in \mathcal{B}} A|$.
The {\bf Hall marriage theorem} is
\index{Hall marriage problem}
\index{marriage condition}

\satz{$\mathcal{A}$ has a transversal $\Leftrightarrow$ $\mathcal{A}$ satisfies marriage condition.}

The theorem was proven by Philip Hall in 1935. It implies for example that if a deck of cards with 52
cards is partitioned into 13 equal sized piles, one can chose from each deck a card so that the 13
cards have exactly one card of each rank. The theorem can be deduced from a result in graph geometry:
if $G=(V,E)=(X,\emptyset)+(Y,\emptyset)$ is a bipartite graph,
then a {\bf matching} in $G$ is a collection of edges which pairwise have no common vertex.
For a subset $W$ of $X$, let $S(W)$ denote the set of all vertices adjacent to some
element in $W$. The theorem assures that there is an {\bf X-saturating matching}
(a matching that covers $X$) if and only if $|W| \leq |S(W)|$ for every $W \subset X$.
The reason for the name ``marriage" is the situation that $X$ is a set of men and $Y$
a set of women and that all men are eager to marry. Let $A_i$ be the set of women which
could make a spouse for the $i$'th man, then marrying everybody off is an
$X$-saturating matching. The condition is that any set of $k$ men has a combined list
of at least $k$ women who would make suitable spouses. See \cite{Brualdi2004}.
\index{bipartite graph}
\index{marriage theorem}

\section{Mandelbulb}

The {\bf Mandelbrot set} $M=M_{2,2}$ is the set of vectors $c \in \mathbb{R}^2$ 
for which $T(x)=x^2+c$ leads to a bounded orbit starting at $0=(0,0)$, where
$x^2$ has the {\bf polar coordinates} $(r^2,2\theta)$ if $x$ has the polar
coordinates $(r,\theta)$. (The map $T$ is just a real reformulation of the 
complex map $T(z)=z^2+c$ in $\mathbb{C}$ and written in the real so
that the construction can be done in arbitrary dimensions.)
The {\bf Mandelbulb set} $M_{3,8}$ is defined as the set of vectors $c \in \mathbb{R}^3$
for which $T(x)=x^8+c$ leads to a bounded orbit starting at $0=(0,0,0)$, where $x^8$
has the {\bf spherical coordinates} $(\rho^8,8 \phi, 8\theta)$ if $x$ has the
{\bf spherical coordinates} $(\rho,\phi,\theta)$. Like the Mandelbrot set, it is a compact
set (just verify that for $|x|>2$, the orbits go to infinity). 
The topology of $M_8$ is unexplored. Also like in the complex plane, one could look
at the dynamics of a polynomials $p=a_0 + a_1 x + \cdots + a_r x^r$ 
in $\mathbb{R}^n$. If $(\rho, \phi_1, \dots, \phi_{n-1})$
are spherical coordinates, then $x \to x^m = (\rho^m, m \phi_1, \dots, m \phi_{n-1})$ 
is a higher dimensional ``power" and allows to look at the dynamics of $T_{n,p}(x) = p(x)$. 
This defines then a corresponding {\bf Mandelbulb} $M_{n,p}$. As with all celebrities, there is a scandal:
\index{Mandelbulb set}
\index{Mandelbrot set}

\satz{There is no theorem about the Mandelbulb $M_{n,m}$ for $n>2$.} 

Except of course the just stated theorem. But you decide whether it is true of not. 
The Mandelbulb set has been discovered only recently. An attempt to 
trace some of its history was done in \cite{KnillSIC}:
already {\bf Rudy Rucker} had experimented with a variant of
$M_{3,2}$ in 1988. {\bf Jules Ruis} wrote me to have written a computer program in ``Basic" 
in 1997.  The first person we know who wrote down the formulas used today is {\bf Daniel White} ,
mentioned in a 2009 fractal forum. Jules Ruis 3D printed the first 
models in 2010.  See also \cite{Bourke2017} for some information on generating the
graphics. 

\section{Banach Alaoglu}

A {\bf Banach space} $X$ is a linear space equipped with a {\bf norm} $| \cdot |$
defining a metric $d(x,y)=|x-y|$ with respect to which the space $X$ is complete. The {\bf unit ball}
in $X$ is the closed ball $\{ x \in X \; | \; |x| \leq 1 \}$. 
The {\bf dual space} $X^*$ of $X$ is the linear space of {\bf linear functionals}
$f:X \to \mathbb{R}$ with the norm $|f| = \sup_{|x| \leq 1, x \in X} |f(x)|$. 
It is again a Banach space. The {\bf weak* topology} is the smallest topology on $X^*$
which makes all maps $f \to f(x)$ continuous for all $x \in X$. 
\index{Banach space}
\index{weak* topology}
\index{linear functionals}

\satz{The unit ball in a dual Banach space $X^*$ is weak* compact.}

The theorem was proven in 1932 in the separable case by Stefan Banach and 
in 1940 in general by Leonidas Alaoglu. The result essentially follows from Tychonov's theorem
as $X^*$ can be seen as a closed subset of a product space. Banach-Alaoglu therefore
relies on the axiom of choice. A case which often appears in applications is when $X=C(K)$ 
is the space of continuous functions on a compact Hausdorff space $K$. In that case $X^*$ is
the space of {\bf signed measures} on $K$. One implication is that the set of 
{\bf probability measures} is compact on $K$.  An other example are $L^p$ spaces ($p \in [1,\infty)$, 
for which the dual is $L^q$ with $1/p+1/q=1$ (meaning $q=\infty$ for $p=1$) and showing
that for $p=2$, the Hilbert space $L^2$ is self-dual. 
In the work of Bourbaki the theorem was extended from 
Banach spaces to {\bf locally convex spaces} (linear spaces equipped with a family 
of semi-norms). Examples are {\bf Fr\'echet spaces} (locally convex
spaces which are complete with respect to a translation-invariant metric).
See \cite{ConwayFunctionalAnalysis}.
\index{Locally convex}
\index{Fr\'echet space}

\section{Whitney trick}

Let $M$ be a smooth orientable simply connected $d$-manifold and two smooth connected
sub-manifolds $K,L$ of dimension $k$ and $l$ such that $k+l=d$ which have the property
that  $K$ and $L$ {\bf intersect transversely} in points $x,y$ in the sense 
that the tangent spaces at the intersection points span $T_xM$ and $T_yM$ and that
they have opposite intersection sign.
The two manifolds $K,L$ can be {\bf isotoped from each other along a disc}
if there exists a smooth 2-disk embedded in $M$ such that $M \cap K$ and $M \cap L$ are single points. 
The disk is called a {\bf Whitney disk}. The {\bf Whitney trick}  or {\bf Whitney lemma} is: 

\satz{Any transverse $K,L$ of $\geq 3$ manifolds in $M$ has a Whitney disk. }

See \cite{DonaldsonKronheimer}. In \cite{Lackenby} there are counter examples in $d \leq 4$.
The author writes there ``A hypothesis of algebraic topology given by the signs of the intersection points 
leads to the existence of an isotopy". The failure of the Whitney trick in smaller dimensions
is one reason why some questions in manifold theory appear hardest in three or four dimension.
There is a variant of the Whitney trick which works also in dimensions 5, where $K$ has dimension $2$
and $L$ has dimension $3$. 

\section{Torsion groups}

An {\bf elliptic curve} $E$ over $\mathbb{Q}$ is also called a
{\bf rational elliptic curve}. The curve $E$ carries an Abelian group
structure where every addition of a point $x \to x+y$ is a morphism. The {\bf torsion
subgroup} of $E$ is the subgroup consisting of elements which all
have finite order in $E$. The Mordell-Weil theorem (which applies more
generally for any Abelian variety) assures
that $E = \mathbb{Z}^r \oplus T$, where $T$ is a finite group and
$r$ is a finite number called the {\bf rank} of $E$.
{\bf Mazur's torsion theorem} states that
the only possible finite orders in $E$ are $1,2,3, \dots, 9,10$ and $12$.
Only 15 different torsion subgroups appear in rational elliptic curves:
$Z_1, \dots, Z_{10}, Z_{12}$ or $Z_2 \times Z_2, Z_2 \times Z_{4},Z_2 \times Z_6$
and $Z_2 \times Z_8$. Lets call this collection of groups the {\bf Mazur class}.
The theorem is:
\index{Mordell-Weil theorem}

\satz{The torsion group of a rational elliptic curve is in the Mazur class.}

The theorem was proven by Barry Mazur in 1977.  \cite{Silverman1986}.
\index{Mazur torsion theorem}

\section{Coloring}

A graph $G=(V,E)$ with vertex set $V$ and edge set $E$ is 
called {\bf planar} if it can be embedded in the Euclidean plane 
$\mathbb{R}^2$ without any of the edges intersecting. 
By a theorem of Kuratowski, this is equivalent to a graph theoretical statement: 
$G$ does not contain a homeomorphic image of neither the complete graph $K_5$ 
nor the bipartite utility graph $K_{3,3}$. A {\bf graph coloring} with $k$ colors
is a function $f: V \to \{1,2, \dots, k\}$ with the property 
that if $(x,y) \in E$, then $f(x) \neq f(y)$. In other words, 
adjacent vertices must have different colors. The {\bf 4-color theorem} is:
\index{planar graph}
\index{graph coloring}
\index{coloring}

\satz{Every planar graph can be colored with $4$ colors.}

Some graphs need 4 colors like a wheel graph having an odd number of spikes.
There are planar graphs which need less. The 1-point graph $K_1$ needs only one color, 
trees needs only 2 colors and the graph $K_3$ or any wheel graph with an even number
of spikes only need 3 colors. The theorem has an interesting history: 
since August Ferdinand M\"obius in 1840 spread a precursor problem given to him by Benjamin Gotthold Weiske, 
the problem was first known also as the {\bf M\"obius-Weiske puzzle} \cite{Soifer}. The actual problem was first posed
in 1852 by Francis Guthrie \cite{MaritzMouton}, after thinking about it with his brother Frederick, 
who communicated it to his teacher Augustus de Morgan, a former teacher of Francis who told
William Hamilton about it. Arthur Cayley in 1878 put it first in print, (but it was still not in the language of graph theory).
Alfred Kempe published a proof in 1879. But a gap was noticed by
Percy John Heawood 11 years later in 1890. There were other unsuccessful attempts like one by Peter Tait in 1880.
After considerable theoretical work by various mathematicians including Charles Pierce, 
George Birkhoff, Oswald Veblen, Philip Franklin, Hassler Whitney, Hugo Hadwiger, Leonard Brooks, 
William Tutte, Yoshio Shimamoto, Heinrich Heesch, Karl D\"urre or Walter Stromquist,
a computer assisted proof of the 4-color theorem was obtained by Ken Appel and Wolfgang Haken
in 1976. In 1997, Neil Robertson, Daniel Sanders, Paul Seymour, and Robin Thomas wrote a new computer
program.  Goerge Gonthier produced in 2004 a fully machine-checked proof of the four-color theorem \cite{RobinWilson}. 
There is a considerable literature like \cite{Ore,Heesch,FritschFritsch,SaatyKainen,ChartrandZhang2,RobinWilson}. 
\index{4 color theorem}

\section{Contact Geometry}

Assume $M$ is a smooth compact orientable $(2n-1)$-manifold equipped with an auxiliary Riemannian metric $g$.
A {\bf $1$-form} $\alpha \in \Lambda^1(M)$ defines a {\bf field of hyperplanes} $\xi = {\rm ker}(\alpha) \subset TM$.
Conversely, given a field of hyperplanes, one can define $\alpha = g(X,\cdot)$, where
$X$ is a local non-zero section of the line bundle $\xi^{\perp}$.
A {\bf contact structure} is a hyperplane field $\xi=d\alpha$ for which the volume form
$\alpha \wedge (d\alpha)^n$ is nowhere zero. The 1-form $\alpha$ is then called
a {\bf contact form} and $(M,\xi)$ is called a {\bf contact manifold}. The {\bf Reeb vector field}
$R$ is defined by $d\alpha(R,\cdot)=0$, $\alpha(R)=1$. The {\bf Weinstein conjecture} is a theorem
in dimension 3:

\satz{On a 3-manifold, the Reeb vector field has a closed periodic orbit.}

The theorem was proven by Clifford Taubes in 2007 using Seiberg-Witten theory.
Mike Hutchings with Taubes established 2 Reeb orbits under the condition that all Reeb orbits $R$ are
{\bf non-degenerate} in the sense that the linearized flow does not have an eigenvalue $1$. Hutchings
with Dan Cristofaro-Gardiner later removed the non-degeneracy condition
\cite{HutchingsTaubes,CristofaroHutshings} and also showed that if the product of the
{\bf actions} $\mathcal{A}(\gamma) = \int_\gamma \alpha$ of the two orbits is larger than the {\bf volume}
$\int_M \alpha \wedge d\alpha$ of the contact form, then there are three. To the history:
Alan Weinstein has shown already that if $Y$ is a convex compact hypersurface in $\mathbb{R}^{2n}$, then there is a periodic orbit.
Paul Rabinovitz extended it to star-shaped surfaces. 
Weinstein conjectured in 1978 that every compact hypersurface
of contact type in a symplectic manifold has a closed characteristic.
Contact geometry as an odd dimensional brother of symplectic geometry
has become its own field. Contact structures are the opposite of
integrable hyperplane fields: the {\bf Frobenius integrability} 
condition $\alpha \wedge d\alpha=0$ defines an {\bf integrable
hyperplane field} forming a co-dimension $1$ foliation of $M$. Contact geometry is therefore a
``totally non-integrable hyper plane field". \cite{GeigesContact}.
The higher dimensional case of the Weinstein conjecture is wide open \cite{HutchingsReview}.
Also the symplectic question whether every compact and regular energy surface $H=c$ for 
a Hamiltonian vector field in $\mathbb{R}^{2n}$ has a periodic solution is open. One knows
that there are for almost all energy values in a small interval around $c$. 
\cite{HoferZehnder1994}. 

\section{Simplicial spheres}

A {\bf convex polytope} $G$ is defined as the convex hull of $n$ points in $\mathbb{R}^{d}$ such
that all vertices are {\bf extreme points} called {\bf vertices}.
(Extreme points are points which do not lie in an open line segment of $G$.)
This definition of \cite{gruenbaum} is also called a {\bf polytopal sphere}. 
A {\bf simplicial sphere} is a geometric realization of a
simplicial complex that is homeomorphic to the standard (d-1)-dimensional spheres in $\mathbb{R}^{d}$.
For a polytopal sphere, the boundary of $G$ is made up of $(d-1)$-dimensional polytopes called {\bf $(d-1)$-faces}.
A {\bf cyclic polytope} $C(n,d)$ can be realized as the convex hull of the $n$ vertices
$\{ (t,t^2,t^3, \cdots t^d) \; | \; t =1,2, \dots, n\} \subset \mathbb{R}^{d}$.
Let $f_k(G)$ denote the number of $k$-dimensional
faces in $G$. So, $f_0(G)$ is the number of vertices, $f_1(G)$ the number of line segments and $f_{d-1}$ the
number of facets, the highest dimensional faces in $G$. Extending the definition to $f_{-1}=1$ (counting the
empty complex, which is a $(-1)$-dimensional complex), the vector $f=(f_{-1},f_0,f_1, \cdots f_{d})$ is
called the {\bf extended $f$-vector} of $G$. The {\bf upper bound theorem} is
\index{convex polytop}
\index{extreme point}
\index{cyclic polytop}
\index{simplicial spheres}
\index{f-vector}
\index{h-vector}

\satz{For simplicial spheres with $f_0(G)=n$, then $f_k(G) \leq f_k(C(n,d))$.}

This had been the {\bf upper bound conjecture} of Theodore Motzkin from 1957 which was proven by Peter McMullen in 1970
who reformulated it $h_k(G) \leq {n-d+k-1 \choose k}$ for all $k < d/2$ as the other numbers are determined by
{\bf Dehn-Sommerville conditions} $h_k=h_{d-k}$ for $0 \leq k \leq d$.
The {\bf $h$-vector} $(h_0, \dots h_{d})$ and {\bf $f$-vector} $(f_{-1},f_0, \dots, f_{d-1})$ determine each other via
$\sum_{k=0}^d f_{k-1}(t-1)^{d-k}=\sum_{k=0}^d h_k t^{d-k}$.
Victor Klee suggested the upper bound conjecture to be true for simplicial spheres, which was then
proven in by Richard Stanley in 1975 using new ideas like relating $h_k$ with intersection cohomology of
a projective {\bf toric variety} associated with the dual of $G$. (A {\bf toric variety} is an algebraic variety
containing an algebraic torus as an open dense subset such that the group action on the torus extends to the variety.)
The result for {\bf simplicial spheres} implies the result for convex polytopes because a subdivision of faces
of a convex polytope into simplices only increases the numbers $f_k$.
The {\bf g-conjecture} of McMullen from 1971 gives a complete characterization of $f$-vectors of simplicial spheres.
Define $g_0=1$ and $g_k=h_k-h_{k-1}$ for $k \leq d/2$. The $g$-conjecture claims that $(g_0, \dots g_{[d/2]})$
appears as a g-vector of a sphere triangulation if and only if there exists a multicomplex $\Gamma$ with exactly
$g_k$ vectors of degree $k$ for all $0 \leq i \leq [d/2]$. (A {\bf multi-complex} $\Gamma$ is a set of non-negative integer
vectors $(a_1, \dots, a_n)$ such that if $0 \leq b_i \leq a_i$, then $(b_1, \dots b_n)$ is in $\Gamma$. The degree
of a multicomplex is $\sum_i a_i$.) The {\bf g-theorem} proves this for polytopal spheres (Billera and Lee in 1980 sufficiency)
and (Stanley 1980 giving necessity). The g-conjecture is open for simplicial spheres.
\cite{Ziegler,Stanley1996,ChenStanley}
\index{g-conjecture}
\index{Dehn-Sommerville conditions}
\index{upper bound conjecture}
\index{multi-complex}
\index{toric variety}

\section{Bertrand postulate}

A basic result in number theory is

\satz{For $n>1$, there always exists a prime $p$ between $n$ and $2n$.}

As the theorem was conjectured in 1845 by Joseph Bertrand, it is still called {\bf Bertrand's postulate}.
Since Pafnuty Tschebyschef's (Chebyshev) proof in 1852, it is a theorem.
For a proof, see \cite{NivenZuckermanMontgomery}  page 367.
Srinivasa Ramanujan simplified Chebyshev's proof considerably in 1919 and strengthened
it: if $\pi(x) = \sum_{p \leq x, p \; {\rm prime}} 1$ is the {\rm prime counting function}, then
Bertrand's result can be restated as $\pi(x)-\pi(x/2) \geq 1$ for $x \geq 2$.
Ramanujna shows that $\pi(x)-\pi(x/2) \geq k$, for large enough $x$ (larger or equal than $p_k$).
The primes $p_k$ giving the lower bound for $x$ solving this are called {\bf Ramanujan primes}.
Simple proofs like one of Erd\"os from 1932 are given in Wikipedia or
\cite{Hua1982} page 82, who notes "it is not a very sharp result. Deep analytic
methods can be used to give much better results concerning the gaps between successive primes".
There is a very simple proof assuming the {\bf Goldbach conjecture} (stating that
every even number larger than $2$ is a sum of two primes): \cite{RicardoGoldbach}
if $n$ is not prime, then $2n=p+q$ is a sum of two primes, where one
is larger than $n$ and one smaller than $2n$; on the other hand, if $n$ is prime, 
then $n+1$ is not prime and $2n+2=p+q$ is a sum of two primes, where one, say
$q$ is larger than $n$ and smaller than $2n+2$. But $q$ can not be $2n+1$ (as
that would mean $p=1$), nor $2n$ (as $2n$ is composite) so that $n<q<2n$. 
There are various generalizations like Mohamed El Bachraoui's 2006 theorem that there are primes
between $2n$ and $3n$ or Denis Hanson from 1973 \cite{Hanson1973} that there are primes
between $3n$ and $4n$ for $n \geq 1$. Mohamed El Bachraoui asked in 2006 whether
for all $n>1$ and all $k \leq n$, there exists a prime in $[kn,(k+1)n]$ which is
for $k=1$ the Bertrand postulate. A positive answer would give that there is always 
a prime in the interval $[n^2,n^2+n]$. Already the {\bf Legendre conjecture}, asking
whether there is always a prime $p$ satisfying $n^2<p<(n+1)^2$ for $n \geq 1$ is open.
The {\bf Legendre's conjecture} is the fourth of the super famous great problems of Edmund Landau's 1912 list:
the other three are the {\bf Goldbach conjecture}, the {\bf twin prime conjecture} and then the
{\bf Landau conjecture} asking whether there are infinitely many primes of the form $n^2+1$.
Landau really nailed it. There are 4 conjectures only, but all of them can be stated in half
a dozen words, are completely elementary, and for more than 100 years, 
nobody has proven nor disproved any of them. 
\index{prime counting function}
\index{Goldbach conjecture}
\index{Legendre conjecture}
\index{Landau conjecture}
\index{Twin prime conjecture}
\index{Bertrand postulate}
\index{Ramanujan primes}

\section{Non-squeezing theorem}

The Euclidean space $M=\mathbb{R}^{2n}$ carries the standard symplectic
$2$-form $\omega(v,w) = (v, J w)$ with the skew-symmetric matrix
$J=\left[ \begin{array}{cc} 0 & I \\ -I & 0 \end{array} \right]$.
A linear transformation $f: M \to M, x \to Ax$ is called {\bf symplectic}, if
$A$ satisfies $A^T J A = J$. A smooth transformation $f: M \to M$
is called a {\bf symplectomorphism} if it is a diffeomorphism and if
the derivative $df$ is a symplectic map from $T_xM \to T_{f(x)} M$ at every point $x \in M$.
Any smooth map for which $df$ is symplectic is automatically a diffeomorphism as
symplectic matrices have determinant $1$ and are so invertible.
Let $B(r)=\{ x \in M \; | \; x \cdot x \leq r^2 \}$ denote the {\bf round solid ball of radius $r$}
and $Z(r) = \{ x \in M \; | \; x_1^2+y_1^2 \leq r^2 \}$ the {\bf solid cylinder
of radius $r$}. Given two sets $A,B$, one says there is a symplectic embedding of
$A$ in $B$, if there exists a symplectomorphism $f$ such that $f(A) \subset B$.
As symplectic maps are volume preserving, a necessary condition is
${\rm Vol}(A) \leq {\rm Vol}(B)$. Is this the only constraint? Yes, for $n=1$,
where the cylinder and the ball are the same as defined $B(r)=Z(r)$.
But no in higher dimensions $n \geq 2$ by the {\bf Gromov non-squeezing theorem}:
\index{diffeomorphism}
\index{symplectomorphism}

\satz{A symplectic embedding $B(r) \to Z(R)$ implies $r \leq R$.}

The theorem has been proven in 1985 by Michael Gromov.
It has been dubbed as the {\bf principle of the symplectic camel} by Maurice de Gosson
referring to the ``eye of the needle" metaphor. A reformulation of the Gosson allegory
\cite{GossonCamel} after encoding ``camel" = ``ball in the phase space", ``hole = ``cylinder", and ``pass"=``symplectically
embed into", ``size of the hole" = ``radius of cylinder" and ``size of the camel" = ``radius of the ball"
is: ``There is no way that a camel can pass through a hole if the size of the hole
is smaller than the size of the camel". See \cite{McDuff2009,HutchingsFun2016} for expositions.
The non-squeezing theorem motivated also the introduction of {\bf symplectic capacities}, 
quantities which are monotone $c(M) \leq c(N)$ if there is a symplectic embedding of 
$M$ into $N$, which are conformal in the sense that if $\omega$ is scaled by $\lambda$, 
then $c(M)$ is scaled by $|\lambda|$ and such that $c(B(1)) = c(Z(1))=\pi$. 
For $n=1$, the area is an example of a symplectic
capacity (actually unique). The existence of a symplectic capacity obviously proves the
squeezing theorem. Already Gromov introduced an example, the {\bf Gromov width}, which is the
smallest. More are constructed in using calculus of variations. 
See \cite{HoferZehnder1994,McDuffSalamon1998}.
\index{symplectic capacity}
\index{symplectic embedding}

\section{K\"ahler Geometry}

A {\bf K\"ahler manifold} is a complex manifold $(M,J)$ together with a 
Hermitian metric $h$ whose associated {\bf K\"ahler form} $\omega$ is closed. 
(The manifold can be given by a Riemannian metric $g$ 
{\bf compatible} with the complex structure $g(JX,JY) = g(X,Y)$. 
The {\bf K\"ahler form} $\omega$ is then a 2-form
$\omega(X,Y)=g(JX,Y)$ satisfying $d\omega=0$ and the metric $h=g+i\omega$ is 
the Hermitian metric. $(M,\omega)$ is then also a symplectic manifold.) 
As $\omega$ is closed, it represents an element 
in the cohomology class $H^2(M)$ called {\bf K\"ahler class}. 
The {\bf Calabi inverse problem} is: given a compact K\"ahler manifold 
$(M,\omega_0)$ and a $(1,1)$-form $R$ representing $2\pi$ times the {\bf first
Chern class} of $M$, find a metric $\omega$ in the K\"ahler class of $\omega_0$ 
such that ${\rm Ricci}(\omega)= R$. In local coordinates, one can write
${\rm Ricci}(\omega) = -i \partial \overline{\partial} \log {\rm det}(g)$.
For compact $M$:
\index{K\"ahler manifold}
\index{K\"ahler class}
\index{Chern class}
\index{Ricci curvature}

\satz{ The Calabi inverse problem has a unique solution $\omega$. }

This was conjectured in 1957 by Eugenio Calabi and proven in 1978 by 
Shing-Tung Yau by solving nonlinear {\bf Monge-Amp\`ere} equations using
analytic Nash-Moser type techniques. 
The theorem implies that if the first Chern class of $M$ is zero, 
then $(M,\omega_0)$ carries has a unique {\bf Ricci-flat} K\"ahler 
metric $g$ in the same K\"ahler class than $\omega_0$. 
K\"ahler geometry deals simultaneously
with Riemannian, symplectic and complex structures: $(M,g)$ is a Riemannian,
$(M,\omega)$ is a symplectic and $(M,J)$ is a complex manifold. 
The inverse problem of characterizing geometries from curvature 
data is central in all of differential geometry. Here are some examples:
a) $M=\mathbb{C}^n$ with Euclidean metric $g$ is K\"ahler with 
$\omega=(i/2) \sum_k dz^k \wedge d\overline{z}^k$
but it is not compact. 
But if $\Gamma$ is a lattice, then the induced metric on the torus 
$\mathbb{C}^n/\Gamma$ is K\"ahler. 
b) Because complex submanifolds of a K\"ahler manifold are K\"ahler ,
and the complex projective space $\mathbb{C}P^n$ with the 
{\bf Fubini-Study metric} is K\"ahler 
(with $\omega=i \partial \overline{\partial} \rho$, where
$\rho=\log(1+\sum_k |z_k|^2/2)$ is the {\bf K\"ahler potential}), 
any {\bf complex projective variety} is K\"ahler. 
d) For the complex hyperbolic case where $M$ is the unit ball in 
$\mathbb{C}^n$, the K\"ahler potential is $\rho=1-|z|^2$. 
By Kodeira, K\"ahler forms representing an integral cohomology class
correspond to projective algebraic varieties. 
c) {\bf Calabi-Yau} manifolds are complex K\"ahler manifolds with 
zero first Chern classes. Examples are {\bf K3} surfaces. 
The existence theorem assures that they carry a 
{\bf Ricci-flat metric}, which are examples of {\bf K\"ahler-Einstein} metrics. 
Also Hodge theory works well for K\"ahler manifolds. In the complex, the 
{\bf Dolbeault operators} $\partial, \overline{\partial}$ and $d=\partial+\overline{\partial}$ 
lead to {\bf Hodge Laplacians} $\Delta_\partial, \Delta_{\overline{\partial}}$ and $\Delta_d$,
and so to {\bf harmonic forms} $H^{p,q}$ for differential forms of type $(p,q)$ and harmonic 
$r$-forms for $\Delta$. In the K\"ahler case, 
$H^r = \sum_{p+q=r} H^{p,q}$. An example result due to Lichnerowicz is that if 
${\rm Ricci}(\Omega) \geq \lambda>0$, then the first eigenvalue 
$\lambda_1$ of $\Delta$ satisfies $\lambda_1 \geq 2 \lambda$.  
See \cite{BallmanLectures,Westrich,AubinNonlinear}.
\index{Dolbeault operators}
\index{harmonic forms}

\section{Projective Geometry}

A {\bf conic section} is a curve which is obtained when intersecting a 
{\bf cone} $x^2+y^2=z^2$ with a plane $ax+by+cz=d$. 
A bit more general is a {\bf conic}, an algebraic curve 
$ax^2+bx y + c y^2 + d x+e y+g=0$ of degree $2$. They are either 
{\bf non-singular conics}, classified as {\bf ellipses} 
like $x^2+y^2=1$, {\bf hyperbola} $x^2-y^2=1$ or {\bf parabola} $x^2=y$, 
or then {\bf degenerate conics} like a point $x^2+y^2=0$, the {\bf cross}
$x^2=y^2$, the {\bf line} $x^2=0$ or {\bf pair of parallel lines} $x^2=1$. 
Given $6$ different points $A_1,A_2,A_3,B_1,B_2,B_3$ on a conic, where $A_1,A_2,A_3$ 
are neighboring and $B_1,B_2,B_3$ are neighboring,
a {\bf Pascal configuration} is the set of lines $A_i B_j$ with $i \neq j$. 
The {\bf intersection points} of this Pascal configuration is the set
of three intersections of $A_i B_j$ with $A_j B_i$, where $\{i,j\}$ runs
over all three 2-point subsets of $\{1,2,3\}$. 
\index{conic section}
\index{conic}
\index{hyperbola}

\satz{The intersection points of a Pascal configuration are on a line.}

The theorem was found in 1639 by Blaise Pascal (as a teenager) in the 
case of an ellipse.  A limiting case where we have two crossing lines is the 
{\bf Pappus hexagon theorem}, which goes back to Pappus of Alexandria 
who lived around 320 AD. The {\bf Pappus hexagon theorem} is one of the first 
known results in {\bf projective geometry}.
\index{projective geometry}
\index{Pascal configurations}
\index{Pascal theorem}
\index{Pappus hexagon theorem}

\section{Vitali theorem}

A {\bf Lebesgue measure} in Euclidean space $\mathbb{R}^n$ is a Borel 
measure which is invariant under Euclidean transformations. It is the Haar
measure of the locally compact group $\mathbb{R}^n$ and unique if one normalizes
is so that the unit cube has measure $1$. In dimension $n=1$, the Lebesgue measure 
of an interval $[a,b]$ is $b-a$.  In dimension $n=2$, the Lebesgue measure of 
a measurable set is the area of the set. In particular, a ball of radius $r$ 
has area $\pi r^2$. When constructing the measure one has to specify a $\sigma$-algebra,
which is in the Lebesgue case the Borel $\sigma$-algebra generated by 
the open sets in $\mathbb{R}^n$. One has for every $n \geq 1$: 

\satz{There exist sets in $\mathbb{R}^n$ that are not Lebesgue measurable.}
\index{Vitali theorem}
\index{Lebesgue measure}
\index{Haar measure}

The result is due to Giuseppe Vitali from 1905. It justifies why one has
to go through all the trouble of building a $\sigma$-algebra carefully 
and why it is not possible to work with the complete $\sigma$-algebra of 
all subsets of $\mathbb{R}^n$ (which is called the {\bf discrete $\sigma$-algebra}).
The proof of the Vitali theorem shows connections with the foundations of mathematics:
by the {\bf axiom of choice} there exists a set $V$ which represents equivalence
classes in $\mathbb{T}/\mathbb{Q}$, where $T$ is the circle.
For this {\bf Vitali set} $V$, all translates $V_r = V+r$ are all disjoint with $r \in \mathbb{Q}$. 
$\{ r+V ,  \; r \in \mathbb{Q} \} = \mathbb{R}$ and so form
a partition. By the Lebesgue measure property, all 
translated sets $V_r$ have the same measure. As they are a countable set 
and are disjoint and add up to a set of Lebesgue measure $1$, they have to have
measure zero. But this contradicts $\sigma$-additivity. Now lift $V$ to $R$ and then 
build $V \times \mathbb{R}^{n-1}$. More spectacular are decompositions of the unit ball 
into $5$ disjoint sets which are equivalent under Euclidean transformations and which 
can be reassembled to get two disjoint unit balls. This is the 
{\bf Banach-Tarski construction} from 1924. 
\index{Banach-Tarski construction}
\index{discrete $\sigma$ algebra}

\section{Wilson's theorem}

The {\bf factorial} $n!$ of a number  defined as $n!=1 \cdot 2 \cdots n$. 
For example, $5!=120$. 

\satz{$n>1$ is prime if and only if $(n-1)!+1$ is divisible by $n$.}

For $n=5$ for example $(5-1)!+1=25$ is divisible by $5$. 
For $n=6$ we have $(6-1)!+1 = 121$ which is not divisible by $6$. Indeed, 
$6=2*3$ is not prime. 
The theorem is named after John Wilson, who was a student of Edward Waring. It seems
that Joseph-Louis Lagrange gave the first proof in 1771. It is not a practical
way to determine whether a number is prime:
\cite{SteinElementary}: {\it from a computational point of view, it is probably one of the
world's least efficient primality tests, since computing $(n-1)!$ takes so many
steps.}
Also named after Wilson are the {\bf Wilson primes}. These are 
primes for which not only $p$ but $p^2$ divides $(p-1)!+1$. 
The smallest one is $5$. It is not known whether there are infinitely many. 
\index{Wilson theorem}
\index{Factorial}
\index{primality test}
\index{Wilson primes}

\section{Carleson theorem}

If $f \in L^2(\mathbb{T})$, where $\mathbb{T}=\mathbb{R}/(2\pi \mathbb{Z})$ is the circle, then
the Fourier transform $L^2(\mathbb{T}) \to l^2(\mathbb{Z})$ gives a {\bf Fourier series}
$g(x) = \sum_{k \in \mathbb{Z}} c_k e^{i k x}$, where $c=(\dots, c_{-2},c_{-1},c_0,c_1,c_2, \dots) \in l^2(\mathbb{Z})$
is given by $c_k = (2\pi)^{-1} \int_{\mathbb{T}} f(x) e^{i k x} \; dx$.
For smooth $f$, one knows $g=f$ and Parseval's identity $\int_{\mathbb{R}} f^2(x) \; dx = \sum_k c_k^2$ so that
the Fourier transform extends to an unitary operator $L^2(\mathbb{T})  \to l^2(\mathbb{Z})$. This does not say
anything yet about the convergence of the sequence $g_n(x)$.
We say the Fourier series {\bf converges to $f$ at a point $x$}, if the sequence $g_n(x) = \sum_{k=-n}^n c_k e^{i k x}$
converges to $f(x)$ for $n \to \infty$. We say, a sequence $g_n(x)$ converges almost everywhere to $f$,
if there exists a set $Y \subset \mathbb{T}$ of full Lebesgue measure $\mu(\mathbb{T})=1$ such that the
series converges for all $x \in Y$. (The Lebesgue measure is the normalized Haar measure $dx/(2\pi)$ on the circle).
That the question can be subtle is illustrated by the result of Andrey Kolmogorov from 1923 to 1926,
who gave examples of $L^1(\mathbb{T})$ functions for which the Fourier series diverges everywhere.
\index{Fourier series}
\index{Parseval's identity}
\index{Fourier series}
\index{almost everywhere convergence}

\satz{The Fourier series of a $L^2$ function converges almost everywhere.}

The statement had been conjectured by Nikolai Luzin in 1915 and was known as the Luzin conjecture.
The theorem was proven by Lennart Carleson in 1966. 
An extension to $L^p$ with $p \in (1,\infty]$ was proven by
Richard Hunt in 1968. The proof of the Carleson theorem is difficult. 
While mentioned in harmonic analysis texts like \cite{Katznelson} or surveys 
\cite{KhavinNikolskij}, who say about the Carleson-Hunt theorem that
it is {\it one of the deepest and least understood parts of the theory}.
\index{Carleson theorem}

\section{Intermediate value}

Let $(X,\mathcal{O})$ be connected topological space         
and $f: X \to \mathbb{R}$ a continuous map. We say  that
$f$ {\bf reaches both positive and negative signs} if there exists
$a,b \in X$ such that $f(a)<0$ and $f(b)>0$. A {\bf root}
of $f$ is a point $x \in X$ such that $f(x)=0$. Let $C(X)$ denote
the set of continuous functions from $X$ to $\mathbb{R}$.
This means that for $f \in C(X)$ and all open sets $U$ in 
$\mathbb{R}$, one has $f^{-1}(U) \in \mathcal{O}$.

\satz{$f \in C(X)$ reaching both signs on a connected $X$  has a root.}

The theorem was proven by Bernard Bolzano in 1817 for functions 
from the interval $[a,b]$ to $\mathbb{R}$. 
The proof follows from the definitions: as $P=(0,\infty)$ is open, also $f^{-1}(P)$ is
open. As $N=(-\infty,0)$ is open, also $f^{-1}(N)$ is open. If
there is no root, then $X=N \cup P$ is a disjoint union of two open
sets and so disconnected. This contradicts the assumption of $X$ being 
connected. A consequences is the {\bf wobbly table theorem}: 
given a square table with parallel equal length legs and a ``floor" given
by the graph $z=g(x,y)$ of a continuous $g$ can be rotated and possibly
translated in the $z$ direction so that all 4 legs are on the table.
The proof of this application is seen as a consequence of the intermediate
value theorem applied to the height function $f(\phi)$ of the fourth leg
if three other legs are on the floor. 
A consequence is also {\bf Rolle's theorem}, assuring that if
a continuously differentiable function $[a,b] \to \mathbb{R}$     
with $f(a)=f(b)$ has a point $x \in (a,b)$ with $f'(x)=0$. Tilting
Rolle gives the {\bf mean value theorem} assuring that for a 
continuously differentiable function $[a,b] \to \mathbb{R}$,
there exists $x \in (a,b)$ with $f'(x)=f(b)-f(a)$. 
The general theorem shows that it is the connectedness and not the 
completeness of $X$ which is the important assumption. 

\section{Perron-Frobenius}

A $n \times n$ matrix $A$ is {\bf non-negative} if $A_{ij} \geq 0$ for all
$i,j$ and {\bf positive} if $A_{ij} >0$ for all $i,j$. The {\bf Perron-Frobenius}
theorem is:
\index{Perron-Frobenius theorem}
\index{positive matrix}
\index{non-negative matrix}

\satz{A positive matrix has a unique largest eigenvalue.}

The theorem has been proven by Oskar Perron in 1907 \cite{Perron1907}
and by Georg Frobenius in 1908 \cite{Frobenius1908}.
When seeing the map $x \to Ax$ on the projective space, this is in suitable
coordinates a contraction and the Banach fixed point theorem applies. 
This is the proof of Garret Birkhoff who used the {\bf Hilbert metric} \cite{KohlbergPratt}.
The Brouwer fixed point theorem only gives existence, not uniqueness, but the
Brouwer fixed point applies for non-negative matrices.
This has applications in graph theory, Markov chains or Google page rank.
The {\bf Google matrix} is defined as $G=d A + (1-d)E$, where
$d$ is a {\bf damping factor} and $A$ is a Markov matrix defined by
the network and $E$ is the matrix $E_{ij}=1$. Sergey Brin and Larry Page write
``the damping factor $d$ is the probability at each page the random surfer
will get bored and request another random page". The {\bf page rank equation}
is $Gx=x$. In other words, the Google Page rank vector (the one billion dollar
vector), is a Perron-Frobenius eigenvector. It assigns page rank values to the
individual nodes of the network. See \cite{LangvilleMeyer}. For the linear algebra
of non-negative matrices, see \cite{MincNonnegative}.
\index{page rank}
\index{Hilbert metric}
\index{Markov matrix}
\index{Google matrix}
\index{damping factor}

\section{Continuum hypothesis}

$\aleph_0$ is the cardinality of the {\bf natural numbers} $\mathbb{N}$.
$\aleph_1$ is the next larger cardinality. The cardinality of the {\bf real
numbers} $\mathbb{R}$ is $2^{\aleph_0}$. The statement $2^{\aleph_0}=\aleph_1$ is the
{\bf continuum hypothesis} abbreviated CH. The {\bf Zermelo-Fraenkel axiom system} ZFC
of set theory is the most common foundational axiomatic framework of mathematics. The
letter $C$ refers to the {\bf axiom of choice}.
\index{Continuum hypothesis}
\index{Axiom of choice}
\index{Zermelo Fraenkel}
\index{cardinality}

\satz{Neither $2^{\aleph_0}=\aleph_1$ nor $2^{\aleph_0}\neq \aleph_1$ can be proven in ZFC.}

This result combines a result of Kurt Goedel from 1938 \cite{Goedel1940} (CH is consistent with ZFC)
and Paul Cohen (Negated CH is independent of ZVC) from 1963 \cite{Cohen1963,Cohen1966}. 
Cantor had for a long time tried to prove that the continuum hypothesis holds. 
The Goedel-Cohen's theorem shows that any such effort has been
in vain and illustrates why Cantor was doomed not to succeed.
The problem had then been the first of Hilbert's problems of 1900. For more, see \cite{Shelah} or
\cite{Woodin2001} who summarizes the result in words: {\it G\"odel solved the {\bf substructure problem} 
in 1938.  Over 25 years later Cohen, arguably the Galois of set theory, solved the 
{\bf extension problem}.}
\index{Hilbert's problems}

\section{Homotopy-Homology}

Given a path connected pointed topological space $X$ with base $b$, the $n$'th {\bf homotopy group}
$\pi_n(X)$ is the set of equivalence classes of base preserving maps
from the pointed sphere $S^n$ to $X$. It can be written as the set of homotopy classes of maps
from the {\bf $n$-cube} $[0,1]^n$ to $X$ such that the boundary of $[0,1]^n$ is mapped to $b$.
It becomes a group by defining addition as $(f+g)(t_1,\dots, t_n)=f(2t_1,t_2, \dots t_n)$ for
$0 \leq t_1 \leq 1/2$ and $(f+g)(t_1,\dots, t_n) = g(2t_1-1,t_2, \dots ,t_n)$ for $1/2 \leq t \leq 1$.
In the case $n=1$, this is ``joining the trip": travel first along the first curve with twice
the speed, then take the second curve. The groups $\pi_n$ do not depend on the base point.
As $X$ is assumed to be connected, $\pi_0(X)$ is the trivial group. 
The group $\pi_1(X)$ is the {\bf fundamental group}.
It can be non-abelian. For $n \geq 2$, the groups $\pi_n(X)$ are always Abelian $f+g = g+f$.
The $k$'th {\bf homology group} $H_n(X)$ of a topological
space $X$ with integer coefficients is obtained from the chain complex of the free abelian group
generated by continuous maps from $n$-dimensional simplices to $X$. The {\bf Hurewicz theorem} is
\index{Hurewicz theorem}
\index{homotopy group}
\index{fundamental group}
\index{homology group}

\satz{There exists a homomorphism $\pi_n(X) \to H_n(X)$.}

Higher homotopy groups were discovered by Witold Hurewitz during the years 1935-1936.
The Hurewitz theorem itself has then been established in 1950 \cite{Hurewicz}. In the case $n=1$,
the homomorphism can be easily described: if $\gamma: [0,1] \to X$ is a path, then since $[0,1]$
is a $1$-simplex, the path is a singular $1$-simplex in $X$. As the boundary of $\gamma$ is empty,
this singular $1$-simplex is a cycle. This allows to see it as an element in $H_1(X)$. If two
paths are homotopic, then their corresponding singular simplices are equivalent in $H_1(X)$.
There is an elegant proof using Hodge theory if $X=M$ is a compact manifold: the image $C$ of a map
$\pi_p(M)$ can be interpreted as a Schwartz distribution on $M$. Let $L=(d+d^*)^2$ be the Hodge
Laplacian and let the heat flow $e^{-t L}$ act on $C$. For $t>0$, the image $e^{-t L} C$ is now
smooth and defines a differential form in $\Lambda^p(M)$. As all the non-zero eigenspaces get damped
exponentially, the limit of the heat flow is a {\bf harmonic form}, an eigenvector to the 
eigenvalue $0$. But Hodge theory identifies ${\rm ker}(L|\Lambda^p)$ with $H^p(M)$ and so with 
$H_p(M)$ by Poincar\'e duality. The Hurewitz homomorphism is then even constructive. ``Just heat up the curve to get 
the corresponding cohomology element, the commutator group elements get melted away by the heat."
A space $X$ is called {\bf $n$-connected} if $\pi_i(X)=0$ for all $i \leq n$. So, 
$0$-connected means {\bf path connected} and $1$-connected is {\bf simply connected}.
For $n \geq 2$, one has $\pi_n(X)$ isomorphic to $H_n(X)$ if $X$ is $(n-1)$-connected. 
In the case $n=1$, this can already not be true as $\pi_1(X)$ is in 
general non-commutative and $H_1(X)$ is but $H_1(X)$ is the isomorphic to the {\bf abelianization}
of $G=\pi_1(X)$ which is the group obtained by factoring out the commutator subgroup $[G,G]$
which is a normal subgroup of $G$ and generated by all the commutators $g^{-1} h^{-1} g h$ 
of group elements $g,h$ of $G$. See \cite{Hatcher}.
\index{n-connected space}
\index{Hurewicz homomorphism}

\section{Pick's theorem}

Let $P$ be a {\bf simple polygon} in the plane $\mathbb{R}^2$. This means
that it is given by as finite ordered set of points called {\bf vertices} $P_i=(x_i,y_i)$
$i=0,\dots,n$ such that the line segments $P_i P_{{\rm mod}(i+1,n)}$ 
called {\bf edges} joining neighboring points do not intersect. 
The polygon defines a {\bf polygonal region} $G$ with
area $A$. Assume now that all coordinates $x_i,y_i$ are integers. Let $I$ be the number
of lattice points $(k,l) \in \mathbb{Z}^2$ inside $G$ and $B$ the number
of lattice points at the boundary of $G$. {\bf Pick's theorem} assures:
\index{Pick theorem}

\satz{ $A=I+B/2-1$. }

The result was found in 1899 by Georg Pick \cite{PicksTheorem}. For a triangle for example
with no interior points, one has $0+3/2-1=1/2$, for a rectangle parallel to the coordinate axes
with $I=n*m$ interior points and $B=2n + 2m +4$ boundary points and area $A=(n+1)(m+1)$
also $I-B/2-1=A$. 
The theorem has become a popular school project assignment in early geometry courses
as there are many ways to prove it. An example is to cut away a triangle and use 
induction on the area then verify that if two polygons are joined along a line segment, the
functional $I+B/2-1$ is additive. There are other explicit formulas for
the area like Green's formula $A=\sum_{i=0}^{n-1} x_i y_{i+1}-x_{i+1}y_i$
which does not assume the vertices $P_i=(x_i,y_i)$ to be lattice points.
\index{Green's formula}

\section{Isospectral drums}

On a compact region $G \subset \mathbb{R}^2$ with piecewise smooth boundary $\delta G$ one can
look at the {\bf Dirichlet problem} $-\Delta f = 0$ in the interior of $G$ and
$f=0$ on $\delta G$. The region is considered a ``drum". If hit, one hears the 
spectrum of the  Laplacian $\Delta u = u_{xx} + u_{yy}$. 
There is a sequence of Dirichlet eigenvalues 
$0=\lambda_0 < \lambda_1 \leq \lambda_2 \leq \cdots$, real values which 
solve $-\Delta u_n = \lambda_n u_n$ for some functions $u_n$ which are zero
on the boundary. For example, if $G$ is the
square $[0,\pi] \times [0,\pi]$, then the eigenvalues are $n^2+m^2$ with eigenvectors
$\sin(n x) \sin(m x)$. The eigenvalue $0$ belongs to the constant eigenfunction.
Two drums are called {\bf isospectral}, if they have the same eigenvalues.
Two drums are non-isometric, if there is no transformation generated by rotations,
translation and reflections which maps one drum to the other. 

\satz{There exist non-isometric but isospectral drums.} 

Mark Kac had asked in 1962 ``Can one hear the sound of a drum" \cite{Kac66}". 
Caroline Gordon, David Webb and Scott Wolpert answered this question negatively \cite{GWW}.
In the convex case, the question is still open. 
\index{Isospectral drum}
\index{Laplacian}
\index{Eigenvalues}
\index{Dirichlet problem}

\section{Bertrand theorem}

The path $r(t)$ of a particle in $\mathbb{R}^n$ moving in a {\bf central force potential} $V(x)=f(|x|)$
experiences the {\bf central force} $F=-\nabla V(x)= -f'(|x|) x/|x|$. In the case of the Newton
potential $V(x)=-GMm/|x|$, where the central mass $M$, the body mass $m$ as well as  
the {\bf gravitational constant} $G$ determines the force $F(x)=-x GM m/|x|^3$.
The motion of the particle follows the differential equations $r''(t) = -M G r(t)/|r|^3$, 
which conserve the {\bf energy} $E(r)=m r'^2/2+V(r)$
and {\bf angular momentum} $L=m r \wedge r'$, a $n(n+1)/2$ dimensional quantity. The invariance of 
$L$ assures that $r(t)$ stays in the plane initially spanned by $r(0)$ and $r'(0)$
and that the area of the parallelogram spanned by $r(t)$ and $r'(t)$ is constant. 
To see the natural potential in $\mathbb{R}^n$ is, one has to go beyond Newton and
pass to Gauss, who wrote the gravitational law in the form ${\rm div}(F) = 4 \pi \rho$, where $\rho$ is the 
mass density. It expresses that mass is the source for the force field $F$. To get the
force field in a central symmetric mass distribution, one can use the {\bf divergence theorem} in $\mathbb{R}^n$ 
and relate the integral of $4 \pi \rho$ over a ball of radius $r$ with the flux of $F$ through the sphere $S(r)$
of radius $r$. The former is $4 \pi M$, where $M$ is the total mass in the ball, the later is
$-|S(r)| F(r)$, where $|S(r)|$ is the surface area of the sphere and the negative sign is because for an attractive
force $F(r)$ points inside.  So, in three dimensions, Gauss recovers the Newton gravitational law
$F(r)=-4\pi G M/|S(r)| = -G M/|r|^2$. There is a natural central
force {\bf Kepler problem} in any dimensions: in $\mathbb{R}^n$, we have $F(r)=-C_n r/|r|^n$ 
where $C_n$ is a constant. For $n=1$, there is a constant force pulling the particle towards the center,
for $n=2$, one has a $1/|r|$ force which corresponds to a logarithmic potential, for $n=3$, it is the
Newtonian inverse square $1/r^2$ force, in $n=4$, it is a $1/r^3$ force. For $n=0$, one formally gets the 
{\bf harmonic oscillator} which is {\bf Hook's law}. Which potentials lead to periodic motion? The answer is surprising
and was given by Bertrand: only the harmonic oscillator potential and the Newtonian potential in $\mathbb{R}^3$ work. Let us call
a central force potential {\bf all periodic} if every bounded (position and velocity) solution $r(t)$ of the differential equations
is periodic. Already for the Kepler problem, there are not only motions on ellipses but also 
scattering solutions moving on parabola or hyperbola, or then suicide motions, with $r'(0)=0$, where the particle 
dives into the singularity.
\index{Bertrand's theorem}
\index{Hook's law}
\index{central force}
\index{Newton potential}
\index{Angular momentum}
\index{Kepler problem}

\satz{Only the Newton potentials for $n=-1$ and $n=3$ are all periodic.}

This theorem of Joseph Bertrand from 1873 tells that three dimensional space is special as it in
any other dimension, calendars would be almost periodic as the solutions to the Kepler problem would not
close up. We could live with that but there are more compelling
reasons why $n=3$ is dynamically better: in other dimensions, only very special orbits stay bounded. A small
perturbation leads to the planet colliding with the sun or escaping to infinity.
Gauss's analysis allows also to compute the force $F(r)$ in distance $r$ to the center of a $n$-dimensional
ball with constant mass density. The divergence theorem gives $4\pi \rho |B(r)| = -|S(r)| F(r)$, where 
$|B(r)|$ is the volume of the solid ball of radius $r$ and $|S(r)|$ the surface area of the sphere.
This gives the {\bf Hook law force} $F(x)=-4\pi \rho x/n$, where $n$ is the dimension. 

\section{Catastrophe theory}

Catastrophe theory describes the singularity structure of smooth functions
$f$ on a $n$-manifold $M$ parametrized by some $r$ parameters. A basic assumption is that
{\bf configurations of interest} of the functional $f$ are {\bf critical points} of $f:M \to R$.
Especially interesting are minima, stable configurations. When changing parameters of $f$,
{\bf bifurcations}, structural changes of the critical set can happen.
Especially, minima can change their nature or disappear. In particular, the function $f_t(x_t)$,
where $x_t$ is a local minimum can change discontinuously, even if the function
$(t,x) \to f_t(x)$ is smooth. Such discontinuous changes are called {\bf catastrophes}.
The stage for Thom's theorem is a smooth function $f: \mathbb{R}^n \to \mathbb{R}^r$.
One can think of $f$ as a $r$ parameter family of functions on {\bf space} $\mathbb{R}^n$.
Let $\nabla_x= (\partial_{x_1}, \dots \partial_{x_n})$ is the {\bf gradient operator} with respect to the
space variables and $M_f = \{ (x,y) \in \mathbb{R}^n \times \mathbb{R}^r \; | \; \nabla_x f=0 \}$
is the submanifold on which points are critical. 
The space $X=C^{\infty}(\mathbb{R}^n \times \mathbb{R}^r)$ of smooth functions in space and 
parameter can be equipped with the {\bf Whitney topology}, the topology
generated by a basis which is the union of all the basis sets of $C^k$ Whitney
topologies. A basis for the later is the set of all functions for
which $f^{(j)}(x,y) \in U_j$ for all $0 \leq j \leq k$ and $U_0, \cdots U_k$ are
all open intervals. With the Whitney $C^{\infty}$ topology, $X$ is a Baire space so that
{\bf residual sets} (countable intersections of open dense sets), are dense.
The next theorem works $n=2, r \leq 6$ and for $n \geq 3$ if $r \leq 5$ \cite{Michor}
\index{bifurcation}
\index{catastrophe}
\index{Whitney topology}
\index{Baire space} 
\index{residual}

\satz{For a residual set in $X$, $M_f$ is an $r$-dimensional manifold.}

The theorem was due to Ren\'e Thom who initiated
catastrophe theory in a 1966 article and wrote \cite{ThomMorphogenesis} building on
previous work by Hassler Whitney.  More work and proofs
were done by various mathematicians like John Mather or Bernard Malgrange.
There is more to it: the restriction $X_f$ of the projection of the singularity set $M_f$ onto
the parameter space $\mathbb{R}^r$ can be classified.
Thom proved that for $r=4$, there are exactly
{\bf seven elementary catastrophes}: `fold", ``cusp", ``swallowtail",
``butterfly", ``hyperbolic umbillic", ``elliptic umbillic" and
``parabolic umbillic". For $r=5$, the number of catastrophe types is $11$.
The subject is part of {\bf singularity theory} of differentiable maps,
a theory that started by Hassler Whitney in 1955. The theory of {\bf bifurcations}
was developed by Henri Poincar\'e and Alexander Andronov.
See also \cite{Michor,StewartCatastrophe,ArnoldCatastrophe}.
It is also widely studied in the context of dynamical systems \cite{Misbah}.
\index{dynamical system}
\index{singularity theory}
\index{elementary catastrophe}

\section{Phase transition}

Given a finite simple graph $G=(V,E)$, an {\bf interaction function}
$J: E \to \mathbb{R}$ and a scalar field $h: V \to \mathbb{R}$ defines
a {\bf Hamiltonian} $H(\sigma) = \sum_{(i,j) \in E} J_{ij} \sigma_i \sigma_j
- \mu \sum_{j \in V} h_j \sigma_j$ on the set of all
functions $\sigma: V \to \{ -1,1\}$. The interpretation is that $\sigma_i$
are {\bf spin values}, $h_j$ an {\bf external magnetic field} and $J_{ij}$ is an
{\bf interaction function}. The additional parameter $\mu$ is a {\bf magnetic moment}.
The energy $H$ defines a probability measure $P$ on the set
$\Omega=\{-1,1\}^V$ of all spin configurations. It is the Gibbs-Boltzmann
distribution $P[ \{ \sigma \} ] = e^{-H(\sigma)}/Z$, where $Z$ is
is the normalization constant rendering $P$ a probability measure. One calls $Z$
the {\bf partition function} (as it is usually considered to be a function of some
of the parameters like temperature). Given a random variable =observable
$X: \Omega \to \mathbb{R}$, one is interested in the expectation ${\rm E}[X]$.
An example is $X(\sigma) = \sigma_i \sigma_j$, which leads to the correlation.
When replacing $H$ with $\beta H$, where $\beta=1/(K T)$ is an inverse temperature
parameter ($T$ is the {\bf temperature} and $K$ the {\bf Bolzmann constant}),
one can study the expectation of a random variable $X$ in dependence of $\beta$. One
writes now also ${\rm E}[X] = \langle X \rangle_\beta$ to stress the dependence on
$\beta$. In the case when $G$ is a $d$-dimensional lattice $G=[-L,L]^d$, where
two lattice points $x,y$ are connected if $\sum_k |x_k-y_k|=1$ one look at the
$L \to \infty$ {\bf van Hove limit}, where $G=\mathbb{Z}^d$. In the case $J=1,h=0$ this is the
{\bf Ising model}. As $J$ is positive, this is a {\bf Ferromagnetic situation}.
A parameter value, where a quantity like $Z_{\beta}$ or a derivative of it
changes discontinuously is called a {\bf phase transition}.
\index{Ising model}
\index{phase transition}
\index{lattice gas model}
\index{van Hove limit}
\index{Boltzmann constant}

\satz{The Ising model in two dimensions has a phase transition.}

This was first proven by Lars Onsager in 1944, who in a tour de force
gave analytical solutions. The analysis shows that there is a {\bf phase transition}.
The temperature $T$ at which this happens is called the {\bf Curie temperature}.
The one dimensional case had been solved by Ernst Ising in 1925, who got it as
a PhD project from his adviser Wilhelm Lenz. In one dimensions, there is no
phase transition. In three and higher dimensions, there are no analytical solutions.
The Ising model is only one of many models and generalizations. If the $J_{ij}$ are
random one deals with {\bf disordered systems}. An example is the {\bf Edwards-Anderson}
model, where $J_{ij}$ are Gaussian random variables. This is an example of a
{\bf spin glass model}. An other example is the {\bf Sherrington-Kirkpatrick model}
from 1975, where the lattice is replaced by a complete graph and the $J_{ij}$ define
a {\bf random matrix}. An other possibility is to change the spin to $Z_n$ or the
symmetric group (Potts) or then some other Lie group (Lattice gauge fields) and
then use a character to get a numerical value. Or one replaces the zero-dimensional
sphere $\mathbb{Z}_2$ with a higher dimensional sphere like $S^2$ and takes 
$\sigma_i \cdot \sigma_j$ (Heisenberg model). See \cite{SimonStatMechanics}.
\index{spin gas model}
\index{Curie temperature}
\index{Heisenberg model}
\index{Sherrington-Kirpatrick}
\index{Disordered system}
\index{Edwards-Anderson model}
\index{Potts model}

\section{Ceva theorem}

Given a {\bf triangle} ABC in the Euclidean plane $\mathbb{R}^2$ and a point $O$ in the interior.
For any choice of points $A'$ on the segment $BC$, any point $B'$ on the segment $AC$
and any point $C'$ on the segment $AB$, one can look at the ratios
$r(AB)= AC'/C'B$ and $r(BC)=BA'/A'C$ and $r(CA)=CB'/B'A$ in which the points bisect
the sides of the triangle. The {\bf Ceva theorem} is
\index{Ceva theorem}

\satz{ $r(AB) r(BC) r(CA) = 1$ }

The theorem is called after Giovanni Ceva who wrote it down in 1678. The result is older however:
Al-Mu'taman ibn Hud from Zaragoza  proved it already in the 11'th century. \cite{Hogendijk}.
See \cite{Russell1893}.

\section{Angle theorem}

Given a {\bf circle} $C$ in the plane $\mathbb{R}^2$. Denote by $M$ its center point.
Pick  two points $A,B$ on $C$.
If $P$ is a point on $C$, then $APB$ is constant for all $P$ in $C$
which are on the same side than $M$ with respect to the segment $AB$.
The angle $APB$ is called the {\bf inscribed angle} of the secant $AB$.
The next theorem is also called the {\bf inscribed angle theorem}.

\satz{The angle $APB$ is half the angle $AMB$. }

The theorem is believed to have been known already to Thales of Miletus who is
the first Greek mathematician known by name (624 - 546 BC).
It is usually called {\bf Thales theorem} in the special case is if $A,B$ are on a diagonal.
Then the angle $APB$ is a right angle. A consequence
of the theorem is that the opposite angles of a quadrilateral which is
inscribed in a circle add up to $\pi$.
Unlike the special case of the right angle which immediately follows
from symmetry, the full version of Thales theorem can surprise at first.
\index{Thales theorem}
\index{Inscribed angle theorem}

\section{Total curvature}

A smooth simple closed curve $C$ in $\mathbb{R}^3$ is called a {\bf knot}.
If $r(t)$ is the parametrization of $C$, then $\kappa(t) = |r'(t) \times r''(t)|/|r'(t)|^3$ is
called the {\bf curvature} of the parametrization of $r$ at the point $r(t)$.
The integral $K(C) = \int_0^{2\pi} \kappa(t) \; dt$ is the {\bf total curvature} of $r$.
We say $C$ is {\bf unknotted} if $C$ can be continuously deformed to a circle
$S=\{ x^2+y^2=1,z=0 \} = \{ r_1(t) = (\cos(t),\sin(t),0), t \in [0,2\pi] \}$
meaning that there exists a smooth function $R(t,s)$ such that $R(t,0)=r(t)$ and $R(t,1)=r_1(t)$
such that for any $s$, the curve $C_t: t \to R(t,s)$ is a simple closed curve.
\index{knot}
\index{unknot}
\index{total curvature}
\index{unknotted}

\satz{If $C$ is a knot and $K(C) \leq 4\pi$, then $K$ is unknotted.}

This is the {\bf theorem of Fary-Milnor}, proven by Fary in 1949 and Milnor in 1950.
The theorem follows also from the existence of {\bf quadrisecants}, which are
lines intersecting the knot in 4 points \cite{Denne2004}. The existence of
quadrisecants was proven by Erika Pannwitz in 1933 for smooth knots and generalized
in 1997 by Greg Kuperberg to {\bf tame knots}, knots which are equivalent to polygonal knots. 
\index{quadrisecant}

\section{Morley's Theorem}

An {\bf angle trisector} of an angle $\alpha= \angle(CAB)$ in $\mathbb{R}^2$ 
is a pair of lines $PA$, $QA$ through $A$ such that the angles $\angle(CAP),\angle(PAQ),\angle(QAB)$ 
are all equal.  Given a triangle $ABC$, we can look at the angle trisectors at each point and intersect
the adjacent trisectors, leading to a triangle $PQR$ inside the triangle. The triangle
$PQR$ is called the {\bf Morley triangle} of $ABC$.  Morley's theorem is
\index{Angle trisector}
\index{Morley triangle}

\satz{For any triangle $ABC$, the Morley triangle is equilateral.}

Morley's theorem was discovered in 1899 by Frank Morley. A short proof was given
in 1995 by John H. Conway: assume the triangle $ABC$ had angles $3\alpha,3\beta,3\gamma$
so that $\alpha+\beta+\gamma=\pi/3$. Start with an equilateral triangle $PQR$ of length $1$.
Build three triangles $PQA$ with angles $\beta+\pi/3,\alpha,\gamma+\pi/3$,
$QCA$ with angles $\alpha+\pi/3,\gamma,\beta+\pi/3$ and a triangle
$RBQ$ with angles $\gamma+\pi/3, \beta,\alpha+\pi/3$. Then fill in three other triangles
$ACQ, CBR,BAP$ with angles $\alpha,\gamma,\beta+2\pi/3$ and $\gamma,\beta,\alpha+2\pi/3$
$\beta,\alpha,\gamma+2\pi/3$. These triangles fits together to a triangle of the shape
$ABC$. See \cite{EricksonBeautiful}.

\section{Rising sun lemma}

Given an interval $[a,b]$, the space $C([a,b])$ denotes the vector space
of all continuous functions on $[a,b]$. For $g \in C([a,b])$, we say the
set $E(g)=\{ x \in (a,b) \; | \; g(t)>g(x) \; {\rm for} \;  t>x\}$
has the {\bf rising sun property} if
$E$ is open, and $E$ is empty if and only if $g$ is decreasing and if
not empty, then $E$ can be written as $E=\bigcup_n (a_n,b_n)$
with pairwise disjoint intervals with
$g(a_n) \leq g(b_n)$. See \cite{Berberian1998}.
\index{Rising sun property}

\satz{$f \in C([a,b],\mathbb{R})$ has the rising sun property.}

The theorem is due to F. Riesz.
The name ``rising sun lemma" appeared according to \cite{Berberian1998}
first in  \cite{AsplundBungart}. The picture is to draw the graph of the function $f$.
If light comes from a distant source parallel to the x-axis,
then the intervals $(a_n,b_n)$ delimit the hollows that remain in the shade at the
moment of sunrise.  The lemma is used in real analysis to prove that every monotone
non-decreasing function is almost everywhere differentiable.
\index{Rising sun lemma}

\section{Uniform continuity}

Uniform continuity is a stronger version of continuity. But unlike continuity, which
is defined for maps between topological spaces, uniform continuity needs
more structure like a metric spaces or more generally a topological space 
with a {\bf uniform structure}. Given two metric spaces $X$ and $Y$, a function 
$f:X \to Y$ is called {\bf continuous} if $f^{-1}(U)$ is open for every open $U$ in $Y$.
A function $f$ is called {\bf uniformly continuous} if there exists
a sequence of numbers $M_n \to 0$ such that for every positive $n \in \mathbb{N}$, the
condition $d(x,y) \leq 1/n$ implies that $d(f(x),f(y)) \leq M_n$.

\satz{For compact $X$, continuous implies uniformly continuous.}

The theorem is due to Eduard Heine and Georg Cantor. Heine is known also for the
Heine-Borel theorem which states that in Euclidean spaces, the class of closed and bounded
sets agrees with the class of compact sets. The proof of the {\bf Heine-Cantor theorem} uses the
{\bf extreme value theorem} assuring that a continuous function on a compact space $X$
achieves a maximum.
Look for every $n$ and every $x$ at the minimal $M_n(x)$ such that if $|x-y| \leq 1/n$,
then $|f(x)-f(y)| \leq M_n(x)$. Now $M_n(x)$ is non-negative and finite and depends
continuously on $x$. By the extremal value theorem there is a maximum. We call it $M_n$.
This assures now that if $|x-y| \leq 1/n$, then $|f(x)-f(y)| \leq M_n$.
The {\bf Bolzano-Weierstrass} or sequential compactness theorem assures that
a bounded sequence in $\mathbb{R}^n$ has a convergent subsequence. This is used
in the intermediate value theorem assuring that if $f(a)<0$ and $f(b)>0$, then there is
an $x$ with $f(x)=0$. The Heine-Cantor theorem together with the intermediate value
theorem assures that continuous functions are Riemann integrable.
The additional {\bf uniform structure} or {\bf metric structure} is also necessary when
defining completeness in the sense that every Cauchy sequence converges.
Completeness is not a property of topological spaces: $(0,1)$ is not complete but
$\mathbb{R}$ is complete even so the two spaces are homeomorphic.
\index{uniform structure}
\index{Heine-Cantor theorem}
\index{Bolzano-Weierstrass theorem}

\section{Jordan normal form}

A $n \times n$ matrix $A$ is {\bf similar} to an other $n \times n$ matrix $B$
if there exists an invertible $n \times n$ matrix $S$ such that $B=S^{-1} A S$.
A matrix is in {\bf Jordan normal form} (also called {\bf Jordan canonical form}) if it is
block diagonal, where each block is a {\bf Jordan block}. A $m \times m$ matrix $J$
is a {\bf Jordan block}, if $J e_1=\lambda e_1$, and $J e_k = \lambda e_k + e_{k+1}$ for
$k=2,\dots, m$. An example of a $3 \times 3$ Jordan block matrix is
$J = \left[ \begin{array}{ccc} 3 & 1 & 0 \\ 0 & 3 & 1 \\ 0 & 0 & 3 \end{array} \right]$.
In other words, $A$ is of the form $A= \lambda 1 + N$, where $N$ is nilpotent: $N^m=0$
and more precisely only has $1$ in the super diagonal above the diagonal.

\satz{Every $n \times n$ matrix is similar to a matrix in Jordan normal form.}

Up to the order of the Jordan blocks, the Jordan normal form is unique. If each
Jordan block is a $1 \times 1$ matrix, then the matrix is called {\bf diagonalizable}. The
{\bf spectral theorem} assures that a normal matrix $A A^* = A^* A$ is diagonalizable.
Not every matrix is diagonalizable as the
{\bf shear matrix} $A=\left[ \begin{array}{cc} 1 & 1 \\ 0 & 1 \end{array} \right]$,
a $2 \times 2$ Jordan block, shows. 
The theorem has been stated first by Camille Jordan in 1870. For history, see
\cite{Brechenmacher}. The {\bf Jordan-Chevalley}
generalization states that over an arbitrary perfect field, a matrix is similar to $B+N$,
where $B$ is {\bf semi-simple} and $N$ is {\bf nilpotent} and $BN=NB$. (See \cite{Humphreys} page 17).
A matrix $B$ is called {\bf semi-simple} if every $B$-invariant linear subspace $V$ has a
complementary $B$-invariant subspace. For algebraically closed fields, semi-simple is equivalent
to be conjugated to a diagonal matrix. To the condition on the field: a field $k$ is called {\bf perfect}
if every irreducible polynomial over $k$ has distinct roots.
\index{perfect field}
\index{Jordan-Chevalley decomposition}
\index{Jordan normal form}
\index{Jordan block}
\index{block diagonal}

\section{Hippocrates theorem}

The Hippocrates theorem dealing with the {\bf lunes of Hippocrates}
or the {\bf lunes of Alhazen} is a theorem in planar geometry:
given a triangle $ABC$ in $\mathbb{R}^2$ with right angle $\beta$ at $B$, 
one can draw the circles with diameter $AC,AB$ and $BC$ centered
at the midpoints $(A+C)/2, (A+B)/2$ and $(B+C)/2$. They
define two ``moon-shaped" regions $U,V$ bounded by circles called the 
{\bf lunes}. 
\index{lunes of Hippocrates}
\index{Hippocrates Theorem}

\satz{The area of $U$ plus the area of $V$ is the area of the triangle.}

The proof directly follows from Pythagoras by relating the areas
of half discs and triangle.
The result is remarkable as it was historically the first attempt
for the {\bf quadrature of the circle}. The lunes are bound by circles,
while the triangle is bound by line segments. The theorem does the
{\bf quadrature of the lunes}. Hippocrates of Chios lived from about 470 to 410 BC. 
\index{Hippocrates theorem}
For history see \cite{RouseBallHistory} page 37. 

\section{Fermat-Hamilton principle}

A point $x$ is called a {\bf critical point} of a differentiable function $f: \mathbb{R}^m \to \mathbb{R}$,
if $\nabla f(x)=0$, where $\nabla f$ is the {\bf gradient} of $f$. A point $x_0$ is called a {\bf local maximum}
of $f$ if there exists $r>0$ such that $f(x) \leq f(x_0)$ for all $|x-x_0|<r$. The local maximum does not
have to be isolated. For a constant function for example, every point is a local maximum. The local maximum also does
not have to be a {\bf global maximum}. The function $f(x) = x^4-x^2$ has a local maximum at $x=0$ but this is
not a global maximum because $f(2)>f(0)$. 
\index{Fermat's principle}
\index{Hamilton principle} 
\index{critical point}  
\index{local maximum}
\index{global maximum} 

\satz{If $x_0$ is a local maximum of $f$, then $\nabla f(x_0)=0$. }

This generalizes to the {\bf calculus of variations}, where $\nabla f$ is replaced by the 
{\bf variation}. In the case when $f(x) = \int_a^b L(x(t),x'(t)) \; dt$ is a function 
on the space of curves $[a,b] \to \mathbb{R}^n$ (one calls this then a functional or action
functional) then we an look at the problem to minimize the action. In that case, the 
gradient is $\delta S = L_x(x(t),x'(t)) - \frac{d}{dt} L_{x'}(x(t),x'(t)) = 0$. This so called
{\bf Hamilton principle} can be seen as a generalization of the Fermat principle to infinite dimensions.
The equations $\delta S=0$ are called the {\bf Euler-Lagrange equations} or 
{\bf Lagrange equations of the second kind}. They are the starting point of {\bf Lagrangian mechanics.}
Fermat's original paper deals with the single variable situation but the higher dimensional 
situation is similar. Fermat in some sense already looked at the action principle which 
is the situation to minimize the arc length of a path in a medium with two different properties
like water and air. In that case the shortest path is described by the {\bf Fermat law} or
{\bf Fermat's principle}.
\index{Calculus of variations}
\index{Euler-Lagrange equations}
\index{Fermat law}
\index{Fermat's principle}
\index{Lagrange equations}
\index{Calculus of variations}

\section{Alternating sign}

An {\bf alternating sign matrix} is a square matrix with entries in $\{ 0,1,-1 \}$
such that the sum of each row and column is $1$ and the nonzero entries in
each row and column alternate in sign.

\satz{The number of $n \times n$ alternating sign matrices is $\prod_{k=0}^{n-1} \frac{(3k+1)!}{(n+k)!}$.}

The numbers $\prod_{k=0}^{n-1} (3k+1)!/(n+k)!$ are known as the
{\bf Robbins numbers} or {\bf Andrews-Mills-Robbins-Rumsey numbers} and are the
integer sequence A005130 \cite{A005130}.
The alternating sign conjecture was popularized by David Robbins in \cite{Robbins1991}.
The theorem was proven by Doron Zeilberger in 1994 \cite{Zeilberger1996}.
A short proof was given by Greg Kuperberg in 1996 \cite{Kuperberg1996}. A book about it
is \cite{BressoudProofsConfirmations}.
\index{Robbins numbers}
\index{Alternating sign matrix}
\index{Alternating sign conjecture}

\section{Combinatorial convexity}

A finite set $P$ of points in $\mathbb{R}^d$ is
called {\bf $r$-convex}, if there is a partition of $P$ into $r$ sets
such that their convex hulls intersect simultaneously in a non-empty set.
{\bf Tverberg's theorem} states:

\satz{A set of $(r-1)(d+1)+1$ points in $\mathbb{R}^d$ is $r$-convex.}

The decomposition of $P$ into $r$ subsets is called the {\bf Tverberg partition}.
In the one-dimensional case $d=1$, the theorem assures that $2r-1$ points on 
the line are $r$-convex. For $r=3$ for example, this means that $5$ points are
$3$-convex. If the points are arranged $x_1<x_2<x_3<x_4<x_5$, the Tverberg partition 
$\{ x_1,x_4\},\{x_2,x_5 \},\{x_3\}  \}$.
For $r=2$, it implies {\bf Radon's theorem} which tells that 
$d+2$ points in $\mathbb{R}^d$ can be partitioned into 2 sets
whose convex hulls intersect. For example, $4$ points $\{x_1,x_2,x_3,x_4\}$ 
in $\mathbb{R}^2$ can be partitioned into two sets such that their convex
hull intersect. Indeed, the 4 points define a quadrilateral and the partition
$\{ \{x_1,x_3\}, \{x_2,x_4\} \}$ define the two diagonals of the quadrilateral.
The theorem has been proven by {\bf Helge Tverberg} in 1966.
See \cite{Tverberg,Tverberg50}.  \index{Tverberg's theorem}
\index{Radon theorem}
\index{Tverberg's theorem}
\index{Tverberg partition}
\index{Combinatorial convexity}

\section{The Umlaufsatz}

Let $r$ be a continuously differentiable closed curve in $\mathbb{R}^2$. If $r(t)$ is a
parametrization for which the speed is $1$, we have $r'(t) = (\cos(\alpha(t)),\sin(\alpha(t)))$ and a 
{\bf signed curvature} $\kappa(t) = \alpha'(t)$. If $[0,2\pi]$ is the parameter interval,
then $K=\int_0^{2\pi} \kappa(t) \; dt$ is the {\bf total curvature}.
The Hopf Umlaufsatz is:

\satz{For $r \in C^1$, the total curvature of a plane curve is $2\pi$. }

The paper was proven in 1935 by Heinz Hopf \cite{hopf35} using a
homotopy proof: define $f(s,t)$ as the argument of the line through $r(s)$ and 
$r(t)$ or continuously extend it $s=t$ as the argument of the tangent line. The
direct line from $(0,0)$ to $(1,1)$ in the parameter $st$-plane gives a total
angle change of $n 2\pi$ where $n$ is an integer. Now deform the curve from
$(0,0)$ to $(1,1)$ so that it first goes straight from $(0,0)$ to $(0,1)$,
then straight from $(0,1)$ to $(1,1)$. Both lines produce a deformation 
of $\pi$ and show that $n=1$.
The theorem can be generalized to a Gauss-Bonnet theorem for 
planar regions $G$. The total curvature of the boundary is $2\pi$ times the
Euler characteristic of $G$. For a discrete version, see \cite{elemente11}.

\index{Hopf Umlaufsatz}
\index{Gauss Bonnet}

\section{Frobenius determinant}

The {\bf Frobenius determinant theorem} tells how
the determinant of the ``multiplication table matrix"
factors into irreducible polynomials:
if $G = \{g_1,\dots,g_n\}$ is a {\bf finite group} and $x_i$ is a
variable associated to the group element $g_i$, then
the matrix $A_{ij} = x_{g_i g_j}$ satisfies

\satz{
${\rm det}(A) = \prod_{j=1}^r p_j(x_1, \dots, x_n))^{d_j}$
}

Here, $d_j={\rm deg}(p_j)$ and $r$ is the number of conjugacy classes of $G$.
For an Abelian group $G$, there are $n$ conjugacy classes. The theorem had been
conjectured in 1896 by Richard Dedekind. Frobenius proved it. See
\cite{FrobeniusContext,ConradOriginRepresentationTheory}.
\index{Frobenius determinant}
\index{Frobenius determinant theorem}


\section{K\"onig's theorem}

\paragraph{}
A {\bf matching} $M$ in a {\bf finite simple graph} $G=(V,E)$ with 
{\bf vertex set} $V$ and {\bf edge set} $E$
is a subset $M$ of the edges $E$ 
in which no two edges have a common vertex. A {\bf vertex cover} $C$ is a set of vertices
such that $\bigcup_{x \in C} S(x) = V$, where $S(x)$ is the unit sphere of a vertex $x$. 
A {\bf bipartite graph} is a graph for which $V=V_1 \cup V_2$
can be partitioned into two disjoint sets $V_1,V_2$ such that all edges connect vertices from 
different sets. K\"onig's theorem, from 1931, also known as 
{\bf K\"onig-Egev\'ary theorem} is:

\satz{For bipartite $G$, matching number = vertex cover number.}

The {\bf vertex cover problem} is the problem to find the vertex cover number
is a classical {\bf NP-complete problem}. For example, for a cyclic graph $G=C_{10}$ with $2n$ vertices 
$\{1,2,3, \cdots, 2n\}$ (which is an example of a bipartite graph), 
the set $C=\{ 2,4, \cdots 2n \}$ is a minimal vertex cover. The edges
$M=\{ (1,2), (3,4), \cdots (2n-1,2n)\}$ are a maximal matching. The example of an odd cyclic graph like 
$C_{9}$ (which is not bipartite) already shows that the bipartite condition is necessary:
for $C_9$, the set $\{1,3,5,7,8\}$ is a minimal cover and 
$M=\{ (1,2), (3,4), (5,6), (7,8) \}$ is a maximal matching.  \\

The origin of the theorem is attributed to D\'enes K\"nig, who proved it in 1931 and 
wrote a precursor paper in 1916, where he proved that a regular (constant vertex degree) 
bipartite graph has a perfect matching (a matching which covers all vertices). 
For a proof, see \cite{Diestel} (Chapter 2).
\index{K\"onig's theorem}
\index{vertex cover problem}
\index{vertex Cover}
\index{minimal Matching}
\index{matching}
\index{bipartite graph}
\index{NP complete}

\section{Polynomial Ergodic theorems}

Birkhoff's ergodic theorem stating that $S_{n,f}(x) = \frac{1}{n} \sum_{k=1}^n f(T^k x)$
converges for $n \to \infty$ point-wise for $\mu$ almost every $x$
for an automorphism $T$ of a probability space $(X,\mathcal{A},\mu)$ and a function
$f \in L^p(X)$ with $1 \leq p < \infty$ has been generalized in 1988 by Jean Bourgain
\cite{Bourgain1988} to {\bf polynomial averages}
$S_{P,n,f}(x) = \frac{1}{n} \sum_{k=1}^n f(T^{P(k)} x)$, where
$P$ is a polynomial with integer coefficients.

\satz{$S_{P,n,f}(x)$ converges point-wise almost everywhere if $p>1$.}

Bourgain proves in \cite{Bourgain1988} first a maximal ergodic theorem and extends it also to
$\mathbb{Z}^d$ actions generated by $d$ commuting transformations.
The starting point is that for $f \in L^2(X,\mu)$,
there is for any integer $t$ a bound $|S_{n^t,n,f}|_2 \leq C |f|_{2}$. This implies for
example that $\frac{1}{n} \sum_{k=0}^{n-1} f(x+m^t \alpha) \to \int_0^1 f(x) \; dx$ for
any irrational $\alpha$ and any bounded measurable function. The case $t=2$
leads to results to sums like $\frac{1}{n} \sum_{k=0}^{n-1} e^{\pi i k^2 \alpha}$
which relates to {\bf Gauss sums} $S(q,a) = \frac{1}{q} \sum_{k=0}^{q-1} e^{2\pi i k^2 a/q}$.
One can for example estimate
$\sum_{k=0}^{n-1} e^{2\pi i k^2 \alpha} \leq C (n/\sqrt{q} + \sqrt{n \log(q)} + \sqrt{q \log(q)})$.
\cite{Bourgain1988}.  The case $p=1$ is known to fail \cite{BuczolichMauldin}.
The results have been generalized to {\bf correlation expressions} like
$S_{n,f,a,b}(x) = \frac{1}{n} \sum_{k=0}^{n-1} f(T^{an}x) g(T^{bn}x)$ for integers $a,b$
where $f \in L^p,g \in L^q$ with $1 < p,q \leq \infty$ and $1/p+1/q<3/2$
\cite{Demeter2007,Lacey2000} and to {non-conventional bilinear polynomial averages}
$\frac{1}{n} \sum_{k=0}^{n-1} f(T^n x) g T^{P(n)}x)$ \cite{KrauseMirekTao2020},
where $P$ is an integer polynomial of degree $d \geq 2$ and
$f \in L^p,g \in L^q$ with $1 < p,q \leq \infty$ and $1/p+1/q \leq 1$.
\index{Polynomial averages}
\index{Bourgain's theorem}
\index{Polynomial ergodic theorem}
\index{Gauss sum}

\section{Wantzel's theorem on Angle trisection}

A classical problem in geometry asks to {\bf trisect} an angle using
an unmarked {\bf straightedge} (ruler) and {\bf compass} only. The insistence on
restricting constructions to ruler and compass has been proposed already by Euclid
and Archimedes already knew how to solve the problem using a marked straightedge
meaning that one has additionally to the constructed points also an additional real
number to work with. One can trisect and angle using an additional curve
like an {\bf Archimedean spiral} \cite{OnSpirals} given in polar coordinates as $r=\theta$.
In that case, the trisecting the radius $r = \sqrt{x^2+y^2}$ of a given point
$(x,y)= (r \cos(\theta), r \sin(\theta))$ gives the angle $\theta/3$
by intersecting the circle of radius $r/3$ with the spiral $r = \theta$.
More generally, a curve which can be used to trisect an angle is called a
{\bf trisectrix}.
\index{ruler and compass}
\index{straightedge and compass}

\satz{One can not trisect a general angle with ruler and compass.}

The theorem follows from Galois theory. An angle $\alpha$ can be trisected if and only
if the polynomial $5x^3-3x-\cos(\alpha)$ is reducible over the field $\mathbb{Q}(\cos(\alpha))$.
The angle $\alpha=60^{\circ}= \pi/3$ for example is not trisectable.
The first proof of the impossibility of trisecting an arbitrary angle was given by
Pierre Wantzel in 1837. Wantzel also solved there the problem of doubling the cube
and characterized {\bf constructable regular $n$-gons} as the ones with
$n=2^k p_1 \cdots p_k$ with distinct {\bf Fermat primes} $p_k = 2^{2^{m_k}}+1$.
Bieberbach realized in 1932 that every cubic construction can be traced back to the
trisection of an angle and the extraction of the third root. \cite{Biberbach1932}.
This has been formulated more precisely by Gleason in 1988 \cite{Gleason1988}
who states in in that article as Theorem 1:
a real cubic equation can be solved geometrically using ruler
and compass and angle-trisector if and only if its roots are all real.
Gleason shows from this also that a {\bf regular $n$-gon} can be constructed
by ruler, compass and angle-trisector if and only if the prime factorization
of $n$ has the form $2^r 3^s p_1 p_2 \cdots p_k$ with $k \geq 0$,
where all primes $p_k>3$ are distinct and have the form $2^n 3^m + 1$.
An example is $p=13=2^2  3 + 1$. The corresponding 13-gon is called the
{\bf triskaidecagon} for which Gleason gives a concrete construction using
that $2 \cos(2\pi k/13)$ are the roots of the polynomial
$x^6+x^5-5x^4-4x^3+6x^2+3x-1$ which factors over $\mathbb{Q}(\sqrt{13})$
because with $\lambda = (1-\sqrt{13})/2, \overline{\lambda} = (1+\sqrt{13})/2$
one can write it as $(x^3-x-1+\lambda (x^2-1)) (x^3-x-1-\overline{\lambda}(x^2-1))$,
where the first factor has the root $2 \cos(2\pi/13)$.
For more on angle trisector and especially many failed attempts, see \cite{Trisectors}.
\index{Triskaidecagon}
\index{Angle trisection}
\index{trisectrix}
\index{Archimedean spiral}

\section{Preissmann's Theorem}

Let $\mathcal{M}_{-}$ denote the class of compact negatively curved Riemannian manifolds $M$.
{\bf Negative curvature} means that all sectional curvatures of $M$ are negative
everywhere. Let $\pi_1(M)$ denote the {\bf fundamental group} of $M$. 
For positively curved manifolds, the theorem of Synge shows that the fundamental
group $\pi_1(M)$ is finite; it can be trivial like for a sphere $\mathbb{S}^d, d >1$ or
be a finite group like $\pi_1(M)=\mathbb{Z}_2$ for the projective space $M=\mathbb{RP}^d$ 
for $d>1$. For a flat manifold like the torus $\mathbb{T}^d$, the fundamental group can 
already be the infinite group $\mathbb{Z}^d$. This changes for negative curvature. Preissmann
showed that if $\pi_1(M)$ is cyclic, then there is only one closed geodesic and
that there is maximally one geodesic in each homotopy class of closed curves in $M$. Here
is Preissmann's theorem which deals with non-trivial subgroups $G$ of $\pi_1(M)$ meaning
that $G$ should not be the trivial $1$-point group.

\satz{
If $G \subset \pi_1(M)$ for $M \in \mathcal{M}_-$ is Abelian then $M=\mathbb{Z}$.
}

A consequence is that the torus $\mathbb{T}^n$ can not admit a Riemannian metric
of negative sectional curvature. Preissmann gives in his paper also the corollary that the 
product of two negatively curved Riemannian manifolds can not carry a metric with negative curvature. 
An analogue result for positive curvature is not known. The famous product conjecture of Heinz Hopf asks
whether the product manifold $S^2 \times S^2$ can carry a metric of positive curvature
(see \cite{YauSeminar1982}). Preissmann who was born at Neuch\^atel in Switzerland
in 1916, went to school at La Chaux-de-Fonds. He studied mathematics from 1934 to 1938 at
the ETH and worked there until 1942 as an assistant to Kollros and Gonseth, writing his thesis under
the guidance of Heinz Hopf, where the theorem appears \cite{Preissmann1942}.
Preissmann later later got interested in hydraulic computations given the Swiss boom of               
hydro-power developments. After having been an actuary in a life insurance from 1942-1946, he
joined VAWD until 1958, then led the Department of Mathematical Methods of the
hydraulics laboratory SOGREAH in Grenoble from 1958-1972, retiring in 1981. 
See \cite{PreissmannBio}.
\index{Preissmann's theorem}
\index{Negative curvature manifolds}
\index{Riemannian manifolds}

\section{Killing-Hopf theorem}

A {\bf space form} $M$ is a quotient $A/G$, where $A$ is a sphere, an
Euclidean space or hyperbolic space and $G$ is a group acting freely 
$(g x = x$ is only possible for g=1) and discontinuously.
The later means that for any compact $K$ in $M$, and any $g \in G$
the set $g K  \cap K$ is finite). A Riemannian manifold has 
{\bf constant curvature} if all sectional curvatures are the same
everywhere. The {\bf Killing-Hopf theorem} is:
\index{space form}
\index{constant curvature}
\index{Killing-Hopf theorem}

\satz{Constant curvature manifolds are space forms.}

The theorem is due to Wilhelm Killing from 1891 \cite{Killing1891} 
and Heinz Hopf 1926 \cite{Hopf1926}.
See \cite{Wolf2011} for the topic of constant curvature manifolds. 

\section{Ballot theorem}

\paragraph{}
Let $X_j$ be independent identically distributed random variables taking values
$e_k=$ $(0,\cdots,0$, $1$, $0,\cdots 0)$ in $\mathbb{Z}^d$ with probability $p_k$.
If $p_1 > \cdots > p_d$, we can look at the multi-dimensional {\bf random walk} 
$S_n = \sum_{k=1}^n X_k$. What is the probability that the walk starting at $0$
remains in open cone $Q=\{ x_1>x_2 \cdots > x_d \}$ at all positive times? The answer is
given by the {\bf Ballot theorem}. It expresses the probability as a
{\bf van der Monde determinant}:
\index{Ballot theorem}
\index{Van der Monde determinant}

\satz{ ${\rm P}[ S_n \in Q, \forall n>0 ] = \prod_{i<j} (p_i-p_j)$. }

The case $d=2$ is the classical result is due to Joseph Bertrand
\cite{Bertrand1887} and appears in virtually every probability textbook like
\cite{Feller} who also points out that the theorem has been proven earlier by
William Whitworth \cite{Whitworth1886} who looked at the problem in a different context like
the problem of counting the number of weak orderings. 
The historical context is voting and explains the etymology of the theorem \cite{AddarioBerryReed}:
if candidate $A$ gets $m$ votes and candidate $B$ gets $n$ votes,
then the probability that during the counting process $A$ always has more votes than
$B$ is $(n-m)/(n+m)$. If $P_{n,m}$ counts the number of paths always favorable for $A$, then
the recursion $P_{n+1,m+1}=P_{n+1,m} + P_{n,m+1}$ holds. As
Binomial coefficients $B_{n,m}$ and so $D_{n,m}=B_{n,m}-P_{n,m}$ satisfies the same
recursion, it can be shown by induction that $D_{n,m} = 2m B_{n,m}/(n+m)$, leading to the result.
The multidimensional result has been studied in \cite{Zeilberger1983,GesselZeilberger1992}.

\section{Poincar\'e-Hopf}

If $F$ be a smooth vector field on a compact $n$-manifold $M$ with
finitely many equilibrium points $F(x)=0$. The {\bf index} $i_F(x)$ of $F$ at
such an equilibrium point $x_k$ is defined as the {\bf degree} of the map
$u \in S(x) \to  F(u)/|F(u)| \in \mathbb{S}^{n-1}$, where
$S(x)$ is the boundary of a small enough ball containing $x$ in the interior.
Let $\chi(M)$ denote the {\bf Euler characteristic} of $M$.
The {\bf Poincar\'e-Hopf index formula} links the topological quantity
$\chi(M)$ with the analytic index sum:

\satz{ 
$\sum_{x, F(x)=0} i_F(x) = \chi(M)$
}

The formula can be used to compute the Euler characteristic of a
manifold $M$: just construct a smooth vector field $F$ with finitely many equilibrium
points and add up their indices. For example, on the $n$-torus $M=\mathbb{T}^n$,
there is the constant vector field $F(x) = v$ without equilibrium points.
Therefore $\chi(M)=0$. On a $2n$-sphere embedded as $\{ |x|=1 \}$ in $\mathbb{R}^{2n+1}$
there are circles in $SO(2n+1,\mathbb{R})$ that have two fixed points of index $1$,
the Euler characteristic is $2$. On a $2n+1$ sphere $M$,
there are circles in $SO(2n+2,\mathbb{R})$ without fixed points so that $\chi(M)=0$.
A special case is if $f$ is a Morse function on $M$, where $F=\nabla f$, the
equilibrium points of $F$ are the critical points of $f$. In that case
$i_F(x)= (-1)^{m(x)}$, where $m(x)$ is the Morse index, the number of negative
eigenvalues of the Hessian $d^2 f(x)$.
Poincar\'e wrote the first article in 1885 \cite{Poincare1885}.
Then appeared Hopf's articles \cite{HopfCurvaturaIntegra} for hypersurfaces and
\cite{HopfVectorfields} for vector fields.
\index{Poincar\'e-Hopf}
\index{index}
\index{Euler characteristic}
\index{degree}

\section{Sampling theorem}

\paragraph{}
Let $S$ be the {\bf Schwartz space} of complex-valued functions in $C^{\infty}(\mathbb{R},\mathbb{C})$
such that $||f||_{m,n} = \sup_{x \in \mathbb{R}} |x^m f^{(n)}(x)| < \infty$. 
The {\bf Fourier transform} $\hat{f}$ of $f \in S$ is defined as
$\hat{f}(k)$ = $\frac{1}{\sqrt{2\pi}}$ $\int_{\mathbb{R}} f(x) e^{-ix k} \; dx$. 
The {\bf Nyquist-Shannon sampling theorem} tells that if
$\hat{f}$ supported on $[-\pi,\pi]$. Then $\{ f(n), n \in \mathbb{Z} \}$ determines $f$:
\index{Schwartz space}

\satz{$f(t) = \frac{1}{\pi} \sum_{n=-\infty}^{\infty} f(n) {\rm sinc}(\pi (n-t))$}

It uses the {\bf sinc} function ${\rm sinc}(x)=\sin(x)/x$.
The explicit reconstruction formula is also known as the 
{\bf Whittaker-Shannon} interpolation formula as
the formula appeared in the book \cite{Whittaker35}. Whitaker has already found
that formula in 1915 while \cite{Shannon48} which is the start of information theory. 
The result was also spearheaded by Nyquist 1928. We followed \cite{Sternberg2019}.
\index{Shannon sampling theorem}
\index{Nyquist-Shannon sampling theorem}

\section{Peter Weyl theorem}

Let $G$ be a compact {\bf topological group}
and let $C(G,\mathbb{C})$ denote the Banach space of continuous complex-valued   
functions on $G$, equipped with the uniform norm
$|f|_{\infty} = \max_{x \in G} f(x)$.
Denote by $\pi: G \to Gl(V)$ a {\bf group representation}
of $G$, where $V$ is a complex vector space.
This means $\pi(g h) = \pi(g) \pi(h)$ for any $g,h \in G$.
A {\bf matrix coefficient} of $G$ is a map $\phi: G \to \mathbb{C}$
which has the form $L(\pi(x))$, where $\pi$ is a representation of $G$ and
where $L$ is a {\bf linear functional} $Gl(V) \to \mathbb{C}$.
An example of a linear functional on $Gl(V)$ is the trace ${\rm tr}(A)$
or an other linear combination of matrix entries $A_{ij}$ (explaining
the name). The {\bf Peter-Weyl theorem} is:

\satz{The set of matrix coefficients is dense in $C(G,\mathbb{C})$.}

This implies that the matrix coefficients are also dense in
the Hilbert space $L^2(G,\mu)$ defined by the Haar measure
$\mu$ on $G$. If $\pi$ is a unitary representation on a Hilbert space
$H=(V,( \cdot, \cdot))$, one can write $\pi$ as a direct sum of irreducible unitary
representations and the matrix elements give an explicit
{\bf orthornormal basis} in $L^2(G)$: make a list of
representatives of the isomorphism classes $\pi$ of irreducible
unitary representations of $G$, then take the basis elements
$\sqrt{d(\pi)} \pi(g)_{ij}$, where $d(\pi)$ is the degree of the
representation. The theorem was proven by Fritz Peter and Herman Weyl 
in 1927 \cite{PeterWeyl1927}.
The result follows from the Stone-Weierstrass theorem if $G$ is a matrix
group and especially for Lie groups which are known to be matrix groups.
Not much seems to be known about Fritz Peter (1899-1949)
whose residence is in the paper \cite{PeterWeyl1927} given as Karlsruhe and 
to whom Weyl refers as ``his student". The book \cite{Hawkins2000}
states that Peter got a doctorate in G\"ottingen in 1923
(with the title: \"Uber Brechungsindizes und Absorptionskonstanten des
Diamanten zwischen {$\lambda$} 644 und 266), under the guidance of
Max Born \cite{Schirrmacher2019}. A conference proceeding lists him later as a 
teacher at a school in Schloss Salem near \"Uberlingen in Germany.
See \cite{Hawkins2000}
\index{Peter-Weyl theorem}
\index{Irreducible representation}
\index{Matrix elements}

\section{Kruskal-Katona theorem}

\paragraph{}
A {\bf finite abstract simplicial complex} $G$ is a finite set of non-empty sets
which is closed under the operation of taking finite non-empty subsets. The {\bf dimension}
of a set $x$ is the $|x|-1$, where $|x|$ is the {\bf cardinality} of $x \in G$. The
{\bf $f$-vector} $f=(f_0,f_1, \cdots, f_{d}) \in \mathbb{N}^{d+1}$ counts the number $f_k$ of
$k$-dimensional sets $x$ in $G$. If $n=B(n_i,i) +B(n_{i-1},i-1) + \cdots + B(n_j,j)$ 
is the {\bf Binomial development} of $n$ at level $i$, 
define $n^{(i)} = B(n_i,i+1)+ \cdots + B(n_j,j+1)$.
The {\bf theorem of Kruskal-Katona} characterizes the possible $f$-vectors which 
simplicial complexes can have:

\satz{$f$ is the $f$-vector of a complex if and only if $f_i \leq f_{i-1}^{(i)}$.}

The theorem was found by Joseph Kruskal (1963) 
(a brother of Martin Kruskal known in the context of solitons)
and Gyula Katona (1968). See \cite{Frankl1983}.
Because the result is sharp, it is often mentioned in the context of 
{\bf extremal set theory}.
The result implies the {\bf Erdoes-Ko-Rado theorem} \cite{ErdoesKoRado1961}.
The later is the result about a finite set $G$ of sub-sets of $\{1, \dots, n\}$ 
of cardinality $k$ such that each pair has a non-empty 
intersection and $n>2k$, then the number of sets in $G$ is less or equal
than the Binomial coefficient $B(n-1,k-1)$. A bit easier to state is 
the following special case of the Kruskal-Katona theorem formulated by Lovasz: 
if $f_k = B(m,i)$, then $f_{k-r} \geq B(m,i-r)$ for any $r \geq 0$. The fact that these
statements are sharp can be seen when looking at the 
complete complex $G$ consisting of all non-empty subsets of                   
$\{1,2, \dots, n\}$, where $f_k(G) = B(n,k-1)$ which means $m=n,i=k-1$ in the
above notation. More specifically, if $G=\{\{1\},\{2\},\{3\},\{4\}$, 
$\{1,2\},\{1,3\},\{1,4\},\{2,3\},\{2,4\},\{3,4\}$,
$\{1,2,3\},\{1,2,4\},\{1,3,4\},\{2,3,4\}$, $\{1,2,3,4\}\}$, where the $f$-vector
is $(4,6,4,1)$ we have the situation of Lovasz. 
\index{Kruskal-Katona Theorem}
\index{Erdoes-Ko-Rado Theorem}
\index{Extremal set theory}

\section{Computational complexity}

An {\bf NP decision problem} has a {\bf probabilistically checkable proof} (PCP) if
given any probability $p<1$, there exists a polynomial $f$ such 
that every mathematical proof of length $n$ can be     
rewritten with a proof of length $f(n)$ and that can be formally verified 
with accuracy $p$. The later means that one can formally verify         
$p*f(n)$ letters of the proof of an NP decision problem. 
Examples of NP hard decision problems are the {\bf traveling salesperson problem},
the {\bf knapsack decision problem}, {\bf clique problems in graphs}.
The PCP theorem is:

\satz{Every NP decision problem has probabilistically checkable proof.}

To cite \cite{Dinur2007}: {\it
"Every language in NP has a witness format that can be checked probabilistically
by reading only a constant number of bits from the proof. 
The celebrated equivalence of this theorem and inapproximability of 
certain optimization problems, due to Feige et al. 1996, has placed the 
PCP theorem at the heart of the area of inapproximability."}

The theorem has been proven by various mathematicians starting 
with 1990 by Laszlo Babai, Lance Fortnow and Carsten Lund. 
More work was done by Sanjeev Arora and Shmuel Safra from 1998.
The theorem is considered one of the most important results in complexity theory
as it shows that certain problems can have
no polynomial-time approximation schemes. See \cite{Wegener2005}.
\index{Traveling salesperson problem}
\index{knapsack decision problem}
\index{PCP theorem}
\index{NP decision problem}
\index{inapproximability}

\section{Fenchel duality theorem}

In the theory of {\bf convex analysis}, one can look at convex
bounded continuous functions $f: X \to \mathbb{R}$, $g:X \to \mathbb{R}$
on a Banach space $X$ and at a bounded linear map $A: X \to Y$ from $X$
to an other Banach space to compute $p^* = \inf_{x \in X} f(x) + g(Ax)$ and
$d^* = \sup_{y \in Y^*} -f^(A^* y^*) + g(-y^*)$.
If $X^*,Y^*$ are the dual Banach spaces of $X,Y$ and $A^*: Y^* \to X^*$ is
the adjoint map $(z^*,Ay) = (A^* z^*, y)$ for the pairing of $Y$ with $Y^*$,
then the {\bf strong duality theorem of Fenchel} states:

\satz{$p^* = d^*$}

The theorem is due to Werner Fenchel \cite{BonnesenFenchel}.
It can be generalized, allowing for milder regularity and even unbounded
functions $f,g$ but then only the {\bf weak duality} result
$p^*  \geq  d^*$ holds.
\index{convex analysis}
\index{Fenchel duality}

\section{Legendre transform duality}

The {\bf Legendre transform} of a convex function $f:X \to \mathbb{R}$ defined
on a convex set $X$ in $\mathbb{R}^n$ with inner product $(x,y)$ is
defined as the function $f^*(x^*) = {\rm sup}{x \in X} (x^*,x)-f(x)$ on
$X^* = \{ x^*, {\rm sup}{x} (x^*,x)-f(x)< \infty \}$.
The convex function $f^*$ on $X^*$ is also called the {\bf convex conjugate} of $f$.
One has the following duality result:

\satz{ $f^{**}=f$. }

In the simplest one-dimensional case, convexity means $f''(x)>0$.
The derivative
$( (x^*,x) - f(x))'=0$ means $x^* = f'(x)$ so that $g'(x^*) = f'(x)^{-1}(x^*)$.
For $f(x)=e^x$, one has $f^*(x^*) = x^* \log(x^*)-x$ and for $f(x)=x^2$ one has
$f^*(x^*)=x^2/4$. For the function $f(x)=e^{x-1}=y$ one has $y=f'(x)=e^{x-1}$
and $x= 1+\log(y)$ so that $f^*(x^*) = x^* - x x^* = x^*-(1+\log(x^*)) x^* 
= -x^* \log(x^*)$ which is the function appearing when defining
{\bf entropy}. See \cite{Rockafellar}.
\index{Legendre transform}
\index{convex conjugate}

\section{Gershgorin Circle Theorem}

If $A$ be a complex $n \times n$ matrix, denote by $\lambda_j(A)$ the
{\bf eigenvalues} of $A$. These are the solutions to the polynomial equation 
$p_A(\lambda) = {\rm det}(A - \lambda I) = 0$ of degree $n$. 
By the fundamental theorem of algebra, there are exactly $n$ eigenvalues, 
counted with multiplicity. If $R_i = \sum_{j \neq i} |A_{ij}|$ is the $l^{\infty}$
norm of the $i$'th row vector with the diagonal entry $|A_{ii}|$ missing, 
the disk $G_{ij} = B_{R_i}(A_{ij})$ is called a
{\bf Gershgorin disc}. The {\bf Gershgorin circle theorem} is a result in 
matrix theory.  
\index{Gershgorin disc}
\index{Gerschgorin disc}
\index{Gershgorin circle theorem}

\satz{Every eigenvalue $\lambda_j$ lies in at least one Gershgorin disk}

The theorem can also be seen as a {\bf perturbation result} because if
$A$ is a permutation matrix multiplied with a diagonal matrix, then the 
Gershogorin discs have radius $0$.  The result can be used to estimate how
much the eigenvalues can deviate if such a matrix is perturbed. The result can 
also be used to estimate the determinant ${\rm det}(A) = \prod_j \lambda_j$
of $A$. A special case, attributed by Gershgorin to Bendixson and Hirsch is
that $|\lambda_j| \leq n {\rm max}_{1 \leq i,j \leq n} |A_{ij}|$. 
The result can also be used to estimate the error when computing 
solutions $Ax = b$ of linear equations. This is useful in numerical methods
like when expressing the error of $x$ in terms of the error in $A,B$
using the {\bf condition number} $||A^{-1} || ||A||$ of $A$.
Gershgorin also mentions the corollary that if $|A_{ii}| > \sum_{j \neq i} |A_{ij}|$
for all $i$, then the matrix $A$ is invertible. The result was found by
Semyon Aranovich Gershgorin in
1931. See \cite{GershgorinAndHisCircles}. This book contains also a
copy of Gershgorin's paper from 1931. (In that original article, Gershgorin writes
his name as Gerschgorin.)

\section{The Canada Day Theorem}

For any symmetric $n \times n$ matrix $A$ , the sum
of all $k \times k$ minors ${\rm det}(A_{I \times J})$ with $|I|=|J|=k$  of $A$ is equal to
the sum of the {\bf principal $k \times k$ minors} ${\rm det}((TA)_{I \times I})$ of
the matrix $T A$ , where $T$ is the lower triangular $n \times n$ matrix that is
$T_{kk} = 1$ in the diagonal, and $T_{kl}=2$ for $k>l$ and $T_{kl}=0$ for $k<l$. 
The notation is that if $I,J$ are subsets of $\{1,\dots, n\}$ with cardinality $k$, 
then $P=I \times J$ is the product set which defines the $k \times k$ matrix 
$A_{I \times J}$ in which only the elements in the pattern $P$ appear. 
The {\bf minor} is then defined as the determinant ${\rm det}(A_{I \times J})$ of that
sub-matrix. 

\satz{$\sum_{|I|,|J|=k} {\rm det}(A_{I \times J}) = \sum_{|I|=k} {\rm det}( (TA)_{I \times I})$. }

For example, if $A=\left[ \begin{array}{cc} a & b \\ b & d \\ \end{array} \right]$ and
$T=\left[ \begin{array}{cc} 1 & 0 \\ 2 & 1 \end{array} \right]$, then for $k=2$,
this means ${\rm det}(A)  = {\rm det}(TA)$. For $k=1$, the theorem can be verified by
computing $TA= \left[ \begin{array}{cc} a & b \\ 2a+b & 2b+c \end{array} \right]$ and
checking $a+b+b+c = a + (2b+c)$.
The paper appeared first in \cite{HoneLundmarkSzmigielski2009} and was published in  
\cite{HoneLundmarkSzmigielski2013}.  Since the peak of the discovery
appeared on a July 1 2008 which is {\bf Canada Day}, the name stuck. The proof of the result
uses the {\bf Cauchy-Binet theorem} which reduces it to show              
$$ \sum_{|I|,|J|=k} {\rm det}(A_{I \times J}) 
 = \sum_{|I|,|J|=k} {\rm det}(T_{I \times J}) {\rm det}(A_{I \times J})  \; . $$
Now, ${\rm det}(T_{I \times J})$ is $2^{p(J,I)}$ if $J<I$ and $0$ otherwise, where
$p(J,I) = |J \setminus I \cap J|$. 
\index{Canada Day theorem}
\index{Cauchy-Binet theorem}

\section{Nash embedding theorem}

A Riemannian $m$-dimensional manifold $(M,g)$ is {\bf isometrically embedded}
in $\mathbb{R}^n$ if there is an injective smooth map $\phi: M \to \mathbb{R}^n$
that is {\bf an isometry}. This means that
$g(u,v) = (d\phi(u),d\phi(v))$ for all $u,v \in T_xM$ for the Riemannian metric
$g$ of $M$. Let us say {\bf the $(m,n)$-embedding problem can be solved} for $M$
or an {\bf $(m,n)$-embedding is possible},
if an isometric smooth embedding into $\mathbb{R}^n$ can be achieved for every
compact Riemannian manifold $(M,g)$ of dimension $m$. The {\bf Nash embedding theorem} is

\satz{ An $(m,n)$-embedding is possible for $n \geq 1.5 m^2 + 5.5 m$.}

For non-compact manifolds $(M,d)$, an isometric embedding needs the dimension $n$ of the Euclidean space
to be a bit larger. It is possible if $n \geq 1.5 m^3 + 7m^2 + 5.5 m$.  These constants
appeared in the original 1955 paper of Nash
(reprinted in \cite{EssentialNash} Chapter 11). The embedding cannot not work for $n<0.5 m^2 + 0.5m$
because the right hand side is the number of freedoms of the Riemannian tensor at a point.
Nash's paper includes also some history: Ludwig Schl\"afli in 1871 conjectured 
an embedding in $n \geq 0.5 m^2+0.5m$ but Hilbert in 1901 showed that a constant 
negative curvature manifold can not be embedded in $\mathbb{R}^3$.
Chern and Kuiper in 1952 showed that the flat torus $(\mathbb{T}^n,d)$ can not be embedded
in $\mathbb{R}^{2n-1}$. This is sharp because for even $n$, the {\bf Clifford torus}
is an embedding in $\mathbb{R}^{2n}$ using that
$\mathbb{T}^1$ has an isometric embedding in $\mathbb{R}^{2}$.
For local embeddings, \'Elie Cartan was able to verify in 1927
(following work of M. Janet in 1926) that the Schl\"afli constant
works. A modern proof of the Nash-embedding theorem uses the
{\bf Nash-Moser inverse function theorem} (combining the method of Nash from 1955 and from 
a paper of J. Moser of 1966, who fashioned it into an abstract theorem in  
functional analysis \cite{Hamilton1982,KrantzParksImplicitFunction}).
The Nash embedding theorem is much harder than the
{\bf Whitney embedding theorem} which solves the embedding problem 
without insisting that $\phi$ is an isometry. In that case, $n \geq 2m$ is possible.  
For a more recent simplification of the proof improving also the constant
to $n \geq {\rm max}(0.5 m^2+2.5m,0.5 m^2+1.5m +5)$, see \cite{Guenther1991}. 
The local embedding is first solved based on the {\bf Cauchy-Kowalevski theorem} for
partial differential equations in an analytic setting. The problem is considerably harder 
in the smooth case and this is where already an iterative smoothing process is needed.
\index{Nash embedding theorem}
\index{Whitney embedding problem}
\index{Nash-Moser inverse function theorem}

\section{Erd\"os Straus relation}

The {\bf Diophantine equation} $4/n = 1/x+1/y+1/z$ for unknown positive integers $x,y,z,n$
is called the {\bf Erd\"os Strauss relation}. It is equivalent to $4x y z = n (xy + xz + yz)$.
One only needs to study this in the case when $n$ is {\bf prime} because if $4/p = 1/a+1/b+1/c$ is
solved, then $4/(pq) = 1/(aq) + 1/(bq)+ 1/(cq)$.
As the equation can be solved modulo any prime, by the {\bf Hasse principle} one should
be able to get solutions for any $n$; but this is still unknown. It can appear silly to put the
following as a ``theorem" because it is ``obvious" (or ``trivial" to use a curse word), 
once one sees it, but it illustrates that the difficulty of a Diophantine problem can be hard to 
judge, if one sees it for the first time.

\satz{If $n+1$ is divisible by $3$, then the Erd\"os Straus equation is solvable.}

There is an easy explicit solution formula which one can look up, but
which can be fun to search for, but only if one has not seen it yet. 
The {\bf Erd\"os-Straus conjecture} or {\bf $4/n$ problem} states that for all 
integers $n$ larger than $1$, the rational number $4/n$ can be expressed as the sum of 
three positive unit fractions. Paul Erd\"os and Ernst G. Straus formulated the conjecture in
1948. The problem is still open.  Related is a {\bf conjecture of Sierpinski}, the conjecture
that $5/n=1/x+1/y+1/z$ can be solved. \cite{Guy}.
These problems have appeal because they tap into an 
old theme of {\bf Egyptian fractions} which already appear on the {\bf Rhynd papyrus} 
from around 1650 BC.  On that document, numbers $2/n$ were written as Egyptian fractions 
for all odd numbers $n$ between $5$ and $101$.  
An other interesting problem is to count or estimate the number $f(n)$ of solutions 
of the $4/n$ problem. In \cite{ElsholtzTao}, the sum 
$S(n)=\sum_{p \leq n, p \; {\rm prime}} f(p)$ is bound both from below and above
by $n \log^2(n) \leq S(n) \leq n \log^2(n) \log \log(n)$.

\index{4/n problem}
\index{Erd\"os Straus equation}
\index{Erd\"os Straus conjecture}
\index{Rhynd papyrus}
\index{Egyptian fractions}

\section{Dieudonn\'e Determinant}

If $A$ is a $n \times n$ matrix with entries in a not necessarily commutative
ring like the quaternions $\mathbb{Q}$, one can still look at the {\bf Leibniz determinant}
${\rm det}(A) = \sum_{\sigma} {\rm sign}(\sigma) A_{1 \sigma(1))} \cdots A_{n \sigma(n)}$.
This is a sum over all permutations $\sigma$ of $\{1, \dots, n\}$, where ${\rm sign}(\sigma)$ 
is the {\bf signature} of $\sigma$. It does not satisfy the {\bf Cauchy-Binet} identity
${\rm det}(AB) = {\rm det}(A) {\rm  det}(B)$ in general. There are two ways to get a determinant
which satisfies the later: the first one is called the {\bf Study determinant} \cite{Study1920}. It is
a real-valued determinant defined if $R$ is a {\bf real normed division ring}, meaning $|ab| = |a| |b|$. The 
second is the {\bf Dieudonn\'e determinant} \cite{Dieudonne1943} which takes values in the 
{\bf Abelianization} $R/[R,R]$ 
of the division ring (this is the unique largest subring of $R$ that is Abelian. It is obtained
by factoring out all elements of the commutator form $ab a^{-1} b^{-1}$). The Dieudonn\'e determinant 
has the property that it agrees with the Leibniz determinant in the commutative case 
like $R=\mathbb{R}$ or $R=\mathbb{C}$, the 
Study determinant is a bit easier to compute because we do not bother with commutators 
and allows directly go to the norm. Both determinants rely on the ability to make
{\bf row reduction} which requires that one can divide from the left or from the right. 
They work especially in all normed real division algebras 
$\mathbb{R},\mathbb{C},\mathbb{H},\mathbb{O}$, where in the quaternion and octonion case, 
the Study and the Dieudonn\'e determinant agree. The axiomatic definition of the 
Dieudonn\'e determinant is by asking it take values in the Abelianization $\overline{R}$ and 
demanding for example ${\rm det}(A) {\rm det}(B)= {\rm det}(AB)$ and 
${\rm det}(A)=\prod_{i} \overline{A}_{ii}$ if $A$ is upper triangular.  
\index{Abelianization}
\index{quaternion}
\index{signature}
\index{division ring}

\satz{For a division ring, there is a unique Dieudonn\'e determinant. }

It follows from the axioms that ${\rm det}(A)=1$ and that 
${\rm det}(1+E_{ij})=1$, if $E_{ij}$ is the elementary $0-1$ matrix which is $0$ everywhere
except in the diagonal and the entry $ij$, where it the value is $1$. It also follows
that ${\rm det}(\lambda A) = \overline{\lambda} {\rm det}(A)$ so that row reduction 
allows to compute the determinant depending on whether $\overline{(-1)}=1$ or not.
It also follows that ${\rm det}(A)=0$ if and only if $A$ is singular because that is 
equivalent to having $A$ row reduce to a triangular matrix with a zero in the diagonal.
For quaternions for example $\overline{(-1)} =1$ because $i j i^{-1} j^{-1}=k k=-1$.
Because $SU(2)$ has a trivial Abelianization, one has $\overline{q}=|q|$ for quaternions. 
In order to show the existence of the determinant, one can use row reduction and note
that for $n \geq 2$, any diagonal entry $a b a^{-1} b^{-1}$ can be morphed into $1$
using row reduction steps. One can verify the product property by writing the matrix $A$
as a product of elementary matrices and abelianized ring elements. 
The Dieudonn\'e determinant is treated in 
\cite{ArtinGeometricAlgebra,Brenner1968,RosenbergKTheory}.
\index{Determinant}
\index{row reduction}
\index{Dieudonn\'e determinant}
\index{Noncommutative determinant}
\index{Study determinant}

\section{Centroid theorem}

The {\bf surface area} of a surface $S$ or revolution in $\mathbb{R}^3$ obtained by rotating
a piecewise smooth curve $T$ around the axis of symmetry $L$ is equal to the 
arc length $|T|$ times the length $|C|$ of the circle which is
traced by the geometric centroid of $T$. This is the {\bf Pappus surface centroid theorem} 
and it can be written as $|S| = |T| |C|$. Similarly, the {\bf volume} $|E|$ of a solid 
of revolution $E$ obtained by taking the unions of all projection lines from $S$ to $L$ is equal 
to the area $|A|$ of the flat lamina $A$ between $L$ and $T$, multiplied by the arc length 
$|C|$ of the circle which the centroid of $A$ traces when rotated around $L$.
The {\bf Pappus solid centroid theorem} is then the formula $|E|= |A| |C|$.
This can be generalized: let $C$ be a finite curve connecting two points $P$ and 
$Q$ and let $A$ be a bounded closed region with smooth boundary $T$ contained
in the plane perpendicular to the curve at $A$ such that $P$ is the centroid of $A$. 
The region $A$ can be transported along $C$ using the Fr\'enet frame and 
defines a solid $E$ with boundary $S$. We assume that the tube $S$ remains 
smooth and is a smooth embedding of a $2$-dimensional 
cylinder in $\mathbb{R}^3$. 
\index{Pappus centroid theorem}
\index{Centroid theorem}
\index{surface area}
\index{volume}

\satz{ For surface area $|S| = |T| |C|$, for volume $|E|= |A| |C|$. }

For example, if $T$ is a half circle of radius $r$ in the $xz$-plane
connecting $P=(0,0,-r)$ with $Q=(0,0,r)$ and $L$ is the $z$-axes, then
$|\gamma|=\pi r$ and $|C|= 2 \pi (2r\/\pi)=4r$ so that the surface area of the
sphere $S$ is $4\pi r^2$. A lamina $A$ is a half disc in the $xz$-plane of radius
$r$ which has area $|A|=\pi r^2/2$. The centroid of $A$ has distance $d=4r/(3\pi)$ 
from $L$ moving on a circle of length $|C|=2\pi d = 8r/3$. The volume of the 
sphere of radius $r$ therefore is $|A| |C| = (\pi r^2/2)  (8r/3) = 4\pi r^3/3$. 
The result of Pappus is also used to compute the surface area and volume of
{\bf tubes}. Here is an other example: if $C$ is a smooth closed curve 
in $\mathbb{R}^3$ such that the {\bf tube} 
$\bigcup_{x \in C} B_r(x)$ forms a solid $E$ with piecewise smooth boundary 
surface $S$ that does not have any self intersection, then 
the surface area
is $|S|= |T| |C| = 2 \pi r |C| + 4\pi r^2$ and the volume is
$|E| = |A| |C| + 4\pi r^3$ (the additional terms come from the sphere ``roundings at the
end points"). In this case, the lamina $A$ are disks of radius $\pi r^2$ 
and the curves $T$ are circles of arc-length $2\pi r$. 
Even more general versions have been discussed in detail in \cite{GoodmanGoodman1969}.
For tube methods in differential geometry also in higher dimensions (which are certainly 
also inspired by the Pappus centroid theorem), see \cite{Gray}.
Herman Weyl used {\bf tubes} as a powerful tool in differential geometry \cite{Weyl1939}.
\index{Tubes}

\section{The Borsuk antipodal theorem}

Let $M=\mathbb{S}^n$ denote the $n$-sphere $\{ |x|^2 = 1 \} \subset \mathbb{R}^{n+1}$
equipped with metric induced from open sets in the Euclidean space $\mathbb{R}^{n+1}$.
Let $A_0, \dots, A_n$ be {\bf cover} of $M$ by {\bf closed sets}. This means that
$\bigcup_{k=0}^n A_k = S$. We say, a subset $A$ of $M$
contains an {\bf antipodal pair}, if there is a pair of points
$\{ x,-x \} \in \mathbb{M}$ which are both in $A$.

\satz{A cover of the $n$-sphere by $n+1$ sets contains an antipodal pair.}

The theorem is also known as the {\bf Lusternik-Schnirelman-Borsuk antipodal theorem} 
(already called so by \cite{Hopf1944}),
much of the literature just calls it the Borsuk theorem, maybe because of simplicity.
The theorem is equivalent to the {\bf Borsuk-Ulam theorem} stating 
that every map $f$ from $M$ to $\mathbb{R}^n$ 
has the property that some antipod pair $x,y$ has the property that $f(x)=f(y)$. 
Stan Ulam was credited in the Borsuk paper as the originator of the problem.
In \cite{Bollobas2006} section 41 contains elegant proofs that
the statements in Borsuk's theorem and in the Borsuk-Ulam are equivalent.
The result generalizes to the situation when $M$ with a manifold homeomorphic 
to $\mathbb{S}^n$ equipped with an involution $T: M \to M$ which is conjugated to the 
antipodality on $\mathbb{S}^n$.
See \cite{Bollobas2006} Section 150.
For $n=1$, the theorem is equivalent to the {\bf intermediate value theorem}: $M$ is a circle
and the function $f(x)-f(x')$, if not constant $0$, takes both positive and negative values
so that there must be a point where $f(x)=f(x')$ with antipodal points $x'$. 
For $n=2$, if we cover the $2$-sphere with $3$ open sets, there is one of the sets
which contains an antipode. The more surprising equivalent Borsuk-Ulam statement is then
that there are two anti-podes on earth, where both the temperature and the pressure 
are the same. The theorem appeared first in 1930 in a paper by Lusternik and Schnirelman and
then more generally in 1933 by Karol Borsuk \cite{Borsuk1933}. 
The fact that there is a general theorem on Lusternik-Schnirelman 
category by Lusternik and Schnirelman is a reason to stick to Borsuk for the antipodal theorem. 
Heinz Hopf generalized in 1944 the theorem as follows: if $A_0, \cdots, A_n$ are
$n$ closed sets covering the unit sphere $\mathbb{S}^n$ in $\mathbb{R}^{n+1}$
and $0<d \leq 2$ is a distance, then there exists a set $A_k$ in which 
there exists two points of distance $d$. The special case $d=2$ is the Borsuk theorem.  
Hopf notes that this implies that if the $n$-sphere is covered by $n+2$ non-empty closed
sets such that none of them contains a antipodal pair, then every collection of $n+1$
sets has a non-empty intersection and states in a footnote that this means that the
{\bf nerve} of the cover $F_0, \dots, F_{n+1}$ is then isomorphic to the boundary complex
of a $(n+1)$-dimensional simplex. 
\index{intermediate value theorem}
\index{Borsuk-Ulam theorem}
\index{Borsuk theorem}
\index{Lusternik Schnirelman Borsuk antipodal theorem}
\index{Antipode}
\index{Antipodal theorem}
\index{nerve}

\section{Zagier's inequality}

Assume $f,g$ are non-negative and decreasing functions on $[0,T]$. They are then
automatically integrable. Denote by 
${\rm E}[f] = \frac{1}{T} \int_0^T f(x) \; dx$ the {\bf average} of $f$.

\satz{If $f,g$ are non-negative and decreasing, then
      ${\rm E}[f g] \geq {\rm E}[f] {\rm E}[g]$.}

\cite{Zagier1995} formulates this more generally as follows: 
if $f,g$ are decreasing and non-negative on $[0,\infty)$ and 
$F,G \in L^1([0,\infty))$ take values in $[0,1]$, then 
 $(f,g) \geq (f,F) (g,G)/{\rm max}(I(f),I(g))$, where $I(F)=\int_0^{\infty} F(x) dx =|F|_1$. 

The Zagier inequality has also been called a {\bf anti-Cauchy-Schwarz} inequality 
\cite{Bollobas2006} because in {\bf Cauchy-Schwarz} $|f \cdot g| \leq |f| |g|$ in a {\bf Hilbert space},
the inequality works in the opposite direction. In \cite{Bollobas2006}, the 
inequality on finite intervals is called {\bf Chebychev's inequality} but the 
later should maybe be reserved for the inequality
${\rm P}[|X-{\rm E}[X]|>\epsilon] \leq {\rm Var}[X]/\epsilon^2$ for a         
random variable variable $X \in L^2(\Omega,\mathcal{A},P)$ 
on a probability space $(\Omega,\mathcal{A},{\rm P})$. 
The Zagier inequality also works for {\bf decreasing sequences} $f_n$,
where ${\rm E}[f] = \frac{1}{n} \sum_{k=0}^{n-1} f_k$ is the Birkhoff average.
Now, the same statement ${\rm E}[f g ] \geq {\rm E}[f] {\rm E}[g]$ holds.
In the simplest case, for $f=(a,b)$ and $g=(c,d)$, this is equivalent to
$2 (ac+bd)  \geq (a+b)(c+d)$ for $a \geq b, c \geq d$ which is already
not totally obvious as it is equivalent to $ac+bd \geq ad+bc$.
\index{Zagier inequality}
\index{anti Cauchy-Schwarz inequality}
\index{Chebychev inequality}
\index{Cauchy-Schwarz}

\section{Gini coefficient}

If $x_1, \cdots, x_n$ are non-negative real numbers with {\bf mean}
$m = \frac{1}{n} \sum_{k=1}^n x_k$, the
number $G = \frac{1}{2n^2 m} \sum_{i=1}^n \sum_{j=1}^n |x_i-x_j|$ is called
the {\bf Gini coefficient} of the data. Using $|a-b|=a+b-2 {\rm min}(a,b)$ it can be rewritten
as $G=1-\frac{1}{n^2 m} \sum_{i=1}^n \sum_{i=1}^n {\rm min}(x_i,x_j)$.
A common interpretation is that $x_k$ is the {\bf income} of person $k$ in a
population $X=\{1, \dots, n\}$ of $n$ people. The number $m$ is then the
{\bf mean income} of the population. If a population $X$ of $n$ people
is split into smaller groups $X_k, k=1, \dots, r$ of size $n_k$
and have mean income $m_k$, then $\sum_{k=1}^r n_k = n, \sum_{k=1}^r n_k m_k = n m$.
If $G(X)$ is the Gini coefficient of $X$ and $G(X_k)$ the Gini coefficient of the
sub-population $X_k$, then

\satz{ $n G(X) \geq \sum_{k=1}^r n_k G(X_k)$ }

This and many more inequalities relating $G(X)$ with $G(X_k)$ appear in
\cite{Zagier1983}. There is also a continuum analog:
for a probability density function $f$ on $[0,\infty)$ with
$\int_0^{\infty} f(x) \; dx=1, m=\int_0^{\infty} x f(x) \; dx$, the
{\bf continuum Gini coefficient} is defined as
$G=\frac{1}{2m} \int_0^{\infty} \int_0^{\infty} |x-y| f(x) f(y) dx dy$
which is equivalent to
$G=1-\frac{1}{m} \int_0^{\infty} \int_0^{\infty} {\rm min}(x,y) f(x) f(y) dx dy$.
The Gini coefficient is twice the area between the {\bf Lorenz curve}
and the diagonal $G=2 \int_0^1 p - L(p) \; dp$, where $p = \int_0^x f(t) \; dt$ is the cumulative
distribution value and $L(p) = \frac{1}{m} \int_0^x t f(t) \; dt$.
In the context of {\bf income inequality}, where the subject has come up in
economics, $L(p)$ represents the fraction of the total income which is
earned by the poorest $n p$ people.
The graph of $L(p)$ is a convex curve from $(0,0)$ to $(1,1)$, the slope
$L'(p)$ being the {\bf relative income} in the corresponding percentile
of the population. The Gini coefficient is also called {\bf Gini index}.
It has been introduced by Corrado Gini in 1912. It is a natural
quantity because on the real line the {\bf Green's function} of the Laplacian
$-\Delta/2$ with $\Delta f=f''$ one has
$g(x,y) = |x-y|$. The potential $V(x)=|x|$ is the
natural ``Newton potential". For $M=\mathbb{R}^d$ in dimension $d \neq 2$ it is
$g(x,y) = |x|^{2-d}$ for the Laplacian $-\Delta/|S_{d-1}|$, where $|S_k|$ is
the volume of the $k$-dimensional unit sphere;
it is the logarithmic potential $\log|z|/(2\pi)$ in dimension $d=2$.
The most familiar case is the 3-dimensional Euclidean space $\mathbb{R}^3$, where the Newton
potential $1/|x|$ appears in {\bf electro magnetism} and {\bf gravity}. 
The {\bf Gini potential} $|x|$ is roughly the force between two planar parallel mass sheets
like two galaxies rotating around the same axis. 
In general, for any Riemannian manifold $M$ with {\bf Greens function} $g(x,y)$
(the inverse of the Laplacian) and measure $\mu$ (mass distribution)
the integral $I(\mu) = \int_{M} \int_{M} g(x,y) d\mu(x) d\mu(y)$ is the
{\bf potential theoretical energy} of the measure $\mu$. The Gini index therefore is
proportional to the potential theoretical energy for a mass
distribution with density $\mu = f(x) dx$ on $[0,\infty)$. The above inequality
could therefore be interpreted as an inequality for the potential energy of particles
which are partitioned into non-interacting groups. Switching off energies between
non-interacting parts lowers the energy. 
\index{Gini coefficient}
\index{Gini potential}
\index{Lorenz curve}
\index{Income curve}
\index{Laplacian}
\index{Newton potential}

\section{Denjoy-Koksma theorem}

If $T: X \to X$ be an ergodic automorphism of a probability space $(\Omega,\mathcal{A},\mu)$.
(Automorphism means $\mu(T(A))=\mu(A)$ for all $A \in \mathcal{A}$ and ergodic means that 
$T(A)=A$ implies $\mu(A) \in \{0,1\}$.)
The {\bf Birkhoff ergodic theorem} assures that for all $g \in L^1(\Omega)$ and almost every $x \in \Omega$ we have
$S_n(x)/n \to {\rm E}[g] = \int_{\Omega} g(x) \; d\mu(x)$ with the Birkhoff sum $S_n = \sum_{k=0}^{n-1} g(T^kx)$.
An example dynamical system is the {\bf irrational rotation} $T:x \to x+\alpha$ on the circle
$\mathbb{T}^1=\mathbb{R}^1/\mathbb{Z}^1$ equipped with the Lebesgue measure $\mu=dx$.
{\bf Denjoy-Koksma theory} estimates the growth of $S_n(x)$ depending on {\bf Diophantine properties} 
of $\alpha$ and {\bf regularity properties} of $g$.
A real number $\alpha$ is called {\bf Diophantine}, if there exists a constant $C$ such that
$|p \alpha-q| \leq C q$, for all integers $p,q$.
A function $g$ has {\bf bounded variation} if ${\rm Var}(g) = \sup_{P} \sum |g(x_{i+1})-g(x_i)|$ is finite,
where the supremum is over all finite sets $P = \{x_1,\dots ,x_n=x_0 \}$ in $\mathbb{T}^1$. In the simplest case,
the {\bf  Denjoy theorem} says $S_n \leq C \log(n) {\rm Var}(g)$ for all $n$ and that there is a 
sequence of integers $q_n$, for which $S_{q_n}(x) \leq {\rm Var}(f)$, the periodic approximations $p_n/q_n \to \alpha$.
For $r \geq 1$, a real number $\alpha$ is called {\bf $r$-Diophantine}, 
if $|q \alpha-p| \leq C q^r$ for all integers $p,q$. The Denjoy-Koksma theorem 
was generalized in 1999 by Svetlana Jitomirskaja to

\satz{If $\alpha$ is $r$-Diophantine, then $|S_n| \leq  C n^{1-1/r} \log(n) {\rm Var}(g)$.}

For a {\bf periodic approximation} $p/q$ of $\alpha$ \cite{Chinchin92}                
one has $|S_q| \leq {\rm Var}(f)$: to see this divide $\mathbb{T}^1$ into $q$ intervals centered at
$y_k=k p/q$. The intervals have length $1/q \pm O(1/q^2)$ and each
contains exactly one point. Renumber the points to have $y_k$ in $I_k$.
By the {\bf intermediate value theorem}, there exists a Riemann sum
$\frac{1}{q} \sum_{i=0}^{q-1} f(x_i) = \int f(x) \; dx =0$ for which every
$x_i$ is in an interval $I_i$. Choosing $x_i=\min_{x \in I_i} f(x)$ gives
an lower and $x_i = \max_{x \in I_k} f(x)$ gives an upper bound.
Now, $\sum_{j=0}^{q-1} f(y_j) - f(x_j)  \leq \sum |f(y_j) - f(x_j)| + |f(x_j) - f(y_{j+1})| \leq {\rm Var}(f)$.
Therefore, if $q_k \leq n \leq q_{k+1}$ and
$n = b_k q_k + b_{k-1} q_{k-1} + \cdots + b_1 q_1 + b_0$, then
$S_n \leq  \sum_{i=0}^{n} (b_0+ \cdots + b_n) {\rm Var}(f)$.
where $b_k \leq q_{i+1}/q_i$. So,
$S_n \leq  \sum_{i=0}^{n} \frac{q_{i+1}}{q_i} {\rm Var}(f)$.
If $\alpha$ is $r$-Diophantine, then $|q \alpha| \leq c/q^r$
and $q_{i+1} \leq q_i^r/c$. Because $n \leq q_{k+1} \leq q_k^r/c$, we have
$q_k \geq (c n)^{1/r}$ and $n/q_k \leq c^{-1/r} k^{1-1/r}$. 
Because $k \leq 2 \log(q_k)/\log(2)$, the claim follows.
For $r=1$, see \cite{CFS} (page 84). In general, see \cite{Jitomirskaya1999}.
\index{Denjoy theorem}
\index{Denjoy Koksma theorem}
\index{Diophantine}
\index{r-Diophantine}
\index{Bounded variation}
\index{Birkhoff sum}
\index{irrational rotation}
\index{ergodic}

\section{Quadrilateral theorem}

Let $ABCD$ denote a convex {\bf quadrilateral} in $\mathbb{R}^2$. Alternatively,
the four arbitrary points $A,B,C,D$ in $\mathbb{R}^3$ define a tetrahedron. Assume the 
side lengths are $a=|AB|,b=|BC|,c=|CD|,d=|DA|$ and that the diagonal lengths are $e=|AC|,f=|BD|$.
Let $M=(A+C)/2$ and $N=(B+D)/2$ be the midpoints of the 
diagonals and $g=2|MN|$. The {\bf Euler law on quadrilaterals} is

\satz{$a^2+b^2+c^2+d^2=e^2+f^2+g^2$.}

One can verify this by just expanding out what one gets when 
writing the condition in coordinates. The proof then shows that the
statement gives also a statement about lengths of a {\bf tetrahedron}
in space $\mathbb{R}^3$: if $a,b,c,d,e,f$ are the side lengths of an arbitrary tetrahedron
in space and the edges $L,M$ belonging $e,f$ have no common point, and $g$ is twice
the length between the midpoints of the two segments $L$ and $M$, then the same 
relation holds in space. This has been noted in \cite{Kandall2002}.
In the case of a rectangle, where $a=c,b=d,g=0,e^2=f^2=a^2+b^2, g^2=0$
one has the {\bf Pythagorean theorem}. In the case of a parallelogram, where
$a=c,b=d,g=0$ one has $2a^2+2b^2=e^2+f^2$, it is the
{\bf parallelogram law}. 
Some other themes of Euler come to mind too like Diophantine equations:
if the points $A,B,C,D$ have integer coordinates and all
distances between points are integers, one has a problem in number theory. 
For rectangles, this leads to {\bf Pythagorean triples}.
The problem of {\bf perfect Euler bricks} comes then to mind, which asks for a
cuboid with integer side and diagonal lengths.
\index{Pythagorean triples}
\index{Parallelogram law}
\index{Euler law on quadrilaterals}
\index{quadrilateral theorem}
\index{Euler-Pythagoras theorem}
\index{Perfect Euler brick}

\section{Reeb sphere theorem}

Let $M$ be a closed, compact $d$-dimensional differentiable manifold. 
{\bf Closed} means that the boundary of $M$ is empty. If $f: M \to \mathbb{R}$
is a smooth real-valued function, then points $x \in M$ with vanishing {\bf gradient}
$\nabla f(x)=0$ are called 
{\bf critical points}. A critical point $x$ is called {\bf non-degenerate}, 
if the {\bf Hessian $d \times d$ matrix} $H(f)(x)$ is invertible at $x$.
Let $c(M)$ denote the minimal number of non-degenerate 
critical points which a function $f$ on $M$ can have. We say $M$ {\bf is a $d$-sphere}, if
there is a homeomorphism of $M$ to the standard unit sphere $\{ |x|=1 \}$ 
in $\mathbb{R}^{d+1}$.

\satz{$c(M)=2$ if and only if $M$ is a $d$-sphere for some $d \geq 0$.}

The level curves $f^{-1}(c) = \{ f = c \}$ of $f$ form then a foliation of $M$
which are $(d-1)$-dimensional spheres which only degenerate to points at
the critical points. The proof of the theorem goes by showing that a manifold
which admits exactly $2$ critical points can be covered by $2$ balls,
then use that this characterizes spheres.
The Reeb sphere theorem was proven in 1952 \cite{Reeb1952}. It is 
referred to and generalized in \cite{McAuley} who generalizes and improves on 
results by Milnor and Rosen. 
The assumption of $f$ has two critical points does not imply that $M$ is
diffeomorphic to the standard unit sphere. There are {\bf exotic spheres}
which are homeomorphic to the standard unit sphere but not diffeomorphic to it.
The Reeb theorem is covered in \cite{Mil65}. In the first proof of the existence of
exotic 7-spheres, \cite{Milnor1956}, the Reeb Sphere theorem was used as
hypothesis H.
\index{Exotic spheres}
\index{Milnor sphere}
\index{Reeb sphere theorem}
\index{critical points}
\index{non-degenerate critical points}

\section{Hausdorff distance}

Let $(X,d)$ be a metric space. Given a compact subset $U$ of $X$,
let $B_r(U)$ the set of all points that are in distance $\leq r$ from a point
of $U$. In other words $B_r(U)=\bigcup_{x \in U} B_r(x)$, where $B_r(x)$ is
the ball $\{ y \in X, d(x,y) \leq r\}$ in $X$. The {\bf Hausdorff
distance} $\delta$ between two non-empty compact subsets $U,V$ of $X$
is defined as the infimum over all $r \geq 0$ such that
$U \subset B_r(V)$ and $V \subset B_r(U)$. It is a metric on the set of
all compact subsets. 
This space $(\chi,\delta)$ is a new metric space. It is again compact:

\satz{If $(X,d)$ is compact, then $(\chi,\delta)$ is again compact.}

The process could therefore be iterated and produce a sequence of compact metric spaces, where
in each step the Hausdorff metric is used on the previous one. 
For Hausdorff distance, see \cite{Falconer}, Chapter 9, in the context of {\bf iterated function systems}
in the {\bf theory of fractals}. A sequence of contractions defines an {\bf attractor} which can be 
seen as a limit of a sequence of compact sets. In the simplest situations, one can then use the 
{\bf Banach fixed point theorem} to establish the existence of a limit.  The distance has been
used by Maurice Fr\'echet in 1906 to measure the distance between curves. The distance
was introduced by Felix Hausdorff in 1914 \cite{Hausdorff1914} (page 303). 

The Hausdorff distance allows also to define a distance between compact metric spaces
$(X_1,d_1)$, $(X_2,d_2)$. {\bf The Gromov-Hausdorff distance} of two compact metric
spaces is defined as the infimum over all possible Hausdorff distances
$\delta (\phi_1(X_1),\phi_2(X_2))$, where $\phi_i:X_i \to X$ are isometric embeddings
of $(X_i,d_i)$ into a third metric space $(X,d)$. This metric space $(\mathcal{X},\mathcal{D})$
of all compact metric spaces has a dense set of finite metric spaces so that it is separable.
It is also complete, from which one can deduce that it is connected.
David Edwards \cite{Edwards1975} called this ``superspace".
\index{Super space}
\index{Hausdoff distance}
\index{Gromov Hausdorff distance}
\index{compactness}
\index{iterated function system}
\index{fractals}

\section{Grove-Searle theorem}

The set of compact even-dimensional Riemannian $2d$-manifolds which admit a positive
curvature metric contains spheres $\mathbb{S}^{2d}$,   
projective spaces
$\mathbb{RP}^{2d},\mathbb{CP}^d,\mathbb{HP}^d,\mathbb{OP}^2$ over the
division algebras $\mathbb{R},\mathbb{C},\mathbb{H},\mathbb{O}$,
the three {\bf Wallach flag manifolds} $W^6,W^{12},W^{24}$ \cite{Wallach1972} and the
{\bf Eschenburg manifold} $E^6$  \cite{Eschenburg1982}.
No other example is known \cite{Ziller2}.  
All these manifolds admit a positive metric with a continuum 
isometry group. In particular they admit a metric which allows for an isometric circle action.
The fixed point set $N=\phi(M)$ of such an action is never empty \cite{Berger2002}.
By a theory started by Conner and Kobayashi it is again a positive curvature manifold $N$
that is totally geodesic and of even co-dimension. The components of $N$
can have different dimension but by Lefschetz, the Euler characteristic
of $N$ is the Euler characteristic of $M$ \cite{Conner1957,Kobayashi1958}. 
Lets call a manifold with circular symmetry {\bf Grove-Searle} if the fixed
point set $N$ has a connected component of co-dimension $2$.
The {\bf Grove-Searle theorem} \cite{GroveSearle} now tells:
\index{Positive curvature manifolds}
\index{Wallach manifolds}
\index{Eschenburg manifold}
\index{Projective spaces}

\satz{ 
If $M$ is Grove-Searle, then $M=\mathbb{S}^{2d}, \mathbb{RP}^{2d}$ or $\mathbb{CP}^d$.
}

In odd dimensions, there is beside $M=\mathbb{S}^{2d+1}$ or $M =\mathbb{RP}^{2d+1}$
also the possibility of space forms $\mathbb{S}^{2d+1}/\mathbb{Z}_m$.   
An application of the theorem is that all $2d$-dimensional positive curvature
manifolds admitting a circular symmetry have positive Euler characteristic if $2d \leq 8$.
Proof: $N$ is not empty by Berger and $\chi(N)=\chi(M)$. 
$N$ has a co-dimension $2$ component, Grove-Searle forces $M$ to be in
$\{ \mathbb{RP}^{2d},\mathbb{S}^{2d},\mathbb{CP}^d \}$. By Frankel \cite{Frankel1961}, 
there can be not two co-dimension $2$ connected components.
In the remaining cases, Gauss-Bonnet-Chern \cite{Chern1966} 
forces all to have positive Euler characteristic.
There is huge interest in even-dimensional positive curvature manifolds because of the
open {\bf Hopf conjecture}  \cite{Hopf1932,Hopf1946,Hopf1953,BishopGoldberg,BergerPanorama}
asking whether every even-dimensional compact positive
curvature manifold has positive Euler characteristic. The above corollary of Grove-Searle
assures that the Hopf conjecture with circle symmetry holds in dimension $\leq 8$.    
It is also known for $2d=10$: \cite{PuettmannSearle,RongSu2005,Wilking2003}.
See also $2d=6$ in \cite{Petersen3} (2. Edition, Cor. 8.3.3).
While in dimension $2$ and $4$ the classification of positive metric manifolds with  
circular symmetry is known like $\{ \mathbb{S}^4, \mathbb{RP}^4, \mathbb{CP}^2 \}$ in 
dimension $4$ \cite{HsiangKleiner}, in dimension $6$ one knows so far the cases      
$\{ \mathbb{S}^6, \mathbb{RP}^6, \mathbb{CP}^3, E^6,W^6 \}$ and it is not known whether
they are all. 
\index{Grove-Searle theorem}
\index{Hopf conjecture}
\index{Gauss-Bonnet-Chern}

\section{Radon-Nikodym theorem}

A {\bf measurable space} $(\Omega,\mathcal{A})$ is a set equipped with
a {\bf $\sigma$-algebra} $\mathcal{A}$. This means that $\mathcal{A}$ is a
set of subsets of $X$ containing $X$, that is closed under forming complements
and the operation of taking countable unions. A non-negative valued function $f:\Omega \to [0,\infty)$
is called {\bf measurable} if $f^{-1}(B) \in \mathcal{A}$ for every $B$ in the
{\bf Borel $\sigma$-algebra} on $[0,\infty)$, the smallest $\sigma$-algebra 
containing the open sets. Given two {\bf $\sigma$-finite measures} $\mu,\nu$,
(meaning that $\Omega$ is in each case a countable union of sets of finite measure),
on $(\Omega,\mathcal{A})$ one calls $\mu$ {\bf absolutely continuous} with respect to $\nu$, if
$\nu(A)=0$ implies $\mu(A)=0$. An example is if there exists a function $f \in L^1(\omega,\mathcal{A},\nu)$
such that $\mu(A) = \int_A f(x) \; d\nu(x)$, then $\mu$ is absolutely continuous with respect to $\nu$
and the function $f$ is called the {\bf Radon-Nikodym derivative} of $\mu$ with respect to $\nu$,
as $d\mu(x)=f(x) d\nu(x)$ suggests to write $d\mu/d\nu = f$. The {\bf Theorem of Radon-Nikodym} assures
that this situation is the general case. 
Let us abbreviate $\mu << \nu$ if $\mu$ is absolutely continuous with respect to $\nu$.

\satz{ If $\mu<<\nu$, there exists $f \in L^1(\Omega,\mathcal{A},\nu)$ with $\mu= f \nu$. }

The theorem is important in {\bf probability theory}, where the measures under consideration
are usually {\bf probability measures}, meaning $\mu(\Omega)=1$. If $\mu$ is absolutely
continuous to $\nu$ then every set of zero probability with respect to $\nu$
has zero probability with respect to $\mu$.
An example of a measure $\mu$ on the Lebesgue space $([0,1], \mathcal{A},\nu=dx$ 
is a Dirac point measure $\delta_x$ for a point in $[0,1]$.
An application of the Radon-Nikodym theorem is the {\bf Lebesgue decomposition} of a measure.
One can split every $\sigma$-finite measure into an absolutely continuous, a singular continuous
and a pure point part. This is important in spectral theory of mathematical physics \cite{ReedSimon}.
For measure theory and real analysis in general, see for example \cite{SteinShakarchi2005}.
For the history, \cite{Simon2017} (page 257): the theorem was first proven by Radon in 1913 in
$\mathbb{R}^n$ and then by Nikodym in 1930.
\index{Radon-Nikodym theorem} 
\index{Lebesgue decomposition}
\index{$\sigma$-algebra}
\index{measurable}
\index{absolutely continuous}
\index{Borel $\sigma$ algebra}
\index{measurable space}
\index{Radon-Nikodym derivative}
\index{singular continuous}

\section{Crofton formula}

If a needle of length $l<1$ is thrown randomly into a periodic grid of lines
spaced distance $1$ apart, the probability of hitting a grid line is $2l/\pi$.    
This method of computing $\pi$ is an example of a {\bf Monte-Carlo method}.
A probability space of needle configurations can be given as 
$(\Omega,\mathcal{A},\mu) = ([0,1/2] \times [-\pi/2,\pi/2],\mathcal{A},2 d\theta dr/\pi)$
with product Lebesgue measure, where $r$ is the minimal distance of the center
of the needle to a grid line and $\theta$ is the polar angle. The needle 
obviously hits if and only if $r \leq (l/2) \cos(\theta)$. The probability therefore
is obtained by integrating the density $2/\pi$ over this region. It gives
$\int_{-\pi/2}^{\pi/2} \int_0^{(l/2) \cos(\theta)} 2/\pi dr d\theta = 2l/\pi$.
This can now be generalized for any rectifiable curve of length $l$. One has only to look
at the {\bf random variable} $X$, which counts the number $X$ of intersections of 
the randomly placed curve with a grid.
The {\bf Crofton formula} in the plane is now ${\rm E}[X]=2l/\pi$:
(to see this, approximate the curve by a polygon and look at each segment $l_i$ as a ``needle"
of length $l/n$. Then $X=X_1 + \cdots + X_n$ where $X_j$ counts the number of intersections
with $L_j$. By linearity of expectation and additivity of length, the Crofton formula follows.)
One can look at the problem also in $\mathbb{R}^n$, where one has a system of parallel hyperplanes spaced
a unit apart and a rectifiable curve of length $l$. 
Now, the volume $|B^{n-1}|$ of the $(n-1)$-dimensional unit ball 
and the volume $|S^{n-1}|$ of the $(n-1)$-dimensional sphere matters. 
Again, $X$ is the number of intersections of the curve with the 
periodic plane grid.

\satz{${\rm E}[X] = 2l |B^{n-1}|/|S^{n-1}|$.}

In the case $n=2$, this was $|B^1|=2, |S^1|=2\pi$ and the original Buffon formula follows.
The {\bf Buffon needle problem} is the fist connection between probability theory 
and geometry. It appeared first in 1733 and was reproduced again in 1777 by Buffon.
Morgan Crofton extended this in 1868 \cite{Crofton1968}.
The mathematical field of integral geometry started to blossom with Blaschke 
\cite{blaschke} in the late 1930ies.
Probability spaces can be used to study more geometrical quantities like surface area, or
curvature \cite{Banchoff67,Banchoff70,MilnorKnot}. General references 
are \cite{Santalo,Santalo1,KlainRota,Schneider1}. The $n$-dimensional Crowfton formula can
be found in \cite{KlainRota}.
\index{Integral geometry}
\index{Crofton formula}
\index{Buffon needle problem}
\index{Geometric probability}

\section{Desnanot-Jacobi identity}

If $A$ is a $n \times n$ matrix, the matrix entries are accessed as $A_{ij}$.
Call $A_i^j$ the matrix obtained by deleting row $i$ and column $j$ in $A$.
The expression $(-1)^{i+j} {\rm det}(A_i^j)$ is also known as a {\bf cofactor}
of the {\bf minor} ${\rm det}(A_i^j)$.
Similarly, let $A_{ij}^{kl}$ be the matrix in which rows $i,j$ and
columns $k,l$ are deleted.  The {\bf Desnanot-Jacobi identity} is the following
relation between sub-determinants of a matrix:

\satz{${\rm det}(A) {\rm det}(A_{1n}^{1n}) = {\rm det}(A_1^1) {\rm det}(A_n^n) - {\rm det}(A_1^n) {\rm det}(A_n^1)$.}

It allows to write ${\rm det}(A)$ in terms of the $(n-2) \times (n-2)$ matrix in which 
the boundary rim is removed and all the four possible $(n-1) \times (n-1)$ matrices, 
where one boundary row and boundary column is removed from the matrix.
In the case when $n=2$ the identity still works if one interprets 
${\rm det}(A_{1n}^{1n})= {\rm det}(A_{12}^{12})$
as $1$, which is usually assumed the value for the determinant of the empty matrix. 
In that case, the Desnanot-Jacobi identity is just
${\rm det}( \left[ \begin{array}{cc} a & b \\ c & d \end{array} \right]) = ad-bc$.
The identity is also called the {\bf Desnanot-Jacobi adjoint matrix theorem}.
A generalization is called the {\bf Sylvester determinant identity}.
The Desnanot-Jacobi identity leads to a process called {\bf Dodgons condensation} or
{\bf Alice in Wonderland condensation} because Charles Lutwidge Dodgson is
also known as Lewis Carroll, the author of ``Alice in Wonderland" \cite{AlicesAdventuresInWonderland}.
The condensation method was described in \cite{Dodgson1866} in 1866. The Desnanot result
appears first in 1819, in the book \cite{Desnanot} on page 152 and in \cite{Jacobi1827} in 1827.
Like for Cauchy-Binet, it is historically remarkable that this identity was found 
before matrices were formed.  Indeed, the word ``matrix", related to the latin word 
``mater" for mother was later used in a more generalized sense as ``womb".
The word ``matrix" therefore appeared because
 matrices are devices which bear determinants.    
There are more references in \cite{Knuth1995}.

\index{Alice in Wonderland}
\index{Desnanot-Jacobi}
\index{Dodgons condensation}
\index{Desnanot-Jacobi adjoint matrix theorem}
\index{matrix}
\index{determinant}

\section{Existence of Minimal surfaces}

A $2$-dimensional {\bf surface} $S$ in $\mathbb{R}^n$ is the image of a
parametrization $r(u,v): R \to \mathbb{R}^n$,
where $R \subset \mathbb{R}^2$ is the parameter domain, an open,
simply connected region in the plane which one can assume
to be the unit disc $R$ with circle $C$ as boundary.
The surface is called {\bf minimal surface} if very component of $r$ is
{\bf harmonic} $\Delta r = 0$ and furthermore $E-F=|r_u|^2-|r_v|^2=0$ and
$F=r_u \cdot r_v=0$, expressing that the Riemannian metric
$g$ on $S$ is conformal. The {\bf Plateau problem} is to find for 
a given simple closed curve $\Gamma$ in $\mathbb{R}^n$, a minimal 
surface $S$ which has $\Gamma$ as the boundary. One wants the map $r$
to be smooth in $R$ and continuous up to the boundary $C$. The surface
$S$ does not necessarily have to be embedded ($r$ is not necessarily injective), 
it can just be {\bf immersed}.

\satz{There is a solution to the Plateau problem.}

Note that this does not mean that the solution is unique. 
Indeed, in general there are multiple solutions even-so generically 
only finitely many. In general, solutions also can have branch points, 
self-intersections or can be physically unstable and so would be
difficult to observe in soap bubble experiments. 
The problem was solved first in 1931 by Jesse Douglas and Tibor Rado in 1930.
If more generally, the region $R$ has larger genus and so several boundary 
curves, the problem is called the {\bf Douglas problem}. 
When looking at how soap films change in dependence of parameters, 
huge changes like catastrophes can happen. For example,
in that if $\Gamma$ is changed, suddenly, solutions to
a genus one Douglas problem appear as it has
lower energy.\cite{FomenkoPlateau}
In order to solve the Plateau problem
one is led to the variational problem of extremizing
the {\bf Dirichlet integral} $\mathcal{L}(r)=\iint_R |r_u|^2+|r_v|^2 dudv$. 
The harmonicity condition $\Delta r = 0$ is the {\bf Euler equation}
of the variational problem. This is a special case of a
{\bf Dirichlet principle}. The problem was raised by Joseph-Louis Lagrange
in 1760 and named after the physics and anatomy professor Joseph Plateau who made experiments.
Poisson realized that soap films are surfaces of constant mean curvature.
In higher dimensions, the problem has led to {\bf geometric measure theory}.
We followed partly \cite{Courant1940}. More information is in 
\cite{Struwe1989}, where also the history of soap films and soap bubbles is
described as one of the oldest objects in mathematical analysis and pointed out
that for a long time, since Lagrange's derivation of the minimal surface equation,
the analysis was too difficult even for mathematicians like
Riemann, Weierstrass or Schwarz. In \cite{FomenkoPlateau} (part I)
there is more history and many pictures and relations where minimal films in nature
as the most economical surfaces forming skeletons of
{\bf radiolarians}, tiny marine organisms.
\index{Plateau problem}
\index{Douglas problem}
\index{Minimal surface}
\index{Soap bubbles}

\section{Fermat's right angle theorem}

A positive integer is a {\bf congruent number} if it is the area of a right triangle
with rational sides. The 3-4-5 triangle for example has the area $n=6$ so that 6 is a
congruent number. The 3/2,20/3,41/6 triangle has area $n=5$. The example $n=5$ shows that
one have to use rational numbers in general. If $x,y,z$ are the lengths of the triangle, then
the condition is $x^2+y^2=z^2, xy = 2n$. {\bf Rational Pythagorean triples} can
be generated with $x=u^2-v^2, y=2uv, z=u^2+v^2$. This leads to
congruent numbers $n = u v (u^2-v^2)$. For $u=3,v=2$ for example, one gets the
12-5-13 triangle with with area $30$. Fermat showed:

\satz{No square number can be a congruent number}

Fermat's proof from 1670 using decent can be found in a self-contained way in \cite{ConradCongruent}.
While integer solutions $(x,y,z)$ can be done by finite search
for a fixed $n$, the task to find rational solutions $x,y,z$ for a given $n$
can be difficult. For example, the smallest example for n=101 found by Bastien in 1914 is
$x=711024064578955010000/q$, $y=3967272806033495003922/q$, $z=4030484925899520003922)/q$
with $q=118171431852779451900$ \cite{Chandrasekar1998} shows that already for smaller $n$,
the smallest rational numbers $x,y,z$ solving the problem can become complicated.
Arabic mathematicians have known that numbers like $5,6,14,15,21,30,34,67,70,110,154,190$
were congruent numbers. Leonardo Pisano (Fibonacci) established that $n=7$ is a congruent
number with $(x,y,z) = (35/12,24/5,337/60)$ and conjectured that no square can be a congruent
numbers. Fermat then with his method of infinite descent proved that no square is a
congruent number. Already $n=1$ is interesting as it illustrates the {\bf decent method}:
if $n=1$ is congruent then $x^4=y^4+z^2$ has a non-trivial solution.
Let $a$ be a rational number such that $a^2+n,a^2-n$ are squares of rational numbers.
Then $x=\sqrt{a^2+n}+\sqrt{a^2-n}, y = \sqrt{a^2+n} + \sqrt{a^2-n}$, $z=2a$ is a solution
as $xy/2 = (a^2+n) - (a^2-n))/2 = n$. Work of 1922 by 
Louis Mordell related the congruent number problem to elliptic curves. If $u$ is so
that $u^2+n,u^2-n$ are rational squares, then $u^4-n^2$ is a rational square $v^2$
so that $u^6-n^2 u^2 = u^2 v^2$, with $x=u^2, y=uv$ this gives $y^2=x^3-n^2x$.
So, if $n$ is a congruent number, there is a rational point on the curve $y^2=x^3-n^2x$.
Kurt Heegner proved in his 1952 paper that if a prime is congruent to $5$ or $7$ modulo $8$, then $p$ is a
congruent number and that if a prime is congruent to $3$ or $7$ modulo $8$
then $2p$ is a congruent number \cite{Birch2004}.
Jerold Tunnell (a student of Tate) showed in 1983 \cite{Tunnell1983}
that the congruent number problem would have a full solution under
the {\bf Birch and Swinnerton-Dyer conjecture}, one of the Millenium problems.
Having that established would allow to test in finitely
many steps whether a given integer $n$ is a congruent number or not. 
\index{Congruent number}
\index{Fermat's right angle theorem}
\index{Pythagorean triples}
\index{Birch and Swinnerton-Dyer conjecture}

\section{Stark-Heegner theorem}

A {\bf imaginary quadratic field} $K=\mathbb{Q}[\sqrt{-n}]$ has {\bf class number} $1$
if there is a unique prime factorization in $K$. Carl Friedrich Gauss found already
9 cases $\{1,2,3,7,11,19,43,67,163\}$. These cases turned out to be all
and are now called {\bf Heegner numbers}.
\index{Heegner number}

\satz{There are exactly 9 imaginary quadratic fields of class number $1$.}

The theorem is now known as the {\bf Stark-Heegner theorem}.
Kurt Heegner proved this in 1952 \cite{Heegner1952}. The proof was for more than a decade
labeled to ``have a gap", but it got rehabilitated by 
Harold Stark in 1969 \cite{Stark1969} thanks also
to \cite{Birch2004} who was one of the first to 
recognize Heegner's achievement in his 1952 paper.
The introduction of Heegner's paper is a master piece, skillfully
pointing out how the {\bf class number theory} has relations to the 
congruent number problem that historically has led Fermat to his 
{\bf descent method}. As Stark pointed
out, the dismissal of Heegners proof must also have been due to 
{\bf professional bias} as there was just one step missing showing 
that a concrete equation $x^{24}-a x^8-16=0$ has a six degree factor 
whose coefficients are algebraic integers of degree $1$, to
which he refers to Weber's textbook \cite{Weber1898}. About the 
origin of bias \cite{Birch2004}: {\it Heegner was a fine
mathematician, with a rather low-grade post in a gymnasium in East Berlin.}
It was a widely held view that the trouble in Heegners proof should 
be traced to Weber. Today thanks to Stark, it is now clear that
the gap was actually not existent and Heegners proof correct.
Stark also points out that Weber's part is correct but 
could have given more details, if he had seen any need to do so. 
One can justify the name Stark-Heegner theorem because Stark 
not just gave a new clarified proof but took the trouble to investigate
whether there was indeed mistake in Heegner's proof. 
Bryan Birch certainly also played an important role in discovering
Heegner as also ``Heegner's numbers" got into the spot light in the 
context of the {\bf Birch and Swinnerton-Dyer conjecture} and the 
{\bf Gross-Zagier theorem}.
The largest Heegner number got a bit of a ``cult status" 
as it appears in {\bf Ramanujan's constant}
$e^{\pi \sqrt{163}}$ that is less than $10^{-12}$ close to the 
integer $640320^3+744$. This can be justified by the fact that 
if $n$ is a Heegner number, then the {\bf $j$-invariant} 
of $(1+\sqrt{-n})/2)$ is an integer and a $q$-expansion 
gives then a theoretical error is of the order $O(e^{-\pi \sqrt{163}})$.  
\index{Ramanujan constant}
\index{Stark-Heegner theorem}
\index{decent method}

\section{Equichordal point theorem}

If $C$ is a smooth convex curve in the plane a point
$P$ in its interior is called an {\bf equichordal point}
if all the line segments through $P$ have the same length.
For the circle $C$, this happens at the center. For the
polar curve $r(t) = 2+\sin(t)$, the center is an equichordal
point.

\satz{A convex curve can not have two equichordal points.}

The problem had been posed by Fujiwara in 1916 \cite{Fujiwara1916} and
appeared in a problem section of Blaschke, Rothe and Weitzenb\"ock:
\cite{BRW1917}. It seems that also Erd\"os was independently conjecturing
this as Gabriel Andrew Diracs work of 1952 indicates \cite{Dirac1952}.
The conjecture was proven by Marek Rychlik in 1997 \cite{Rychlik1997}
who established it more generally for star-like curves. The proof uses methods
from dynamical systems, complex analysis and algebraic geometry. 
\index{Equichordal point}
\index{Rychlik's theorem}

\section{Lucas fundamental theorem}

The {\bf Fibonacci sequence} $F(n)$ is defined by the second order
recursion $F(0)=F(1)=1$ and $F(n+1)=F(n)+F(n-1)$. 
When looking at the prime factorizations one can notice that 
the even terms $F(2n)$ have lots of prime divisors while
the odd terms $F(2n+1)$ have only a few.
Indeed, it follows from Lucas work that all primes appear as
factors of the even Fibonacci numbers.  Let $GCD$ denote
the greatest common divisor. 
 
\satz{$GCD(F(m),F(n)) = F(GCD(m,n))$.}

This {\bf fundamental theorem of Lucas} of 1878 \cite{Lucas1878} tells that 
the sequence $F(n)$ is a {\bf strong divisibility sequence}. 
Together with Lucas {\bf law of apparition} and {\bf Lucas law of repetition},
it implies that every integer divides infinitely many Fibonacci numbers. 
\cite{Lagarias2014}.
In the context or primality testing, Lucas also looked that the 
{\bf Lucas numbers}, $L(n)$ which satisfy the same recursion 
but have a different initial condition $L(1)=1,L(2)=3$. One has then 
$F(2n)=F(n) L(n)$. Lagarias proved in 1985 an anlogue of the
{\bf Chebotarev Density Theorem} using a method of Hasse.
He showed that the
density of prime divisors of the Lucas sequence is $2/3$ \cite{Lagarias1985}.
That article mentions that it is believed that the set of primes dividing
the terms $U(n)$ of any non-degenerate second order linear 
recurrence has a positive density and that this is conditionally true 
under the assumption of the {\bf generalized Riemann hypothesis}.
A bit about the history (see \cite{Lagarias1985}): 
the Fibonacci sequences appeared first in the third book
of ``Liber Abbaci" of Leonardo Pisano from 1227, a book that contains 90 sample
problems, with 50 from Arabic sources. It also contains 
the famous rabbit problem.  \'Edouard Lucas 
had been an artillery officer in the Franco-Prussian war and then was a 
high school teacher in Paris, who also was interested in recreational mathematics
and invented the tower of Hanoi problem \cite{Lucas1891}.
\index{Liber Abbaci}
\index{Leonardo Pisano}
\index{Fundamental theorem of Lucas}
\index{strong divisibility sequence}
\index{Fibonacci sequence}
\index{Lucas sequence}
\index{Tower of Hanoi problem}

\section{Hilbert distance}

The {\bf Hilbert distance} $d(x,y)$ is defined for points $x,y$ a bounded convex domain $X$
in a Hilbert space: construct the line through $x,y$. It intersects the boundary of $X$ in    
exactly two points $p,q$. The Hilbert distance is now defined as $d(x,y) = \frac{1}{2} \log(C(x,y,p,q))$,
where $C(x,y,p,q)=(|x-p| |y-q|)/(|y-p| |x-q|)$ is the {\bf cross ratio} between these four points.
Due to its projective invariance, the Hilbert distance defines then also
a {\bf Hilbert distance} on the {\bf projective space} $\mathbb{RP}^{n-1}$
which has the property that positive $n \times n$ matrices are contractions.
Lets call a metric on the projective space {\bf Perron-Frobenius} if it has 
this property. 

\satz{The Hilbert metric is the unique Perron-Frobenius metric.}

In the simplest case $\mathbb{P}^1$, elements are described as
${\bf t} = [1,t]$ with $t \in \mathbb{R} \cup \{\infty\}$. The Hilbert metric then 
is $d({\bf t},\bf{s}) = |\log(t/s)|$.
A positive matrix $A=\left[ \begin{array}{cc} a & b \\ c & d \end{array} \right]$
maps $[1,t]$ to $[1,(c+dt)/(a+bt)]$. 
David Hilbert defined the Hilbert metric in 1895 in a letter to Felix Klein \cite{Hilbert1895}.
The Hilbert metric between two points depends on the domain in which the points
are considered. The larger the domain, the smaller the distance. Also, if $z$ is 
in the line segment $[x,y]$, then $d(x,y) = d(x,z) + d(z,y)$. For strictly convex
region, there is a unique geodesic (with respect to this metric) connecting 
two points. It was Garret Birkhoff \cite{Birkhoff1957} 
and Hans Samelson \cite{Samelson1957}, who independently first suggested 
to use the Banach fixed point theorem to prove the 
{\bf Perron-Frobenius theorem} \cite{Perron1907,Frobenius1908,Frobenius1912} stating that 
a positive matrix has a unique maximal eigenvalue.  
\cite{LemmensNussbaum}.
Birkhoff called it the {\bf projective metric}.
For that, one only needs the mere existence of a Hilbert metric and not the
uniqueness.  Uniqueness is shown in \cite{KohlbergPratt}.
\index{Hilbert distance}
\index{Hilbert metric}
\index{projective metric}
\index{Perron Frobenius}

\section{Gross-Zagier}

The {\bf projective special linear group} $G=PSL(2,\mathbb{Z})$ is the group of integer
matrices $A$ of determinant $1$ for which the matrices $A$ and $-A$ are identified. 
It is also called the {\bf modular group} as its 
elements act as {\bf M\"obius transformations} $z \to T_{a,b,c,d}(z) = (az+b)/(cz+d)$ 
on the {\bf upper half plane} $H \subset \mathbb{C}$. A {\bf congruence subgroup} $\Gamma$ of $G$ 
is a subgroup of $G$ which has a {\bf principal congruence subgroup} $\Gamma(N)$, a set
of matrices in $G$ congruent to the identity matrix modulo $M$. The smallest $N$
for which this happens, is called the {\bf level} of $\Gamma$.
An important example is the {\bf Hecke congruence group} $\Gamma_0(N) =\{ T_{a,b,c,d}, N|c \}$. 
A {\bf modular elliptic curve} $E$ is a quotient $H/\Gamma$, where $\Gamma$ is a 
congruence subgroup of the modular group. Elliptic curves are the simplest positive-dimensional 
projective algebraic curves that carry a commutative algebraic group structure.
The set of rational points $E(\mathbb{Q})$ is finitely generated by the {\bf Mordell-Weil theorem},
so that $E(\mathbb{Q})$ is isomorphic to $\mathbb{Z}^r \times T$, where $r  \geq 0$ is called 
the {\bf rank of $E$} and $T$ is a finite Abelian group called the {\bf torsion subgroup of $E$}. 
The {\bf Birch and Swinnerton-Dyer conjecture} claims that $r$ is the order ${\rm ord}_{s=1}(L(E,s)$ 
of $L(E,s)$ at $s=1$, where the {\bf L-function} $L(s)$ for an elliptic curve $E$ over $K$ is an explicitly
given Dirichlet series $L(s)=\sum_{n=1}^{\infty} a_n n^{-s}$.
[It can defined as follows: for a prime $p$ let $\mathbb{F}_{p^e}$ denote
the field with $p^e$ elements. Define $t_1(E)=2$, $t_p^e(E)=p^e+1-|E(\mathbb{F}_{p^e})|$ and
the {\bf counting zeta function} $\zeta_p(z) = \exp(\sum_{e \geq 1} \frac{t_{p^e}(E)}{e} z^e)$ at $p$
and $1_E(p) = 1$ if $p$ does not divide $N$ and $1_E(p)=0$ if $p|N$. Then 
$\zeta_p(z)=(1-t_p(E) z + 1_E(p) p z^2)^{-1}$. The {\bf $L$-function} is then defined as
the {\bf Euler product} $L(s)=\prod_{p \; {\rm prime}} \zeta_p(p^{-s})$. While the Dirichlet series only 
converges for ${\rm Re}(s)$ larger than the abscissa of convergence, one knows from work in the 1970ies 
like Shimura that in the modular case, $L$ has an analytic continuation to all of $\mathbb{C}$.]
The {\bf j-invariant} $j(\tau)$ is a modular function of weight zero on $G$. It can be explicitly written down
and was originally used to represent isomorphism classes of elliptic curves. It is known that
the {\bf field of modular functions} is $\mathbb{C}(j)$. 
If $\tau$ is an element of an imaginary quadratic field with positive imaginary 
part, then $j(\tau)$ is an algebraic integer by a result of Theodor Schneider from 1937.
Now, a modular elliptic curve can be parametrized as $r(z)=(j(z),j(N z)) \in \mathbb{C}^2$,
where $N$ is the level of $\Gamma$.
If $\omega \in H$ is a quadratic irrational number (a number of the form $a+b\sqrt{D} \in H$ 
with rational $a,b$) which solves $N A \omega^2+B \omega+C=0$
then $\omega$ and $N \omega$ both have the same discriminant $D=B^2-4NAC$ so that 
$P=r(\omega) \in E(\mathbb{Q}(D))$. This $P$ is called a {\bf Heegner point} on $E$ \cite{Birch2004}.
The {\bf global canonical height} function $h: E \to \mathbb{R}$ is a function on $E$ with the
property that $h(Q)=0$ if $Q$ is a torsion point and such that the 
{\bf parallelogram law} $h(P+Q)+h(P-Q)=2h(P)+2h(Q)$ holds for all pair of points $P,Q$ on $E$. 
It is difficult to compute but the {\bf Gross-Zagier formula} \cite{GrossZagier1986} 
relates it in an explicit way with the order of the root at $1$ of the function $L$: 

\satz{The height $h(P)$ of a Heegner point is a non-zero multiple of $L'(1)$.}

This implies that if $L'(1)=0$, then $P$ is a torsion point and
that if $L'(1) \neq 0$, then the rank $r$ of $E$ is positive. 
Heegner points have been used to construct a rational point on the 
curve of infinite order. The theorem was later used to prove much of the Birch and Swinnerton-Dyer
conjecture for rank 1 elliptic curves. \cite{Birch2004} illuminates the history of the theorem.
\index{Heegner point}
\index{canonical height}
\index{modular elliptic curve}
\index{Birch and Swinnerton-Dyer}
\index{modular function}
\index{parallelogram law}
\index{Gross-Zagier formula}

\section{Schur determinant identity}

The {\bf Schur determinant identity} is an identity for partitioned matrices
$M=\left[ \begin{array}{cc} A & B \\ C & D \end{array} \right]$, where $A,B,C,D$
are all $n \times n$ matrices. Assume $A$ is invertible, one can write
$M=\left[ \begin{array}{cc} A & 0 \\ 0 & 1 \end{array} \right]$
                    $\left[ \begin{array}{cc} 1 & 0 \\ C & 1 \end{array} \right]$
                    $\left[ \begin{array}{cc} 1 & A^{-1} B \\ 0 & D-C A^{-1} B \end{array} \right]$.
Using the Cauchy-Binet product formula, one gets the {\bf Schur identity}

\satz{${\rm det}(M) = {\rm det}(A) {\rm det}(D-C A^{-1} B)$}

The matrix $D-C A^{-1} B$ is called the {\bf Schur complement}.
Given two $n \times m$ matrices $F,G$
one can compare the determinant of $A B = \left[ \begin{array}{cc} 1 & -F \\ G & 1 \end{array} \right]
             \left[ \begin{array}{cc} 1 & F  \\ 0 & 1 \end{array}  \right]$
with the determinant of $B A$ to get the {\bf Weinstein-Aronszajn identity}
${\rm det}(1+F^T G) = {\rm det}(1+G^T F)$.  
See \cite{Deift1978,TaoMatrixIdentities}. 
This identity also follows from the formula
${\rm det}(1+F^T G)  = \sum_P {\rm det}(F_P) {\rm det}(G_P)$ 
involving the summation
over all minors \cite{CauchyBinetKnill}. 
(Compare that the classical Cauchy-Binet formula for $n \times m$ matrices $F,G$
states ${\rm det}(F^T G) = \sum_P {\rm det}(F_P) {\rm det}(G_P)$ which is a sum over all 
$m \times m$ minors. For $n=m$, it becomes the product formula ${\rm det}(F G) = {\rm det}(F) {\rm det}(G)$.)
In \cite{TaoMatrixIdentities} many more identities are listed like 
${\rm det}(A+BC) = {\rm det}(A) {\rm det}(1+C A^{-1} B)$ if $A$ is invertible 
(which means especially ${\rm det}(A+B) = {\rm det}(A) {\rm det}(1+A^{-1} B)$
which is special case of the Schur identity for $C=D=1$) or 
${\rm det}(A+B) {\rm det}(A-B) = {\rm det}(B) {\rm det}(A B^{-1} A - B)$, if $B$ is invertible.
\index{Schur complement}
\index{Caucny-Binet}
\index{Minor}

\section{Herman's subharmonicity theorem}

If $(\Omega,\mathcal{A},\mu)$ is a probability space and $T$
an automorphism and $A \in L^{\infty}(\Omega,SL(2,\mathbb{C}))$
define the {\bf non-abelian Birkhoff product}
$A^n(x)= A(T^{n-1}x)) A(T^{n-2}(x)) \cdots A(Tx) A(x)$. 
An example is when $\Omega$ is a 2-manifold and $T$ an area-
preserving diffeomorphism $\Omega \to \Omega$ and $A(x)=dT(x)$
is the Jacobian. An other example is when $(Lu)(n) =u(n+1)+u(n-1)+V(T^nx) u(n)$
where the time equation $Lu = Eu$ leads to the transfer matrix 
$A(x) = \left[ \begin{array}{cc} E-V(x) & 1 \\ -1 & 0 \end{array} \right]$. 
Define $A^n(x) = A(T^{n-1}x) \cdots A(T(x)) A(x)$. The {\bf Lyapunov exponent}
$\lambda(A) = \lim_{n \to \infty} \frac{1}{n} \int_{\Omega} \log||A^n(x)|| d\mu(x)$ 
exists because it is a limit of a sub-additive sequence.
Assume $z \in \mathbb{C}^d \to SL(2,\mathbb{C})$ is analytic in the sense
that each matrix entry is an analytic function in each of the variables.
Assume also that $T: \mathbb{D}_r^d \to \mathbb{D}_r^d$ is analytic in a neighborhood
of the polydisc $\mathbb{D}_r^d$ and maps the boundary $\Omega=\mathbb{T}^d$ into itself
and that $T(0)=0$ and $T$ preserves the Haar measure on $\Omega$. Herman's theorem \cite{Her83} is

\satz{$\lambda(A) \geq \lambda(A(0)) = \log({\rm max}(\sigma(A(0))))$}

The reason is that $z \to \log(||A^n(z)||)$ is {\bf pluri-subharmonic} so that the integral
over the torus is bounded below by the Lyapunov exponent value at $0$. 
For example, if $p(z)=c(z+z^{-1})/2$ and $T(z) = w z$ with $w=e^{i \alpha}$
induces the dynamical system $T(\theta) = \theta+ \alpha \; {\rm mod} \; 2\pi$ 
on the boundary $\mathbb{T}^1$, then the Lyapunov exponent of 
$A(\theta) = \left[ \begin{array}{cc} c \cos(\theta) & -1 \\ 1 & 0 \end{array} \right]$
over the dynamical system is then larger or equal than $\log(c/2)$. 
The reason is that the Lyapunov exponent of 
$B(z) = z A(z) = \left[ \begin{array}{cc} c (z^2+1)/2 & -z \\ z & 0 \end{array} \right]$
is bounded below by the logarithm of the spectral radius of $B(0) = 
\left[ \begin{array}{cc} c/2 & 0 \\ 0 & 0 \end{array} \right] = \log(c/2)$. 
An other application is if $A \in L^{\infty}(\Omega,SL(2,\mathbb{C}))$ is 
arbitrary and $T: (\Omega,\mathcal{A},\mu) \to (\Omega,\mathcal{A},\mu)$ is 
a dynamical system, then for $A_\beta(x) = A(x) \left[ \begin{array}{cc} \cos(\beta) & -\sin(\beta) \\
                                                                            \sin(\beta) &  \cos(\beta) \end{array} \right]$,
the Lebesgue measure of values $\beta$ with $\lambda(A(\beta))>0$ is positive
if $A \notin SU(2,\mathbb{C})$ on some positive measure.
This can be used to show that the set of 
$A \in L^{\infty}(\Omega,SL(2,\mathbb{C}))$ with $\lambda(A)>0$ is dense \cite{knill_1992}.
The method of Herman has been extended in various way: \cite{SoSp91} use the 
{\bf Jensen inequality in complex analysis} to show that
for a non-constant real analytic $f$ and 
$A(x) = \left[ \begin{array}{cc} E-c f(x) & 1 \\ -1 & 0 \end{array} \right]$
and dynamical system $T(x) = x+\alpha \; {\rm mod} \; 2\pi$ with irrational $\alpha$, the Lyapunov exponent of $A$ is
positive for all $E$, if $c$ is large enough. 
Herman's and the Soret's Spencer theorem are the starting point in \cite{Bourgain2005}.
\index{subharmonic}
\index{polydisc}
\index{Lyapunov exponent}
\index{Haar measure}

\section{Gabriel's theorem}

A {\bf quiver} $(V,E)$ is an other word for a {\bf multidigraph}, a directed graph in which multiple 
directed connections = {\bf arrows} and self connections = {\bf loops} are allowed. 
The {\bf graph} defined by $(V,E)$ is the multigraph one obtains 
if the directions of the arrows are ignored.  A {\bf representation} $V$ of a quiver assigns
vector space over an algebraically closed field 
to each node $x \in V$ and a linear map $V( x \to y ): V(x) \to V(y)$ 
attaching to each arrow $x \to y$ a linear map. 
It is indecomposable if it can not be written as the direct sum of smaller positive
dimensional representations. 
A quiver is of {\bf finite type}, if it has only finitely many isomorphism classes of 
indecomposable representations. The {\bf Quiver diagrams}
are formed by the {\bf simply laced Dynkin diagrams} $A_n,D_n,E_6,E_7,E_8$. 
Gabriel's theorem classifies the connected quivers of finite type.

\satz{Connected quivers of finite type correspond to quiver diagrams}.

The theorem was proven by Peter Gabriel in 1972 \cite{Gabriel1972}. Written in German, the article uses the
word ``K\"ocher" is used there for quiver. Peter Gabriel (1933-2015) was a French
and Swiss mathematician also known as Pierre Gabriel. On Wikipedia, he is listed as a student of Alexander
Grothendieck with a thesis done in 1960 on Abelian categories 
\cite{Gabriel1962} (on his personal website which is still active, Henri Cartan was listed as the Jury,
and Jean Pierre Serre as the rapporteur, on the Mathematics Genealogy page, Jean-Pierre Serre is listed as
the advisor. [According to Serre, Gabriel wrote an independent thesis and pointed out that in 1960, the advisor status
had not been yet as formal as today. In the published article, it is also not visible
who the formal advisor was.] Remarkably,  Gabriel was doing his military service 1960-1962 just after finishing his thesis
and the Abelian category paper was submitted in 1961. Gabriel worked at the University of Z\"urich from 
1974-1998.
\index{quiver}
\index{ADE diagrams}
\index{Gabriel's theorem}

\section{Zeckendorf representation}

Let $F(n)$ denote the $n$'th {\bf Fibonacci number}. It is defined by
the recursion $F(n+1) = F(n) + F(n-1)$ and $F(0)=0,F(1)=1, F(2)=1$.
Given a positive integer $n$, a representation
$n=\sum_{k=0}^m F(c(k))$ with $c(k) \geq 2$ and $c(k+1)>c(k)+1$
is called a {\bf Zeckendorf representation}. The
finite sequence $n_F = (c(0),c(1),\dots,c(m))$ a notation of Knuth,
this is called the {\bf Fibonacci coding} of $n$.
For example, $11=(1010000)_F$ and $13=(10000000)$.            

\satz{Every positive integer has a unique Zeckendorf representation.}

Edouard Zeckendorf published this in 1972 and mentions to have proven it
already in 1939. Lederkerker independently found the result in 1952 \cite{Lekkerkerker}.
The proof of existence and uniqueness can both be done by induction. 
As Donald Knuth realized \cite{Knuth1988}, the Zeckendorf representation of an integer 
leads to an {\bf associative multiplication}
$x \circ y = \sum_{i=0}^{m(x)} \sum_{j=1}^{m(y)} F(c_i(x) + c_j(y))$ 
for positive integers $x,y$. This is called the {\bf Fibonacci product}.
The proof of associativity is the realization that $(x \circ y) \circ z$ is equal to
$\sum_{i=0}^{m(x)} \sum_{j=1}^{m(y)} \sum_{k=1}^{m(y)} F(c_i(x) + c_j(y) + c_k(z))$. 
Knuth mentions that the Fibonacci product asymptotically satisfies $x \circ y \sim \sqrt{5} x y$
and that the multiplication $x * y = xy + [\phi x] [\phi  x]$ by Porta and Stolarsky
is asymptotically $(1+\phi^2) mn \sim 3.62 mn$, where $\phi=(1+\sqrt{5})/2$ is the
{\bf golden ratio}.
\index{Fibonacci coding}
\index{Zeckendorf representation}
\index{Golden ratio}
\index{Zeckendorf multiplication}

\section{Turan's theorem}

A finite simple graph $G=(V,E)$ has $n=|V|$ vertices and $m=|E|$ edges.
A {\bf $p$-clique} is a complete subgraph of $G$ with $p$ vertices.
The $1$-cliques can be identified with $V$ and the $2$ cliques can
be identified with $E$. {\bf Tur\'an's graph theorem} \cite{Turan1941} is

\satz{If $m>\frac{p-2}{p-1} \frac{n^2}{2}$, then $G$ has a $p$-clique.}

It assures that a triangle free graph can have at most $n^2/4$ edges
so that if a graph has more than a quarter of all edges connected,
there must be a triangle in it. This is called Mantel's theorem from 1907.
The {\bf Turan graphs} are graphs of the form $P_{n_1} + ... + P_{n_k}$
where for all $n_j$, we have $n_j \in \{a,a+1\}$ for some integer $a$.
For $n_j=n/(p-1)$ are constant, these are graphs without $p$-cliques and
$B(p-1,2) (n/(p-1))^2 = \frac{p-2}{p-1} n/2$. This shows that the result
is sharp.  \cite{AigZie} contains four short proofs, the first one doing
induction with respect to $n$. See also \cite{Aigner1995} who states that the
theorem of Tur\'an initiated extremal graph theory and that the theorem
had been rediscovered man8y times.
\index{Turan graphs}
\index{Turan graph theorem}

\section{The Szpilrajn-Marczewski theorem}

A finite simple graph $\Gamma=(V,E)$ is {\bf represented} by a set of sets $G$
if $V=G$ and $E=\{ (x,y) \; | x \neq y, x \cap y \neq \emptyset\}$.
The graph $\Gamma$ is the {\bf connection graph} of the set of sets $G$.

\satz{Every graph is the connection graph of a set of sets.}

An arbitrary set of sets is sometimes also called a {\bf multigraph}.
The theorem shows that from the point of view of connectivity, a multigraph
can be studied by its connection graph. It does not encode other properties like
subset property. The set of sets $G=\{ \{ 1,2\}, \{2,3\} \}$ and the set of sets
$H=\{ \{ 1,2\}, \{2\} \}$ both have the same connection graph $K_2$.
The theorem was shown by Edward Szpilrajn-Marczewski (1907-1976) in 1945 \cite{Szipilrajn-Marczewski}. 
The Polish mathematician was born Szpilrajn but changed his name while hiding from
Nazi persecution.
Erd\"os, Goodman and Posa showed in 1964 that one can realize any
graph of $n$ vertices as a set of subsets of a set with $[n^2/4]$ elements and that for 
$n \geq 4$, one can even require all sets to be distinct. The result is sharp for $n \geq 4$
The smallest number $d(n)$ of sets needed to represent every graph with n vertices
satisfies $d(2)=2, d(3)=3)$ and $d(n) = [n^4/4]$ for $n \geq 4$. For example $d(4)=[4^2/4]=4$
and $d(5)=[5^2/4]=6$.  The Erd\"os-Goodman-Posa proof is done by induction $n \to n+2$ and
by first establishing the cases $4$ and $5$ which can be done by looking at all cases.
The Szpilrajn-Marczewski theorem has been abbreviated SM theorem in
\cite{ErdoesGoodmanPosa} and is a much referenced theorem in {\bf intersection
graph theory}. The theorem does not assume the graph to be finite. `
\index{connection graph}
\index{Szpilrajn-Marczewski theorem}

\section{Sakai theorem}

Let $\mathcal{B}(H,\mathbb{C})$ denote the Banach algebra of all
bounded linear operators on a Hilbert space $H$. The {\bf commutant} $X'$ of a subset $X \subset \mathcal{B}(H)$
is the set of all elements in $\mathcal{B}(H)$ that commute with every element in $X$.
Because of the contra-variance condition $X \subset Y \Rightarrow Y' \subset X'$, the {\bf bicommutants}
satisfy $X'' \subset Y''$ so that, using $\mathcal{B}(H)'=\mathbb{C},\mathbb{C}'=\mathcal{B}(H)$,
any subset $X$ is contained in the {\bf bicommutant} $X''$.
A subalgebra $X$ satisfying $X=X''$ is called a {\bf von-Neumann algebra}. It is called a {\bf factor}
if its {\bf center} $X \cap X'$ is $\mathbb{C}$. Von Neumann showed the {\bf bicommutant theorem} stating
that $X''=X$ is equivalent to $X$ being weakly closed. (The {\bf weak operator topology} means pointwise
convergence in the sense that $A_n \to A$ in the weak operator topology if and only if
for every pair $f,g \in H$ one has $(g,A_n f) \to (g,A f)$, meaning that given a
basis in $H$ that the matrix elements of operators converge pointwise.)
The bicommutant theorem is remarkable as it equates the algebraic bicommutant condition with
the topological weak-closed condition. Von Neumann algebras can also be defined more abstractly using $C^*$ algebras
without referral to operator algebras but the {\bf GNS construction} justifies the more intuitive operator algebra
definition. Like the bicommutant theorem, there are other characterizations of von Neumann algebras. One of them
is {\bf Sakai's theorem}
\index{von Neumann algebra}
\index{Sakai theorem}
\index{commutant}
\index{bicommutant}
\index{weakly closed}
\index{type I factor}
\index{type II factor}
\index{type III factor}
\index{factor (von Neumann algebra)}

\satz{A $C^*$ algebra is von Neumann if and only if has a pre-dual.}

Sakai's theorem was proven in 1956 \cite{Sakai1956}.
Examples of von Neumann algebras are $X=\mathcal{B}(H)$, any finite dimensional subalgebra $X$ of the algebra
of operators $\mathcal{B}(H)$ or any algebra $X=(S \cup S^*)''$ generated by an arbitrary subset $S$ of $\mathcal{B}(H)$.
For example, every commutative von-Neuman algebra is of the form $L^{\infty}(\Omega,\mathcal{A},\mu)$;
the predual is then $L^1(\Omega,\mathcal{A},\mu)$. Since $L^{\infty}(\Omega,\mathcal{A},\mu)$ (acting as
multiplication operators on $H=L^2(\Omega,\mathcal{A},\mu)$)
for a measure $\mu$ completely encodes the {\bf measure theory of} $(\Omega,\mathcal{A},\mu)$, 
the theory of von Neumann algebras has been seen as {\bf non-commutative measure theory}.
This is the picture of Alain Connes \cite{Connes}. Von Neumann algebras are pretty well understood: each is a direct
integral of factors. Factors are classified as type I (meaning that it has a non-zero minimal projection like 
operator algebras on Hilbert spaces), type II (meaning that there is a non-zero finite projection) or then type $III$
(meaning that it contains no non-zero finite projection). There are other characterizations of von Neumann algebras:
the {\bf Kaplanski density theorem states} that if $A$ is a $C^*$ subalgebra of an operator algebra $\mathcal{B}(H)$
then the unit ball of $A$ is strongly dense in the unit ball of the weak closure of $A$.
This implies that a subalgebra $M$ of $\mathcal{B}(H)$ containing $1$ is a 
von Neumann algebra if and only if the unit ball of $M$ is weakly closed. More references are
\cite{ReedSimon,Blackadar2005,Dixmier,Jones2015}.
\index{non-commutative measure theory}

\section{Takens's theorem}

Let $M$ be a $d$-dimensional manifold and $T:M \to M$ be a smooth map from $M$ to $M$.
A compact $T$-invariant set $A \subset M$ is called an {\bf attractor} for $T$
if there there is an open neighborhood $N$ of $K$ such that 
$\bigcap_{n \geq 0} T^n(N)= A$. It is called a {\bf minimal attractor} if no
proper sub attractor exists. 
The map $T$ is called {\bf partially hyperbolic} if the {\bf Lyapunov exponent}
$\lambda(\mu) = \lim_{n \to \infty} n^{-1} \int_A  \log|dT^n(x)| \; d\mu(x)$ is non-zero
for some $T$-invariant measure $\mu$ on $A$, where $dT(x)$ is the Jacobian matrix.
The partial hyperbolic attractor is called {\bf strange} if it is not a countable union
of lower dimensional sets homeomorphic to varietes in $M$. 
This happens for example if $A$ is a {\bf fractal}, meaning that the 
{\bf Hausdorff dimension} of $A$ is not an integer. 
A {\bf Takens embedding} of $M$ is given by a transformation $T$ and a smooth $C^2$ function
$f:M \to \mathbb{R}$ and an integer $k$ and defined as the 
{\bf time series} $x \to (f(x),f(T(x)),\dots,f(T^{k-1}x)) \subset \mathbb{R}^k$.  
One can often reconstruct $M$ and so also the attractor $A$ from such measurements.
This happens Bair generically in $C^2(M,M) \times C^2(M, \mathbb{R})$.     

\satz{For a Bair generic set of pairs $T,f$, a Takens embedding exists.}

One can therefore use a dynamical system $T:M \to M$ to embed $M$ into some Euclidean      
space $\mathbb{R}^k$. This {\bf Takens's embedding theorem} is           
analogue to the {\bf Whitney embedding theorem} which assures
that if $f$ is allowed to be $\mathbb{R}^m$ valued, then $A$ can be embedded in 
$\mathbb{R}^{m}$, even for a time series with $k=1$ observation so that 
no dynamics is needed: $f: M \in \mathbb{R}^m$ embeds $M$ and so $A$ into a Euclidean space. 
The significance of the Takens's theorem is that one can ``see" $M$ or the attractor $A$
using a time series of a single real {\bf observable} $f:M \to \mathbb{R}$ and
then use {\bf time}, that is the dynamical system, to generate the {\bf coordinates} of the embedding.
This is extremely practical. One can for example observe the times, when a drop leaves a faucet and
use the differences of the times between two drops to create an attractor without having any model
of drop formation. 
The time series of course does work in general as the functions $f$ and the transformation $T$ must
be interesting enough. For the identity $T$ for example, the time series does not give enough information. 
A special case is if $A$ consists of a single point $a$ which is hyperbolic in the sense that        
all eigenvalues of the Jacobian matrix $dT(a)$ are smaller than $1$ in absolute value. In that case, 
the manifold $M$ is the {\bf stable manifold of $a$} and that remains true for an open set of 
transformations near $T$. Takens theorem then implies that for a generic $C^2$ function $f$, 
one can chose a $k$ such that the time series reconstructs the manifold $M$. The same works if $A$
is a {\bf hyperbolic attractor}, because the {\bf structural stability} of $T$ allows then to
restrict the genericity statement to the function $f$. 
Floris Takens 1940 -2010 was a Dutch mathematician. Together with David Ruelle, he introduced the
notion of {\bf strange attractor}. See \cite{DGS} for dynamical systems in general.
Takens's article is in \cite{RandYoung1981} (p 366-381).
\index{Takens's embedding theorem}
\index{time series}
\index{Attractor}
\index{strange attractor}
\index{hyperbolic attractor}
\index{partially hyperbolic attractor}

\section{Perfect graphs}

A finite simple graph is called {\bf perfect} if every induced subgraph
has a {\bf chromatic number} (minimal number of colors needed for a vertex coloring)
which is equal to the {\bf clique number} the graph. (The clique number is
maximal number $n$ of vertices for which there exists a complete subgraph $K_n$ with that number of vertices).
A finite graph satisfies the {\bf Berge condition}, it none of the induced subgraphs
are cyclic graphs $C_{2n+1}$ with $n \geq 2$ nor that it is the complement of such a cyclic graph.
The {\bf strong perfect graph theorem} states:

\satz{The set of perfect graphs is the set of Berge graphs.}

Because the odd cycle condition is invariant under graph complement formation,
the following {\bf weak perfect graph theorem} follows:
if $G$ is perfect, then its graph complement is perfect.
Examples of perfect graphs are trees, bipartite graphs,
wheel graphs with even boundary length or Barycentric refinements of graphs
(the graph in which the cliques are the vertices and two cliques are connected if
one is contained in the other, where obviously the dimension function is a coloring
and agrees with the clique number). The strong perfect graph conjecture had been 
conjectured by Berge in 1961 \cite{Berge1961}.  Maria Chudnovsky, 
Neil Robertson, Paul Seymour  and Robin Thomas proved the theorem in 2006 
\cite{StrongPerfectGraph}.
\index{perfect graphs}
\index{Berge graphs}
\index{strong perfect graph conjecture}
\index{strong perfect graph theorem}

\section{Kochen-Specker theorem}

Let $H$ be a {\bf Hilbert space} and let $\mathcal{X}$ denote the set of {\bf self-adjoint}
operators on $H$. These operators $A$ are also known as {\bf quantum mechanical observables}.
The mathematical frame work of quantum mechanics considers a time evolution $\psi = i L \psi$
with a Hamiltonian $L$ and then does for $A \in \mathcal{X}$
produce data $\langle \psi(t),A \psi(t) \rangle$ (Schr\"odinger picture) or 
$\langle \psi, A(t) \psi \rangle$ with $A(t) = U(t)^* A U(t)$
with unitary $U(t) = \exp(i L)$  (Heisenberg picture). Since $\mathcal{X}$ is non-commutative, one
can not expect to do measurements as in the classical calculus. The non-commutativity is illustrated
best with the famous {\bf anti-commutation relation} $[P,Q]= i$ which holds for the self-adjoint operators
$Pf(x)=i f'(x), Q(x) = x f(x)$ on $L^2(\mathbb{R})$ which represent momentum and position of a
particle on the real line. Before John Bell and Simon Kochen and Ernst Specker, it was not excluded that one could
use some hidden variables and still be close to a classical theory. By formulating this precisely,
one can also produce theorems. A function $v: \mathcal{X} \to \mathbb{R}$ is called a
{\bf classical value function}, if it is linear $\mathcal{X}$ and satisfy
$f(v(A)) = v(f(A))$ for all continuous functions $f$ as well as $v(A B) = v(A) v(B)$.
In other words, $v$ is a multiplicative linear functional
on $\mathcal{X}$, honoring the functional calculus and being compatible with multiplication.
For a continuous real function, the value $f(A)$ is defined by the functional calculus
which exists by the spectral theorem for any self-adjoint operator. 
The Kochen-Specker theorem is a {\bf no-go theorem}:

\satz{If ${\rm dim}(H) \geq 3$, there is no classical value function.}

It was proven by Simon Kochen and Ernst Specker in 1967 \cite{KochenSpecker}
even in the case when the dimension is $3$ or higher
and complements {\bf Bells theorem} on ``hidden variables". An important precursor was
Gleason's theorem. Kochen and Specker show more generally
that there is no partial Boolean algebra $D$ has no homomorphism into $\mathbb{Z}_2$.
It is refreshingly simple and elegant especially, considering the
difficulties that surround interpretations of quantum mechanics.  A bit simpler
is the argument if the dimension of $H$ is assumed to be $4$ or higher:
let $u_1,u_2,u_3,u_4$ be four orthogonal vectors in $H$ and let $P_k$ be the projection 
operators onto the line spanned by $u_k$. They satisfy $P_1+P_2+P_3+P_4=1$ so that by 
linearity, $v(P_1)+v(P_2)+v(P_3)+v(P_4)=1$. 
The condition $v(A B) = v(A) v(B)$ implies for {\bf projections} $P$ 
(elements in $\mathcal{X}$ satisfying $P^2=P$) that $v(P^2)=v(P)=v(P) v(P)$ so
that $v(P)=0$ or $1$. The linearity condition now implies that exactly one value is $1$. 
\cite{Kernaghan1994} simplifies \cite{Peres1991} uses the following list of 11 inconsistent
equations for 20 vectors which can not be satisfied because each vector appears 2 or four times
but on the left one has column sum which is 11 and so odd.
\begin{tiny}
\begin{eqnarray*}
    1&=&v([1,0,0,0]) + v([0,1,0,0]) + v([0,0,1,0]) + v([0,0,0,1]) \\
    1&=&v([1,0,0,0]) + v([0,1,0,0]) + v([0,0,1,1]) + v([0,0,1,-1]) \\
    1&=&v([1,0,0,0]) + v([0,0,1,0]) + v([0,1,0,1]) + v([0,1,0,-1]) \\
    1&=&v([1,0,0,0]) + v([0,0,0,1]) + v([0,1,1,0]) + v([0,1,-1,0]) \\
    1&=&v([-1,1,1,1]) + v([1,-1,1,1]) + v([1,1,-1,1]) + v([1,1,1,-1]) \\
    1&=&v([-1,1,1,1]) + v([1,1,-1,1]) + v([1,0,1,0]) + v([0,1,0,-1]) \\
    1&=&v([1,-1,1,1]) + v([1,1,-1,1]) + v([0,1,1,0]) + v([1,0,0,-1]) \\
    1&=&v([1,1,-1,1]) + v([1,1,1,-1]) + v([0,0,1,1]) + v([1,-1,0,0]) \\
    1&=&v([0,1,-1,0]) + v([1,0,0,-1]) + v([1,1,1,1]) + v([1,-1,-1,1]) \\
    1&=&v([0,0,1,-1]) + v([1,-1,0,0]) + v([1,1,1,1]) + v([1,1,-1,-1]) \\
    1&=&v([1,0,1,0)]) + v([0,1,0,1])  + v([1,1,-1,-1]) + v([1,-1,-1,1])  \; .
\end{eqnarray*}
\end{tiny}
\index{Kochen-Specker theorem}
\index{Quantum mechanics}
\index{Hidden variables}

\section{Perfect difference sets}

A subset $D$ of $\mathbb{Z}_m$ is called a {\bf perfect difference set} if
every nonzero number in $\mathbb{Z}_m$ can be written uniquely as $a-b$ for $a,b \in D$.
An example for $m=13$ is $D=\{1,2,5,7 \} \subset \mathbb{Z}_{13}$.
For $D$ to exist we need $m = n^2+n+1$ and $|D|=n+1$. The number $n$ is called the
{\bf order} of the perfect difference set. Any perfect difference set $D$ produces a
{\bf finite projective plane} $P(2,n)$ with $m=n^2+n+1$ lines. Singer showed in 1938 \cite{ Singer1938}
that perfect difference sets exist if $n = p^k$ is a prime power:
\index{perfect difference set}
\index{order of a perfect difference set}
\index{finite projective plane}

\satz{For every prime power $n=p^k$ there exists a finite projective plane}

Singer obtained the perfect difference set in the following way:
Let $\zeta$ be generator of the multiplicative group in the Galois field $G_3=F_{q^3}^n$
which is a Galois extension of $G_1=\mathbb{F}_{q^n}$, then $\zeta$ is the root of an
irreducible cubic polynomial in $G_1$ so that every element can be written
as $a+b \zeta + c \zeta^2, a,b,c \in G_1$. Every element different from $0$ in $G_3$ can be
written as $\zeta^k$. Look at all elements
$D = \{ k ,  \zeta^k = a+b \zeta$ for $a,b \in G_1 \} \cup 0$. Two such elements are called equivalent
if one is the multiple of the other. The equivalence classes partition all numbers
into $n+1$ equivalence classes. If they are written as $a_i + b_i \zeta = \zeta^{k_i}$, then
the set of exponents $k_i$ is a perfect difference set.
The {\bf prime power conjecture} claims that for any finite projective plane the order
is a prime power. One already does not know
whether there exists a projective plane of order $n=12$.
The prime power conjecture has been verified for all $n \leq 20 \cdot 10^9$ by Gordon.
Sarah Peluse recently showed \cite{Peluse2020} that the number of positive integers $n<N$ such that
$Z_{n^2+n+1}$ contains a perfect difference set is asymptotically $N/\log(N)$
giving more evidence for the prime power conjecture.
Perfect difference can be used to define {\bf Sidon sets} 
if $a+b=c+d$ for $a,b,c,d \in D$, then $\{a,b\} = \{c,d\}$. 
Small sets typically are Sidon sets. Sidon
sets $D$ can not be too large as $|D| (|D|+1)/2 < 2n$ implies $|D| < 2 \sqrt{n}$. 
The set $D=\{ (x,x^2), x \in \mathbb{Z}_p\}$ is a Sidon set in $\mathbb{Z}_p^2$
\index{Sidon sets}
\index{prime power conjecture}

\section{Trace Cayley-Hamilton theorem}

For a $n \times n$ matrix $A$, let $p_A(x) = {\rm det}(A-x)=\sum_{k=0}^n c_{n-k} x^k$
denote its {\bf characteristic polynomial}. The {\bf Cayley-Hamilton theorem} $p_A(A)=0$ assures that
$\sum_{k=0}^n c_{n-k} A^k=0$. While obvious for matrices which allow diagonalization (like normal operators),
the Cayley-Hamilton theorem is remarkably non-shallow.
The {\bf trace Cayley-Hamilton theorem} is

\satz{$k c_k + \sum_{j=1}^k {\rm tr}(A^j) c_{k-j} = 0$}

This implies that if all trace powers are zero, then $p_A(x)=(-x)^n$.
The reason for the name trace-Cayley-Hamilton theorem is that for $k \geq n$, the result
can be obtained from the Cayley-Hamilton theorem $\sum_{j=0}^n c_{n-j} A^j$ by multiplying with
$A^{k-n}$ and taking traces. The trace Cayley Hamilton theorem 
implies also that if two $n \times n$ matrices have the same traces 
${\rm tr}(A^k)={\rm tr}(B^k)$ for $k=1, \dots, n$, then $A,B$
have the same characteristic polynomial and so are isospectral. 
This is extremely useful as computing the traces of $n$ matrices can be 
more convenient than computing the characteristic polynomial. One can use the theorem especially
in theoretical settings better.
For {\bf normal matrices} one can conclude that $A$ is the zero matrix if
${\rm tr}(A^k)=0$ for $k=1, \dots, n$. See \cite{Grinberg2019,Zeilberger1993}. 
For moment problems see \cite{Schmuedgen2017}.
The Cayley-Hamilton theorem was first tackled in 1984 by William Rowan Hamilton in
the context of quaternions, meaning for $n=2$ complex or $n=4$ real matrices.
Arthur Cayley stated the theorem in 1858 for $n \leq 3$ but only proved $n=2$. 
In 1878, the general case was proven by Ferdinand Georg Frobenius. 

\index{Cayley-Hamilton}
\index{trace Cayley Hamilton}
\index{matrix computation}
\index{trace powers}

\section{Maximal permanent}

The {\bf permanent} of a $n \times n$ matrix $A$ is ${\rm per}(A) = \sum_{\pi} \prod_{i=1}^n A_{i,\pi(i)}$, where
the sum is over all permutations of $\{1,2, \dots, n\}$. It takes the Leibniz definition
 ${\rm det}(A) = \sum_{\pi} {\rm sign}(\pi) \prod_{i=1}^n A_{i,\pi(i)}$ of the determinant
determinant but ignores the signatures ${\rm sign}(\pi)$ of the permutations.
Unlike determinants which can be computed in polynomial times using row reduction, there is no
polynomial way known to compute permanents in polynomial time.
A {\bf probability vector} $p=(p_1, \dots, p_n)$ is an element in $\mathbb{R}^n$ for which all entries are in $[0,1]$
and add up to $1$. A $n \times n$ matrix is {\bf doubly stochastic}, if each row and each column of $A$ are
probability vectors. In 1926 Bartel van der Waerden conjectured that the maximal permanent which a doubly stochastic
$n \times n$ matrix can have, is obtained if all entries are $1/n$. These are the matrices with {\bf maximal entropy}
in the sense that the {\bf Shannon entropy} $S(p)=-\sum_{k=1}^n p_k \log(p_k)$ is maximal for each column or row of the matrix.

\satz{Doubly stochastic maximal permanent $\Leftrightarrow$ maximal entropy.}

The van der Waerden conjecture was proven in 1980 by B\'ela Gyires \cite{Gyires1980}
and in 1981 by G.P. Egorychev and by D.I. Falikman.  In \cite{Gyires1996} it was pointed out
that the conjecture had already been proven in 1977 \cite{Gyires1977}. For permanents, see
\cite{MincPermanents}. B\'ela Gyires was a Hungarian mathematician who lived from 1909 to 2001.
In his last paper \cite{Gyires2001}, Gyires gives an other account on the proof of the van 
der Waerden conjecture and two proofs.

\index{permanent}
\index{probability vector}
\index{stochastic matrix}
\index{doubly stochastic matrix}
\index{Shannon entropy}
\index{van der Waerden conjecture}

\section{Billiards in polygons}

A convex compact polygon  in $\mathbb{R}^2$ defines a {\bf billiard dynamical system}.
Parametrize the boundary by $x \in \mathbb{T} = \mathbb{R}/\mathbb{Z}$.
Given $(x_1,x_2) \in \mathbb{T}^2$ where both points are not at a vertex of the polygon,
we get a new point $x_3$ such that the path $x_1,x_2,x_3$ satisfies the law of reflection
at $x_2$. The set of points in $\mathbb{T}^2$ for which no future point $x_k$ is a vertex
has full measure. A point $x_0$ is called a {\bf periodic point} if $x_n=x_0$
for some $n>0$ and $x_k$ are all points not on vertices of the polygon. 
It is unknown already in the case of an obtuse triangle, whether a
periodic point exists. Fagnano already observed in 1775 that any acute triangle has a periodic
trajectory, the {\bf orthopic triangle}. A polygon is called {\bf rational} if all
angles $\alpha_j$ have the property that the angles $\alpha_j/\pi$ are rational.
\index{orthopic triangle}
\index{Fagnano triangle}

\satz{A rational polygon has a periodic orbit.}

Actually, there is a dense set of directions $\theta$ for which there is a periodic orbit.
This is called the {\bf Masur theorem} named after Howard Masur who proved this in 1986
\cite{Masur1986} by reducing the problem to flows defined by $e^{i \theta} \phi$
where $\phi$ is a holomorphic $1$-form on a compact Riemann surface $R$ of genus $ \geq 2$.
More generally, if $q$ is a holomorphic quadratic differential on such an $R$, there exists
a dense set of $\theta$ such that $e^{i \theta} q$ has a closed regular vertical trajectory.
The existence theorem uses {\bf Teichm\"uller theory}. The basic questions about billiards
in polygons has been raised by Carlo Boldrighini, Michael Kean and
Federico Marchetti in 1978 \cite{BoldrighiniKeaneMarchetti}.
\index{Masur theorem}

Billiards in polygons are also interesting from an ergodic point of view. A Bair generic
polygon produces an ergodic flow. For rational polygons, this is not the case as
the directions of the flow stay in a finite set generated by the rational angles $\pi n_i/m_i$
at the vertices. There is then an interval $[0,\pi/n)$ which parametrizes invariant hypersurfaces
in the phase space. One knows that for Lebesgue all directions $\theta \in [0,\pi/n)$ the flow is uniquely
ergodic, even weakly mixing but not mixing and has zero entropy.
This implies that there exists a generic set of ergodic and even weakly mixing (non mixing)
polygons (they are then non-rational) with $n$ vertices.
For more on billiards in polygons, see \cite{Gut86,Tab95,HalbeisenHungerbuehler2000}.

\section{Elasticity}

If $G \subset \mathbb{R}^n$ be an open and connected domain. For a {\bf vector field} $v: G \to \mathbb{R}^n$ and $x \in G$
denote with $v(x)=(v^1(x), \dots, v^n(x))$ its coordinates. Let $dv^{i}_j(x) = \partial_j v^i(x)$ denote the {\bf Jacobian
matrix} of $v$ at $x$. Let $||v||_{H^1}$ denote the {\bf Sobolev norm} obtained from
the {\bf inner product} $\langle v,w \rangle_{H^1} = \int_G v(x) \cdot w(x) + {\rm tr}(dv^T dw) \; dx$
on smooth vector fields and let $H^1(G)$ be the Hilbert space obtained by completing this set of vector fields
with respect to that norm. Let $(\partial_i v^j + \partial_j v^i)/2$ be abbreviated as $dv^s(x) = (dv^T(x)+dv(x))/2$
and denote the {\bf symmetric part of Jacobian} matrix at $x$. Let $||v||_{H^1_S}$ denote the {\bf symmetrized Sobolev norm}
obtained from the inner product $\langle v,w \rangle_{H^1_S} = \int_G v(x) \cdot w(x)+ {\rm tr}(d^sv^T(x) d^sw(x)) \; dx)$.
This  means that the {\bf Hilbert-Schmidt product} ${\rm tr}(A^T B)$ of the Jacobian matrices $A=dv$ and $B=dw$ is replaced
by the Hilbert-Schmidt product of the symmetrized Jacobian matrices $d^sv$ and $d^sw$.

\satz{There exists $C=C(G)$ such that $||v||_{H^1} \leq C ||v||_{H^1_S}$.}

This inequality is called the {\bf Korn inequality}. It is used in linear elasticity and continuum mechanics.
The constant $C$ is called the {\bf Korn constant} of $G$. The inequality had first been established by Arthur Korn in 1909
\cite{Korn1909} in the case $G=\mathbb{R}^n$, where for smooth vector fields $v$, we have using integration by parts
$\int_G |d^s f(x)|^2 dx = \int_G |df(x)|^2/2 + \int_G ({\rm div}(f))^2 \; dx$
so that the constant $C=2$ would do. See \cite{Korninequality2020}.
The inequality has been generalized to $W^{1}(G)$ if the region is bounded with
Lipschitz boundary. It has also been generalized to other Sobolev spaces $W^{1,p}(G)$ for $p \in (1,\infty)$ if the boundary
is smooth enough. It fails for $p=1,\infty$.
Arthur Korn was a German physicist born in 1870. He was also an inventor, involved in the
development of the fax machine and  Bildtelegraph which were early television systems, as well as a mathematician working
on partial differential equations. He had been dismissed from his post in 1935 and left Germany to the US, 
working at the Stevens Institute of Technology in Hoboken. For more on the inequality, see \cite{Ciarlet2010}.
\index{Korn inequality}
\index{Elasticity}
\index{Hilbert-Schmidt}
\index{Symmetrized Sobolev norm}
\index{Sobolev norm}

\section{Twin primes}

A pair $(p,q=p+2)$ of two {\bf rational primes} is called a {\bf prime twin}.
Examples are $(p,q)=(5,7)$. One might have wondered since antiquity about the infinitude of prime twins.
The {\bf twin prime conjecture} claims that infinitely many prime twins exist. The first
known source about the conjecture is Alphonse de Polignac in 1849 so that the conjecture
is sometimes also called the {\bf Polignac conjecture}. 
Let $\pi_2(x)$ denote the number of twin primes up to $x$. The {\bf sieve bound} has
first been established by Viggo Brun who showed $\pi_2(x) = O(x (\log \log(x)/\log(x))^2$.
Let $Li_2(x)= \int_2^x \frac{dt}{\log(t)}^2$ and $\mathcal{S}=2 \prod_{p {\rm prime} \; \geq 3} (1-2/p) (1-1/p)^{-2}$.
The sieve bounds theorem is

\satz{There is a constant $C$ with $\pi_2(x) \leq C \mathcal{S} Li_2(x)$.}

The constant $\mathcal{S}$ has a probabilistic background. It is 
$\mathcal{S} = \prod_{p \; {\rm prime}} \mathcal{S}_p$
where $\mathcal{S}_2=(1-1/2) (1-1/2)^{-2}=2$ and $\mathcal{S_p}=(1-2/p) (1-1/p)^{-2}$ for $p \geq 3$.
One expects then from a probabilistic point of view a prime twin density of $\mathcal{S} x/\log^2(x)$.
The sieve bound implies that the sum $\sum_{p,q \; {\rm prime}} \frac{1}{p} + \frac{1}{q} = 
(1/3+1/5) + (1/5 +1/7) + \cdots \sim 1.902$
of all reciprocals $1/p$ of all twin primes converges. The constant limit is called
the {\bf Brun's constant}. In the context of the {\bf twin prime conjecture}
there is {\bf Chen's theorem} telling that there are infinitely many 
primes $p$ such that $p+2$ has at most $2$ prime factors. 
{\bf Zhang's theorem} from 2014 about the existence of infinitely many 
bounded gaps has been pushed further: there are infinitely many pairs $(p,q)$ 
of distinct primes such that $|p-q| \leq 246$. 
See \cite{Maynard2019} for a recent review. 

\index{Sieve bound for prime twins}
\index{Brun's constant}
\index{prime twin}
\index{Chen's theorem}
\index{Zhang's theorem}

\section{Auction theory}

A real $n \times m$ {\bf signal matrix} $S=S_{ik}$ for $n$ {\bf buyers=bidders} 
and $m$ {\bf goods=merchandise} encodes {\bf real signal values} $S_{ik}$ which buyer $i$ 
can observe about the good $k$. Fixed also before hand is $T_{i}$, the set of $S$ matrix 
entries which buyer $i$ can see for good $k$. 
Buyers have {\bf private values} if they do not see what others do and
 {\bf common values} if they do see all what others can observe about $k$.
A {\bf valuation matrix} $V$ evaluates the signals relevant to the $i$'th buyers evaluation 
$V_{ik}$ for good $k$. It defines a {\bf welfare} of {\bf value system}
$V_i(S) = \sum_{j \in T_i} V_{ij}$. Given a {\bf payment} $P_i(S)$, the {\bf utility} is the difference
$U_i(S) = V_i(S)-P_i(S)$ of {\bf value minus payment}. A {\bf strategy} $\Sigma_i$
of buyer $i$ consists of defining $V(S)$ given the constraint $T$ of what they can see.
A {\bf pure strategy} is a deterministic choice of $V$, meaning that
buyers do not randomize. Given $P$ defined by the auction, its {\bf expected utility} is denoted by 
$U_i(\Sigma)$. A strategy $\Sigma^*$ is a {\bf Nash equilibrium} if all buyers optimize 
their own utility, meaning that $U_i(\Sigma^*) \geq U_i(\Sigma)$ for all $\Sigma$. 
An auction with a Nash equilibrium is called {\bf effective} if it is a Nash equilibrium 
for which $U$ is a {\bf global maximum} $U$. The problem is to find conditions and mechanisms which lead 
to Nash equilibria or even effective equilibria. The {\bf auction process} consists of
a {\bf bidding} that allows buyers to form a strategy $\Sigma$ to find the value $V$,
an {\bf allocation process} assigning goods to buyers according to $V$ and then define a {\bf payment} $P$ 
leading to the utility $U$. The goal is to find an auction process which leads to an effective Nash equilibrium. 
A {\bf Vickrey auction} is an auction process for private values and one good, 
a {\bf Vickrey-Clarke-Groves auction} (VCG) extends this to several goods.

\satz{There is a VCG bidding leading to an effective Nash equilibrium}

Auction theory is a chapter in game theory and is part of mathematical economics. 
It deals with the problem to use a bidding setup to allocate goods among buyers who bid for a fair prize.
It is a way to discover a correct price for a good.
Game theory started with von Neumann's paper of 1928. Von Neumann and Morgenstern 
\cite{NeumannMorgenstern} developed it in their book
in 1944. The concept of Nash equilibrium was introduced by John Nash in 1950 (see e.g. \cite{EssentialNash,Maskin2011}). 
In game theoretical settings, this means that players choose strategies from which unilateral deviations from the strategy do not pay better.
The Vickrey auction from 1961 in which "the highest bidder wins but pays the second highest bid", 
is a private auction where each person's bid only depends on its own value. The theory has shown to be so valuable
that Vickrey was awarded a Nobel prize in economics for his work. See \cite{KlempererAuctions,Milgram2004}.
\index{Auctions}
\index{Vickrey-Clarke-Groves auction}
\index{Nash equilibrium}
\index{bidding}
\index{valuation matrix}
\index{signal matrix}

\section{Wiener's 1/f theorem}

The {\bf Wiener algebra} $A(\mathbb{T})$ is the set of continuous $2\pi$-periodic functions $f$
with absolutely convergent {\bf Fourier series} $f(x) = \sum_{n \in \mathbb{Z}} c_n e^{i n x}$.
Equipped with the norm $||f||=\sum_{n \in \mathbb{Z}} |c_n|^2 < \infty$, it is a commutative Banach algebra,
meaning $||f \cdot g|| \leq ||f|| \cdot  ||g||$. It is not a $C^*$ algebra although, one would have to
change the norm to the supremum norm which then completes to the larger set $C(\mathbb{T})$.
The algebra consists of mildly regular continuous
functions because $C^{\alpha}(\mathbb{T}) \subset A(\mathbb{T}) \subset C(\mathbb{T})$
for all $\alpha>1/2$. The Fourier transform $f \in A(\mathbb{T}) \to \hat{f} \in l^1(\mathbb{Z})$ is
an isomorphism of Banach algebras. {\bf Wiener's 1/f theorem} is 

\satz{$f \in A(\mathbb{T})$ and $f(x) \neq 0, \forall x$ $\Rightarrow 1/f \in A(\mathbb{T})$. }

Wiener proved this in $1932$ (\cite{WienerTauberian}, Lemma IIe). 
In \cite{Katznelson} the theorem is called ``one of the nicest applications of
the theory of Banach algebras to {\bf harmonic analysis}". It is also known as the 
{\bf Wiener-L\'evy theorem} as Paul L\'evy extended the result showing that for any function $\phi$
that is analytic on the image of $f$, the function $\phi(f(x))$ is in $\A(\mathbb{T})$ \cite{Levy1935}.
L\'evy gives the example $f(x)=1/\log|\sin(x)/2|$ which has Fourier coefficients $c_n \sim 4/(n \log^2(n))$
and so has an absolutely convergent Fourier series. Its derivative $f'(x)=-f^2(x) \cot(x)$ is no 
more continuous. The 1/f theorem was proven by Israel Gelfand in 1939 using the structure theorem for 
commutative Banach algebras \cite{GelfandCollectedPapers,Zimmer1990}: 
it uses the fact (actually a lemma in Gelfand's
first paper on normed rings in 1939) that in order that an element in a normed ring has an inverse, it is 
necessary and sufficient that it is not in a maximal ideal: (if $f$ has an inverse then it does not
belong to any maximal ideal $I$ as then $1 \in I$ and so every $g \in I$; on the other hand if
$f$ does not have an inverse then $I=\{ gf, f \in A(\mathbb{T})$ \} is an ideal which is not the entire ring
and so is contained in a maximal ideal different from the ring). It also uses that the maximal ideals in $A(\mathbb{T})$ 
are the set of functions $f$ which vanish at some point $t_0 \in \mathbb{T}$. So, functions which do not vanish
are not contained in a maximal ideal and so are invertible. 
A short direct proof is given by Donald Newman in \cite{Newman1975}
using the inequality $|f|_{\infty} \leq ||f|| \leq |f|_{\infty} + 2 |f'|_{\infty}$ (which holds for
differentiable $f$ and in particular for finite partial Fourier sums): if $f \in A(\mathbb{T})$
is given which is nowhere $0$, scale it so that $|f(x)| \geq 1$ and take a partial sum $P$ such that $||P-f|| \leq 1/3$.
Now look at the geometric sum $S(x)=\sum_{n=1}^{\infty} (P(x)-f(x))^{n-1}/P^n$ which converges because $P(x) \geq 2/3$ and
$||1/P^n|| \leq (3 |P'|_{\infty}+1) (3/2)^n$. Because the geometric series converges to $S(x)=1/P (1/(1-(P-f)/P) = 1/f(x)$,
the theorem is proven. 

\index{harmonic analysis}
\index{Banach algebra}
\index{1/f theorem}
\index{Wiener 1/f theorem}
\index{Wiener algebra}

\section{Well ordering theorem}

A set $X$ is called {\bf well-ordered} if there is a {\bf total order} on $X$ such that
every non-empty subset $Y \subset X$ has a least element. [A total order is a binary relation $\leq$
that satisfies {\bf antisymmetry} ($a \leq b$ and $b \leq a$ implies $a=b$), {\bf reflexivity} $x \leq x$) 
{\bf transitivity} ($a \leq b$ and $b \leq c$ implies $a \leq c$) and {\bf connexity} ($a \leq b$ or $b \leq a$).
Without connexity, one only has a {\bf partial order}.
The {\bf least element} of a set $Y$ in a totally ordered set is an element $y \in Y$ such that $y \leq z$
for all $z \in Y$.] The {\bf well ordering theorem} is like {\bf Zorn's lemma} or {\bf Tychonov's theorem}
equivalent to the {\bf axiom of choice} and leads to seemingly paradoxa like the
{\bf Banach-Tarsky paradox} telling that one can partition the unit ball in $\mathbb{R}^3$ into 5 disjoint
sets such that three of them can be translated and rotated to become the unit ball again and the 2 remaining
can be translated and rotated to become the unit ball again which is a paradox because the doubling of the ball
is incompatible with volume.

\satz{Every set can be well ordered.}

The theorem was suggested by Georg Cantor in 1883 \cite{CantorUnendlichePunktmannigfaltigkeiten} 
and proven first by Ernst Zermelo in 1904 who 
called it the ``true fundament of the whole theory of number". The integers $\mathbb{Z}$ can be well 
ordered (we write $<<$ to distinguish from the $\leq$) for example with 
$x << y$ if $|x|<|y|$ or $|x|=|y|$ and $x=y$ or $x<0,y=-x>0$.
By using a bijection from $\mathbb{Q}$ to $\mathbb{Z}$, also $\mathbb{Q}$ can be well ordered as such. 
But this does not work for $\mathbb{R}$ any more. 
Koenig in 1904 at the Heidelberg Congress, gave on August 9th a wrong proof that the real numbers can 
not be well ordered. Already on August 10th, Zermelo pointed out an error in K\"onig's argument.
It was Felix Hausdorff who found an essential problem on September 1904 in a letter to Hilbert 
and also Cantor pointed to a problem. Hilbert, Hensel, Hausdorff and Schoenfliess had met in
Wengen (in the Swiss alps) at a successor congress. K\"onig then in October 1904 revoked his Heidelberg proof.            
On September 24 1904, Zermelo found the proof of the well-ordering theorem in M\"unden, near G\"ottingen and
acknowledges Erhard Schmidt for the idea, to base it on the axiom of choice. The letter
is printed in \cite{EbbinghausPeckhaus}, where it also pointed out that the proof 
was object of intensive criticism which only ebbed after decades.
While the fact that $\mathbb{R}$ can be well ordered is a consequence of the well ordering theorem,
it is impossible to explicitly construct such an ordering without assuming an 
{\bf axiom of constructibility}. See \cite{EbbinghausPeckhaus}  (section 2.5, 2.6) or 
\cite{Ebbinghaus2007}. The well ordering theorem is now known in the  mathematics history as one of the greatest mathematical controversies of all times. 
\index{well ordered}
\index{Well ordering theorem}
\index{axiom of choice}
\index{Banach-Tarsky}
\index{Zorn's lemma}
\index{Tychonov theorem}

\section{Caristi-Kirk-Ekeland theorem}

Let $(X,d)$ be a {\bf complete metric space} and $T: X \to X$ an arbitrary map,
not necessarily continuous. The map $T$ satisfies an {\bf inward condition} if there
exists a {\bf lower semi-continuous} function $f(x) \geq 0$ such that
$d(x,T(x)) \leq f(x) - f(T(x))$. [A function $f$ is {\bf lower semi-continuous} if
limits only can ``jump down" that is if $f(a) \leq \lim_{x \to a} f(x)$ for every $a$.
For example $f(x)=-1$ for $x \leq 0$ and $f(x)=1$ for $x>1$ is lower semi-continuous but
not continuous. $f$ is lower semi-continuous  if and only if $-f$ is {\bf upper semi continuous}.
If $f$ is both lower and upper continuous, then $f$ is {\bf continuous}.]

\satz{If $T$ satisfies an inward condition, then $T$ has a fixed point.}

An example is if we take $y \in X$ and $f(x)=d(x,y)$.
The condition then means $d(x,T(x)) \leq d(x,y) - d(T(x),y)$, implying
$d(T(x),y) < d(x,y)$ if $x \neq T(x))$, justifying the name ``inward condition".
In general, the condition means $f(T(x)) \leq f(x) - d(x,T(x)) < f(x)$ as long as $x \neq T(x)$.
The sequence $x_n = T^n(x)$ has the property that $y_n=f(x_n) \geq 0$ is decreasing, and so
some sort of Lyapunov function. By completeness, the sequence $y_n$ then 
must have an accumulation point $y$, which is a fixed point of $T$.
The theorem is easier to see if $f$ is continuous because there exists then by completeness
of $X$ an element $x$ with $f(x) = y$ and $f(T(x))=f(x)$ so that $f(T(x))=f(x)$.
James Caristi \cite{Caristi1976} (Theorem 2.1') mentions that the theorem was
suggested by Felix Browder and that I. Ekeland has proven an equivalent
theorem in 1972 (\cite{Ekeland1972} Theorem 1) as an abstraction of a lemma of Bishop and Phelps.
W.A. Kirk, who was the PhD advisor of Caristi, proved already a related theorem in 1965 
\cite{Kirk1965}: Kirk assumed that $X$ is a bounded convex subset of a reflexive Banach space
with a {\bf normal structure} [for every convex subset $H$ of $X$ with more than one point, there
is a point that is not a {\bf diametral point}. A diametral point in $H$ is a 
point $x$ which appears as in the supremum ${\rm sup}_{y \in H} ||x-y||$ being the diameter of $H$] 
and that $T$ does not increase distances. Caristi's statement is more general and 
elegant in comparison with the results of Ekeland and Kirk who were more concerned
with convex analysis \cite{EkelandTemanConvexAnalysis}.
\index{inward condition}
\index{lower semi-continuous}
\index{fixed point theorem}
\index{convex analysis}

\section{Shapley-Folkman theorem}

The {\bf Minkowski addition} of two subsets $A,B$ in $V = \mathbb{R}^d$
is defined as the set $\{a+b \; | \; a \in A, b \in B\}$.  A set $A$ in $V$
is called {\bf convex}, if for any two points $x,y \in A$ also the {\bf connecting
interval points} $\{ x+t (y-x), t \in [0,1] \}$ is part of $A$. The {\bf convex hull} of a set
$A$ is the smallest subset $c(A)$ in $B$ which is convex. A set $A$ is convex if
and only if the {\bf Minkowski distance} $d(A,c(A))$ of $A$ to its
{\bf convex hull} $c(A)$ is zero. One has in general the relation $c(A+B)=c(A)+c(B)$.
Let us call a sequence of sets $A_n$ {\bf uniformly bounded} if there is a ball
$B=B_r$ such that $A_n$ are all subsets of $B$. Define the {\bf Minkowski average}
$S_n = \frac{1}{n} \sum_{k=1}^n A_k$. The {\bf Shapley-Folkman theorem} is:

\satz{If $A_n$ are uniformly bounded then $d(S_n,c(S_n)) \to 0$. }

This is some sort of a {\bf law of large sets} in the sense that the Minkowski average
converges to the ``average" which is the convex hull.
There are uniform bounds for the distance which do not depend on $n$ as long as $n \geq d$.
For convex analysis, see \cite{EkelandTemanConvexAnalysis}.
For convexity in economics, see \cite{GreenHeller1981}.
\index{Minkowski addition}
\index{Minkowski sum}
\index{convex}
\index{Shapley-Folkman theorem}

\section{Dirichlet's unit theorem}

An element $r$ in a ring $R$ is called a {\bf unit} if it has an inverse.
The units form a group called the {\bf group of units}. In a division ring $R$,
it is the multiplicative group $R \setminus \{0\}$.
[ In a normed division algebra $R$ one sometimes calls the elements of norm $1$ ``units" in the sense that they are
elements of norm 1. Units here are all invertible elements in a ring. ]
A {\bf algebraic number field} is an algebraic field extension of the field of rational
numbers $\mathbb{Q}$. Let $O_K$ be the {\bf ring of integers} of the {\bf number field} $K$.
Its {\bf degree} is the dimension of $K$ as a vector field over $\mathbb{Q}$.
For {\bf quadratic fields} $\mathbb{Q}(\sqrt{d})$ for example, the degree is $2$ if the integer
$d$ is not a square integer. The field $K$ is the field of fractions of $O_K$.
{\bf Dirichlet's unit theorem} tells:

\satz{The group of units in a ring of integers is finitely generated.}

The {\bf rank} $r$ of a ring is the maximal number of multiplicative independent elements
in the group of units. It is $r=r_1+r_2-1$, where $r_1$ is the number of real embeddings (the number of real
conjugates of a primitive element) and $2r_2$ is the number of conjugates which are complex.
For a ring of integers in a {\bf real quadratic field} like $\mathbb{Z}[\sqrt{5}]$, the rank is $1$.
In an {\bf imaginary quadratic field} like the ring of Gaussian integers $\mathbb{Z}[\sqrt{-1}]$,
the rank is $0$. For all other fields, the rank is larger than $1$.
For algebraic number theory \cite{Lang1994,Neukirch1999,StewartTall}. A relatively short proof of the unity theorem 
can be found in \cite{Stein2012}. 
\index{Rank of a ring of integers}
\index{Dirichlet's unit theorem}
\index{Group of units}
\index{Ring of integers}
\index{number field}

\section{Spectrum of a countable theory}

{\bf Model theory} investigates how a formal theory build by sentences in a formal language is {\bf modeled}
and interpreted in concrete structures. A {\bf theory} $T$ a set of sentences in a language $L$.
A model $M$ of $T$ is a set with interpretations of functions, relations, symbols in that language $L$.
A model is {\bf complete} if every substructure of a model of $T$
which is a itself a model of $T$ can be axiomatized in first order logic.
The {\bf spectrum} of a {\bf complete theory} $T$ is the number $I(T,k)=I(k)$ of isomorphism          
classes of {\bf models} as a function of the {\bf cardinality} $k$. 
The {\bf L\"owenheim-Skolem theorem} tells that if $I(T,k)>0$ for some countable $k$, then $I(k)>0$ for 
all cardinalities $k$. First order logic theories can therefore not control the cardinality of their models. 
The theorem also shows that a theory with arbitrary large finite models must have an infinite model. 
A theory $T$ is called {\bf $k$-categorical} if $I(T,k)=1$ has only
one model up to isomorphism. L\"owenheim-Skolem shows that a first order theory with an infinite model 
is not $k$-categorical. This also follows from G\"odels incompleteness theorem.     
Michael Morley conjectured in 1961 that $I(T,k)=1$ for some uncountable     
$k$ then $I(T,k)=1$ for all uncountable $k$. In other words, if a theory is $k$-categorical for an 
uncountable power $k$, then it is $k$-categorical for every uncountable power $k$.
It was proven in 1965 by Morley \cite{Morley1965} and called the {\bf Morley's categoricity theorem}.

\satz{ $I(T,k)=1$ for $k$ uncountable $\Rightarrow$  $I(T,k)=1$ $\forall$ uncountable $k$. }

The theorem is remarkable in comparison with the L\"owenheim-Skolem theorem which       
tells that a theory in a countable language has an countably infinity model, then it has a model 
of any infinite cardinality. The categoricity theorem is considered the beginning of 
modern model theory. Michael Morley who died on October 16, 2020 had won the 2003 Steele prize 
for seminal contributions to research for his paper \cite{Morley1965} which had been initiated when writing this PhD thesis
in 1962. See \cite{Marker2002} chapter 6 and \cite{Shelah1990}. 
\index{Categoricity Theorem}
\index{Loewenheim-Skolem theorem}
\index{complete theory}

\section{Peano axioms}

The {\bf Dedekind-Peano axioms} (PA) formalize the arithmetic of the natural numbers $\mathbb{N}$.
The axiom system first lists five axioms that are already true in first order logic with
equality. The next three axiomatize the {\bf successor function} $S$: 1) for every $n \in \mathbb{N}$, there is a
successor $S(n)$. 2) $S$ is injective and 3) there is no $n$ with $S(n)=0$. And then there is the
{\bf axiom of induction:} 4) if $K$ is a set such that $0 \in K$ and $n \in K$ implies $S(n) \in K$, then $\mathbb{N} \subset K$.
Not all statements which are true for integers can be proven by the Peano axioms. Already
Kurt G\"odel established the existence of statements in PA that are true but unprovable within PA.
An accessible and natural example has been given by Jeff Paris and Harrington \cite{ParisHarrington}:
a finite set $H \subset \mathbb{N}$ is called {\bf relatively large} if ${\rm card}(H) \geq {\rm min}(H)$.
Given a finite set $M \subset \mathbb{N}$ and $e,r,k \in \mathbb{N}$ let $F(M,k,r,e)$ denote the statement: 
for every coloring map $P: M^e \to \{1, \dots, r\}$ (producing a partition of $M^e$), there is a relatively large 
$H \subset M$ with ${\rm card}(H) \geq k$ on which $P$ is constant (there is only one color on $H$).
The {\bf extended finite Ramsey theorem} is ``for all $e,r,k \in \mathbb{N}$, there exists $M$
such that $F(M,k,r,e)$ holds".

\satz{The extended finite Ramsey theorem is not provable in PA.  }

Paris and Harrington point out that when working with natural numbers, working in PA
amounts of replacing the axiom of infinity by its negation in ZF. They then give first a proof of the
extended finite Ramsey theorem as follows. (We write $\mathbb{N}$ for $\omega$, the order type of $\mathbb{N}$): 
fix $e,r,k$ and assume there is no such $M$. Let $P:M \to \{1, \dots ,r\}$
be the counter example map. There is no relatively large homogeneous set of size at least $k$.
The set of counter examples is a graph where $(P,M),(Q,N)$ are connected if $M \subset N$ and 
$P$ is the restriction of $Q$ to $M$. This is an infinite tree with finite vertex degree at every point.
By K\"onig's lemma, there is $P: \mathbb{N}^e \to \{1, \dots, r\}$ such that for every $M \subset \mathbb{N}$
the restriction of $P$ to $M^e$ is a counter example for $M$. By the infinite Ramsey theorem, there is an
infinite $H \subset \mathbb{N}$ that is homogeneous for $P$. By choosing $M$ large enough 
(compared to $k$, ${\rm min}(H)$) $H \cap M$ is a relatively large homogeneous set for $P|M^e$ 
of size at least $k$. This finishes the proof of the extended finite Ramsey theorem.
The proof of the Paris-Harrington theorem uses model theoretic techniques and the G\"odel's incompleteness
theorem. Paris and Harrington define a ``beefed up" theory $T$ and show that the consistency of $T$ implies the 
consistency of $PA$ using $PA$ only. Then they show that the extended finite Ramsey theorem implies the
consistency of $T$ and so the consistency of $PA$. This contradicts G\"odel's incompleteness theorem: one can 
not prove the consistency of $PA$ within $PA$.                           
Laurence Kirby and Jeff Paris have produced even more accessible examples \cite{KirbyParis1983}, especially 
the {\bf Hydra game}: which is a game in which the player has to cut off heads of a tree to which the tree
reacts by growing a multiple copies of branches. The theorem is that the player always wins. The surprise
is that one can not prove this within PA. More examples are in \cite{Spencer1983}. 
\index{Parris-Harrington theorem}
\index{Peano axioms}
\index{Inconsistency}
\index{Ramsey theory}
\index{Extended finite Ramsey theorem}
\index{Hydra game}

\section{Simplicial sets}

The {\bf simplex category} $\Delta$ has {\bf simplices} $[n]=\{ 0,1, \dots, n\}$ (non-empty totally ordered finite sets)
as {\bf objects} and {\bf order preserving maps} between them as {\bf morphisms}.
A {\bf simplicial set} is a contravariant functor $\Delta \to Set$. Simplicial sets form a category called $sSet$.
More generally, a {\bf simplicial object} is a contravariant functor from $\Delta$ to an other {\bf category} $\mathcal{C}$.
A {\bf coface map} $d^i: \{0, \dots, n\} \to \{0,\dots,n+1 \}$ is the unique order preserving bijection for which
element $i$ is omitted in the codomain. The {\bf codegeneracy map} $s^i: \{0,\dots,n+1) \to \{0,\dots ,n\}$
duplicates the element $i$ meaning that $s^i(j)=j$ if $0 \leq j \leq i$ and $s^i(j)=j-1$ if $i < j \leq n$.
There is now a {\bf decomposition lemma}:

\satz{$\Delta$ morphisms are a composition of coface and codegeneracy maps.}

This simple lemma is important to appreciate the axiomatic description of simplicial sets given first by May in 1967.
See \cite{Hatcher,RiehlLeisurly}.
It tells that every morphism: $f: \{0,\dots,n\} \to \{0,\dots,m\}$ has a unique
representation $f=d^{i_k} \cdots d^{i_1} s^{j_1} \cdots s^{j_h}$ with $n+k-h=m$ and
$m \geq i_k \geq \dots \geq i_1 \geq 0$ and $0 \leq j_1 \leq \dots \leq j_h < n$.
One has the relations $d^i d^j = d^{j+1} d^i$ and
$s^j s^i = s^i s^{j+1}$ for $i \leq j$ and
$s^i d^j = d^i s^{j-1}$ for $i<j$, $s^j d^i=1$ if $i=j,j+1$ and $s^i d^i = d^{i-1} s^j$.
In the {\bf opposite category} $\Delta^{op}$ (the category with the same objects but reversed morphisms),
the morphisms are denoted by $d_i,s_j$ and called {\bf face and degeneracy maps}. All morphisms in $\Delta^{op}$ are
now generated by composites of $d_i,s_j$. It follows that a contra-variant functor $X$ from $\Delta$ to $\mathcal{C}$
is determined by the images $X \{1,\dots,n\}$ of the simplices and if the face and degeneracy maps $d_i$ and $s_j$ are known.
A simplicial set therefore is a set of sets $X_n$ together with functions $d_i: X_n \to X_{n-1}$ and $s_i: X_n \to X_{n+1}$ satisfying
the composition relations for $d_i,s_i$. The elements $X_0$ are called the vertices, the elements $X_k$ are
called the $k$-simplices. The image of some $s_j$ is called a {\bf degenerate simplex}.
An advantage of looking at simplicial sets rather than the simplicial complexes is that one can use the frame work
in any category and that the Cartesian product of simplicial sets is a simplicial set.
The covariant geometric realization functor $X \to |X|$ from $\Delta$ to $Top$ is
the right adjoint to the {\bf singular homology} of the theory of simplicial sets. An other example is that the
nerve $N(\mathcal{C})$ of a small category $\mathcal{C}$ is a simplicial set constructed from the objects and
morphisms of $C$. A functor $f: \mathcal{C} \to \mathcal{D}$ between two categories induces then a map of
the corresponding simplicial sets and a natural induced transformation between two functors induces a homotopy between the
induced maps. The geometric realization of $N \mathcal{C}$ is called the classifying space of $\mathcal{C}$.
In general, for locally finite simplicial sets one has $|X \times Y| =|X| \times |Y|$ in $Top$
and the geometric realization $|X|$ of a simplicial set $X$ in Euclidean space is a CW complex.
In 1950, Eilenberg and Zilber introduced semisimplicial complexes (see \cite{EilenbergZilber}),
a terminology which later morphed into {\bf simplicial sets}.
According to \cite{Hatcher}, every simplicial complex can be subdivided to become a
simplicial complex so that every simplicial set is homeomorphic to a simplicial complex. But similarly as with CW complexes,
simplicial sets allow computations with fewer simplices. The category of topological spaces of
homotopy type of a CW complex are equivalent to the category of simplicial sets which satisfy an extension condition.
For more literature \cite{May1967,GoerssJardine,Kerodon2020}.
\index{simplicial set}
\index{coface map}
\index{codegeneracy map}
\index{simplicial category}

\section{Ostrowski theorem}

An {\bf absolute value} on $\mathbb{Q}$ is a norm
function from $x \in \mathbb{Q} \to |x| \in \mathbb{R}^+$ with the property that
$|x|=0$ is equivalent to $x=0$, and which is compatible with multiplication
$|x y| = |x| |y|$ and satisfying the {\bf triangle inequality}. The {\bf trivial norm}
is $|x|_1=1$ for $x \neq 0$ and $|0|_1=0$ is not considered. An example is the
usual absolute value $|x|_{\infty}$ or the {\bf $p$-adic absolute value}
$|x|_p=|p^n \frac{u}{v}|_p = p^{-n}$ if $p,u,v$ are all coprime numbers and where $p$ is a rational
prime. The {\bf p-adic norm} $|x|_p$ is a {\bf non-Archimedean} absolute value in
the sense that the {\bf ultra metric property} $|x+y| \leq {\rm max}(|x|,|y|)$ holds.
Two different absolute values are equivalent if $|x|_1 = |x|_2^{c}$ for some positive constant $c$. An
equivalence class of non-trivial absolute values on $\mathbb{Q}$ is a {\bf place}.

\satz{Every place $|x|$ is either $|x|=|x|_\infty$ or $|x|=|x|_p$ with prime $p$.}

Alexander Ostrowski proved this theorem in 1916 \cite{Ostrowski1916}.
It shows that every field containing $\mathbb{Q}$ which is complete with respect to an
{\bf Archimedean absolute value} is either $\mathbb{R}$ or $\mathbb{Q}$. An other curious consequence is
the {\bf product formula} $\prod_{p \leq \infty} |x|_p = 1$ for $x \in \mathbb{Q} \setminus \{0\}$
which combines all possible norms $|x|_p$ as well as $|x|_{\infty}$. The completion of $\mathbb{Q}$
with respect to $|x|_p$ is the space $\mathbb{Q}_p$ of {\bf $p$-adic numbers}. Each number in $\mathbb{Q}_p$ 
can be written in a unique way as {\bf half infinite Laurent series} $x = \sum_{k={-\infty}}^{\infty} a_k p^k$,
where the $a_k \in \{0,\dots,p-1\}$ are zero for $k < n(x)$. The norm is then $|x|_p=p^{-n(x)}$.
The field $\mathbb{Q}$ contains the sub-ring $\mathbb{Z}$ of {\bf p-adic integers}
$x = \sum_{k=0}^{\infty} a_k p^k$ which is the unit ball in $\mathbb{Q}_p$.
The ring $\mathbb{Z}$ of $p$-adic integers has no zero divisors so that
$\mathbb{Q}_p$ is the {\bf field of fractions} of $\mathbb{Z}_p$ (the smallest field containing $\mathbb{Z}_p$).
The $p$-adic integers $G=\mathbb{Z}_p$ with addition and metric coming from the norm
forms a commutative compact topological group with respect to addition. It is a totally disconnected
perfect metric space and so a {\bf Cantor space}, a topological space homeomorphic to the Cantor set.
Its {\bf Pontryagin dual group} $\hat{G}$ is the {\bf $p$-Pr\"ufer group} $\mathbb{Q}_p/\mathbb{Z}_p$,
the $p$-adic rational numbers modulo the integers which is
$\hat{G} = \{ e^{2\pi i k/p^l}, k,l \in \mathbb{Z} \}$.
While $\mathbb{R}$ and $\mathbb{Q}_p$ are all locally compact, 
the circle $\mathbb{T}=\mathbb{R}/\mathbb{T}$ is a compact, connected metric space which is the dual to
the integers $\mathbb{Z}$ which is non-compact and completely disconnected. The
$p$-adic integers $\mathbb{Z}_p$ are all compact, completely disconnected spaces, dual to
the {\bf Pr\"ufer p-group} $\mathbb{P}_p=\mathbb{Q}_p/\mathbb{Z}_p$.
In both cases, one has Haar measures $\mu$ on $G$ and a Fourier
isomorphism $f \in L^2(G,\mu) \to L^2(\hat{G},\hat{\mu})$ which is the Pontryagin involution.
All group translations and multiplications by some integer preserve the Haar measure $\mu$. While for $\mathbb{T}$,
the translations $x \to x+\alpha$ preserving $dx$ can be arbitrarily close to the identity, there is a smallest group translation
$x \to T(x)=x+1$ on the $p$-adic integers. It is called the {\bf adding machine} and is ergodic. The eigenvalues
of the unitary Koopman operator $f \to f(T)$ on $L^2(\mathbb{Z}_p)$ coincides as a set with the Pr\"ufer group
$\mathbb{P}_p \subset \mathbb{T}= \{ z \in \mathbb{Z}, |z|=1 \}$.
Besides group translations, there are also Bernoulli shifts in both cases which preserve the Haar measure.
On the compact topological group $\mathbb{T}=\mathbb{R}/\mathbb{Z}$,
the map $x \to n x$ is a Bernoulli shift for every $n >1$ with entropy $\log(n)$.
On the compact topological group $\mathbb{Z}_p$, the map $x \to p x$ is a Bernoulli shift with 
Kolmogorov-Sinai entropy $\log(p)$.
About the life of Ostrowski see \cite{GautschiOstrowski}. For $p$-adic analysis \cite{KatokPAdic2007,Gouvea1997}.
\index{place}
\index{p-adic numbers}
\index{p-adic integers}
\index{Archimedean}
\index{non-Archimedean}
\index{p-adic norm}
\index{ultra metric}
\index{Pruefer group}
\index{Pontryagin dual}
\index{Fourier transform}

\section{Clifford algebras}

If $V$ is a finite dimensional {\bf vector space} over a field $k$ and $q$ is a {\bf quadratic form}
of {\bf signature} $(p,q)$, its {\bf Clifford algebra}
$Cl(V,q)$ is the quotient $T(V)/I$ of the {\bf tensor algebra} $T(V) = \bigoplus_{n=0}^{\infty} \bigotimes_{k=1}^n V$ by the ideal
$I$ generated by elements $v \otimes v - q(v) 1$.
[We use the sign convention of \cite{Chevalley1995, Snygg1997}. Other authors, 
like \cite{Garling2011,LawsonMichelsohn},
prefer to take the ideal $v \otimes v + q(v) 1$.] 
The Clifford algebra is a unital associative algebra. If the underlying field $k$ is $\mathbb{R}$, it is called
a {\bf geometric algebra}.
For $q=0$, one obtains the {\bf exterior product} with exterior multiplication, where $v \otimes w=-w \otimes v$.
As turning on $q$ deforms the anti-commutativity relation, one sees $Cl(V,q)$ as a {\bf quantization} of
the exterior algebra $Ext(V)$. (One can also see the process of going from $(V,q)$ to $Cl(V,q)$ as a second quantization
if one interprets the tensor algebra as a many-body {\bf Fock space}.)
Examples of Clifford algebras are the {\bf complex numbers} $Cl_{0,1}(\mathbb{R}) = \mathbb{C}$
the {\bf quaternions} $Cl_{0,2}(\mathbb{R})$ and {\bf split complex numbers} $Cl_{1,0}(\mathbb{R})$
or {\bf split quaternions} $CL_{2,0}(\mathbb{R})$. Notable in relativity is the {\bf space time algebra}
$Cl_{1,3}(\mathbb{R})$. Let $i: V \to Cl(V,q)$ be the {\bf inclusion map} which embeds $V$ into $Cl(V,q)$ and which
satisfies $i(v)^2 = q(v)$. The algebra $CL(V,q)$ enjoys now a {\bf universal property}
if given any associative algebra $A$ and any linear map $j:V \to A$
obeying $j(v)^2=q(v)$: there exists a unique algebra homomorphisms $f: Cl(V,q) \to A$
such that $j(v) = f(i(v))$. The {\bf fundamental theorem of Clifford algebras} is that
$CL(V,q)$ is unique. One can speak of {\bf ``the" Clifford algebra} $Cl(V,q)$ defined by $V$ and $q$.

\satz{The Clifford algebra construct satisfies the universal property.}

In category theory one sees $Cl$ as a {\bf functor} 
from the category of finite dimensional pseudo Hilbert spaces $(V,q)$
to the category of unital associative algebras. 
The universal property generalizes the process of
getting from an algebra to the {\bf free algebra} $F(A)$ 
or to get the tensor algebra $T(M)$ from a module $M$ over a ring
\cite{Chevalley1995}. If $V$ has dimension $n$, 
the dimension of the Clifford algebra $Cl(V,q)$
is $2^n$. As a vector space, the Clifford algebra is isomorphic to the exterior algebra $Ext(V)$ which is like $Cl(V,q)$
a {\bf super algebra}. This follows from the universal property. As the involution $v \to -v$ does not change
the quadratic form $q$, it lifts to an involution $\alpha$ on $Cl(V,q)$. This produces a splitting
$Cl(V,q) = Cl(V,q)^{{\rm even}} \oplus Cl(V,q)^{{\rm odd}}$, where $Cl(V,q)^{{\rm even}} = \{ x \in Cl(V,q),  \alpha(x)=x \}$
and $CL(V,q)^{{\rm odd}} = \{ x \in CL(V,q),  \alpha(x) = -x \}$. Multiplication honors this grading.
The quadratic form $q$ can be extended from $V$ to $Cl(V,q)$:
first define the {\bf transpose} $x^T$ which reverses the order $x=v_1 \otimes \cdots \otimes v_k \to v_k \otimes \cdots \otimes v_1$,
and the {\bf scalar part} $x_0$ of an element $x = \sum_n \bigotimes_{|k|=n} x_k  v_{k_1} \otimes \cdots \otimes v_{k_n}$.
The symmetric, bilinear form $q_1$ on $V_1=Cl(V,q)$ is then defined as $x \cdot y = (x^T \cdot y)_0$ and continues to be non-degenerate
if $q$ was. In the case when $q$ is positive definite, where $(V,q)$ is a Hilbert space, the operation $(V,q) \to (V_1,q_1)$ can
now be iterated and produces a sequence $(V_n,q_n)$ or Hilbert spaces.
Clifford algebras have many applications, like in algebraic geometry (starting with Grassmann who introduced exterior algebras),
in representation theory of classical Lie groups (it was \'Elie Cartan who discovered in 1913 first
unknown representations of the orthogonal group and called the elements on which the matrices
operated ``spinors" \cite{Cartan1981}), in physics or in differential geometry.
To the later, Cartier writes in the introduction to
\cite{Chevalley1995} that {\it since the 1950's, spinors and the associated Dirac equation have developed 
into a fundamental tool in differential geometry}. Indeed, on has at every
point $x \in M$ of a Riemannian manifold a Clifford algebra $Cl(T_xM,g(x))$ defined by $g(x)$, the quadratic form in
the tangent space $V=T_xM$. This produces a {\bf Clifford bundle}. One can then ask, under which conditions a {\bf spin
structure} exists on $M$. It is the case if and only if the second Stiefel-Whitney class $w_2(M)$ is zero. This
topological obstruction for the existence of {\bf spin structures} on an orientable Riemannian manifold $(M,g)$
was found by Andr\'e Haefliger in 1956. Haefliger defined the spin structure on $(M,g)$ as
a lift of the principal orthonormal {\bf frame bundle} $F_{SO}(M) \to M$ to $F_{Spin}(M)$. Not every Riemannian manifold
is {\bf spin}. While spheres $\mathbb{S}^n$ are spin, the $2n$-manifolds $\mathbb{CP}^{2n}$ (complex projective spaces)
are not spin. The space of {\bf spinors} of $(V,q)$ is the fundamental representation of a Clifford algebra $Cl(V,q)$.
Spinors belong also to vectors in a representation of the {\bf double cover} Lie algebra $Spin(p,q)$ of the special orthogonal
group $SO(p,q)$ of signature $p,q$. Representation theory of classical groups like $SO(n)$ or $Spin(n)$, a subgroup
of the group of invertible elements in a Clifford algebra of a Hilbert space, are a major motivator for Clifford algebras.
\index{Clifford algebra}
\index{Universal property}
\index{quadratic form}
\index{Spinors}
\index{spin structure}
\index{quaternion}
\index{quantization}

\section{Transcendental number theory}

A complex number is called {\bf algebraic} if it is the root of a polynomial $a_0 + a_1 x + \cdots + a_n x^n \in \mathbb{Z}[x]$
with integer coefficients $a_0,\dots,a_n \in \mathbb{Z}$. 
The algebraic numbers $\mathbb{A}$ form a {\bf field}. They can be enumerated and so are a countable set in $\mathbb{C}$.   
As a consequences, almost all real numbers are not algebraic. This argument of Cantor is a non-constructive but elegant 
proof of the existence of non-algebraic numbers, numbers in the complement $\mathbb{C} \setminus \mathbb{A}$ which are
also called {\bf transcendental numbers}. Let us call a {\bf Gelfond-Schneider pair} a pair of algebraic numbers
$\alpha,\beta$ for which $\alpha \neq 0,1$ and for which $\beta$ is rational. The Gelfond-Schneider theorem is:       

\satz{$(\alpha,\beta)$ Gelfond-Schneider $\Rightarrow$ any choice of $\alpha^{\beta}$ is transcendental.}

The theorem is named after Alexander Osipovich Gelfond and Theodor Schneider.
Gelfond proved a special case in 1929 ($\beta$ imaginary quadratic) 
and the full version in 1934. Schneider proved the same in 
his PhD thesis in 1934 written under the advise of Carl Siegel, who already proved 
it for real quadratic $\beta$. An example is
the {\bf Gelfond-Schneider constant} $2^{\sqrt{2}}$. An other example is
{\bf Gelfond constant} $e^{\pi} = (e^{i\pi})^{-i} = (-1)^{-i}$.   
One has to say ``any choice" because $(-1)^{-i}$ invokes the complex
logarithm which has many branches. Any other branch like $(-1)^{-i} = e^{-i \log(-1)} = 
e^{-\pi + 2k \pi} = e^{-(1+2k)\pi}$ is also transcendental. A third example is 
the ``eye for an eye" number $i^i = (e^{i \pi/2})^i = e^{-\pi/2}$ which is already transcendental as a consequence
of the Gelfond theorem because $i$ is algebraic, solving the equation $1+x^2=0$. 
The problem whether $\alpha^{\beta}$ is transcendental for a Gelfond-Schneider pair had been asked
by David Hilbert and got to be known as {\bf Hilbert's seventh problem} in 1900. 
Questions about transcendental numbers are difficult. For example, one still does not know, whether 
$\pi^e$ is transcendental or not.  See \cite{Gelfond1960} and especially chapter 4 and 5 of \cite{BurgerTubbs}.
\index{Gelfond theorem}
\index{algebraic number}
\index{transcendental number}
\index{Gelfond constant}
\index{Gelfond-Schneider constant} 

\section{Conformal maps}

Let $D_r(w) = \{ |z-w|<r \}$ denote a disk of radius $r$ centered at $z \in \mathbb{C}$.
If $f: \mathbb{D}_r(w) \to \mathbb{C}$ is a holomorphic function with $f'(w) \neq 0$,
then by the implicit function theorem, the map is invertible and by Bloch's theorem,
$f(D_r(w))$ contains a disk $D(f(w),c |f'(w)| r)$ for some constant $c$.
The best constant $c$ for which this works is is called the {\bf Bloch constant}.
An analytic, injective function $f$ on $D_r(w)$ is also called {\bf univalent} and
$f: D_r(w) \to f(D_r(w))$ is called a {\bf conformal mapping}.
{\bf Koebe's quarter theorem} is

\satz{If $f$ is univalent on $D_r(w)$, then $D_{|f'(w)|r/4}(w) \subset f(D_r(w))$. }

The result had been conjectured in 1907 by Paul Koebe and was first proven by Ludwig Bieberbach
in 1916 \cite{Bieberbach1916}. The {\bf Koebe function} $f(z) = z/(1-z)^2 = \sum_{n=1}^{\infty} n z^n$
shows that the constant $1/4$ can not be improved upon.
In \cite{Tao11282020}, the result is stated and proven that for polynomials of degree $n$,
the image $f(D_r(w))$ contains the disk $D_{|f'(w)| r/n}(f(w))$. (The blog cites \cite{Miller1993} but
this the disc result is not that obvious there).
For more information on Koebe, see \cite{Carlson}.
Paul Koebe's is famous also for a his theorem generalizing
the {\bf Riemann mapping theorem}. It states that
any finitely connected domain is conformally equivalent to a circle domain unique up to
M\"obius transformations. (A {\bf circle domain} is an open subset of $\mathbb{C}$ such that every
connected components of its boundary is either a circle or a point.) Koebe's Kreisnormierungsproblem
from 1909 asks whether every domain in $\mathbb{C}$ is conformally equivalent to a circle domain
unique up to a M\"obius transformation. The problem is open.
\index{Koebe function}
\index{Koebe one quarter theorem}
\index{Bloch constant}
\index{holomorphic function}
\index{Moebius transformation}
\index{Kreisnormierungsproblem}

\section{Shannon capacity}

The {\bf independence number} $\alpha(G)$ of a finite simple graph $G=(V,E)$
is the maximum number of independent points in $G$ (a set of vertices is 
{\bf independent} if the members of the set
are pairwise not adjacent). The
{\bf strong product} $G*H$ of two graphs $G,H$ has as the vertex set the Cartesian
product $V(G) \times V(H)$ of vertices in $G$ and $H$ and as edges all connections which when projected on
any of the graphs gives either a vertex or edge. In communication theory, where $V$ is
the alphabet and $E$ gives letters which can be confused, then
$\alpha(G^k)$ the maximal number of $k$ letter messages which can be sent without
the danger of confusion. The limit $\Theta(G) = \lim_{k \to \infty} \alpha(G^k)^{1/k}$ is called
the {\bf Shannon capacity} of the graph. One has clearly $\Theta(G) \geq \alpha(G)$ because there
are at least $\alpha(G)^k$ words which can not be confused. The extreme cases is $P_n$, the graph
with $n$ vertices and no edges and $K_n$, the graph with $n$ vertices and all edges present. In
these cases $\Theta(P_n)=n$ and $\Theta(K_n)=1$.

\satz{The Shannon capacity of $G=C_5$ is $\sqrt{5}$.}

The Shannon capacity was introduced by Claude Shannon in 1956 \cite{Shannon1956} who wrote:
{\it The zero error capacity of a noisy channel is defined as the least upper bound of rates
at which it is possible to transmit information with zero probability of error.}
Shannon took the logarithm and called
$C_0 = \log(\Theta(G)) = \lim_{k \to \infty} \frac{1}{k} \log(\alpha(G^k))$
the {\bf zero-error capacity} which reminds of a Lyapunov exponent measuring the exponential 
growth of a cocycle. Shannon computed the capacity for all graphs with $n=1,2,3,4,5$ nodes and the 
pentagon had been the smallest, where he had been
unable to determine the value, he only established $\sqrt{5} = 2.236.0 \leq \Theta(G) \leq 5/2=2.5$.
This true value for the pentagon was then computed in \cite{Lovasz1979}, where also the notation
$\Theta(G)$ appears. An exposition about Shannon capacity appears in \cite{Matousek} (Miniature 28 and 29).
The problem of computing $\Theta(G)$ is formidable. One does not even know $\Theta(C_7)$. 
\index{Shannon capacity}
\index{Shannon zero error capacity}
\index{communication theory}
\index{Strong product}
\index{independent set}

\section{Outer billiards}

A convex curve $C$ in $\mathbb{R}^2$ defines an area-preserving map
$T: X \to X$, where $X$ is the unbounded region outside of the table. A point $(x,y)$
is mapped into $T(x,y)$ which is the point reflection at the point
$(p,q)$ which is the midpoint of the interval $I$ obtained by intersecting
the counter clockwise tangent from $(x,y)$ to $C$. The map can be extended to $C$ by
defining $T(x,y)=(x,y)$ there. For most points $(x,y)$, the
interval $I$ is a single point but already for polygons, we want to have $T$
defined everywhere, even so it is not continuous. The map $T$ is called the
{\bf outer billiard map} or {\bf dual billiard map} defined by $C$. The {\bf Penrose polygon}
is the {\bf quadrilateral} $ABPQ$ defined by the 5-gon $A,B,C,D,E$, where
$P=(AD) \cap (BE)$ and $Q=(BD) \cap (CD)$, with $(AD),(BE),(BD),(CD)$ denoting diagonal segments
in the pentagon. A table is called {\bf unstable} if there exists $(x,y)$ such that $|T^n(x,y)|$ is 
unbounded. 

\satz{Outer billiard at the Penrose kite is unstable.}

The outer billiard $T$ is smooth if $C$ is a smooth and strictly convex. 
The dynamical system had been introduced by B.H. Neumann in the Manchester University Mathematics
students journal of 1959 \cite{Neumann1959} and was popularized in 
\cite{MoserStableRandom,Mos78}. The question whether there exists a convex table for which
an unbounded orbit of the map $T$ exists is known as the {\bf Moser-Neumann question}
\cite{Schwartz2007,Schwartz2009}. If $C$ is smooth and strictly convex, then KAM theory
establishes invariant curves for $T$ and so stability of the table \cite{Dou82}.
For a class of tables called quasi-rational polygons, which includes rational polygons and 
regular $n$-gons, all orbits are bounded \cite{ViSh87,Kolodziej,GutkinSimanyi}. Also trapezoids
lead to bounded orbits \cite{Li2009}.

\index{Outer billiards}
\index{Dual billiards}
\index{Penrose polygon}
\index{Quasi-rational polygons}
\index{Moser-Neumann problem}

\section{Sandwich theorem}

Let $G=(V,E)$ denote a finite simple graph. In information theory, where
$V$ is an alphabet of symbols, the graph is the {\bf confusion graph} where
connecting symbols which can be confused. A function $f: V \to S^(n+1)$
is called an {\bf orthonormal representation} if the orthogonality condition
$\langle f(u),f(v) \rangle=0$ holds if the vertices are not adjacent.
The {\bf Lovaz number} is defined as
$$  \theta(G) = \min_{c,U} {\rm max}_{v \in V} \langle c,U(v) \rangle^{-2}  \; , $$
where $c$ is a unit vector and $U$ is an orthonormal representation.
This corresponds to minimizing the half-angle $\alpha$ of a rotational cone
as $\theta(G) = 1/\cos^2(\alpha)$, where $c$ is the symmetry axes of
the cone. The {\bf Lovasz number} is multiplicative in the graph product because one
can build for every power $G^n$ also an umbrella $U^n$.
Let $\alpha(G)$ be the {\bf independence number} of $G$. It is the clique number $c(G)$
of the graph complement $\overline{G}$. Let $\chi(G)$ denote the {\bf chromatic
number}, the minimal number of colors which one can use to color the graph. It is
the {\bf clique covering number} $\beta(\overline{G})$ of the graph complement. The following
sandwich identity is the key to estimate the {\bf Shannon capacity}
$\Theta(G) = \lim_{n \to \infty} \alpha(G^n)^{1/n}$ \cite{Shannon1956}.

\satz{$c(\overline{G^n}) = \alpha(G^n) \leq \Theta(G) \leq \theta(G) \leq \beta(G)=\chi(\overline{G})$.}

The Lovasz number $\theta(G)$ can be computed in polynomial time in the number of vertices.
The Shannon capacity is sandwiched between the independence number of any power of $G$
and the Lovasz number. See \cite{KnuthSandwich}. An example is
$\alpha(C_5^2)=5$ where $(1,1),(2,3),(3,5),(5,4),(4,2)$ is an independent set in the Shannon
product $G^2$, we have $\Theta(C_5) \geq \sqrt{5}$. The Lovasz umbrella
$U=\{ u_1,u_2,u_3,u_4,u_5 \}$ with
$u_k = [\cos(t) \sin(s), \sin(t) \sin(s), \cos(s) ]$ with
$\cos(s) = 1/5^(1/4), t=2 \pi k/5$ gives $\theta(C_5) \leq \sqrt{5}$.
Therefore, the Shannon capacity of the pentagon is $\sqrt{5}$.
One does not know the Shannon capacity of the heptagon. In \cite{Lovasz1979}, where also the notation
$\Theta(G)$ appears. An exposition about Shannon capacity appears in \cite{Matousek}.
\index{Sandwich theorem}
\index{Lovasz umbrella}
\index{Shannon capacity}
\index{independence number}
\index{chromatic number}
\index{clique covering number}
\index{clique number}
\index{confusion graph}


\section{Shannon capacity theorem}

For a communication with bandwidth $B$, the {\bf signal to noise ratio} $S/N$ (also abbreviated
SNR) has maximal capacity $C$. These quantities are related by 

\satz{$C= B \log_2(1+S/N)$.}

This is also called the {\bf Shannon capacity theorem}.
The units are $C$ as bits per second. The bandwidth is in given Herz,
$S$ is the average received signal power measured in Watts and $N$ is
the average power of the noise measured in Watts. The number $\log_2(1+S/N)$
is the spectral efficiency.

In 1993, {\bf turbo codes} appeared \cite{TurboCodes}. 
These were first practical codes to get to the Shannon limit.
These codes were already patented by
Claude Berrou in 1991. These codes are used in modern 3G, 4G mobile telephony
standards. In 5G wireless communication other codes like {\bf Polar codes} 
are used, which reach Shannon channel capacity \cite{PatilPawarSaquib}.
\index{Capacity theorem}
\index{signal to noise ratio}
\index{wireless communication}

\section{Differential Galois theory}

A {\bf differential ring} is a field $R$ equipped with a {\bf derivation} $D: R \to R$
which is linear and satisfies the {\bf Leibniz rule} $D(fg)=D(f) g + f D(g)$.
The {\bf field of fractions} of an integral domain $R$ 
(a ring $R$ for which the product of two non-zero
elements is not zero) is the smallest field containing $R$.
A {\bf differential ring extension} $R<S$ has the ring $R$ as a sub-ring and the derivation of $S$
on restricted to $R$ agreeing with the derivation on $R$.
A {\bf differential ideal} is an ideal $I \subset R$ that is invariant under $D$. If defines
the quotient ring $R/I$ with derivation $D(a+I) = D(a) + I$.
The {\bf ring of differential polynomials} over $R$ is the polynomial ring $R[ Y_1,Y_2, \dots ]$
with a countable set of variables in which $D$ is extended as $D Y_i = Y_{i+1}$.
If $F$ is a differential field and $K$ a field extension, then $t \in K$ is called {\bf elementary}
if it is generated by algebraic, a logarithm or an exponential functions.

\satz{$f =e^{x^2}$ can not be integrated in elementary terms.}

After differentiation, there would have to exist a function $f$ with $1=f'+2fx$.
\cite{Magid1999}. For a book, \cite{MagidDiffGalois} or the lectures \cite{Singer2006}.

\index{Leibniz rule}
\index{Derivation}
\index{differential ring}
\index{elementary function}
\index{differential Galois theory}

\section{Non-linear Schroedinger equation}

The {\bf non-linear Schr\"odinger equation} (NLSE)
$i u_t = -\Delta u + |u|^p u$ is an example of a nonlinear
{\bf partial differential equation} for $u(t,x)$ with $x \in \mathbb{R}^d$.
It is an example of a classical field equation which can be used to
describe {\bf Langmuir waves} in hot plasmas or wave 
propagation in fiber optics in which
the non-linearity comes from self-phase modulation. It also appears to have relevance in
understanding the formation of {\bf rogue waves} in the ocean.
The later are unexpectedly large waves that can endanger ships. 
In dimension $d=1$, the differential equation is an example of an
{\bf integrable system} featuring non-linear phenomena like {\bf solitons}.
The $L^2$-norm square of $u$ is called the
{\bf mass} $|u|_2$ of $u$. It is preserved under the evolution.
For any $\lambda$, the function $u_\lambda(t,x) = \lambda^{2/p} u(\lambda^2t,\lambda x)$ is also a solution
and its mass is $M(u_{\lambda}) = \lambda^{-d+4/p} M(u)$. The mass {\bf subcritical case} is $p<4/d$.
The {\bf Sobolev space} $H^{k}(\mathbb{R}^d)$ is the space of functions $f$
such that $f$ as well as all its weak derivatives up to order $k$ have finite $L^2$ norm.

\satz{Global solutions of a subcritical NLSE exist in $H^1(\mathbb{R}^d)$. }

The problem is ill-posed in $H^{2/d-2/p}(\mathbb{R}^d)$. The {\bf mass-critical case} is when $p=4/d$.
In that case there is a minimal mass $m_0$ for solutions to blow up.  See \cite{TaoVisanZhang}.

\index{Integrable system}
\index{nonlinear partial differential equation}
\index{nonlinear Schroedinger equation}
\index{soliton}
\index{plasma}
\index{rogue waves}
\index{mass subcritical}
\index{mass critical}
\index{blow up of solutions}

\section{Menger's theorem}

Let $G=(V,E)$ be a finite simple graph. For two disjoint subsets $A,B$, a
{\bf minimal AB separator} is the minimal number of
vertices disjoint from $A,B$ which when removed disconnects $A$ from $B$. A {\bf maximal AB-connector}
is the maximal number of pairwise disjoint paths connecting $A$ with $B$. Let us
denote by $|Minimal AB-separators|$ the number of minimal AB-separators and similarly for 
the maximal AB-separators. The result is:

\satz{|Minimal AB-separators| = |maximal AB-connectors|.}

If $A$ and $B$ have an intersection, both numbers are just the cardinality of $A \cap B$ as
zero length paths $\{x \} \subset A \cap B$ are considered connectors.
An other special case is $G$ is $2$-connected with cut $\{x\}$ and where
$A \cup \{x\} \cup B$ is a disjoint union. Now, $\{x\}$ is a minimal $AB$-separator.
Since every path from $A$ to $B$ crosses $x$, a maximal $AB$-connector consists of only
one path. More generally, if $G$ is {\bf $k$-connected} meaning that we need to remove a vertex
set $X$ of cardinality $k$ to make it disconnected, then if $V=A \cup X \cup B$ is a disjoint
union, the set $X$ is a minimal $AB$ separator and a maximal $AB$ connector consist of $|X|$
paths. The proof is done with respect to the number of edges in $G$. 
Menger proved this theorem in 1927 \cite{Menger1927}. Menger did not use the language of graphs but 
proved it for {\bf curves} which are compact connected topological spaces 
for which the boundary of arbitrary small neighborhoods is disconnected. 
He considered them as one-dimensional continua. Menger's research was 
part of a program about dimension which works for general topological spaces independent
of metric. The graph theoretical version is a special case as a
geometric realization of the one-dimensional skeleton complex $V \cup E$ of a graph $G=(V,E)$ 
defines a curve in Menger's sense. 
\index{Menger's theorem}
\index{minimal separator}
\index{maximal connector}
\index{k-connected}

\section{Ap\'ery's theorem}

The {\bf Ap\'ery constant} $\zeta(3) = \sum_{n=1}^{\infty} \frac{1}{n^3}$
is a special value of the {\bf Riemann zeta function}.
R. Ap\'ery proved in 1979 that

\satz{
The Ap\'ery constant is irrational.
}

While one knows that all $\zeta(2n)$ are irrational for $n \geq 1$ starting with
$\zeta(2) = \pi^2/6, \zeta(4) = \pi^4/90$, the odd numbers are not yet known for $2n+1>3$.
One does not know for example whether $\zeta(5)$ is irrational. 
The problem of whether the Ap\'ery constant is irrational is in \cite{Nahin2021} the
``most mysterious unsolved math problem".  One only knows that infinitely many
of the numbers $\zeta(2n+1)$ are irrational. 
To the history: Euler, who gained fame with the computation of $\zeta(2) = \pi^2/6$ already
computed $\zeta(3)$ to several digits. The entire book \cite{Nahin2021} is dedicated
to Zeta-3. 
\index{Ap\'ery constant}
\index{Zeta function}
\index{Ap\'ery's theorem}
\index{Zeta-3}

\section{Vietoris theorem}

\paragraph{}
A topological space $(X,\mathcal{O})$ is {\bf compact} if every open cover 
(a subset $\mathcal{F}$ of $\mathcal{O}$ whose union is $X$) has a finite sub-cover
(a finite subset of $\mathcal{F}$ whose union is $X$). If $A \in \mathcal{O}$ then 
$X \setminus A$ is called {\bf closed}. 
A topological space is {\bf normal} if any two {\bf closed sets} $A,B$ in $X$
have disjoint open neighborhoods $U,V$.  Non-normal topological spaces are relevant
in mathematics: the Zariski topology on the spectrum of a ring for example is non-normal. 
The {\bf Vietoris theorem} is

\satz{
A compact topological space is normal.
}

Leopold Vietoris proved this in 1921 and considered it his most important result
even so he lived 110 years and wrote his last paper on trigonometric sums with 103
and more than half of his papers were written after his sixties birthday
\cite{Reitberger2002}. Normality is also called {\bf Tietze's normality
condition}. In modern topology books the normality condition is called Axiom $T_4$
but one has to be careful, as sometimes, it also assumes Hausdorff (any two points in $X$
can be separated by open neighborhoods). Normality $T_4$ does not imply Hausdorff $T_2$: an example is
the topological space $X=(\{a,b\}, \mathcal{O}=\{ \emptyset,X \})$ which has only $\emptyset,X$ 
as closed sets and both sets are both open and closed. They
also have disjoint open neighborhoods as they themselves are open neighborhoods (
they are both clopen sets). Indeed, any {\bf indiscrete topological space} is normal. 
But if $X$ contains at least two points, it is not Hausdorff. There are two points $a,b$ 
that can not be separated by open neighborhoods. 
Any indiscrete topological space with at least two points 
is also an example of a {\bf compact non-Hausdorff space}. 
The Theorem of Vietoris assures that the seemingly stronger condition
of normality holds for all compact topological spaces while the Hausdorff property does not
always hold. Vietoris is the father of modern convergence concepts like {\bf filter base} or 
{\bf nets} and modern notions of compactness.
Normality is important because of the {\bf Tietze's extension theorem} stating that a
continuous function on closed subset of a normal topological space can be extended to the
entire space. The Tietze theorem was proven by Brouwer and Lebesgue for Euclidean spaces,
extended by Tietze to metric spaces and by Urysohn for normal space.
\index{Vietoris theorem}
\index{Thietze extension theorem}
\index{Normal topological space}
\index{Indiscrete topological space}
\index{Compact non-Hausdorff space}

\section{Whitney extension theorem}

\paragraph{}
Given a $C^m$ function on $\mathbb{R}^n$, {\bf Taylor's theorem} assures that
$f(x) = \sum_{|k| \leq m} f^{(k)}(y)/k! (x-y)^k + \sum_{|k|=m} R_k(x,y) (x-y)^m/m!$ with $R_k(x,y) \to 0$
uniformly as $x,y \to a$. 
This gives relations $f^{(r)}(x) = \sum_{|k| \leq m-|r|} f^{(k+r)}(y) (x-y)^k/k! + R_r(x,y)$. 
A set $F=\{ f_k\}$ of functions on $\mathbb{R}^n$ 
with multi-index $|k| \leq m$ satisfying these {\bf Taylor compatibility conditions} is called a {\bf Taylor compatible set}. 
A {\bf closed} subset $A$ of $\mathbb{R}^n$ has the {\bf Whitney extension property} if there is a function 
$f \in C^m$ such that $f^{k}(x)=f_k(x)$ for $x \in A$ 
and such that $f$ is real analytic on every point $\mathbb{R}^n \setminus A$. 

\satz{A Taylor compatible $F,A$ has the Whitney extension property. }

Hassler Whitney proved this result at Harvard in 1934 \cite{Whitney1934} just two years after getting his 
PhD there (remarkably in the completely different field of graph theory). He cites \cite{Beikowitsch1924} 
who proved a special case.  Chapter 12 in \cite{KendigDullMoment} gives an exposition of Whitney's work on the theorem. 
To cite from this book: {\it Hass found it a real challenge to go beyond the first dimension. He drew picture after
picture, but the problem seemed stubbornly intent on putting up a succession of frustrating barriers.
His eventual success in 1933 was a real tour de force for the 26-year old}.
\index{Whitney extension theorem}
\index{Taylor's theorem}
\index{Whitney extension property}
\index{Taylor compatibility}

\section{Markov's inequality}

\paragraph{}
Let $X: \Omega \to \mathbb{R}$ a random variable on a probability space
$(\Omega,\mathcal{A},{\rm P})$ and let $a>0$ be a real constant. In all what follows, using ${\rm E}[f(X)]$ for
some function $f$ assumes that this expectation is finite, meaning that $f(X) \in L^1(\Omega,{\rm P})$.
The {\bf Markov inequality} is ${\rm P}[|X| \geq a] \leq {\rm E}[|X|]/a$.
More generally, if $f:[0,\infty) \to [0,\infty)$ is a strictly
monotonically increasing function with $f(0)=0$ and $a>0$, then

\satz{
${\rm P}[|X| \geq a] \leq {\rm E}[f(|X|)]/f(a)$.
}

The proof is done by defining for all $a>0$ a new random variable $A(x)=a$ if $X(x) \geq a$
and $A(x)=0$ else. Then $0 \leq f(A(x)) \leq f(X(x))$ and
${\rm E}[f(X)] = \int_{\Omega} f(X(x)) \; dP(x)
\geq \int_{\Omega} f(A(x)) \; dP(x) = f(a) {\rm P}[X \geq a]$. This gives
$P[X \geq a] \leq {\rm E}[f(X)]/f(a)$. 
For example, if $f(x)=x^2$, applying the inequality to $X=Y-{\rm E}[Y]$
such that $f(X)={\rm Var}[Y]$, one has the {\bf Chebyshev's inequality}
${\rm P}[|Y-{\rm E}[Y]| \geq a] \leq \frac{{\rm Var}[X]}{a^2}$. 
For $f(x)=e^{x}$ one gets the {\bf Chernoff inequality}
${\rm P}[X \geq a] = {\rm P}[e^X \geq e^a] \leq \frac{{\rm E}[e^X]}{e^a}$. 
Since this works also for every $f(x)=e^{tx}$ with $t \geq 0$. One has the
{\bf Chernoff bound}
${\rm P}[X \geq a]  \leq {\rm inf}_{t \geq 0} \frac{{\rm E}[e^{t X}]}{e^{t a}}$ 
which is of interest as ${\rm E}[E^{t X}]$ is the {\rm moment generating function} of
the random variable $X$. 

\index{Chebychev's inequality}
\index{Markov's inequality}
\index{Chernoff bound}
\index{moment generating function}
\index{probability space}

\section{Magnus Freiheitssatz}

Let $G=(X,R)$ be a finitely presented group with {\bf generators} $X=\{ x_1,\dots ,x_q \}$ and
{\bf relations} $R=\{ r_1, \dots r_q \}$. The group $G$ is called a {\bf 1-relator group}
or {\bf Magnus group} if $q=1$ and if the relation $r$ has the property that $r$ is
{\bf cyclically reduced} and that all generators appear in $r$.  
[A word is called cylically reduced if every cylic permutation of the 
word is reduced. A word is called {\bf reduced} if it does not
contain subwords of the type $x x ^{-1}$ or $x^{-1} x$. For example, the 
word $aba^{-1}$ is not cyclically reduced. A cylic permutation can 
be reduced to $b$.]

\satz{$Y \subset X, Y \neq X$ generates a free group in a Magnus group $G$.}

This is a result of Wilhelm Magnus of 1930. It means that given say the generators
$x_1, \dots, x_{q-1}$, the only relations involving them are the trivial ones.
It is also called the {\bf Independence Theorem}. In \cite{FreiheitssatzMagnus1994} the
Freiheitssatz is considered a non-commutative analog of a similar result in a commutative
algebraic structure. If $V$ is a $n$-dimensional linear space over a field and $W \subset V$
is  a linear subspace given by a single equation $\sum_i a_i x_i = 0$, then $W$ has
dimension $n-1$ meaning that $W$ is a free Abelian subgroup of $V$. An other analogy
in that overview paper is to compare it with a situation in algebraic geometry:
{\it an irreducible algebraic equation of $n$ complex variables in which all $n$ variables
appear can not be used to derive any irreducible algebraic equation in which not all
of these variables appear.}
The theory of finitely presented groups
was initiated by Max Dehn in 1912. Magnus wrote his thesis in 1931 under the guidance
of Dehn. Dehn had raised the word problem, the conjugacy problem
and the isomorphism problem. Dehn also proposed that the Freiheitssatz could hold true. 
As Magnus pointed out, the Freiheitssatz assures that the 1-relator
group has a positive word problem solution. In 1954, Novikov came up with the 
first finitely presented group with insoluble word problem.  
See \cite{FreiheitssatzMagnus1994,MusingsBaumslag,Baumslag}.
\index{Freiheitssatz}
\index{Independence Theorem}
\index{Magnus theorem}
\index{finitely presented group}

\section{Martin's axiom}

Let $k$ be a {\bf cardinal}, the {\bf Martin condition} $M(k)$ is the statement that if $P$ is a
{\bf partial order} satisfying the {\bf countable chain condition} and family $D$ of dense sets in $P$
with cardinality less or equal than $k$, there is a filter $F$ on $P$ such that $F$ intersects
every element in $D$. The Martin axiom states
that $M(k)$ holds for every cardinality smaller than $2^{\alpha_0}$.  One has

\satz{In ZFC, $M(\aleph_0)$ holds but $M(2^{\alpha_0})$ fails.}

The first statement is the {\bf Rasiowa-Sikorski lemma}. The second statement
is proven by an example: the set $[0,1]$ with usual topology is separable and so
satisfies the countable chain condition. An individual pint is nowhere dense
but the union has $2^{\alpha_0}$ points.
Martin's axiom was introduced in 1970 by Tony Martin and Robert Solovay 
\cite{MartinSolovay}.
The {\bf continuum hypothesis} CH implies Martin's axiom but it is also 
consistent with $ZFC$ and the negation of CH. 

\index{ZFC}
\index{Martin Axiom}
\index{Rasiowa Sikorski theorem}
\index{Continuum hypothesis}


\section*{Epilogue: Value}

Which mathematical theorems are the most important ones? This is a complicated
variational problem because it is a general and fundamental problem in economics 
to define {\bf ``value"}. The difficulty with the concept is that ``value"
is often a matter of taste or fashion or 
social influence and so an {\bf equilibrium of a complex social system}. 
Value can change rapidly, sometimes triggered by small things. 
The reason is that the notion of value like in game theory depends on 
how it is valued by others. A fundamental principle of 
{\bf catastrophe theory} is that maxima of a functional can depend discontinuously 
on parameter. As value is often a social concept, this can be especially
brutal or lead to unexpected {\bf viral effects}. First of all, value is often 
linked to {\bf historical or morale considerations}. We tend more and more to link
artistic and scientific value also to the person. In mathematics, the work of 
Oswald Teichm\"uller or Ludwig Bieberbach for example are linked to their 
political view and so devalued despite their brilliance \cite{Segal2003}.
This happens also outside of science, in art or in industry. 
The value of a company now also depends on what ``investors think" or what 
analysts see for potential gain in the future. Social media try to measure 
value using ``likes" or ``number of followers". A majority vote is a measure
but how well can it predict correctly what be valuable in the future? Majority votes 
taken over longer times would give a more reliable value functional. Assume one could persuade every 
mathematician to give a list of the two dozen most fundamental theorems and do that 
every couple of years, and reflect the ``wisdom of an educated crowd", one could probably get a 
pretty good value functional. Ranking theorems and results in mathematics are a mathematical 
optimization problem by itself. One could use techniques known in the ``search industry".
One idea is to look at the finite graph in which the theorems are the nodes and where
two theorems are related to each other if one can be deduced
from the other (or alternatively connect them if one influences the
other strongly). One can then run a {\bf page rank algorithm} \cite{LangvilleMeyer}
to see which ones are important.
Running this in each of the major mathematical fields could give
an algorithm to determine which theorems deserve the name ``fundamental".
Now, there was also a problem with publishing the page rank as people tried to manipulate
it using search engine optimization tricks. Google now does no more give the page rank
of a website, simply to avoid such manipulations. The story illustrates that reflecting
about algorithms that measure value can influence the algorithm itself and even destroy it.
Similarly as in quantum mechanics, the measurement process can influence the
experiment to the point that it is no more reliable. 
\index{Page rank}
\index{value}
\index{viral effects}
\index{catastrophe theory}
\index{page rank}

\section*{Opinions} 

It had been a course ``Math from a historical perspective" taught a couple of times at
the Harvard extension school has motivated to write up the present document.
As part of a project it was often asked to to write about some theorems or mathematical 
fields or a mathematical person and try to rank it. 
The present document benefits from these writings as it is interesting to see 
what others consider important. Sometimes, seeing different opinions can change your 
own view. I was definitely influenced by students, teachers, colleagues and literature
as well of course by the limitations of my own understanding. 
My own point of view has already changed while writing the actual theorems down and will
certainly change more. Value is more like an equilibrium of many different factors.
In mathematics, values have changed rapidly over time. And mathematics
can describe the rate of change of value \cite{StewartCatastrophe}.
Major changes in the appreciation for mathematical topics came throughout the history.
Sometimes with dramatic shifts like when mathematical notations started to appear, at the time of 
Euclid, then at the time when calculus was developed by Newton and Leibniz. Also the 
development of more abstract algebraic constructs or topological notions, 
like for example the start of set theory changed things considerably. In more modern times,
the {\bf categorization of mathematics} and the development of rather general and
abstract new objects, (for example with new approaches taken by Grothendieck) 
changed the landscape. In most of the new development, I remain the puzzled tourist
wondering how large the world of mathematics is. It has become so large that continents
have emerged: we have {\bf applied mathematics}, {\bf mathematical physics}, {\bf statistics}, 
{\bf computer science} and {\bf economics} which have drifted away to independent subjects
and departments. Classical mathematicians like Euler would now be called applied mathematicians, 
de Moivre would maybe be stamped as a statistician, {\bf Newton} a mathematical physicist
and {\bf Turing} a computer scientist and {\bf von Neuman} an economist or physicist.

\section*{Search} 

A couple of months before starting this document in 2018, when looking online for ``George Green",
the first hit in a search engine would be a 22 year old soccer player. 
(This was not a search bubble thing \cite{Filterbubble} as it was
tested with cleared browser cache and via anonymous VPN from other locations, where the search engine 
can not determine the identity of the user). Now, I love soccer, 
played it myself a lot as a kid and also like to watch it on screen, 
but is the English soccer player George William Athelston Green 
really more ``relevant" than the British mathematician George Green, who made fundamental 
break through discoveries which are used in mathematics and physics? 
Shortly after I had tweeted about this strange ranking on December 27, 2017, 
the page rank algorithm must have been adapted, because already on January 4th, 2018,
the Mathematician George Green appeared first (again not a search bubble phenomenon, where
the search engine adapts to the users taste and adjusts the search to their preferences). 
It is not impossible that my tweet has reached, meandering
through social media, some search engine engineer who was able to 
rectify the injustice done to the miller and mathematician George Green. 
The theory of networks shows ``small world phenomena" \cite{SixDegrees,Linked,SmallWorld} 
can explain that such influences or synchronizations are not that impossible \cite{Sync}. But
coincidences can also be deceiving. Humans just tend to observe coincidences even so there
might be a perfectly mathematical explanations. This is prototyped by the {\bf birthday paradox}
\cite{Fluke}. But one must also understand that search needs to serve the majority. 
For a general public, a particular subject like mathematics is not that important. 
When searching for ``Hardy" for example, it is not Godfrey Hardy who is 
mentioned first as a person belonging to that keyword 
but Tom Hardy, an English actor. This obviously serves most of the searches better. 
As this might infuriate particular groups (here mathematicians), search engines have started to adapt the 
searches to the user, giving the search some {\bf context} which is an important 
ingredient in artificial intelligence. The problem is the search bubble phenomenon which runs hard against 
objectivity. Textbooks of the future might adapt their language, difficulty and 
even their citations or the historical credit on who reads it. Novels might adapt the language to
the age of the user, the country where the user lives, and the ending might depend
on personal preferences or even the medical history of the user (the medical history of course being 
accessible by the book seller via `big data" analysis of user behavior and tracking which is not SciFi
this is already happening): even classical books are cleansed for political correctness, 
many computer games are already customizable to the taste of the user. 
A person flagged as sensitive or a young child might be served a happy ending in a novel rather 
than a conclusion of the novel in an ambivalent limbo or
even a disaster. \cite{Filterbubble} explains the difficulty. The issues have amplified
even more in more recent times. The phenomenon of filter bubble  even influences elections 
and polarizes opinions as one does not even hear any more alternate arguments.

\section*{Beauty} 

In order to determine what is a ``fundamental theorem", also 
aesthetic values matter. But the question of ``what is beautiful" is even trickier.  
Many have tried to define and investigate the mechanisms of beauty:
\cite{HardyApology,Wells1,Wells2,RotaBeauty,SinclairMathBeauty,AharoniBeauty,MontanoExplainingBeauty}. 
In the context of mathematical formulas, the question has been investigated within the field of
{\bf neuro-aesthetics}. Psychologists, in collaboration with mathematicians have measured
the brain activity of $16$ mathematicians with the goal to determine what they consider
beautiful \cite{ZRBA}. The {\bf Euler identity} $e^{i \pi} + 1 = 0$ was rated high with a value
0.8667 while a formula for $1/\pi$ due to Ramanujan was rated low with an average rating
of -9.7333. Obviously, what mattered was not only the complexity of the formula but also how
much {\bf insight} the participants got when looking at the equation. The authors of
that paper cite Plato who wrote once {\it "nothing without understanding would ever be 
more beauteous than with understanding"}. Obviously, the formula of Ramanujan is much
deeper but it requires some background knowledge for being appreciated.
But the authors acknowledge in the discussion that that correlating 
``beauty and understanding" can be tricky. 
Rota \cite{RotaBeauty} notes that the appreciation of mathematical beauty in some statement
requires the ability to understand it. And \cite{MontanoExplainingBeauty} notices that
``even professional mathematicians specialized in a certain field might find results or proofs in 
other fields obscure" but that this is not much different from say music, where ``knowledge about
technical details such as the differences between things like cadences, progressions or chords changes 
the way we appreciate music" and that ``the symmetry of a fugue or a sonata are simply invisible without 
a certain technical knowledge". As history has shown, there were also
always ``artistic connections" \cite{Gamwell,BruterArt} as well as ``religious influences"
\cite{Livio2009,Skinner}. The book \cite{Gamwell} cites Einstein who defines 
``mathematics as the poetry of logical ideas". It also provides many examples and
illustrations and quotations. And there are various opinions.  Rota argues that
beauty is a rather objective property which depends on historic-social contexts. 
And then there is {\bf taste}: what is more appealing, the element of {\bf surprise} like the Birthday paradox
or Petersburg paradox in probability theory, the {\bf Banach-Tarski
paradox} in measure theory which obviously does not trigger any {\bf enlightenment} nor
{\bf understanding} if one hears the first time:
one can disassemble a sphere into 5 pieces, rotate and translate these pieces in space
to build up two spheres. Or the surprising fact that the infinite sum $1+2+3+4+5+ \dots$ 
is naturally equal to $-1/12$ as it is $\zeta(-1)$ (which is a value defined by 
{\bf analytic continuation} and can hardly be understood without training in complex analysis).
The role of aesthetic in mathematics is especially important in education, where
mathematical models \cite{FischerModels}, mathematical visualization \cite{Banchoff1990},
artistic enrichment \cite{Fomenko}, 
surfaces \cite{AnalyticSurfaces}, or 3D printing \cite{Segerman2016,CFZ} 
can help to make mathematics more approachable. 
Update 2019: as reported in Science Daily a study of the university of Bath concludes that 
people appreciate beauty in complex mathematics \cite{JohnsonSteinerberger}.
The results which had been chosen in that study had been rather simple however: 
the infinite geometric series formula, the Gauss's summation trick for positive integers,  
the Pigeonhole principle, and a geometric proof of a {\bf Faulhaber formula} for 
the sum the first powers of an integer. 
When judging the mathematics describing physical models, Paul Dirac was probably the most outspoken
advocate for beauty. He stated in \cite{Dirac1963} for example:
{\it It seems to be one of the fundamental features of nature that fundamental
physical laws are described in terms of a mathematical theory of great beauty and power, needing quite a high
standard of mathematics for one to understand it}.

\section*{Deepness}

A {\bf taxonomy} is a way to place objects like theorems in an multi-dimensional cube of numerical 
attributes. Besides the {\bf ugly-beauty} parameter, one can think of all kind of {\bf taxonomies} 
to classify theorems. There is the {\bf simplicity-complexity axes}, which could be measured by the
number of mathematicians who can understand the proof, the {\bf boring-interesting axes} which
measures the entertainment value or potential for pop culture appearances, 
the {\bf useless-applicable axes} which measures how many applications the theorem has 
in engineering, economics or other sciences, the {\bf easy-hard} which could be measured 
in the amount of time one needs to understand the proof. And then there is the 
{\bf shallow-deepness} axes, which is even more subjective but which could be quantified too. 
One could look for example, how long a proof path is from basic axioms to the theorem and 
weight each path with how many other interesting theorems have been visited along. 
Also of benefit are how many different areas of mathematics have been visited along the proof. 
A deep theorem could be obtained by proving it with different long paths, each reaching 
other already established deep results. One can now
argue how to average all these paths, whether one should take the minimum or maximal deep proof
path. The later point was addressed in \cite{Lange2014}. \\

Maybe unlike with other parameters, the antipode ``trivial" of ``deepness" has a positive side too: it is maybe
not ``shallow" but what we call ``fundamental". Fundamental theorems are not necessarily deep. 
The Pythagorean theorem for example or Zorn's lemma are not deep but they are fundamental. 
Basic logical identities based on Boolean algebra which are used in almost every proof step are of 
fundamental importance but not deep. 
One could still go back and measure how fundamental something is by how many deep theorems
can be proven with it. \\

\cite{Urquhart2015} points out that the adjective ``deep" is used for all kind of mathematical
objects: theorems, proofs, problems, insights or concepts can be described as deep and that
often the theorem is called deep if its proof is deep. Urquhart points out however that
``if a simple proof is discovered later, perhaps the result might be reclassified as not deep at all"
and that so, the difficulty of the concept ``mathematical depth" is not so well defined. The author
then mentions the {\bf graph minor theorem} (in every infinite set of graphs there are two for which one is the
minor of the other), which Diestel \cite{Diestel} calls ``one of the deepest theorems that mathematics
has to offer. Some justification for the deepness of the result is that it has made impact also
outside graph theory and that its proof takes well over 500 pages. \\

\cite{Urquhart2015} also collects opinions of philosophers and mathematics about deepness.
Cited is for example \cite{HardyApology} as Hardy gives an extended discussion on depth and
sees mathematical ideas ``arranged somehow in strata, each stratum being linked by a complex
relation both among themselves and with those above and below, the lower the stratum, the deeper
the idea. Also cited is the book of Penelope Maddy \cite{Maddy2011} which expresses 
doubt that that mathematical depth really can be accounted for productively because 
it is a ``catch-all" for the various kinds of virtues
and often used as a term of approbation, but always in an informal context without giving a precise
meaning. Also cited in \cite{Urquhart2015} are present day mathematicians like Gowers
\cite{Gowers2008} who links ``deep" with ``hard" 
and contrasts it with ``obvious". If a proof requires a non-obvious
idea, then it is considered deep. Also cited is a later statement of Gowers telling that
``The normal use of the word `deep' is something like this: a theorem is deep if
it depends on a long chain of ideas, each involving a significant insight".
Finally mentioned is Tao \cite{Tao2007} who lists over twenty meanings to ``good mathematics":
(be a breakthrough for solving a problem, masterfully using technique,
building theory, having insight, discovering something unexpected, having
application, clear exposition, good pedagogy enabling understanding,
long-range vision, good taste, public relations, advancing foundations,
rigorous, beautiful, elegant, creative, useful, sharp to known
counterexamples, intuitive and visualisable, being definitive like a
classification result and finally {\bf deep} which Tao defines as ``manifestly
non-trivial, for instance by capturing a subtle phenomenon beyond the
reach of more elementary tools".)  \cite{Urquhart2015}  also illustrates
the concept of deepness with moves one sees in chess: a combination of moves which are not obvious
and have an element of surprise like in the Byrne-Fischer game of 1963-1964. \\

In a talk ``Mathematical Depth Workshop" of April 11,12, 2014 John Stillwell gave the following
examples of deep theorems: Dirichlet's theorem on primes in an arithmetic progression,
Perelman's theorem on Poincar\'e's conjecture, Fermat's last theorem and then
the classification of finite simple groups. A deep theorem should be
difficult, surprising, important, fruitful, elegant and fundamental.
As less deep but accessible, he gives the independence of the 
parallel postulate, the fundamental theorem of algebra,
the existence of division algebras, the Riemann integrability of continuous functions, the uncountability
of $\mathbb{R}$.  Robert Geroch told in that same workshop that deep theorems
should be detached from connections with people, or then have connections with physics:
examples are representations of the Lie group $SL(2,\mathbb{C})$,
the {\bf TCP theorem} or the appearance of symmetric hyperbolic partial differential equations.
Jeremy Gray stressed then the importance of multiple proofs, to give more
reasoning, show different methodologies, see new routes or produce more purity.
He said that the difference between deep and difficult is that deep things should be more hidden.
Deep according to Gauss has to be ``difficult". The result may be elegant or beautiful,
but the proof needs to be difficult.
Marc Lange \cite{Lange2014} argues to assign the attribute deep to the proof of a theorem and not
the theorem itself. The reason is that there could be multiple proofs, where one proof is deeper 
than the other. This could mean for example that a theorem which is considered deep, remains to have
a deep proof even in the case if it turns out to be provable in a very simple and dull way. 
\index{Deepness}
\index{Mathematical depth}

\section*{The fate of fame} 

Aesthetics is a fragile subject. If something beautiful has become too popular and so 
entered {\bf pop-culture}, a natural aversion against it can develop. The feeling is
justified that popular things are often frivolous. It is also in danger to become a 
{\bf clish\'e} or even become {\bf kitsch} (which is a word used to tear down popular 
stuff or to label poor taste). The Mandelbrot
set for example is just marvelous, but it does hardly does excite anymore because it is so commonly
known. The {\bf Monty-Hall problem} which became famous by Gardner columns in the early 1990'ies
(see \cite{Snell95,Rosenhouse}) was cool to teach in 1994, three years after the infamous
``parade column" of 1991 by Marilyn vos Savant which blew it into the spot light. 
But especially after a cameo in the movie ``21", 
the theorem has become part of {\bf mathematical kitsch}. I myself love 
mathematical kitsch. A topic that gained that status must have been nice and innovative to obtain
that label. Kitsch becomes only tiresome however if it is not presented in a new and original form. 
The book \cite{PR}, in the context of complex dynamics, remains a master piece still 
today, even-so the picture have become only too familiar, but rendering the Mandelbrot set today in
that same way hardly does the rock the boat any more. 
Still, it remains fascinating and more and youtube allows to see sophisticated 
zooms down to the size of $10^{-200}$. In that context, it appears strange that mathematicians do not
jump on the ``Mandelbulb set" $M$, a three dimensional version of the Mandelbrot set
which is one of the most beautiful mathematical objects.
The reason could be that as a ``youtube star" it is not worthy yet any serious academic consideration;
more likely however is that the object is just too difficult for a serious study, as we lack the mathematical 
analytic tools which for example would just to answer a basic question like whether $M$ is connected.
A second example is {\bf catastrophe theory} \cite{StewartCatastrophe,ArnoldCatastrophe} a beautiful part
of {\bf singularity theory} which started with Hassler Whitney and was then developed by Ren\'e Thom
\cite{ThomMorphogenesis}.
It was hyped to much that it fell into a deep fall from which it has not yet fully recovered. 
This happened despite the fact that Thom himself already pointed out the limits, as well as the controversies 
of the theoryy \cite{Boutot1993}. It had to pay a prize for its fame and appears to be forgotten. 
Chaos theory from the 60ies which started to peak with Edward Lorenz and terms like
the ``Butterfly effect" ``strange attractors" 
started to become a {\bf clish\'e} latest after that infamous scene featuring 
the character Ian Malcolm in the 1993 movie Jurassic park. 
It was laughed at already within the same movie franchise, when in
the third Jurassic Park installment of 2001, the kid {\bf Erik Kirby} snuffs on Malcolm's 
``preachiness" and quotes his statement ``everything is chaos" in a condescending way. 
In art, architecture, music, fashion or design also, if something has
become too popular, it is despised by the ``connaisseurs". Hardly anybody would consider a ``lava lamp"
(invented in 1963) a object of taste nowadays, even so, the fluid dynamics and motion is objectively 
rich and interesting, illustrating also geometric deformation techniques in geometry like the Ricci flow. 
The piano piece ``F\"ur Elise" by Ludwig van Beethoven became so popular that 
it can not even be played any more as background music in a supermarket. There is something which 
prevents a ``serious music critic" to admit that the piece is great, genius due to its simplicity.
Such examples suggest that it might be better for an achievement 
(or theorem in mathematics) not to enter pop-culture as this indicates a lack of ``deepness" and
is therefore despised by the elite. The principle of having fame torn down to disgrace is common also outside
of mathematics. Famous actors, entrepreneurs or politicians are not universally admired but sometimes
hated to the guts, or torn to pieces and certainly can hardly live normal lives any more. 
The phenomenon of accumulated critique got amplified with mob type phenomena in social media. 
There must be something fulfilling to trash achievements, the simplest explanation being envy.  
Film critics are often harsh and judge negatively because this elevates their own status as 
they appear to have a ``high standard". Similarly morale judgement is expressed often just to elevate 
the status of the judge even so experience has shown that often judges are offenders themselves and
the critique turns out to be a compensation.
Maybe it is also human ``Schadenfreude", or greed which makes so many to voice critique. 
History has shown however that social value systems do not matter much in the long term. 
A good and rich theory will show its true value if it is appreciated also in hundreds of years, 
where fashion and social influence have no more any impact. The theorem of Pythagoras will be important
independent of fame and even if it has become a clich\'e, it is too important to be labeled as such. It has
not only earned the status of kitsch, it is also a prototype as well as a useful tool. 
 
\section*{Media} 

There is no question that the {\bf Pythagorean theorem}, the {\bf Euler polyhedron formula} $\chi=v-e+f$
the {\bf Euler identity} $e^{i \pi} + 1 =0$, or the 
{\bf Basel problem formula} $1+1/4+1/9+1/16 + \cdots = \pi/6$ 
will always rank highly in any list of beautiful formulas. 
Most mathematicians agree that they are elegant and beautiful. These results will also
in the future keep top spots in any ranking. On social networks, one can find lists of favorite formulas.
On ``Quora", one can find the arithmetic mean-geometric mean inequality
$\sqrt{ab} \leq (a+b)/2$ or the {\bf geometric summation formula} $1+a+a^2 + \cdots =1/(1-a)$
high up. One can also find strange contributions in social media like the identity 
$1=0.99999\dots$ which is used by Piaget inspired educators to probe mathematical 
maturity of kids. Similarly as in Piaget's experiments, there is time of mathematical maturity
where a student starts to understand that this is an identity. A very young student thinks
$1$ is larger than $0.9999...$ even if told to point out a number in between. Such threshold moments
can be crucial for example to mathematical success later. We have a strange fascination with 
``wunderkinds", kids for which some mathematical abilities have come earlier (even so the existence
of each wonder kid produces a devastating collateral damage in its neighborhood as their success sucks
out any motivation of immediate peers). The problem is also that if somebody does not
pass these Piaget thresholds early, teachers and parents consider them lost, they get discouraged 
and become uninterested in math (the situation in other art or sport is similar). 
In reality, slow learners for which the thresholds are passed
later are often deeper thinkers and can produce deeper or more extraordinary results.
At the moment, searching for the ``most beautiful formula in mathematics"
gives the Euler identity and search engines agree. 
But the concept of taste in a time of social media can be confusing. We live in an epoch, where
a 17 year old ``social influencer" can in a few days gather more ``followers"
and become more widely known than {\bf Sophie Kovalewskaya} who made fundamental beautiful
and lasting contributions in mathematics and physics like the Cauchy-Kovalevskaya theorem.
Such a theorem is definitely more lasting than a few ``selfie shots" of a pretty face, but 
measured by a ``majority vote", it
would not only lose, it would completely disappear. One can find youtube videos of kids explaining the
4th dimension, which are watched millions of times, many thousand times more than videos of 
mathematicians who have created deep mathematical new insight about four dimensional space. 
But time rectifies. Kovalewskaya will also
be ranked highly in 50 years, while the pretty face has faded. Hardy put this even more
extremely by comparing a mathematician with a literary heavy weight: 
{\it Archimedes will be remembered when Aeschylus is forgotten, because languages die 
and mathematical ideas do not.} \cite{HardyApology}
There is no doubt that film and TV (and now internet like ``Youtube", social networks 
and ``blogs") has a great short-term influence on value or exposure of a 
mathematical field. Examples of movies with
influence are {\bf It is my turn} (1980), or {\bf Antonia's line} (1995) 
featuring some algebraic topology, {\bf Good will hunting} (1997) 
in which some graph theory and Fourier theory appears,
{\bf 21} from (2008) which has a scene in which the Monty Hall problem has a cameo.
{\bf The man who knew infinity} displays the work of Ramanujan and promotes some combinatorics
like the theory of partitions. There are lots of movies featuring cryptology like {\bf Sneakers} (1992), 
{\bf Breaking the code} (1996), {\bf Enigma} (2001) or {\bf The imitation game} (2014).
For TV, mathematics was promoted nicely in {\bf Numb3rs} (2005-2010). 
For more, see \cite{MathMovies} or my own online math in movies collection. \\

\section*{Professional opinions} 

Interviews with professional mathematicians can also probe the waters.  In \cite{Kondratieva},
Natasha Kondratieva has asked a number of mathematicians: ``What three mathematical
formulas are the most beautiful to you". The {\bf formulas of Euler} or the {\bf Pythagoras 
theorem} naturally were ranked high. Interestingly, Michael Atiyah included even a formula 
{\bf "Beauty = Simplicity + Depth"}. Also other results, like the 
{\bf Leibniz series} $\pi/4 = 1-1/3+1/5-1/7+1/9 - \dots$, the {\bf Maxwell equations} 
$dF=0,d^*F=J$ or the {\bf Schr\"odinger equation} 
$i \hbar u'=(i \hbar \nabla + e A)^2 u + V u$, the {\bf Einstein formula} $E=mc^2$ or the
{\bf Euler's golden key} $\sum_{n=1}^{\infty} 1/n^s = \prod_p (1-1/p^s)^{-1}$  
or the {\bf Gauss identity} $\int_{-\infty} ^{\infty} e^{-x^2} dx = \sqrt{\pi}$ or
the {\bf volume of the unit ball} in $R^{2n}$ given as $\pi^n/n!$ appeared.
Gregory Margulis mentioned an application of the {\bf Poisson summation formula}  
$\sum_n f(n) = \sum_n \hat{f}(n)$
which is $\sqrt{2} \sum_n e^{-n^2} = \sum_n e^{-n^2/4}$
or the {\bf quadratic reciprocity law} $(p|q) = (-1)^{(p-1)/2 (q-1)/2}$, where
$(p|q)=1$ if $q$ is a {\bf quadratic residue} modulo $p$ and $-1$ else. 
Robert Minlos gave the {\bf Gibbs formula}, a {\bf Feynman-Kac formula} or the {\bf Stirling formula}.
Yakov Sinai mentioned the {\bf Gelfand-Naimark realization} of an Abelian $C^*$ algebra as an 
algebra of continuous function or the {\bf second law of thermodynamics}.
Anatoly Vershik gave the generating function $\prod_{k=0}^{\infty} (1+x^k) = \sum_{n=0}^{\infty} p(n) x^n$ 
for the {\bf partition function} $p(n)$ and the {\bf generalized Cauchy inequality} 
between arithmetic and geometric mean. 
An interesting statement of David Ruelle appears in that article who quoted Grothendieck
by `` my life's ambition as a mathematician, or rather my joy and passion, have constantly been 
to discover obvious things \dots". Combining Grothendieck's and Atiyah's quote, 
fundamental theorems should be ``obvious, beautiful, simple and still deep". \\
A recent column ``Roots of unity" in the Scientific American asks mathematicians for their
favorite theorem: examples are {\bf Noether's theorem}, {\bf the uniformization theorem},
the {\bf Ham Sandwich theorem}, the {\bf fundamental theorem of calculus}, {\bf the circumference of the circle},
the {\bf  classification of compact 2-surfaces}, {\bf Fermat's little theorem}, the 
{\bf Gromov non-squeezing theorem}, a theorem about Betti numbers, the {\bf Pythagorean theorem}, the 
{\bf classification of Platonic solids}, the {\bf Birkhoff ergodic theorem}, the {\bf Burnside lemma}, 
the {\bf Gauss-Bonnet theorem}, Conways rational tangle theorem, 
{\bf Varignon's theorem}, an upper bound on {\bf Reidemeister moves in knot theory}, the 
{\bf asymptotic number of relative prime pairs}, the {\bf Mittag Leffler theorem}, 
a theorem about {\bf spectral sparsifiers}, the {\bf Yoneda lemma} and the {\bf Brouwer fixed point theorem}.
These interviews illustrate also that the choices are different if asked for ``personal favorite theorem"
or ``objectively favorite theorem".  \\
\index{Einstein formula}
\index{quadratic reciprocity}
\index{Gibbs formula}
\index{Gauss identity}
\index{partition function}
\index{Maxwell equations}
\index{Stirling formula}
\index{Euler golden key}
\index{Varignon theorem}
\index{Nonsqueezing theorem}
\index{sprasifiers}

\section*{Fundamental versus important} 

Asking for fundamental theorems is different than asking for ``deep theorems" or ``important theorems".
Examples of deep theorems are the {\bf Atiyah-Singer} or {\bf Atiyah-Bott theorems} in differential topology,
the {\bf KAM theorem} related to the strong implicit function theorem, or the {\bf Nash embedding theorem} in 
Riemannian geometry. An other example is the {\bf Gauss-Bonnet-Chern theorem} 
in Riemannian geometry or the {\bf Pesin theorem} in partially hyperbolic dynamical systems. 
Maybe the {\bf shadowing lemma} in hyperbolic dynamics is more fundamental than the much deeper 
{\bf Pesin theorem} (which is still too complex to be proven with full details in any classroom.
Also excellent textbooks like \cite{Pollicott1993,KatokStrelcyn} do not prove the full theorem 
establishing the Bernoulli property on ergodic components). 
One can also argue, whether the {\bf ``theorema egregium"} of Gauss, stating that the curvature of a 
surface is intrinsic and not dependent on an embedding is more ``fundamental" than 
the {\bf ``Gauss-Bonnet"} result, which is definitely deeper. 
In number theory, one can argue that the {\bf quadratic reciprocity formula}
is deeper than the {\bf little Theorem of Fermat} or the {\bf Wilson theorem}. (The
later gives an if and only criterion for primality but still is far less important than the little 
theorem of Fermat which as the later is used in many applications.)
The {\bf last theorem of Fermat} \cite{Boston2003} is an example of an important theorem as it is deep and
related to other fields and culture, but it is not yet so much a ``fundamental theorem". 
Similarly, the {\bf Perelman theorem} fixing the {\bf Poincar\'e conjecture} is important, but it is
not (yet) a fundamental theorem. It is still a mountain peak and not a sediment in a rock. 
Important theorems are not much used by other theorems as they are located at the end of a development. 
Also the solution to the {\bf Kepler problem} on sphere packings
or the proof of the {\bf 4-color theorem} \cite{ChartrandZhang2} or the proof of the 
{\bf Feigenbaum conjectures} \cite{DeMelo,Lanford84} are important results
but not so much used by other results. Important theorems build {\bf the roof} of the building, while 
fundamental theorems form the {\bf foundation} on which a building can be constructed. But this
can depend on time as what is the roof today, might be in the foundation later on, once more floors have
been added.  \\
\index{Atiyah Singer}
\index{Atiyah Bott}
\index{KAM}
\index{shadowing}
\index{last theorem of Fermat}
\index{4-color theorem}
\index{Feigenbaum conjectures}
\index{Gauss-Bonnet}
\index{Perelman theorem}
\index{4-color theorem}
\index{Wilson theorem}
\index{Kepler problem}

\section*{Essential math}

In education it is necessary regularly to reexamine what a student of mathematics needs to know. 
What are essential fields in mathematics? Also here, there are many opinions and things are always in the 
flux. The 7 {\bf liberal arts of sciences} was an early attempt to organize things in a larger scale.
For example, while in the 19th century, quaternions were considered essential, they
fell out of the curriculum and today, it is well possible that a student learns about 
division algebras only in graduate school. 
One of the questions is how to balance {\bf applicability} and {\bf elegance}.
In pure mathematics, one might more focus on beauty and elegance, in applied mathematics, the
applicability is important. As the field of mathematics has expanded enormously, there
is the problem of fragmentation. On the other hand, the mathematical fields have also split. 
Some domains have been ``taken over" by new departments like applied mathematics, computer science 
or statistics. Discrete mathematics courses like graph theory or theory of computation or 
cryptology are now in the hands of computer science, differential equations or numerical analysis by
applied mathematics departments, probability theory courses taught by statistics
departments. Still, there is a core of mathematical content which a mathematician should at least 
have been exposed to. A student studying the subject should probably have an eye on both getting into 
a field which looks promising for research as well as having a broad general education in all
possible fields. One can get an idea what is required in various mathematics departments by 
looking at what are called ``general examinations" or ``qualifying examinations". These are
exams given to first year graduate students which have to be passed. Departments like 
Harvard \cite{HarvardQuals} or Princeton \cite{PrincetonGenerals} have many of these questions 
in the public. Also here, one could go to the AMS classification and grind through all topics.
Instead, let us try an attempt to put it all in one box, being aware that other priorities can
work too:  \\

\begin{tiny}
\begin{tabular}{|l|l|} \hline
Pre-calculus& Algebra, Trig functions,Log and Exp Functions, Graphs, Modeling, Geometry, Solving equations, Inequalities \\
Single variable& Functions, Limits, Continuity, Differentiation, Integration, Series, Differential equations,
Fundamental theorem \\
Multi variable& Vectors, Geometry, Functions, Differentiation, Integration, Vector calculus: Green Stokes and Gauss \\
Linear algebra& Linear equations, Determinants, Eigenvalues, Projection and Data-fitting, Differential equations, Fourier theory \\
Dynamical systems& Iteration of maps, Ordinary and partial differential equations, Bifurcation theory, Integrability, Ergodic Theory \\
Probability& Probability spaces, Random variables, Distributions, Stochastic Processes, Statistics, Data, Estimation \\
Discrete math& Combinatorics, Graphs, Order structures, Counting tools, Theory of computation, Complexity, Game Theory \\
Numerics& Algorithms,Integration, Solving ODE's, PDE's, Approximation techniques, Interpolation, Comput. Geometry \\
Analysis& Functional analysis, Banach algebras, Complex analysis, Harmonic analysis, Fourier theory, Laplace, PDE's \\
Algebra& Groups and Rings, Modules, Vector Spaces, Commutative algebra, Non-commutative Rings, Galois theory \\
Number theory& Primes, Diophantine equations and approximations, Geometry of numbers, Dirichlet Series, Zeta function \\
Geometry& Differential topology, Differential Geometry, Geodesics, Curvature, Invariants, Geometric Measure theory \\
Alg. Geometry& Affine and Projective varieties, Ringed spaces, Schemes, Sheaf Theoretical Methods, Cohomology, Categories \\
Topology&  Set theoretical topology, Fractal Geometry, Differential topology, Homotopy, Algebraic Topology, Topos theory \\
Logic& First/second order Logic, Foundations, Models, Incompleteness, Forcing, Computability, New Axiom systems \\
Real analysis& Foundations, Metric spaces, Measure theory, Theory of integration on delta rings, Non-standard analysis \\
Computer Science& Math software, Programming Paradigms, Computer Architecture, Data structures, Big Data, Machine Learning \\
Connections& History, Big picture, Number systems, Notation, Linguistic, Psychology, Philosophy, Sociology and Pedagogy \\ \hline
\end{tabular}
\end{tiny}

\section*{Open problems} 

The importance of a result is also related to {\bf open problems} attached to 
the theorem. Open problems fuel new research and new concepts. Of course this is a moving target
but any ``value functional" for ``fundamental theorems" is time dependent and a bit
also a matter of {\bf fashion}, {\bf entertainment} (TV series like ``Numbers" or
Hollywood movies like ``good will hunting" changed the value) 
and under the influence of {\bf big shot mathematicians} which serve as ``influencers". 
Some of the problems have {\bf prizes} attached like the {\bf 23 problems of Hilbert}, 
the {\bf 15 problems of Simon} \cite{Simon15Problems}, the {\bf 18 problems of Smale},
the {\bf Yau problems} in geometry \cite{YauSeminar1982},
the {\bf 10 Millenium problems} or the four {\bf Landau problems} (Goldbach conjecture, twin prime conjecture,
the existence of primes between consecutive primes and the existence of infinitely many primes of the 
form $n^2+1$) and then the {\bf oldest problem of mathematics} the existence of {\bf odd perfect numbers}.
\index{perfect numbers}

There are beautiful open problems in any major field and
building a ranking would be as difficult as the problem to rank theorems. It is a bit a personal matter. 
I like the odd perfect number problem because it is the oldest problem in mathematics. 
Also Landau's list of 4 problems are clearly on the top. They are shockingly short and elementary but brutally hard,
having resisted more than a century of attacks by the best minds. There are other problems, where
one believes that the mathematics has just not been developed yet to tackle it, 
an example being the Collatz (3k+1) problem. With respect to the Millenium problems, one could argue that
the Yang-Mills gap problem is a rather vague. The problem looks like ``made by humans" while a problem like the 
odd perfect number problem has been ``made by the gods". 

There appears to be wide consensus that the {\bf Riemann hypothesis} is the most important open problem in 
mathematics. It states that the roots of the Riemann zeta function are all located on the
axes ${\rm Re}(z)=1/2$. In number theory, the {\bf prime twin problem} or the {\bf Goldbach problem} have
a high exposure because they can be explained to a general audience without mathematics background.
For some reason, an equally simple problem, the {\bf Landau problem} asking whether there are infinitely
many primes of the form $n^2+1$ is much less well known. In recent years, due to an alleged proof
by Shinichi Mochizuki of the ABC conjecture using a new theory called 
{\bf Inter-Universal Teichm\"uller Theory} (IUT) which so far is not accepted by the main mathematical
community despite strong efforts. But it has put the ABC conjecture from 1985 in the spot light like
\cite{WolchoverABC}. It has been described in \cite{GoldFeld1996} as the most important problem in Diophantine
equations. It can be expressed using the {\bf quality} $Q(a,b,c)$ of three integers $a,b,c$ which is 
$Q(a,b,c)=\log(c)/\log({\rm rad}(abc))$, where the {\bf radical} ${\rm rad}(n)$ of a number $n$ is the product of the
distinct prime factors of $n$. The ABC conjecture is that for any real number $q>1$ there exist only finitely many
triples $(a,b,c)$ of positive relatively prime integers with $a+b=c$ for which $Q(a,b,c)>q$. The triple with
the highest quality so far is $(a,b,c)=(2,3^{10} 109, 23^5)$; its quality is $Q=1.6299$. 
And then there are entire collections of conjectures, one being the {\bf Langlands program}
which relates different parts of mathematics like number theory, algebra, 
representation theory or algebraic geometry. I myself can not appreciate this program yet because I need first
to understand it. My personal favorite problem is the {\bf entropy problem}
in smooth dynamical systems theory \cite{Katok2007}. The {\bf Kolmogorov-Sinai entropy} of a smooth 
dynamical system can be described using Lyapunov exponents. For many systems like smooth convex billiards,
one measures positive entropy but is unable to prove it. An example is the real analytic $l^4$ table $x^4+y^4=1$ \cite{JeKn96}.
For ergodic theory, see \cite{CFS,DGS,Friedman,Sinai}.

\index{radical}
\index{quality}
\index{ABC conjecture}
\index{Landau problem} 
\index{Open problems}
\index{Millenium problems}
\index{Math in movies}
\index{Hollywood}
\index{prizes}
\index{Riemann hypothesis}
\index{Langlands program}
\index{Inner Universal Teichmuller Theory}
\index{Simon's problems}
\index{Smale's problems}
\index{Millenium problems} 

\section*{Classification results} 

One can also see classification theorems like the above mentioned 
Gelfand-Naimark realization as mountain peaks in the landscape of mathematics. 
Examples of {\bf classification results} are the classification of regular or 
semi-regular polytopes, the classification of discrete subgroups of a Lie group, 
the classification of ``Lie algebras", the classification of 
``von Neumann algebras", the ``classification of finite simple groups", the {\bf classification of
Abelian groups}, or the classification of associative {\bf division algebras} which by Frobenius is 
given either by the real or complex or quaternion numbers. Not only in algebra, also in 
differential topology, one would like to have classifications like the classification of $d$-dimensional 
manifolds. In topology, an example result is that every Polish space is homeomorphic to 
some subspace of the {\bf Hilbert cube}. Related to physics is the question what 
``functionals" are important. Uniqueness results help to render a functional important and fundamental. 
The classification of {\bf valuations} of fields is classified by {\bf Ostrowski's theorem}
classifying valuations over the rational numbers either being the absolute value or the $p$-adic norm. 
The {\bf Euler characteristic} for example can be characterized as the unique {\bf valuation} 
on simplicial complexes which assumes the value $1$ on simplices or functional which is 
invariant under Barycentric refinements. A theorem of 
Claude Shannon \cite{Shannon48} identifies the {\bf Shannon entropy} is the unique functional 
on probability spaces being compatible with additive and multiplicative operations on 
probability spaces and satisfying some normalization condition. \\
\index{Ostrowski theorem}
\index{p-adic valuation}
\index{valuation of a field}
\index{Hilbert cube}
\index{Polish space}
\index{Euler characteristic}
\index{Shannon entropy}
\index{Classification of finite simple groups}
\index{division algebra}

\section*{Bounds and inequalities} 

An other class of important theorems are
{\bf best bounds} like the {\bf Hurwitz estimate} stating that there are infinitely many $p/q$ for which
$|x-p/q| <1/(\sqrt{5} q^2)$. In packing problems, one wants to find the best packing density, like 
for {\bf sphere packing problems}. In complex analysis, one has the {\bf maximum principle}, which 
assures that a harmonic function $f$ can not have a local maximum in its domain of definition. One can
argue for including this as a fundamental theorem as it is used by other theorems like the 
{\bf Schwarz lemma} (named after Hermann Amandus Schwarz) from complex analysis which is used in
many places. In probability theory or statistical mechanics, one often has
thresholds, where some {\bf phase transition} appears. Computing these values is often important. 
The concept of {\bf maximizing entropy} explains many things like why the Gaussian 
distribution is fundamental as it maximizes entropy. Measures maximizing entropy are often special and
often {\bf equilibrium measures}. This is a central topic in statistical mechanics \cite{RuelleStatMech,
RuelleThermo}.  In combinatorial topology, the {\bf upper bound theorem} 
was a milestone. It was long a conjecture of Peter McMullen
and then proven by Richard Stanley that {\bf cyclic polytopes} maximize the volume in the class of polytopes
with a given number of vertices. 
Fundamental area also some {\bf inequalities} \cite{Inequalities} like the 
{\bf Cauchy-Schwarz inequality} $|a \cdot b| \leq |a| |b|$, 
the {\bf Chebyshev inequality} ${\rm P}[|X-[{\rm E}[X]| \geq |a|] \leq {\rm Var}[X]/a^2$. 
In complex analysis, the {\bf Hadamard three circle theorem} is important as gives bounds 
between the maximum of $|f|$ for a holomorphic function $f$ defined on an annulus given 
by two concentric circles. Often inequalities are more fundamental and powerful than 
equalities because they are more widely used. Related to inequalities are {\bf embedding theorems}
like {\bf Sobolev embedding theorems}.  For more inequalities, 
see \cite{DictionaryInequalities}. 
Apropos embedding, there are the important Whitney or Nash
embedding theorems which are appealing.  
\index{Hurwitz estimate}
\index{Sphere packing}
\index{Schwarz lemma}
\index{Phase transition}
\index{Inequalities}
\index{upper bound theorem}
\index{three circle theorem}
\index{maximum principle}
\index{Chebyshev inequality}
\index{Three circle theorem}
\index{Embedding theorems}

\section*{Big ideas} 

Classifying and valuing {\bf big ideas} is even more difficult than ranking individual
theorems. Examples of big ideas are the idea of {\bf axiomatisation} which stated with 
planar geometry and number theory as described by Euclid and the concept of {\bf proof} 
or later the concept of {\bf models}. 
Archimedes idea of {\bf comparison}, leading to ideas like the {\bf Cavalieri principle}, integral geometry 
or measure theory. Ren\'e Descartes idea of {\bf coordinates} which allowed to work on geometry using algebraic tools,
the use of {\bf infinitesimals and limits} leading to calculus, allowing to merge
concepts of rate of change and accumulation, the idea of {\bf extrema} leading to the calculus
of variations or Lagrangian and Hamiltonian dynamics or descriptions of fundamental forces. Maximizing quantities
like {\bf entropy} lead to fundamental distributions like the Gaussian, exponential, Binomial or uniform distributions. 
{\bf Cantor's set theory} allows for a universal
simple language to cover all of mathematics, the {\bf Klein Erlangen program} of 
``classifying and characterizing geometries through symmetry". The abstract idea of a group
or more general mathematical structures like monoids. The concept of extending {\bf number systems}
like completing the real numbers or extending it to the {\bf quaternions} and {\bf octonions} or then
producing {\bf p-adic number} or {\bf hyperreal numbers}. 
The concept of {\bf complex numbers} or more generally the idea of {\bf completion} of a field.
The idea of {\bf logarithms} \cite{StaudacherBuergi}.
The idea of {\bf Galois} to relate problems about solving equations with {\bf field extensions}
and {\bf symmetries}. The idea of {\bf equivalence classes} is used when looking at projective spaces
or {\bf ideals}. The idea of seeing {\bf prime ideals} as a more fundamental replacement for ``maximal ideal" or 
``point", leading to the notion of {\bf spectrum} of a ring and by gluing to the notion of schemes vastly expanding
classical algebraic geometry.  The {\bf Grothendieck program} of ``geometry without points" or ``locales" as topologies without points 
in order to overcome shortcomings of set theory. This lead to new objects like {\bf schemes} or {\bf topoi}. 
Central in algebra, geometry and number theory is the idea of {\bf localilization} which allows to extend a ring
so that one can start ``dividing", the prototype being the {\bf field of fractions} like the construction of rational functions
from polynomials. An other basic big idea is the concept of {\bf duality}, which appears in many places like 
in projective geometry, in polyhedra, {\bf Poincar\'e duality} or {\bf Pontryagin duality} or {\bf Langlands duality}
for reductive algebraic groups. The idea of {\bf dimension} to measure topological spaces numerically leading to 
{\bf fractal geometry}. The idea of {\bf almost periodicity} is an important generalization
of periodicity. Crossing the boundary of integrability leads to the important paradigm of
stability and randomness \cite{MoserStableRandom} 
and the interplay of structure and randomness \cite{TaoStructureRandomness}.
These themes are related to {\bf harmonic analysis} and {\bf integrability} as integrability means that for every 
invariant measure one has almost periodicity. It is also related to spectral properties in solid state physics or via 
{\bf Koopman theory} in ergodic theory or then to fundamental new number systems like the 
{\bf p-adic numbers}: the {\bf p-adic integers} form a compact topological
group on which the translation is almost periodic. It also leads to problems in 
{\bf Diophantine approximation}. The concept of {\bf algorithm} and building the foundation of computation using
precise mathematical notions. The use of algebra to track problems in topology
starting with mathematicians like Kirchhoff, Betti, Poincar\'e or Emmy N\"other.
An other important principle is to reduce a problem to 
a {\bf fixed point problem}. This often leads to {\bf universality} like for the central limit theorem (where the 
Gaussian distribution is the fixed point). The {\bf categorical approach} is not only a unifying language
but also allows for generalizations of concepts allowing to solve problems. Examples are generalizations of 
Lie groups in the form of {\bf group schemes}. 
Then there is the {\bf deformation idea} which was used for example in the Perelman proof of the
{\bf Poincar\'e conjecture}. Deformation often comes in the form of {\bf partial differential equations}
and in particular heat type equations. Deformations can be abstract in the form of {\bf homotopies} or more concrete
by analyzing concrete partial differential equations like the {\rm mean curvature flow} or {\bf Ricci flow}. 
An other important line of ideas is to use {\bf probability theory} to prove results, even in combinatorics. 
A probabilistic argument can often give existence of objects which one can not even construct. Examples are 
to define a sequence of simplicial complexes $G_n$ with $n$ nodes for which the Euler characteristic $\chi(G_n)= \sum_x (-1)^{{\rm dim}(x)}$
is exponentially large in $n$. The idea of {\bf non-commutative geometry} generalizing geometry through functional analysis or the 
idea of {\bf discretization} which leads to numerical methods or computational geometry. 
The power of coordinates allows to solve geometric problems more easily. 
The above mentioned examples have all proven their use. Grothendieck's ideas have lead to the solution of 
the {\bf Weyl conjectures}, fixed point theorems were used in {\bf Game theory} (first by Nash), or be used
to prove uniqueness of solutions of differential equations. It is also used to justify perturbation theory using
renormalization schemes or iterative methods like in the 
{\bf KAM theorem} about the persistence of quasi-periodic motion leading to {\bf hard implicit function theorems}. 
In the end, what really counts is whether the big idea can solve practical problems or that it can be used to new theorems
(or reprove old theorems more elegantly). The history of mathematics clearly shows that abstraction for the sake
of abstraction or for the sake of generalization rarely was able to convince the mathematical community initially.
But it can also happen that the break-through of a new theory or generalization only pays off much later and that a subtle
generalization actually pushes the tool into a realm where it can be used in other contexts. 
A big idea might have to age like a good wine. 
\index{Grothendieck program}
\index{Poincar\'e conjecture}
\index{doformation idea}
\index{Galois theory}
\index{Koopman theory}
\index{Weyl conjectures}
\index{Deformation idea}
\index{homotopy idea} 
\index{duality}
\index{Klein Erlangen program}
\index{Symmetry}

\section*{Paradigms}

There is once in a while an idea which completely changes the way we look at things. 
These are {\bf paradigm shifts} as described by the philosopher and historian Thomas Kuhn
who relates it also to {\bf scientific revolutions} \cite{KuhnScientificRevolutions}. 
For mathematics, there are various places, where such fundamental changes happened:
the introduction of written numbers which happened independently in various different
places. An early example is the {\bf tally mark notation} on tally sticks 
(early sources are the Lebombo bone from 40 thousand years ago or the Ishango bone 
from 20 thousand years ago) or the technology of {\bf talking knots}, the {\bf khipu} 
\cite{UrtonInkaHistoryKnots}, which is a topological writing which flourished in the
Tawantinsuyu, the Inka empire. 
An other example of a paradigm change is the development of {\bf proof}, which required
the insight that some mathematical statements are assumed as {\bf axioms} from which, using 
{\bf logical deduction}, new theorems are proven. Also {\bf proof assistant frameworks} like
SAM \cite{Huet1973},
ACL2 \cite{ACL2CaseStudies}, Coq \cite{BertotCasteran},
Isabelle \cite{Isabelle2002}, Lean \cite{Hales2018} (extended to Xena in an educational setting)
have emerged allowing to build in more {\bf reliability} and {\bf accountability} to proofs. 
The fact that axiom systems can be {\bf deformed} like from Euclidean to non-Euclidean geometry was definitely 
a paradigm change. On a larger scale, the insight that even the axiom systems of 
mathematics can be deformed and extended in various ways came only in the 20th century with G\"odel. 
Before that, one was under the impression that one could  base all of mathematics on a universal axiom system.
This was Hilbert's program \cite{Zach2007}. 
A third example of a paradigm change is the introduction of the {\bf concept of functions}
which came surprisingly late. The modern concept of a function which takes a quantity
and assigns it a new quantity came only late in the 19'th century with the development
of {\bf set theory}, which is a paradigm change too. There had been a long struggle
also with understanding {\bf limits}, which puzzled already Greek mathematicians like Zeno
but which really only became solid with clear definitions like Weierstrass and then with
the concept of topology where the concept of limit is absorbed within set theory, for example
using the notion of {\bf filters}. Related to functions is the use of functions to understand combinatorial
or number theoretical problems, like through the use of {\bf generating functions}, or 
{\bf Dirichlet series}, allowing analytic tools to solve discrete problems like the existence
of primes on arithmetic progressions. 
The opposite, the use of discrete structures like finite groups to understand the continuum
like Galois theory is an other example of a paradigm change. It led to the insight that 
the quadrature of the circle, or angle trisection can not be done with ruler and compass. 
There are various other places, where paradigm changes happened. A nice example is
the axiomatization of probability theory by Kolmogorov or the realization that {\bf statistics
becomes a geometric theory} if random variables are seen as vectors in a vector space:
the correlation between two random variables is the cosine of the angle between centered versions of these random 
variables. Paradigm changes which are really fundamental can be surprisingly simple. 
An example is the Connes formula \cite{Connes} which is based on the simple idea that distance
can be measured by extremizing slope. This allows to push traditional geometry into non-commutative settings
or discrete settings, where a priory no metric (notion of distance) is given. 
An other example is the extremely simple but powerful idea of the {\bf Grothendieck extension}
of a monoid to a group. It has been used throughout the history of mathematics to generate new number
systems starting with getting integers from natural numbers, rational numbers from integers, complex numbers from real numbers
or quaternions from complex numbers, or the construction of surreal numbers or games generalizing numbers.
The idea is also used in dynamical systems theory to generate
from a not necessarily invertible dynamical system an invertible dynamical system by extending time
from a monoid to a group. In the context of Grothendieck, one should mention also that {\bf category theory}
similarly as set theory at the beginning of the last century changed the way mathematics is done
and extended. Like the switch from {\bf relational data bases} to {\bf graph databases}, it is a paradigm change
stressing more the {\bf relations} (arrows) between objects (nodes) and not only the objects (sets) themselves.

\index{Paradigm}
\index{Revolution}

\section*{Taxonomies} 

When looking at mathematics overall, {\bf taxonomies} are important. They not only help to navigate the
landscape, they are also interesting from a pedagogical as well as historical point of view. I borrow here
some material from my course Math E 320 which is so global that a taxonomy is helpful. Organizing a field
using markers is also important when {\bf teaching intelligent machines}, a field which be seen as the 
{\bf pedagogy for AI}. The big bulk of work in \cite{AIEducation} was to teach a bot mathematics, which 
means to fill in thousands of entries of knowledge. It can appear a bit mind numbing as it is a similar task than
writing a dictionary. But writing things down for a machine actually is even tougher than writing things
down for a student. We can not assume the machine to know anything it is not told. This document about
fundamental theorems by the way could relatively easily be adapted into a database of ``important theorems". 
It actually is one my aims to feed it eventually to the Sofia bot. If the machine is asked about
``important theorem in mathematics", it should be well informed, even so it is just a ``stupid"
encyclopedic data entry. Historically, when knowledge was still sparse, one has
classified teaching material using the {\bf liberal arts of sciences}, the 
{\bf trivium}: grammar, logic and rhetoric, as well as the {\bf quadrivium}: 
arithmetic, geometry, music, and astronomy. More specifically, one has built the {\bf eight ancient roots of 
mathematics} which are tied to activities: counting and sorting (arithmetic),
spacing and distancing  (geometry), positioning and locating  (topology),
surveying and angulating  (trigonometry), balancing and weighing (statics),
moving and hitting  (dynamics), guessing and judging (probability) and
collecting and ordering (algorithms). This leads then to topics like
Arithmetic, Geometry, Number Theory, Algebra, Calculus, Set theory, Probability, Topology, Analysis, 
Numerics, Dynamics and Algorithms. The {\bf AMS classification} is much more refined and distinguishes
64 fields. The Bourbaki point of view is given in \cite{DieudonnePanorama}: it partitions mathematics into
algebraic and differential topology, differential geometry, ordinary differential equations, ergodic theory, 
partial differential equations, non-commutative harmonic analysis, automorphic forms, 
analytic geometry, algebraic geometry, number theory, homological algebra, Lie groups, abstract groups, 
commutative harmonic analysis, logic, probability theory, categories and sheaves, commutative algebra 
and spectral theory. What are {\bf hot spots in mathematics}? Michael Atiyah \cite{atiyah2000}
distinguished parameters like {\bf local - global},  {\bf low and high dimensional},
{\bf commutative - non-commutative}, {\bf linear - nonlinear}, {\bf geometry - algebra},
{\bf physics and mathematics}. 
\index{Taxonomy}
\index{AMS classification}

\section*{Key examples} 

The concept of {\bf experiment} came even earlier and has always been part of mathematics. Experiments
allow to get good examples and set the stage for a theorem. 
\footnote{To quote Vladimir Arnold: ``Mathematics is a part of physics where experiments are cheap"}
Obviously the theorem can not contradict
any of the examples. But examples are more than just a tool to falsify statements;
a good example can be the {\bf seed} for a new theory or for an entire subject. 
Here are a few examples: in {\bf smooth dynamical systems} the {\bf Smale horse shoe} comes to mind,
in {\bf differential topology} the {\bf exotic spheres} of Milnor, in one-dimensional dynamics 
the {\bf logistic map}, or {\bf H\'enon map}, in perturbation theory of Hamiltonian systems the {\bf Standard map} featuring
KAM tori or Mather sets, in homotopy theory the {\bf dunce hat} or {\bf Bing house}, 
in combinatorial topology the {\bf Rudin sphere}, the {\bf Nash-Kuiper non-smooth embedding} of a torus into Euclidean space,
in topology there is the {\bf Alexander horned sphere} or the {\bf Antoine necklace}. In complexity theory there is 
the {\bf busy beaver problem} in Turing computation which is an illustration with how small machines one can achieve great things,
in group theory there is the {\bf Rubik cube} which illustrates many fundamental notions for finitely presented groups, 
in fractal analysis the {\bf Cantor set}, the {\bf Menger sponge}, 
in Fourier theory the series of $f(x)=x \; {\rm mod} \; 1$, 
in Diophantine approximation the {\bf golden ratio}, in the calculus of sums the {\bf zeta function}, in 
dimension theory the {\bf Banach Tarski paradox}. In harmonic analysis the {\bf Weierstrass function} 
as an example of a nowhere differentiable function. The case of {\bf Peano curves} giving concrete examples of
a continuous bijection from an interval to a square or cube. In {\bf complex dynamics} not only 
the {\bf Mandelbrot set} plays an important role, but also individual, specific Julia sets can be interesting.
Examples like the {\bf Mandelbulb} have not yet been investigated mathematically.
In mathematical physics, the {\bf almost Matthieu operator} \cite{Cycon} produced a rich theory related to 
spectral theory, Diophantine approximation, fractal geometry and functional analysis. 
Besides examples illustrating a typical case, it is also important to explore the boundary and limitations 
of a theorem or theory by looking at {\bf counter examples}. Collections of counter examples exist in many fields like
\cite{GelbaumOlmsted,SteenSeebach,RajwadeBhandari,Stoyanov,WiseHall,CapobiancoMolluzzo,Klymchuk}.
\index{Smale horse shoe}
\index{Exotic sphere}
\index{Logistic map}
\index{Rudin sphere}
\index{Alexander sphere}
\index{Golden ratio}
\index{Banach-Tarski paradox}
\index{Menger sponge}
\index{Cantor set}
\index{Nash-Kuiper} 
\index{Diophantie approximation}
\index{KAM tori}
\index{Mather set}
\index{Zeta function}
\index{Mandelbrot set}
\index{Mandelbulb set}
\index{Almost Mathieu operator}

\section*{Physics} 

One can also make a list of great ideas in physics \cite{CenturyIdeas} and see the 
relations with the fundamental theorems in mathematics. A high applicability should then
contribute to a {\bf value functional} in the list of theorems. 
Great ideas in physics are {\bf the concept of space and time}, meaning to describe 
physical events using {\bf differential equations}. In cosmology, one of the insights was
to understand the structure of our solar system and getting for a earth centered to a heliocentric system,
an other is to look at {\bf space-time} as a hole and realize the expansion of the universe or that the idea of a
{\bf big bang}. More general is the Platonic idea that {\bf physics is geometry}. Or calculus:
Lagrange developed his {\bf calculus of variations} to find laws of physics. Then there is the idea of
{\bf Lorentz invariance} and symmetries more general which leads to {\bf special relativity}, there
is the idea of {\bf general relativity} which allows to describe gravity through 
geometry and a larger symmetry seen through the {\bf equivalence principle}. 
There is the idea of see elementary particles using {\bf Lie groups}.
There is the {\bf Noether theorem} which is the idea that any {\bf symmetry} is tied
to a {\bf conservation law}: translation symmetry leads to momentum conservation,
rotation symmetry to angular momentum conservation for example. Symmetries also play a role
when {\bf spontaneous broken symmetry} or {\bf phase transitions}. 
There is the idea of quantum mechanics which mathematically 
means replacing differential equations with {\bf partial differential equations} or
replacing commutative algebras of observables with {\bf non-commutative algebras}. 
An important idea is the concept of {\bf perturbation theory} and in particular the
notion of {\bf linearization}. Many laws are simplifications of more complicated laws
and described in the simplest cases through linear laws like Ohms law or Hooks law. 
{\bf Quantization processes} allow to go from commutative to non-commutative structures. 
{\bf Perturbation theory} allows then to extrapolate from a simple law to a more complicated
law. Some is easy application of the {\bf implicit function theorem}, some is harder like
KAM theory. There is the idea of using {\bf discrete mathematics} to describe complicated 
processes. An example is the language of {\bf Feynman graphs} or the language of graph theory
in general to describe physics as in loop quantum gravity or then the language of 
{\bf cellular automata} which can be seen as partial difference equations where also the
function space is quantized. The idea of {\bf quantization}, a formal transition from an
ordinary differential equation like a Hamiltonian system to a partial differential equation
or to replace single particle systems with infinite particle systems (Fock). There are other
quantization approaches through {\bf deformation of algebras} which is related to 
{\bf non-commutative geometry}. There is the idea of using {\bf smooth functions} to 
describe discrete particle processes. An example is the {\bf Vlasov dynamical system} or 
{\bf Boltzmann's equation}  to describe a plasma, or thermodynamic notions to describe large sets of particles like a 
gas or fluid. Dual to this is the use of {\bf discretization} to describe a smooth system 
by discrete processes. An example is {\bf numerical approximation}, like using the Runge-Kutta 
scheme to compute the trajectory of a differential equation. 
There is the realization that we have a whole spectrum of dynamical systems, 
{\bf integrability} and {\bf chaos} and that some of the transitions are 
{\bf universal}. An other example is the {\bf tight binding approximation} in which 
a continuum Schr\"odinger equation is replaced with a bounded {\bf discrete Jacobi operator}.
There is the general idea of finding the {\bf building blocks} or {\bf elementary particles}.
Starting with Demokrit in ancient Greece, the idea got refined again and again. Once, atoms
were detected and charges found to be quantized (Robert Millikan), the structure of the atom was
explored (Rutherford), and then the atom got split (Lisa Meitner, Otto Hahn). The structure of the 
nuclei with protons and neutrons was then refined again using quarks leading the {\bf standard model
in particle physics}. There is furthermore the idea to use statistical methods for complex systems. 
An example is the use of stochastic
differential equations like diffusion processes to describe actually deterministic particle
systems. There is the insight that complicated systems can form {\bf patterns} through
interplay between symmetry, conservation laws and {\bf synchronization}. Large scale patterns
can be formed from systems with local laws. Finally, there is
the idea of solving {\bf inverse problems} using mathematical tools like Fourier theory
or basic geometry (Eratostenes could compute the radius of the earth by comparing the lengths of shadows
at different places of the earth.)
An example is {\bf tomography}, where the structure of some object is explored using 
{\bf resonance} and where the reconstruction solves an {\bf inverse problem}. 
Then there is the idea of {\bf scale invariance} which allows to describe
objects which have {\bf fractal nature}.  \\
\index{conservation law}
\index{Noether theorem}
\index{calculus of variations}
\index{tight binding approximation}
\index{numerical methods}
\index{Standard model}
\index{Eratostenes}
\index{Demokrit}
\index{patterns}
\index{universality}
\index{Runge Kutta} 
\index{value function} 
\index{synchronization}
\index{quantisation}
\index{tomography}
\index{non-commutative geometry}

\section*{Computer science} 

As in physics, it is harder to pinpoint ``big ideas" in computer science as they are in general not 
theorems. But it has been done \cite{KnuthComputerProgramming}.
The initial steps of mathematics was to build a {\bf language}, where {\bf numbers}
represent quantities \cite{conway_guy}. Physical tools which assist in manipulating numbers can
already been seen as a {\bf computing device}. Marks on a bone, pebbles in a clay bag, talking knots in a
Khipu \cite{UrtonInkaHistoryKnots,AscherAscher}, marks on a Clay tablet were the first step. 
Papyri, paper, magnetic, optical and electric
storage, the tools to build {\bf memory} were refined over millenniums. The mathematical language
allowed us to explore topics beyond the finite and also build {\bf data bases}. The Khipu concept was
already an early form of graph database \cite{AggarwalWangGraphData}.
Using a finite number of symbols we can represent and count
infinite sets, have notions of {\bf cardinality}, have {\bf various number systems} and more generally
have {\bf algebraic structures}. Numbers can even be seen as {\bf games} \cite{numbersgames,knuthnumbers}.
A major idea is the concept of an {\bf algorithm}. Adding or multiplying on an {\bf abacus} already
was an algorithm. The concept was refined in geometry, where {\bf ruler and compass} were used as
{\bf computing devices}, like the construction of points in a triangle. 
To measure the effectiveness of an algorithm, one can use notions of {\bf complexity}. 
This has been made precise by computing pioneers like Alan Turing,
as one has to formulate first what a ``computation" is. The concept of the {\bf Turing machine} is 
particularly elegant as it is both a theoretical construct as well as a concrete machine (although 
extremely inefficient). 
In the last century one has seen that computations and
proofs are very similar and that they have similar general restrictions. There are
some tasks which can not be computed with a Turing machine and there are theorems
which can not be proven in a specific axiom system. 
As mathematics is a language, we have to deal with concepts of {\bf syntax}, {\bf grammar}, 
{\bf notation}, {\bf context}, {\bf parsing}, {\bf validation}, {\bf verification}. 
As Mathematics is a {\bf human activity} which is done in our {\bf brains}, it is related to 
psychology and {\bf computer architecture}. 
Computer science aspects are also important also in {\bf pedagogy} and {\bf education} 
how can an idea be communicated {\bf clearly}?
How do we {\bf motivate}? How do we {\bf convince} peers that a result is true? Examples from 
history show that this is often done by {\bf authority} and that the validity of some proofs turned
out to be wrong or incomplete, even in the case of fundamental theorems or when treated by 
great mathematicians. (Examples are the fundamental theorem of arithmetic, the fundamental theorem of
algebra or the wrong published proof of Kempe of the 4 color theorem). On the other hand, there were
also quite many results which only later got recognized. The work of Galois for example only exploded
much later. How come we trust a human brain more than an electronic one? We have to make
some fundamental assumptions for example to be made like that if we do a logical step "if A and B
then ``A and B" holds. This assumes for example that {\bf our memory is faithful}: after having put A and
B in the memory and making the conclusion, we have to assume that we did not forget A nor B! 
Why do we trust this more than the memory of a machine?
As we are also assisted more and more by electronic devices, the question of the validity of
{\bf computer assisted proofs} comes up. The {\bf 4-color theorem} of Kenneth Appel and Wolfgang 
Haken based on previous work of many others like Heinrich Heesch or the proof of the {\bf Feigenbaum
conjecture} of Mitchell Feigenbaum first proven by Oscar Lanford III 
or the proof of the Kepler problem given by Thomas Hales are examples. 
A great general idea is related to the representation of {\bf data}. This can be done using
matrices like in a {\bf relational database} or using other structures like {\bf graphs} leading
to {\bf graph databases}. The ability to use computers allows mathematicians to do {\bf experiments}. A branch
of mathematics called {\bf experimental mathematics} \cite{ArnoldExperimental,BBG}
relies heavily on experiments to find new
theorems or relations. Experiments are related to {\bf simulations}. We are able, within a computer to build
and explore new worlds, like in {\bf computer games}, we can enhance the physical world using {\bf virtual reality} or 
{\bf augmented reality} or then {\bf capturing a world} by {\bf 3D scanning} and 
{\bf realize} a world by {\bf printing the objects} \cite{CFZ}. 
A major theme is {\bf artificial intelligence} \cite{RussellNorvig,JacksonAI}. It is related to 
optimization problems like optimal transport, neural nets as well as {\bf inverse problems} like 
{\bf structure from motion problems}. An intelligent entity must be able to take information, build
a model and then find an optimal strategy to solve a given task. A self-driving car for example has
to be able to translate pictures from a camera and build a map, then determine where to drive.
Such tasks are usually considered part of {\bf applied mathematics} but they are 
very much related with pure mathematics because computers also start to learn how to read mathematics, 
how to {\bf verify proofs} and to {\bf find new theorems}. 
Artificial intelligent agents \cite{Wei65} were first developed in the 1960ies learned also some mathematics. 
I myself learned about it when incorporated computer algebra systems into a chatbots in \cite{AIEducation}. 
AI has now become a big business as {\bf Alexa}, {\bf Siri}, {\bf Google Home}, {\bf IBM Watson} or {\bf Cortana} 
demonstrate. But these information systems must be taught, they must be able to rank alternative answers,
even inject some humor or opinions. Soon, they will be able to learn themselves and 
answer questions like ``what are the 10 most important theorems in mathematics?"
\index{graph database}
\index{memory}
\index{notation}
\index{context}
\index{pedagogy} 
\index{algorithm}
\index{inverse problem}
\index{3D scanning}
\index{Khipu}
\index{artificial intelligence}

\section*{Brevity}

We live in a instagram, snapchat, twitter, microblog, vine, tiktok, watch-mojo,
petcha-kutcha time. Many of us multi task, read news on
smart phones, watch faster paced movies, read shorter novels and feel that a million word
Marcel Proust's masterpiece ``a la recherche du temps perdu" is ``temps perdu". Even
classrooms and seminars have become more aphoristic.
Micro blogging tools are only the latest incarnation of ``miniature stories".
They continue the tradition of older formats like "mural art" by Romans
to modern graffiti or ``aphorisms" 
\cite{KrantzApocrypha1,KrantzApocrypha2}),
poetry, cartoons, Unix fortune cookies \cite{KenArnold}.
Shortness has appeal: aphorisms, poems, ferry tales,
quotes, words of wisdom, life hacker lists, and tabloid top 10 lists
illustrate this. And then there are books like ``Math in 5 minutes", ``30 second math",
``math in minutes" \cite{BehrendsFuenfMinutenMath,MathInMinutes,30SecondMath},
which are great coffee table additions. Also short proofs are appealing like ``Let epsilon be smaller than zero"
which is the shortest known math joke, or ``There are three type of mathematicians,
the ones who can count, and the ones who can't." Also short open problems are attractive, like
the {\bf twin prime problem} ``there are infinitely many twin primes" or 
the {\bf Landau problem} ``there are infinitely many primes of the form $n^2+1$, or the {\bf Goldbach} problem
``every $n>2$ is the sum of two primes". 
\index{Prime twin conjecture}
\index{Goldbach conjecture}
\index{Landau problem}
For the larger public in mathematics shortness has appeal:
according to a poll of the Mathematical Intelligencer from 1988, the most favorite theorems
are short \cite{Wells1,Wells2}. Results with longer proofs can make it
to graduate courses or specialized textbooks but still then, the results are often
short enough so that they can be tweeted without proof.
Why is shortness attractive? Paul Erd\"os expressed short elegant proofs as ``proofs from the book" 
\cite{AigZie}.  Shortness reduces the possibility of error as complexity is always a stumbling block
for understanding. But is beauty equivalent to brevity?
Buckminster Fuller once said: ``If the solution is not beautiful, I know it is wrong."
\cite{AharoniBeauty}. Much about the aesthetics in mathematics is investigated
in \cite{MontanoExplainingBeauty}.
According to \cite{RotaBeauty}, the beauty of a piece of mathematics is frequently
associated with the shortness of statement or of proof:
{\it beautiful theories are also thought of as short, self-contained chapters fitting
within broader theories. There are examples of complex and extensive theories which
every mathematician agrees to be beautiful, but these examples are not the one
which come to mind}.
Also psychologists and educators know that simplicity appeals to children: From \cite{SinclairMathBeauty}
{\it For now, I want simply to draw attention to the fact that even for a young, mathematically naive child,
aesthetic sensibilities and values (a penchant for simplicity, for finding the building blocks of more
complex ideas, and a preference for shortcuts and ``liberating" tricks rather than cumbersome recipes)
animates mathematical experience.} It is hard to exhaust them all, even not with tweets: there are  more than 
${\rm googool}^2 = 10^{200}$ texts of length 140. This can not
all ever be written down because there are more than what we estimate the number of elementary particles. 
But there are even short story collections.
Berry's paradox tells in this context that the shortest non-tweetable text
in 140 characters can be tweeted: "The shortest non-tweetable text".
Since we insist on giving proofs, we have to cut corners. Books containing lots of elegant
examples are \cite{AlsinaNelson,AigZie}.
We should add that brevity is not a new thing. J.E. Littlewood has raised the question how 
short a dissertation can be and proves in an example, that two sentences are enough 
and gives a one-sentence proof of the fact that bounded entire functions are constant 
by using Cauchy's integral theorem. It has been refined a bit in \cite{Zagier2014}.

\section*{Twitter math}

The following 42 tweets were written in 2014, when twitter
still had a 140 character limit. Some of them were actually tweeted.
The experiment was to see which theorems are short enough so that
one can tweet both the theorem as well as the proof in 140 characters.
Of course, that often required a bit of cheating. See \cite{AigZie}
for proofs from the books, where the proofs have full details. 

\tweet{
{\bf Euclid:} The set of primes is infinite.
Proof: let $p$ be largest prime, then $p!+1$ has
a larger prime factor than $p$. Contradiction. }

\tweet{
{\bf Euclid:} $2^p-1$ prime then $2^{p-1} (2^p-1)$ is perfect.
Proof. $\sigma(n)=$ sum of factors of $n$, 
$\sigma(2^n-1) 2^{n-1})  = \sigma(2^n-1) \sigma(2^{n-1}) = 2^n (2^n-1) = 2 \cdot 2^n (2^n-1)$
shows $\sigma(k) = 2k$.}

\tweet{
{\bf Hippasus:} $\sqrt{2}$ is irrational.
Proof. If $\sqrt{2}=p/q$, then $2 q^2 = p^2$.
To the left is an odd number of factors 2,
to the right it is even one. Contradiction.
}

\tweet{
{\bf Pythagorean triples:} all $x^2+y^2=z^2$ are of form $(x,y,z) =(2st,s^2-t^2,s^2+t^2)$.
Proof: $x$ or $y$ is even (both odd gives $x^2+y^2=w^k$ with odd $k$).
Say $x^2$ is even: write $x^2= z^2-y^2 = (z-y) (z+y)$.
This is $4s^2 t^2$. Therefore $2s^2 = z-y, 2t^2 = z+y$. Solve for $z,y$.}

\tweet{
{\bf Pigeon principle:} if $n+1$ pigeons live in $n$ boxes, 
there is a box with 2 or more pigeons. 
Proof: place a pigeon in each box until every box is filled. 
The pigeon left must have a roommate. }

\tweet{
{\bf Angle sum in triangle:} $\alpha + \beta + \gamma = K A + \pi$ 
if $K$ is curvature, $A$ triangle area. Proof: Gauss-Bonnet for
surface with boundary. $\alpha,\beta,\gamma$ are Dirac measures
on the boundary. 
}

\tweet{
{\bf Chinese remainder theorem:} a(i) x = b(i) mod n(i)
has a solution if gcd(a(i),n(i))=0 and gcd(n(i),n(j))=0
Proof: solve eq(1), then increment x by n(1) to solve eq(2),
then increment x by n(1) n(2) until second is ok. etc. 
}

\tweet{
{\bf Nullstellensatz:} algebraic sets in $K^n$ 
are 1:1 to radical ideals in $K[x_1...x_n]$. 
Proof: An algebra over K which is a field is finite 
field extension of K.
}

\tweet{
{\bf Fundamental theorem algebra:} a polynomial of degree $n$ has exactly $n$ roots.
Proof: the metric $g=|f|^{-2/n} |dz|^2$ on the Riemann sphere has curvature 
$K=n^{-1} \Delta \log|f|$. Without root, K=0 everywhere contradicting
Gauss-Bonnet. \cite{AlmiraRomero}:}

\tweet{
{\bf Fermat:} p prime $(a,p)=1$, then $p| a^p-a$
Proof: induction with respect to $a$. Case $a=1$ is trivial
$(a+1)^p -(a+1)$ is congruent to $a^p - a$  modulo  $p$ 
because Binomial coefficients $B(p,k)$ are divisible by 
$p$ for $k=1,\dots p-1$. }

\tweet{
{\bf Wilson:} $p$ is prime iff $p | (p-1)!+1$
Proof. Group $2,\dots p-2$ into pairs $(a,a^-1)$
whose product is $1$ modulo $p$. Now $(p-1)! = (p-1)=-1$ modulo $p$. 
If $p=ab$ is not prime, then $(p-1)! = 0$ modulo $p$ and $p$ does 
not divide $(p-1)!+1$.  }

\tweet{
{\bf Bayes:} $A,B$ are events and $A^c$ is the complement.
$P[A|B] = P[B|A] P[A]/(P[B|A] P[A] + P[B|A^c] P[A^c]$
Proof: By definition $P[A|B] P[B] =P[A \cap B]$.
Also $P[B] = (P[B|A] P[A] + P[B|A^c] P[A^c]$.
}

\tweet{
{\bf Archimedes:} Volume of sphere $S(r)$ is $4 \pi r^3/3$
Proof: the complement of the cone inside the cylinder
has at height $z$ the cross section area $r^2-z^2$, the
same as the cross section area of the sphere at height $z$. 
}

\tweet{
{\bf Archimedes:} the area of the sphere $S(r)$ is $4 \pi r^2$
Proof: differentiate the volume formula with respect 
to $r$ or project the sphere onto a cylinder of height 
$2$ and circumference $2\pi$ and not that this is area 
preserving.  }

\tweet{
{\bf Cauchy-Schwarz:} $|v \cdot w| \leq |v| |w|$.
Proof: scale to get $|w|=1$, define $a=v.w$, so that
$0 \leq (v-a w) \dots (v-a w) =|v|^2-a^2 
= |v|^2 |w|^2 - (v \cdot w)^2$.
}

\tweet{
{\bf Angle formula:} Cauchy-Schwarz defines the angle 
between two vectors as $\cos(A) = v.w/|v| |w|$. 
If $v,w$ are centered random variables, then  
$v \cdot w$ is the covariance, $|v|,|w|$ are standard
deviations  and $\cos(A)$ is the correlation. }

\tweet{
{\bf Cos formula:} $c^2 = a^2+b^2 - a b \cos(A)$ in 
a triangle ABC  (Al-Kashi theorem)
Proof: $v=AB, w=AC$ has length $a=|v|,b=|w|, |c|=|v-w|$.
Now: $(v-w).(v-w) = |v|^2 + |w|^2 - 2 |v| |w| \cos(A)$.
}

\tweet{
{\bf Pythagoras:} $A=\pi/2$, then $c^2=a^2+b^2$.
Proof: Let $v=AB$, $w=AC$, $v-w =BC$ be the sides
of the triangle. Multiply out
$(v-w) \cdot (v-w) = |v|^2 + |w|^2$ and use $v\cdot w=0$.
}

\tweet{
{\bf Euler formula:} $\exp(i x) = \cos(x) + i \sin(x)$.
Proof: $\exp(i x) = 1+(ix) + (ix)^2/2! - ...$
Pair real and imaginary parts and use 
definition $\cos(x) = 1-x^2/2!+x^4/4! ...$
and $\sin(x) = x-x^3/3!+x^5/5!-....$. 
}

\tweet{
{\bf Discrete Gauss-Bonnet} $\sum_x K(x) = \chi(G)$ with
$K(x) = 1-V_0(x)/2 + V_1(x)/3 + V_2(x)/4 ...$ curvature
$\chi(G) = v_0-v_1+v_2-v_3...$ Euler characteristic
Proof: Use handshake $\sum_x V_k(x) = v_{k+1}/(k+2)$. }

\tweet{
{\bf Poincar\'e-Hopf:} let $f$ be a coloring, $i_f(x)=1-\chi(S^-_f(x))$,
where $S^-_f(x)={y \in S(x) | f(y)<f(x) }$
$\sum i_f(x) = \chi(G)$.  Proof by induction. Removing local maximum of $f$
reduces Euler characteristic by $\chi(B_f(x))-\chi(S^-f(x))=i_f(x)$. }

\tweet{
{\bf Lefschetz:} $\sum_x i_T(x) = {\rm str}(T|H(G))$. Proof:
LHS is ${\rm str}(\exp(-0 L) U_T)$ and RHS is ${\rm str}(exp(-t L) U_T)$
for $t \to \infty$. The super trace does not depend on $t$. 
}

\tweet{
{\bf Stokes:} orient edges $E$ of graph $G$. 
$F: E \to R$ function, $S$ surface in G with boundary C. 
$d(F)(ijk) = F(ij)+F(jk)-F(ki)$ is the curl. 
The sum of the curls over all triangles is
the line integral of $F$ along $C$.
}

\tweet{
{\bf Plato:} there are exactly $5$ platonic solids. 
Proof: number $f$ of $n$-gon satisfies $f=2e/n$,
$v$ vertices of degree $m$ satisfy $v=2e/m$
$v-e+f-2$ means $2e/m - e + 2e/n = 2$ or 
$1/m + 1/n = 1/e+1/2$ with solutions: $(m=4,n=3),(m=3,n=5),
(n=m=3),(n=3,m=5),(m=3,n=4)$. }

\tweet{
{\bf Poincar\'e recurrence:} $T$ area-preserving map
of probability space $(X,m)$. If $m(A)>0$ and
$n>1/m(A)$ we have $m(T^k(A) \cap A)>0$ for some $1\leq k \leq n$
Proof. Otherwise $A,T(A),...,T^n(A)$ are all disjoint and
the union has measure $n \cdot m(A)>1$. }

\tweet{
{\bf Turing:} there is no Turing machine which halts if input
is Turing machine which halts: Proof: otherwise
build an other one which halts if the input is a non-halting
one and does not halt if input is a halting one. }

\tweet{
{\bf Cantor:} the set of reals in [0,1] is uncountable. 
Proof: if there is an enumeration $x(k)$, let $x(k,l)$ be 
the $l$'th digit of $x(k)$ in binary form. The number
with binary expansion $y(k)=x(k,k)+1$ mod $2$ is not in the list. }

\tweet{
{\bf Niven:} $\pi \notin Q$: Proof: $\pi=a/b$, $f(x)=x^n (a-b x)^n/n!$
satisfies $f(pi-x)=f(x)$ and $0<f(x)<\pi^n a^n/n^n$
$f^(j)(x)=0$ at $0$ and $\pi$ for $0 \leq j \leq n$ shows
$F(x)=f(x)-f^{(2)}(x)+f^{(4)}(x) \cdots +(-1)^n f^{(2n)}(x)$
has $F(0),F(\pi) \in Z$ and $F+F''=f$.
Now $(F'(x) \sin(x) - F(x) \cos(x)) = f \sin(x)$, so
$\int_0^\pi f(x) \sin(x) dx \in Z$. }

\tweet{
{\bf Fundamental theorem calculus:}
With differentiation $Df(x)=f(x+1)-f(x)$ and integration
$Sf(x)=f(0)+f(1)+...+f(n-1)$ have
$SDf(x)=f(x)-f(0), DS f(x)=f(x)$. }

\tweet{
{\bf Taylor:} $f(x+t) = \sum_k f^{(k)}(x) t^k/k!$. 
Proof: $f(x+t)$ satisfies transport equation $f_t=f_x = Df$
an ODE for the differential operator $D$. Solve 
$f(x+t)=\exp(D t) f(x)$. 
}

\tweet{
{\bf Cauchy-Binet:} $\det(1+F^T G) = \sum_P \det(F_P) \det(G_P)$
Proof: $A=F^TG$. Coefficients of $\det(x-A)$
is $\sum_{|P|=k} \det(F_P) \det(G_P)$.
}

\tweet{
{\bf Intermediate:} $f$ continuous 
$f(0)<0,f(1)>0$, then there exists $0<x<1, f(x)=0$.
Proof. If $f(1/2)<0$ do proof with $(1/2,1)$
If $f(1/2)>0$ redo proof with $(0,1/2)$. 
}

\tweet{
{\bf Ergodicity:} $T(x)=x+a$ mod 1 with irrational a
is ergodic. Proof. $f=\sum_n a(n) \exp(i n x)$
$T f = \sum_n a(n) \exp(i n a) \exp(i n x) = f$
implies $a(n)=0$. 
}

\tweet{
{\bf Benford:} first digit $k$ of $2^n$ appears with probability $\log(1-1/k)$
Proof: $T: x \to x+\log(2) \; {\rm mod} \; 1$ is ergodic. $\log(2^n) {\rm mod} \; 1 =k$ if
$\log(k) \leq T^n(0)<\log(k+1)$. The probability of hitting this interval 
is $\log(k+1)/\log(k)$.
}

\tweet{
{\bf Rank-Nullity:} ${\rm dim}({\rm ker}(A)) + {\rm dim}({\rm im}(A)) = n$ for $m \times n$ matrix $A$. 
Proof: a column has a leading $1$ in $rref(A)$ or no leading $1$. In the first
case it contributes to the image, in the second to a free variable parametrizing
the kernel. 
}

\tweet{
{\bf Column-Row picture:} $A: R^m \to R^n$. The $k$'th column of $A$
is the image $A e_k$. If all rows of $A$ are perpendicular to $x$
then $x$ is in the kernel of $A$.
}

\tweet{
{\bf Picard:} $x'=f(x), x(0)=x_0$ has locally a unique solution if $f \in C^1$.
Proof: the map $T(y) = \int_0^t f(y(s)) \; ds$ is a contraction on $C([0,a])$ 
for small enough $a>0$. Banach fixed point theorem.}

\tweet{
{\bf Banach:} a contraction $d(T(x),T(y)) \leq a d(x,y)$ 
on complete $(X,d)$ has a unique fixed point.
Proof: $d(x_k,x_n) \leq a^k/(1-a)$ using triangle 
inequality and geometric series. Have Cauchy sequence. }

\tweet{
{\bf Liouville:} every prime p=4k+1 is the sum of two squares.
Proof: there is an involution on $S={ (x,y,z) | x^2+4yz=p }$
with exactly one fixed point showing |S| is odd implying
$(x,y,z) -> (x,z,y)$ has a fixed point. \cite{Zagier90}
}

\tweet{
{\bf Banach-Tarski:}
The unit ball in $R^3$ can be cut into 5 pieces,
re-assembled using rotation and translation to
get two spheres. Proof: cut cleverly using axiom
of choice.
}

\section*{Math areas}

We add here the core handouts of Math E320 which aimed to give for each of the
12 mathematical subjects an overview on two pages. For that course, I had
recommended books like \cite{Eves,guinness,MathThroughAges,Stillwell2010,Stillwell2016}. \\

\fontsize{10}{14} \selectfont
\pagebreak {\bf E-320: Teaching Math with a Historical Perspective \hfill O. Knill, 2010-2018}

\chapter{Lecture 1: Mathematical roots}

Similarly, as one has distinguished the {\bf canons of rhetorics}: memory, invention, delivery, style, and arrangement, 
or combined the {\bf trivium}: grammar, logic and rhetorics,
with the {\bf quadrivium}: arithmetic, geometry, music, and astronomy, to obtain the seven 
{\bf liberal arts and sciences}, one has tried to {\bf organize all mathematical activities}. 
\index{trivium}
\index{quadrivium}
\index{canons of rhetorik}
\index{liberal arts and sciences}

\parbox{16.8cm}{
\parbox{5cm}{
Historically, one has distinguished
{\bf eight ancient roots of mathematics}. Each of these 8 activities in turn
suggest a key area in mathematics: 
}
\hspace{5mm}
\fcolorbox{yellow2}{yellow2}{\parbox{11cm}{
\begin{center}
\begin{tabular}{l|l}
  counting and sorting      & {\bf arithmetic} \\
  spacing and distancing    & {\bf geometry} \\
  positioning and locating  & {\bf topology} \\
  surveying and angulating  & {\bf trigonometry} \\
  balancing and weighing    & {\bf statics}   \\
  moving and hitting        & {\bf dynamics} \\
  guessing and judging      & {\bf probability} \\
  collecting and ordering   & {\bf algorithms} \\
\end{tabular}
\end{center}
}}
}

To morph these 8 roots to the 12 mathematical areas covered in this class, we complemented the ancient 
roots with calculus, numerics and computer science, merge trigonometry with geometry, separate arithmetic 
into number theory, algebra and arithmetic and turn statics into analysis.
\index{mathematical roots}
\index{ancient roots}

\parbox{16.8cm}{
\parbox{5cm}{
Let us call this modern adaptation the  \\

\begin{center} \parbox{4cm}{
{\bf 12 modern roots of Mathematics}: 
} \end{center}
}
\hspace{5mm}
\fcolorbox{yellow2}{yellow2}{\parbox{11cm}{
\begin{center}
\begin{tabular}{l|l}
  counting and sorting      & {\bf arithmetic} \\
  spacing and distancing    & {\bf geometry} \\
  positioning and locating  & {\bf topology} \\
  dividing and comparing    & {\bf number theory}\\
  balancing and weighing    & {\bf analysis}  \\
  moving and hitting        & {\bf dynamics} \\
  guessing and judging      & {\bf probability} \\
  collecting and ordering   & {\bf algorithms} \\
  slicing and stacking      & {\bf calculus} \\
  operating and memorizing  & {\bf computer science} \\
  optimizing and planning   & {\bf numerics} \\
  manipulating and solving  & {\bf algebra} \\
\end{tabular}
\end{center}
}}
}

\vspace{2mm}

\parbox{16.8cm}{
\parbox{4cm}{
While relating {\bf mathematical areas} with {\bf human activities} is useful, 
it makes sense to select specific topics in each of this area.
These 12 topics will be the 12 lectures of this course. 
}
\hspace{5mm}
\fcolorbox{yellow2}{yellow2}{\parbox{12cm}{
\begin{center}
\begin{tabular}{l|l}
Arithmetic      &    numbers and number systems \\
Geometry        &    invariance, symmetries, measurement, maps \\
Number theory   &    Diophantine equations, factorizations  \\
Algebra         &    algebraic and discrete structures \\
Calculus        &    limits, derivatives, integrals \\
Set Theory      &    set theory, foundations and formalisms \\
Probability     &    combinatorics, measure theory and statistics \\
Topology        &    polyhedra, topological spaces, manifolds \\
Analysis        &    extrema, estimates, variation, measure   \\
Numerics        &    numerical schemes, codes, cryptology \\
Dynamics        &    differential equations, maps \\
Algorithms      &    computer science, artificial intelligence
\end{tabular}
\end{center}
}}
}

Like any classification, this chosen division is rather arbitrary and a matter of 
personal preferences. The {\bf 2010 AMS classification} distinguishes 64 areas of mathematics.
Many of the just defined main areas are broken off into even finer pieces. Additionally,
there are fields which relate with other areas of science, like economics, biology or
physics:a
\index{AMS classification}

\begin{center}
\begin{tiny}
\parbox{16.8cm}{
\fcolorbox{yellow2}{yellow2}{
\parbox{8cm}{
00 General \\
01 History and biography  \\
03 Mathematical logic and foundations \\
05 Combinatorics  \\
06 Lattices, ordered algebraic structures \\
08 General algebraic systems \\
11 Number theory \\
12 Field theory and polynomials \\
13 Commutative rings and algebras \\
14 Algebraic geometry \\
15 Linear/multi-linear algebra; matrix theory \\
16 Associative rings and algebras  \\
17 Non-associative rings and algebras \\
18 Category theory, homological algebra  \\
19 K-theory  \\
20 Group theory and generalizations \\
}
}
\hspace{5mm}
\fcolorbox{yellow2}{yellow2}{
\parbox{8cm}{
45 Integral equations \\
46 Functional analysis  \\
47 Operator theory \\
49 Calculus of variations, optimization  \\
51 Geometry  \\
52 Convex and discrete geometry \\
53 Differential geometry  \\
54 General topology  \\
55 Algebraic topology \\
57 Manifolds and cell complexes  \\
58 Global analysis, analysis on manifolds  \\
60 Probability theory and stochastic processes  \\
62 Statistics \\
65 Numerical analysis \\
68 Computer science  \\
70 Mechanics of particles and systems \\
}
}
}
\parbox{16.8cm}{
\fcolorbox{yellow2}{yellow2}{
\parbox{8cm}{
22 Topological groups, Lie groups  \\
26 Real functions  \\
28 Measure and integration  \\
30 Functions of a complex variable  \\
31 Potential theory  \\
32 Several complex variables, analytic spaces  \\
33 Special functions  \\
34 Ordinary differential equations \\
35 Partial differential equations \\
37 Dynamical systems and ergodic theory  \\
39 Difference and functional equations \\
40 Sequences, series, summability \\
41 Approximations and expansions  \\
42 Fourier analysis \\
43 Abstract harmonic analysis  \\
44 Integral transforms, operational calculus \\
}
}
\hspace{5mm}
\fcolorbox{yellow2}{yellow2}{
\parbox{8cm}{
74 Mechanics of deformable solids \\
76 Fluid mechanics  \\
78 Optics, electromagnetic theory  \\
80 Classical thermodynamics, heat transfer  \\
81 Quantum theory \\
82 Statistical mechanics, structure of matter \\
83 Relativity and gravitational theory \\
85 Astronomy and astrophysics  \\
86 Geophysics  \\
90 Operations research, math. programming \\
91 Game theory, Economics Social and Behavioral Sciences \\
92 Biology and other natural sciences \\
93 Systems theory and control  \\
94 Information and communication, circuits \\
97 Mathematics education  \\
\vspace{2mm}
}
}
}
\end{tiny}
\end{center}

\parbox{16.8cm}{
\parbox{6cm}{
What are \\
\begin{center}
{\bf fancy developments}  \\
\end{center}
in mathematics today? Michael Atiyah 
\cite{atiyah2000} 
identified in the year 2000 the following
{\bf six hot spots}: \\
}
\hspace{1cm}
\fcolorbox{red1}{red1}{
\parbox{9cm}{
\begin{center}
\begin{tabular}{|lcr|}
local       & and & global \\
low         & and & high dimension \\
commutative & and & non-commutative \\
linear      & and & nonlinear \\
geometry    & and & algebra \\
physics     & and & mathematics \\
\end{tabular}
\end{center}
}
}
}

\vspace{5mm}
Also this choice is of course highly personal. One can easily add 12 other 
{\bf polarizing} quantities which help to distinguish or parametrize different parts of 
mathematical areas, especially the ambivalent pairs which produce a captivating gradient: 

\parbox{16.8cm}{
\fcolorbox{red1}{red1}{
\parbox{8cm}{
\begin{center}
\begin{tabular}{|lcr|}
regularity   & and &  randomness    \\
integrable   & and &  non-integrable \\
invariants   & and &  perturbations \\
experimental & and &  deductive     \\
polynomial   & and &  exponential   \\
applied      & and &  abstract
\end{tabular}
\end{center}
}}
\fcolorbox{red1}{red1}{
\parbox{8cm}{
\begin{center}
\begin{tabular}{|lcr|}
discrete    & and &  continuous   \\
existence   & and &  construction \\
finite dim  & and &  infinite dimensional \\
topological & and &  differential geometric \\
practical   & and &  theoretical \\
axiomatic   & and &  case based   \\
\end{tabular}
\end{center} 
}
}
}

The goal is to illustrate some of these structures from a historical point of view
and show that ``Mathematics is the science of structure".

 \pagebreak
\pagebreak {\bf E-320: Teaching Math with a Historical Perspective \hfill Oliver Knill, 2010-2018}

\chapter{Lecture 2: Arithmetic}
The oldest mathematical discipline is {\bf arithmetic}. It is the theory of the construction and manipulation of numbers.
The earliest steps were done by {\bf Babylonian}, {\bf Egyptian}, {\bf Chinese}, {\bf Indian} and {\bf Greek} thinkers. 
Building up the number system starts with the {\bf natural numbers}  $1,2,3,4 ...$ which can be added and multiplied.
Addition is natural: join 3 sticks to 5 sticks to get 8 sticks.
Multiplication $*$ is more subtle: $3*4$ means to take $3$ copies of $4$ and 
get $4+4+4=12$ while $4*3$ means to take $4$ copies of $3$ to get $3+3+3+3=12$. The first factor counts the number of 
operations while the second factor counts the objects. To motivate $3*4=4*3$, spacial insight motivates to arrange 
the 12 objects in a rectangle. This commutativity axiom will be carried over to larger number systems. 
Realizing an addition and multiplicative structure on the natural numbers requires to define $0$ and $1$.
It leads naturally to more general numbers. There are two major motivations to {\bf to build new numbers}: we want to

\begin{center} \fcolorbox{yellow1}{yellow1}{ \parbox{17cm}{
1. {\bf invert operations} and still get results. \hfill 2. {\bf solve equations}. 
}} \end{center}

To find an additive inverse of $3$ means solving $x+3=0$. The answer is a negative number.
To solve $x*3=1$, we get to a rational number $x=1/3$. 
To solve $x^2=2$ one need to escape to real numbers. 
To solve $x^2=-2$ requires complex numbers. 

\begin{center}
\fcolorbox{yellow2}{yellow2}{
\parbox{16cm}{
\begin{tabular}{|l|l|l|} \hline
Numbers           &  Operation to complete         &  Examples of equations to solve \\ \hline
Natural numbers   &  addition and multiplication   &  $5+x=9$                      \\
Positive fractions&  addition and division         &  $5 x = 8$                    \\
Integers          &  subtraction                   &  $5 + x=3$                    \\
Rational numbers  &  division                      &  $3 x = 5$                    \\
Algebraic numbers &  taking positive roots         &  $x^2 = 2$ , $2x+x^2-x^3=2$   \\
Real numbers      &  taking limits                 &  $x=1-1/3+1/5-+... $,$\cos(x)=x$ \\
Complex numbers   &  take any roots                &  $x^2 = -2$                   \\
Surreal numbers   &  transfinite limits            &  $x^2 = \omega$, $1/x=\omega$ \\ 
Surreal complex   &  any operation                 &  $x^2 +1 = - \omega$          \\ \hline
\end{tabular}
}
}
\end{center}  
\index{extending mathematical operations}

The development and history of arithmetic can be summarized as follows: humans started
with natural numbers, dealt with positive fractions, reluctantly introduced negative numbers and zero
to get the integers, struggled to ``realize" real numbers, were scared to introduce complex numbers, 
hardly accepted surreal numbers and most do not even know about surreal complex numbers.
Ironically, as simple but impossibly difficult questions in number theory show, the modern point of 
view is the opposite to Kronecker's {\bf "God made the integers; all else is the work of man"}:

\begin{center} \fcolorbox{yellow1}{yellow1}{ \parbox{14cm}{
The {\bf surreal  complex} numbers are the most {\bf natural} numbers; \\
The {\bf natural} numbers are the most {\bf complex, surreal} numbers.
}} \end{center} 
\index{surreal numbers}

{\bf Natural numbers}. Counting can be realized by sticks, bones, quipu knots, pebbles or wampum knots.
The {\bf tally stick} concept is still used when playing card games: where bundles of fives are formed, maybe by 
crossing 4 "sticks" with a fifth. There is a "log counting" method in which graphs are used and vertices and edges count.
An old stone age tally stick, the {\bf wolf radius bone} contains 55 notches, with 5 groups of 5.
It is probably more than 30'000 years old.  \cite{SondheimerRogerson}
The most famous paleolithic tally stick is the 
{\bf Ishango bone}, the fibula of a baboon. It could be 20'000 - 30'000 years old. \cite{Eves} 
\index{tally stick}
\index{Ishango bone}
\index{Wolf bone}
Earlier counting could have been done by assembling {\bf pebbles}, tying {\bf knots} in a string, making 
{\bf scratches} in dirt or bark but no such traces have survived the thousands of years. 
The {\bf Roman system} improved the tally stick concept by introducing new symbols for larger numbers like
$V=5,X=10,L=40,C=100,D=500,M=1000$. \label{Roman system}
in order to avoid bundling too many single sticks.
The system is unfit for computations as simple calculations $VIII + VII = XV$ show.
{\bf Clay tablets}, some as early as 2000 BC and others from 600 - 300 BC are known. They feature
{\bf Akkadian arithmetic} using the base 60. 
The hexadecimal system with base $60$ is convenient because of many factors. It survived: we use
$60$ minutes per hour. {\bf The Egyptians} used the base 10.  The most important source on Egyptian mathematics is the 
{\bf Rhind Papyrus} of 1650 BC.  It was found in 1858 \cite{Katz2007,SondheimerRogerson}. 
Hieratic numerals were used to write on papyrus from 2500 BC on.
{\bf Egyptian numerals} are hieroglyphics. Found in carvings on tombs and 
monuments they are 5000 years old. 
The modern way to write numbers like 2018 is the
{\bf Hindu-Arab system} which diffused to the West only during
the late Middle ages. It replaced the more primitive {\bf Roman system}. \cite{SondheimerRogerson}
\label{Roman system}
\label{Hindu-Arab system}
Greek arithmetic used a number system with no place values: 
9 Greek letters for $1,2,\dots 9$, nine for $10,20,\dots,90$ and nine for $100,200,\dots,900$. \\ 

{\bf Integers}.
{\bf Indian Mathematics} morphed the place-value system into a modern method of writing numbers. 
Hindu astronomers used words to represent digits, but the numbers would be written in the opposite order. 
Independently, also the Mayans developed the concept of $0$ in a number system using base $20$. 
Sometimes after 500, the Hindus changed to a digital notation which included the symbol $0$. 
Negative numbers were introduced around 100 BC in the {\bf Chinese} text "Nine Chapters on the Mathematical art".
Also the {\bf Bakshali manuscript}, written around 300 AD subtracts numbers carried out additions
with negative numbers, where $+$ was used to indicate a negative sign.  \cite{Mathbook}
In Europe, negative numbers were avoided until the 15'th century.  \\
\index{Bakshali manuscript}
\index{Nine chapters on the Mathematical art}

{\bf Fractions}: 
{\bf Babylonians} could handle fractions. The {\bf Egyptians} also used fractions, but wrote every fraction a as a sum of fractions
with unit numerator and distinct denominators, like $4/5 = 1/2 + 1/4+1/20$ or
$5/6 = 1/2+1/3$. Maybe because of such cumbersome computation techniques, Egyptian mathematics failed
to progress beyond a primitive stage. \cite{SondheimerRogerson}.
The modern decimal fractions used nowadays for numerical calculations were 
adopted only in 1595 in Europe.  \\ 

{\bf Real numbers:} As noted by the Greeks already, the diagonal of the square is not a fraction. 
It first produced a crisis until it became clear that "most" numbers are not rational. 
{\bf Georg Cantor} saw first that the cardinality of all real numbers is much larger
than the cardinality of the integers: while one can count all rational numbers
but not enumerate all real numbers. One consequence is that most
real numbers are transcendental: they do not occur as solutions of polynomial equations with 
integer coefficients. The number $\pi$ is an example. The concept of real numbers is related to the
{\bf concept of limit}. Sums like $1+1/4+1/9+1/16+1/25 + \dots$ are not rational.

{\bf Complex numbers:} 
some polynomials have no real root. To solve $x^2=-1$ for example, we need new numbers.
One idea is to use pairs of numbers $(a,b)$ where $(a,0)=a$ are the usual numbers and extend addition and multiplication
$(a,b) + (c,d) = (a+c,b+d)$ and $(a,b) \cdot (c,d) = (ac-bd,ad+bc)$. With this multiplication, the number
$(0,1)$ has the property that $(0,1) \cdot (0,1) = (-1,0)=-1$.  It is more convenient to write $a+ib$ where $i=(0,1)$
satisfies $i^2=-1$. One can now use the common rules of addition and multiplication. 


{\bf Surreal numbers:}
Similarly as real numbers fill in the gaps between the integers, the surreal numbers fill in the 
gaps between Cantors ordinal numbers. They are written as $( a,b,c,...| d,e,f ,... )$
meaning that the "simplest" number is larger than $a,b,c...$ and smaller than $d,e,f,..$. 
We have $(|)=0, (0|)=1, (1|)=2$ and $(0|1)=1/2$ or $(|0)=-1$. Surreals contain already transfinite 
numbers like $(0,1,2,3...|)$ or infinitesimal numbers like $(0|1/2,1/3,1/4,1/5,...)$. They 
were introduced in the 1970'ies by John Conway. The late appearance
confirms the pedagogical principle: {\bf late human discovery manifests in increased difficulty to teach it}.

 \pagebreak
\pagebreak {\bf E-320: Teaching Math with a Historical Perspective \hfill Oliver Knill, 2010-2018}

\chapter{Lecture 3: Geometry}

Geometry is the science of {\bf shape, size and symmetry}. 
While arithmetic deals with numerical structures, geometry handles 
metric structures.  Geometry is one of the oldest mathematical 
disciplines. Early geometry has relations with arithmetic: 
the multiplication of two numbers $n \times m$ as an area of a
{\bf shape} that is invariant under rotational {\bf symmetry}. 
Identities like the {\bf Pythagorean triples} $3^2+4^2=5^2$ 
were interpreted and drawn geometrically. The {\bf right angle} is 
the most "symmetric" angle apart from $0$.
Symmetry manifests itself in quantities which are {\bf invariant}. 
Invariants are one the most central aspects of geometry. 
Felix Klein's {\bf Erlangen program} uses symmetry to classify 
geometries depending on how large the symmetries of the shapes are. 
In this lecture, we look at a few results which can all be stated 
in terms of invariants. 
In the presentation as well as the worksheet part of this lecture, 
we will work us through smaller miracles like
{\bf special points in triangles} as well as a couple of gems:
{\bf Pythagoras}, {\bf Thales},{\bf Hippocrates}, {\bf Feuerbach}, {\bf Pappus}, {\bf Morley},
{\bf Butterfly} which illustrate the importance of symmetry. \\
\index{Erlangen program}
\index{Feuerbach theorem}
\index{Pappus theorem}
\index{Morley theorem}
\index{Butterfly theorem}

Much of geometry is based on our ability to measure {\bf length}, the 
{\bf distance} between two points.  
Having a distance $d(A,B)$ between any two points $A,B$, we can look at 
the next more complicated object, which is a set $A,B,C$ of 3 points, 
a {\bf triangle}. Given an arbitrary triangle ABC, are there relations 
between the 3 possible distances $a=d(B,C),b=d(A,C),c=d(A,B)$? 
If we fix the scale by $c=1$, then $a+b \geq 1, a+1 \geq b, b+1 \geq a$. For any 
pair of $(a,b)$ in this region, there is a triangle. 
After an identification, we get an abstract space, which represent all
triangles uniquely up to similarity. 
Mathematicians call this an example of a {\bf moduli space}. \\ 
\index{moduli space}

A {\bf sphere} $S_r(x)$ is the set of points which have distance $r$ from a given point $x$.
In the plane, the sphere is called a {\bf circle}.
A natural problem is to find the circumference $L=2 \pi$ of a unit circle,
or the area $A = \pi$ of a unit disc, 
the area $F=4 \pi$ of a unit sphere and the volume $V=4=\pi/3$ of a unit sphere. 
Measuring the length of segments on the circle leads to new 
concepts like {\bf angle} or {\bf curvature}.
Because the circumference of the unit circle in the plane is $L=2\pi$, 
angle questions are tied to the number
$\pi$, which Archimedes already approximated by fractions. \\

Also {\bf volumes} were among the first quantities, Mathematicians wanted to measure and compute.
A problem on {\bf Moscow papyrus} dating back to 1850 BC
explains the general formula $h (a^2 + a b + b^2)/3$ for a truncated pyramid with base length $a$, 
roof length $b$ and height $h$. 
Archimedes achieved to compute the  {\bf volume of the sphere}:
place a cone inside a cylinder. 
The complement of the cone inside the cylinder has on each height $h$
the area $\pi - \pi h^2$. 
The half sphere cut at height $h$ is a disc of radius $(1-h^2)$ which has area
$\pi(1-h^2)$ too. Since the slices at each height have the same area, the volume must be the same. 
The complement of the cone inside the cylinder has 
volume $\pi - \pi/3 = 2\pi/3$, half the volume of the sphere.  \\

The first geometric playground was {\bf planimetry}, the geometry in the flat two dimensional space.
Highlights are {\bf Pythagoras theorem}, {\bf Thales theorem}, {\bf Hippocrates theorem}, 
and {\bf Pappus theorem}. Discoveries in planimetry have been made later on:
an example is the Feuerbach 9 point theorem from the 19th century.
Ancient Greek Mathematics is closely related to history. 
It starts with {\bf Thales} goes over Euclid's era at 500 BC and
ends with the threefold destruction of Alexandria 
47 BC by the Romans, 392 by the Christians and 640 by the Muslims.  
Geometry was also a place, where the {\bf axiomatic method} 
was brought to mathematics: theorems are proved from a 
few statements which are called axioms like the 5 axioms of Euclid: 

\begin{center}
\fcolorbox{yellow2}{yellow2}{
\parbox{16.8cm}{
1. Any two distinct points $A,B$ determines a line through $A$ and $B$.  \\
2. A line segment $[A,B]$ can be extended to a straight line containing the segment. \\
3. A line segment $[A,B]$ determines a circle containing $B$ and  center $A$.  \\
4. All right angles are congruent.  \\
5. If lines $L,M$ intersect with a third so that inner angles add up to $<\pi$, then $L,M$ intersect.
}} \end{center}

{\bf Euclid} wondered whether the fifth postulate can be derived from the first four and
called theorems derived from the first four the "absolute geometry". 
Only much later, with {\bf  Karl-Friedrich Gauss} and {\bf Janos Bolyai} and {\bf Nicolai Lobachevsky} in the 19'th century 
in {\bf hyperbolic space} the 5'th axiom does not hold. Indeed, geometry can be generalized to non-flat, or even much more abstract situations.
Basic examples are geometry on a sphere leading to {\bf spherical geometry} or geometry on the Poincare disc, a {\bf hyperbolic space}. 
Both of these geometries are non-Euclidean. {\bf Riemannian geometry}, which is essential for {\bf general relativity theory} generalizes
both concepts to a great extent. An example is the geometry on an arbitrary surface. Curvatures of such spaces can be computed by measuring
length alone, which is how long light needs to go from one point to the next. \\

An important moment in mathematics was the {\bf merge of geometry with algebra}: 
this giant step is often attributed to {\bf Ren\'e Descartes}. 
Together with algebra, the subject leads to algebraic geometry which can
be tackled with computers: here are some examples of geometries which are determined from the
amount of symmetry which is allowed: 

\begin{center}
\fcolorbox{yellow2}{yellow2}{
\parbox{16.8cm}{
\begin{tabular}{ll}
Euclidean geometry    &   Properties invariant under a group of rotations and translations  \\
Affine geometry       &   Properties invariant under a group of affine transformations \\
Projective geometry   &   Properties invariant under a group of projective transformations \\
Spherical geometry    &   Properties invariant under a group of rotations \\
Conformal geometry    &   Properties invariant under angle preserving transformations \\
Hyperbolic geometry   &   Properties invariant under a group of M\"obius transformations \\
\end{tabular}
}} \end{center}

Here are four pictures about the 4 special points in a triangle and
with which we will begin the lecture. 
We will see why in each of these cases, the 3 lines intersect in a common point.
It is a manifestation of a {\bf symmetry} present on the space of all 
triangles. {\bf size} of the distance of intersection points is constant 
$0$ if we move on the space of all triangular {\bf shapes}. It's Geometry! \\

 \pagebreak
\pagebreak {\bf E-320: Teaching Math with a Historical Perspective \hfill Oliver Knill, 2010-2018}

\chapter{Lecture 4: Number Theory}

Number theory studies the structure of integers like prime numbers and 
solutions to Diophantine equations. 
Gauss called it the "Queen of Mathematics". Here are a few theorems and open problems. \\
An integer larger than $1$ which is divisible by $1$ and itself only is called a {\bf prime number}.
The number $2^{57885161}-1$ is the largest known prime number. It has $17425170$ digits.
{\bf Euclid} proved that there are infinitely many primes: 
[Proof. Assume there are only finitely many primes $p_1<p_2< \dots < p_n$. Then $n=p_1 p_2 \cdots p_n + 1$
is not divisible by any $p_1, \dots,p_n$. Therefore, it is a prime or divisible by a 
prime larger than $p_n$.] Primes become more sparse as larger as they get. An important result is
the {\bf prime number theorem} which states that the $n$'th prime number has approximately
the size $n \log(n)$. For example the $n=10^{12}$'th prime is $p(n) = 29996224275833$ and 
$n \log(n) = 27631021115928.545...$ and $p(n)/(n \log(n)) = 1.0856..$.
Many questions about prime numbers are unsettled:
Here are four problems: the third uses the notation $(\Delta a)_n = |a_{n+1}-a_n|$ to get
the absolute difference. For example: $\Delta^2 (1,4,9,16,25 ...) =
\Delta(3,5,7,9,11,...) = (2,2,2,2,...)$. Progress on prime gaps has been done in 2013:
$p_{n+1}-p_n$ is smaller than 100'000'000 eventually (Yitang Zhang).
$p_{n+1}-p_n$ is smaller than $600$ eventually (Maynard). 
The largest known gap is 1476 which occurs after $p=1425172824437699411$. 
\index{Twin prime}

\begin{center}
\fcolorbox{yellow2}{yellow2}{
\parbox{16cm}{
\begin{tabular}{|l|l|} \hline
{\bf Landau}     & there are infinitely many primes of the form $n^2+1$. \\
{\bf Twin prime} & there are infinitely many primes $p$ such that $p+2$ is prime.  \\
{\bf Goldbach}   & every even integer $n>2$ is a sum of two primes.  \\
{\bf Gilbreath}  & If $p_n$ enumerates the primes, then $(\Delta^k p)_1=1$ for all $k>0$. \\
{\bf Andrica}    & The prime gap estimate $\sqrt{p_{n+1}}-\sqrt{p_n} <1$ holds for all $n$. \\ \hline
\end{tabular} }} \end{center}

If the sum of the proper divisors of a $n$ is equal to $n$, then $n$ is called a 
{\bf perfect number}. For example, $6$ is perfect as
its proper divisors $1,2,3$ sum up to $6$. All currently known perfect numbers
are even. The question whether odd perfect numbers exist is probably the oldest
open problem in mathematics and not settled. Perfect numbers were familiar to Pythagoras and his followers already. 
Calendar coincidences like that we have 6 work days and the moon needs "perfect" 28 days 
to circle the earth could have helped to promote the "mystery" of perfect number. 
{\bf Euclid of Alexandria} (300-275 BC) was the first to realize that if $2^p-1$ is prime then
$k=2^{p-1} (2^p-1)$ is a perfect number:
[Proof: let $\sigma(n)$ be the sum of {\bf all} factors of $n$, including $n$. Now
$\sigma( 2^n-1) 2^{n-1})  = \sigma(2^n-1) \sigma(2^{n-1}) = 2^n (2^n-1) = 2 \cdot 2^n (2^n-1)$
shows $\sigma(k) = 2k$ and verifies that $k$ is perfect.]
Around 100 AD, {\bf Nicomachus of Gerasa} (60-120) classified in his work
"Introduction to Arithmetic" numbers on the concept of perfect numbers and lists four perfect
numbers. 
Only much later it became clear that Euclid got all the even perfect numbers:
Euler showed that all even perfect numbers are of the form $(2^n-1) 2^{n-1}$, where $2^n-1$ is prime.
The factor $2^n-1$ is called a {\bf Mersenne prime}. 
[Proof: Assume $N=2^k m$ is perfect where $m$ is odd and $k>0$. Then
$2^{k+1} m = 2N = \sigma(N) = (2^{k+1}-1) \sigma(m)$.
This gives
$\sigma(m) = 2^{k+1} m/(2^{k+1}-1) = m (1+1/(2^{k+1}-1)) = m + m/(2^{k+1}-1)$.
Because $\sigma(m)$ and $m$ are integers, also $m/(2^{k+1}-1)$ is an integer. It must
also be a factor of $m$. The only way that $\sigma(m)$ can be the sum of only two of its
factors is that $m$ is prime and so $2^{k+1}-1=m$.]
The first 39 {\bf known Mersenne primes} are of the form $2^n-1$ with 
n = 2, 3, 5, 7, 13, 17, 19, 31, 61, 89, 107, 127, 521, 607, 1279, 2203, 2281, 3217, 
4253, 4423, 9689, 9941, 11213, 19937, 21701, 23209, 44497, 86243, 110503, 132049, 
216091, 756839, 859433, 1257787, 1398269, 2976221, 3021377, 6972593, 13466917.
There are 11 more known from which one does not know the rank of the corresponding Mersenne prime:
n = 20996011, 24036583, 25964951, 30402457, 32582657, 37156667, 42643801,43112609,57885161,
74207281,77232917. The last was found in December 2017 only. It is unknown whether there are infinitely many. \\
\index{Mersenne primes}

A polynomial equations for which all coefficients and variables are integers is called
a {\bf Diophantine equation}. The first Diophantine equation studied already by Babylonians is
$x^2 +y^2 = z^2$. A solution $(x,y,z)$ of this equation in positive integers is called a 
{\bf Pythagorean triple}. For example, $(3,4,5)$ is a Pythagorean triple.
Since 1600 BC, it is known that all solutions to this equation are of the form $(x,y,z) =(2st,s^2-t^2,s^2+t^2)$ or 
$(x,y,z) =(s^2-t^2,2 s t,s^2+t^2)$, where $s,t$ are different integers. 
[Proof. Either $x$ or $y$ has to be even because if both are odd, then the sum $x^2+y^2$
is even but not divisible by $4$ but the right hand side is either odd or divisible by $4$. 
Move the even one, say $x^2$ to the left and write $x^2= z^2-y^2 = (z-y) (z+y)$, 
then the right hand side contains a factor $4$ and is of the form $4s^2 t^2$. Therefore
$2s^2 = z-y, 2t^2 = z+y$. Solving for $z,y$ gives $z = s^2+t^2,
y=s^2-t^2$, $x=2st$.] \\
Analyzing Diophantine equations can be difficult. Only 10 years ago, one has established that
the {\bf Fermat equation} $x^n+y^n=z^n$ has no solutions with $xyz \neq 0$ if $n>2$. 
Here are some {\bf open problems} for Diophantine equations. 
Are there nontrivial solutions to the following Diophantine equations?
\begin{center}
\fcolorbox{yellow2}{yellow2}{
\parbox{9cm}{
\begin{tabular}{|l|l|} \hline
    $x^6+y^6+z^6+u^6+v^6 = w^6$    &  $x,y,z,u,v,w>0$  \\
    $x^5+y^5+z^5 = w^5$            &  $x,y,z,w>0$      \\
    $x^k + y^k = n!  z^k$          &  $k \geq 2,n>1$   \\
    $x^a + y^b = z^c, a,b,c>2$     &  ${\rm gcd}(a,b,c)=1$ \\ \hline
\end{tabular} }} \end{center}
The last equation is called {\bf Super Fermat}.
A Texan banker {\bf Andrew Beals} once sponsored a prize of $100'000$ 
dollars for a proof or counter example to the statement:
"If $x^p+y^q = z^r$ with $p,q,r>2$, then ${\rm gcd}(x,y,z)>1$." 
\index{Beals conjecture}
Given a prime like $7$ and a number $n$ we can add or subtract multiples
of $7$ from $n$ to get a number in $\{0,1,2,3,4,5,6 \; \}$. We write for example $19 = 12 \; {\rm mod} \; 7$
because $12$ and $19$ both leave the rest $5$ when dividing by $7$. 
Or $5*6 = 2 \; {\rm mod} \; 7$ because $30$ leaves the rest $2$ when dividing by $7$.
The most important theorem in elementary number theory is {\bf Fermat's little theorem} which tells that
if $a$ is an integer and $p$ is prime then $a^{p} - a$ is divisible by $p$.
For example $2^7 -2 = 126$ is divisible by $7$. [Proof: use induction. For $a=0$ it is clear.
The binomial expansion shows that $(a+1)^p-a^p-1$ is divisible by $p$. 
This means $(a+1)^p - (a+1) = (a^p-a) + m p$ for some $m$. By induction, $a^p-a$ is divisible by $p$ and so $(a+1)^p - (a+1)$.] 
An other beautiful theorem is {\bf Wilson's theorem} which allows to characterize primes: 
It tells that $(n-1)!+1$ is divisible by $n$ if and only if 
$n$ is a prime number. For example, for $n=5$, we verify that $4!+1=25$ is divisible by $5$.
[Proof: assume $n$ is prime. There are then exactly two numbers $1,-1$ for which $x^2-1$ is divisible by $n$. 
The other numbers in $1,\dots,n-1$ can be paired as $(a,b)$ with $a b=1$. Rearranging the product 
shows $(n-1)!=-1$ modulo $n$. Conversely, if $n$ is not prime, then $n=k m$ with $k,m<n$ and 
$(n-1)! = ... k m$ is divisible by $n=k m$. ] \\
The solution to systems of linear equations like $x=3 \; ({\rm mod} \; 5), x=2 \; ({\rm mod} \; 7)$ is given by 
the {\bf Chinese remainder theorem}. To solve it, continue adding $5$ to $3$ until we reach a 
number which leaves rest $2$ to $7$: on the list $3,8,13,18,23,28,33,38$, the number $23$ is
the solution. Since $5$ and $7$ have no common divisor, the system of linear equations has
a solution.  \\
For a given $n$, how do we solve $x^2 - y n =1$ for the unknowns $y,x$?
A solution produces a square root $x$ of $1$ modulo $n$. For prime $n$, 
only $x=1,x=-1$ are the solutions. For composite $n=pq$, more solutions $x=r \cdot s$ where
$r^2=-1 \; {\rm mod} \;  p$ and $s^2=-1 \; {\rm mod} \;  q$ appear. Finding $x$ is equivalent to 
factor $n$, because the greatest common divisor of $x^2-1$ and $n$ is a factor of $n$. 
{\bf Factoring is difficult} if the numbers are large. It assures that {\bf encryption algorithms} 
work and that bank accounts and communications stay safe. Number theory, once the least applied 
discipline of mathematics has become one of the most applied one in mathematics. 

 \pagebreak
\pagebreak {\bf E-320: Teaching Math with a Historical Perspective \hfill Oliver Knill, 2010-2018}

\chapter{Lecture 5: Algebra}

Algebra studies {\bf algebraic structures} like "groups" and "rings". 
The theory allows to solve polynomial equations, 
characterize objects by its symmetries and is the heart and soul of many puzzles. 
Lagrange claims {\bf Diophantus} to be the inventor of Algebra, 
others argue that the subject started with solutions of 
{\bf quadratic equation} by {\bf Mohammed ben Musa Al-Khwarizmi}
in the book Al-jabr w'al muqabala of 830 AD. Solutions to equation like
$x^2 + 10 x = 39$ are solved there by {\bf completing the squares}: add 25 on both sides
go get $x^2+10 x + 25 = 64$ and so $(x+5) = 8$ so that $x=3$.


The use of {\bf variables} introduced in school in {\bf elementary algebra}
were introduced later. Ancient texts only dealt with particular examples and calculations were 
done with concrete numbers in the realm of {\bf arithmetic}. 
{\bf Francois Viete} (1540-1603) used first letters like $A,B,C,X$ for variables. \\ 

The search for formulas for polynomial equations of degree $3$ and $4$ lasted 700 years. 
In the 16'th century, the cubic equation and quartic equations were solved. 
{\bf Niccolo Tartaglia} and {\bf Gerolamo Cardano} reduced the cubic to the quadratic:
[first remove the quadratic part with $X=x-a/3$ so that $X^3+aX^2+bX+c$ becomes the 
{\bf depressed cubic} $x^3 + p x + q$. 
Now substitute $x=u-p/(3u)$ to get a quadratic equation $(u^6+qu^3-p^3/27)/u^3=0$ for $u^3$.]
{\bf Lodovico Ferrari} shows that the quartic equation can be reduced to the cubic. 
For the {\bf quintic} however no formulas could be found. 
It was {\bf Paolo Ruffini}, {\bf Niels Abel} and {\bf \'Evariste Galois} who independently 
realized that there are no formulas in terms of roots which allow to "solve" equations 
$p(x)=0$ for polynomials $p$ of degree larger than $4$.
This was an amazing achievement and the birth of "group theory". 

\begin{center} \fcolorbox{yellow1}{yellow1}{ \parbox{16cm}{
Two important algebraic structures are {\bf groups} and {\bf rings}.
}} \end{center}

In a {\bf group} $G$ one has an operation $*$, an inverse $a^{-1}$ and a one-element $1$
such that $a*(b*c) = (a*b)*c, a*1=1*a=a, a*a^{-1} =a^{-1}*a=1$. 
For example, the set $Q^*$ of nonzero fractions $p/q$ with multiplication 
operation $*$ and inverse $1/a$ form a group. The integers with addition and inverse 
$a^{-1} = -a$ and "1"-element $0$ form a group too. 
A {\bf ring} $R$ has two compositions $+$ and $*$, where the plus operation
is a group satisfying $a+b=b+a$ in which the one element is called $0$. The multiplication operation 
$*$ has all group properties on $R^*$ except the existence of an inverse. 
The two operations $+$ and $*$ are glued together by the {\bf distributive law} $a*(b+c) = a*b+a*c$. 
An example of a ring are the {\bf integers} or the {\bf rational numbers} or the
{\bf real numbers}. The later two are actually {\bf fields}, rings for which the 
multiplication on nonzero elements is a group too. 
The ring of integers are no field because an integer like $5$ has no multiplicative inverse. 
The ring of rational numbers however form a field. \\

Why is the theory of groups and rings not part of arithmetic? 
First of all, a crucial ingredient of algebra is the appearance of {\bf variables} and 
computations with these algebras without using concrete numbers. Second,
the algebraic structures are not restricted to "numbers". 
Groups and rings are general structures and extend
for example to objects like the set of all possible symmetries of a geometric object. 
The set of all {\bf similarity operations} on the plane for 
example form a group. An important example of a ring is the {\bf polynomial ring} of all
polynomials. Given any ring $R$ and a variable $x$, the set $R[x]$ consists of all polynomials
with coefficients in $R$. The addition and multiplication is done like in 
$(x^2+3x+1) + (x-7) = x^2+4x-7$.
The problem to factor a given polynomial with integer coefficients into polynomials of 
smaller degree: $x^2-x+2$ for example can be written as $(x+1)(x-2)$ have a number 
theoretical flavor. Because symmetries of some structure form a group, we also have 
intimate connections with geometry. 
But this is not the only connection with geometry. Geometry also enters through the polynomial rings with 
several variables. Solutions to $f(x,y)=0$ leads to geometric objects with shape and 
symmetry which sometimes even have their own algebraic structure. 
They are called {\bf varieties}, a central object in {\bf algebraic geometry}, objects
which in turn have been generalized further to {\bf schemes}, 
{\bf algebraic spaces} or {\bf stacks}. \\

Arithmetic introduces addition and multiplication of numbers. 
Both form a group. The operations can be written additively or multiplicatively. 
Lets look at this a bit closer: for integers, fractions and reals and the addition $+$,
the $1$ element $0$ and inverse $-g$, we have a group.
Many groups are written multiplicatively where the $1$ element is $1$. 
In the case of fractions or reals, $0$ is not part of the multiplicative group
because it is not possible to divide by $0$. The nonzero fractions or the nonzero 
reals form a group. 
In all these examples the groups satisfy the commutative law $g*h = h*g$. \\
Here is a group which is not commutative: let $G$ be the set of all rotations in space, which 
leave the unit cube invariant. There are 3*3=9 rotations around
each major coordinate axes, then 6 rotations around axes connecting midpoints of opposite edges,
then 2*4 rotations around diagonals. Together with the identity rotation $e$, these are
24 rotations. The group operation is the composition of these transformations.  \\
An other example of a group is $S_4$, the set of all permutations of 
four numbers $(1,2,3,4)$. If $g: (1,2,3,4) \to (2,3,4,1)$ is a permutation and 
$h: (1,2,3,4) \to (3,1,2,4)$ is an other permutation, then we can combine the two 
and define $h * g$ as the permutation which does first $g$ and then $h$. We end up 
with the permutation $(1,2,3,4) \to (1,2,4,3)$. 
The rotational symmetry group of the cube happens to be the same than the group $S_4$. 
To see this "isomorphism", label the 4 space diagonals in the cube by $1,2,3,4$. 
Given a rotation, we can look at the induced permutation of the diagonals and
every rotation corresponds to exactly one permutation. 
The symmetry group can be introduced for any geometric object. For shapes like the triangle, the cube,
the octahedron or tilings in the plane. 

\begin{center} \fcolorbox{yellow1}{yellow1}{ \parbox{12cm}{
Symmetry groups describe geometric shapes by algebra. 
}} \end{center}

Many {\bf puzzles} are groups. A popular puzzle, the {\bf 15-puzzle}
was invented in 1874 by {\bf Noyes Palmer Chapman} in the state of New York. 
If the hole is given the number $0$, then the task of the puzzle is to order a given 
random start permutation of the 16 pieces. 
To do so, the user is allowed to transposes $0$ with a neighboring piece. Since every 
step changes the signature $s$ of the permutation and changes the taxi-metric 
distance $d$ of $0$ to the end position by $1$, only situations with even $s+d$ can 
be reached. It was {\bf Sam Loyd} who suggested to start with an impossible
solution and as an evil plot to offer 1000 dollars for a solution. 
The 15 puzzle group has $16!/2$ elements and the "god number" is between $152$ and $208$. 
The {\bf Rubik cube} is an other famous puzzle, which is a group.
Exactly 100 years after the invention of the 15 puzzle, the Rubik puzzle was 
introduced in 1974. Its still popular and the world record is to have it solved in 
5.55 seconds. All Cubes 2x2x2 to 7x7x7 in a row have been solved in a total time of 6 minutes.
For the 3x3x3 cube, the {\bf God number} is now known to be 20: one can always solve it in 
20 or less moves. 
\index{Rubik cube}
\index{God number}

\begin{center} \fcolorbox{yellow1}{yellow1}{ \parbox{12cm}{
Many puzzles are groups.
}} \end{center}

A small Rubik type game is the "floppy", which is a third of the Rubik and which has only 192 elements. 
An other example is the {\bf Meffert's great challenge}. 
Probably the simplest example of a Rubik type puzzle is the {\bf pyramorphix}. 
It is a puzzle based on the tetrahedron. Its group has only 24 elements.
It is the group of all possible permutations of the 4 elements. 
It is the same group as the group of all reflection and rotation symmetries of 
the cube in three dimensions and also is relevant when understanding 
the solutions to the quartic equation discussed at the beginning. The circle is closed.

 \pagebreak
\pagebreak {\bf E-320: Teaching Math with a Historical Perspective \hfill Oliver Knill, 2010-2018}

\chapter{Lecture 6: Calculus}

Calculus generalizes the process of {\bf taking differences} and {\bf taking sums}. 
Differences measure {\bf change}, sums explore how quantities
{\bf accumulate}. The procedure of taking differences has a limit called {\bf derivative}. The 
activity of taking sums leads to the {\bf integral}. Sum and difference are dual to each other 
and related in an intimate way. In this lecture, we look first at a simple 
set-up,  where functions are evaluated on integers and where we do not take any limits.

Several dozen thousand years ago, numbers were represented by units like
$1,1,1,1,1,1, \dots $. The units were carved into sticks or bones like the {\bf Ishango bone}
It took thousands of years until numbers were represented with symbols like 
$0,1,2,3,4, \dots $. Using the modern concept of function, we can say
$f(0)=0,f(1)=1, f(2)=2, f(3)=3$ and mean that the {\bf function} $f$ assigns to an input like
$1001$ an output like $f(1001)=1001$. Now look at $Df(n) = f(n+1)-f(n)$, the {\bf difference}.
We see that $Df(n) = 1$ for all $n$. We can also formalize the summation process. 
If $g(n) = 1$ is the constant $1$ function, then then 
$Sg(n) = g(0) + g(1) + \dots + g(n-1) = 1+1+\cdots +1 = n$. We see that $Df=g$ and $Sg=f$. 
If we start with $f(n)=n$ and apply {\bf summation} on that function
Then $Sf(n) = f(0) + f(1) + f(2) + \cdots + f(n-1)$ leading to the 
values $0,1,3,6,10,15,21, \dots$. The new function $g=Sf$ satisfies 
$g(1)=1,g(2)=3,g(2)=6$, etc. The values are called the 
{\bf triangular numbers}. From $g$ we can get back $f$ by taking difference:
$Dg(n) = g(n+1)-g(n)=f(n)$. For example $Dg(5) = g(6)-g(5) = 15-10 = 5$ which indeed is $f(5)$.
Finding a formula for the sum $Sf(n)$ is not so easy. Can you do it? 
When {\bf Karl-Friedrich Gauss} was a 9 year old school kid,
his teacher, a Mr. B\"uttner gave him the task to sum up the first 100 numbers $1+2+ \cdots + 100$.
Gauss found the answer immediately by pairing things up:
to add up $1+2+3+ \dots +100$ he would write this as $(1+100) + (2+99) + \cdots + (50+51) $ leading to 
$50$ terms of $101$ to get for $n=101$ the value $g(n)=n(n-1)/2 = 5050$. 
Taking differences again is easier $Dg(n) = n(n+1)/2 - n(n-1)/2 = n = f(n)$. 
If we add up he triangular numbers we compute $h=Sg$ which has the first values
$0,1,4,10,20,35, ....$. These are the {\bf tetrahedral numbers} because $h(n)$ balls  are needed
to build a tetrahedron of side length $n$. For example, $h(4)=20$ golf balls are needed to build a 
tetrahedron of side length 4. The formula which holds for $h$ is \fbox{$h(n) = n(n-1)(n-2)/6$}. 
Here is the fundamental theorem of calculus, which is the core of calculus: 

\begin{center}
\fcolorbox{yellow1}{yellow1}{ \parbox{8cm}{
$Df (n) = f(n) - f(0)$, \hspace{6mm}  $D S f(n) = f(n) \; . $
}}
\end{center}

Proof.
$$ SDf(n) =  \sum_{k=0}^{n-1} [f(k+1)-f(k) ]     = f(n)-f(0) \; ,  $$
$$ DSf(n) = [\sum_{k=0}^{n-1} f(k+1)  - \sum_{k=0}^{n-1} f(k) ] = f(n) \; . $$

The process of adding up numbers will lead to the {\bf integral} \fbox{$\int_0^x f(x) \; dx$}.
The process of taking differences will lead to the {\bf derivative} \fbox{$\frac{d}{dx} f(x)$}. 

The familiar notation is

\begin{center}
\fbox{
$\int_0^x \frac{d}{dt} f(t) \; dt = f(x)-f(0), \hspace{1cm}  \frac{d}{dx} \int_0^x f(t) \; dt = f(x)$
}
\end{center}

If we define $[n]^0 = 1, [n]^1=n, [n]^2 = n(n-1)/2, [n]^3 = n (n-1)(n-2)/6$
then $D[n]=[1], D[n]^2 = 2 [n], D[n]^3 = 3 [n]^2$ and in general
\begin{center} \fbox{
$\frac{d}{dx} [x]^n = n [x]^{n-1}$
} \end{center}
The calculus you have just seen, contains the essence of single variable calculus. This core idea will become more
powerful and natural if we use it together with the concept of limit. 

{\bf Problem:} The Fibonnacci sequence $1,1,2,3,5,8,13,21, \dots$ satisfies the rule $f(x)=f(x-1)+f(x-2)$. 
For example, $f(6) = 8$. What is the function $g=Df$, if we assume $f(0)=0$? 
We take the difference between successive numbers and get the sequence of numbers 
$0,1,1,2,3,5,8, ... $ which is the same sequence again. We see that \fbox{$Df(x)=f(x-1)$}. 

If we take the same function $f$ but now but now compute the function $h(n) = Sf(n)$, we get the sequence
$1,2,4,7,12,20,33, ...$. What sequence is that? 
{\bf Solution:} Because $Df(x)=f(x-1)$ we have $f(x)-f(0) = S D f(x) = S f(x-1)$ so that $S f(x)=f(x+1)-f(1)$. 
Summing the Fibonnacci sequence produces the Fibonnacci sequence shifted to the left with $f(2)=1$ is subtracted.
It has been relatively easy to find the sum, because we knew what the difference operation did. This example shows:
we can study differences to understand sums. 

{\bf Problem:} 
The function $f(n) = 2^n$ is called the {\bf exponential function}.
We have for example $f(0) =1, f(1)=2,f(2)=4, \dots $. It leads to the sequence of numbers
\begin{center}
\begin{tabular}{lllllllllll}
n=   &   0 & 1 & 2 & 3 & 4  & 5  & 6  & 7   & 8   & $\dots$ \\
f(n)=&   1 & 2 & 4 & 8 & 16 & 32 & 64 & 128 & 256 & $\dots$ \\
\end{tabular}
\end{center}
We can verify that $f$ satisfies the equation \fbox{$Df(x)=f(x)$}.
because $Df(x) = 2^{x+1} - 2^{x} = (2-1) 2^{x} = 2^{x}$. \\
This is an important special case of the fact that 
\begin{center}
\fbox{ The derivative of the exponential function is the exponential function itself.}
\end{center}
The function $2^x$ is a special case of the exponential function when the Planck constant is equal to $1$. 
We will see that the relation will hold for any $h>0$ and also in the limit $h \to 0$, where it 
becomes the classical exponential function $e^x$ which plays an important role in science. 

Calculus has many applications: computing areas, volumes, solving differential equations.
It even has applications in arithmetic. Here is an example for illustration. 
It is a proof that $\pi$ is irrational
The theorem is due to Johann Heinrich Lambert (1728-1777):
We show here the proof by Ivan Niven is given in a book of Niven-Zuckerman-Montgomery.
It originally appeared in 1947 (Ivan Niven, Bull.Amer.Math.Soc. 53 (1947),509).
The proof illustrates how calculus can help to get results in arithmetic. 

{\bf Proof}. Assume $\pi = a/b$ with positive integers $a$ and $b$.
For any positive integer $n$ define
$$ f(x) = x^n (a-b x)^n/n!  \; . $$
We have $f(x) = f(\pi-x)$ and
$$ 0 \leq  f(x) \leq  \pi^n a^n/n!  (*) $$
for $0 \leq x \leq \pi$. For all $0 \leq j \leq n$, the j-th derivative of $f$
is zero at $0$ and $\pi$ and for $n <= j$, the j-th derivative of $f$
is an integer at $0$ and $\pi$. \\
The function $F(x) = f(x) - f^{(2)}(x) + f^{(4)}(x)  - ...  + (-1)^n f^{(2n)}(x)$
has the property that $F(0)$ and $F(\pi)$ are integers and $F + F '' = f$. Therefore,
$(F'(x) \sin(x) - F(x) \cos(x))' = f \sin(x)$.
By the fundamental theorem of calculus, $\int_0^{\pi} f(x) \sin(x) \; dx$ is an integer.
Inequality (*) implies however that this integral is between 0 and $1$ for large enough $n$.
For such an $n$ we get a contradiction.

 \pagebreak
\pagebreak {\bf E-320: Teaching Math with a Historical Perspective \hfill Oliver Knill, 2010-2018}

\chapter{Lecture 7: Set Theory and Logic}

{\bf Set theory} studies sets, the fundamental building blocks of mathematics. While {\bf logic}
describes the language of all mathematics, set theory provides the framework 
for additional structures like category theory. 
In {\bf Cantorian set theory}, one can compute with subsets of a given set $X$ like with numbers. 
There are two basic operations: the {\bf addition} $A + B$ of two sets is defined as the set of
all points which are in exactly one of the sets. 
The {\bf multiplication} $A \cdot B$ of two sets contains all the 
points which are in both sets. With the symmetric difference as addition and the 
intersection as multiplication, the subsets of a given set $X$ become a 
{\bf ring}. This {\bf Boolean ring} has the property $A+A = 0$ and $A \cdot A = A$ for all sets. 
The zero element is the empty set $\emptyset=\{\}$. 
The additive inverse of $A$ is the complement $-A$ of $A$ in $X$.
The multiplicative $1$-element is the set $X$ because $X \cdot A = A$. 
As in the ring $\mathbb{Z}$ of integers, the addition and multiplication on sets
is commutative. Multiplication does not have an inverse in general. 
Two sets $A,B$ have the {\bf same cardinality}, if there exists a one-to-one map 
from $A$ to $B$. For finite sets, this means that they have the same number of elements. 
Sets which do not have finitely many elements are called {\bf infinite}. 
Do all sets with infinitely many elements have the same cardinality? 
The integers $\mathbb{Z}$ and the natural numbers $\mathbb{N}$ for example are 
infinite sets which have the same cardinality: the map
$f(2n)=n,f(2n+1)=-n$ establishes a bijection between $\mathbb{N}$ and  $\mathbb{Z}$. 
Also the rational numbers $\mathbb{Q}$ have the same cardinality than $\mathbb{N}$. 
Associate a fraction $p/q$ with a point $(p,q)$ in the plane. Now cut out the column 
$q=0$ and run the {\bf Ulam spiral} on the modified plane. This provides a numbering
of the rationals. Sets which can be counted are called of cardinality $\aleph_0$. 
Does an interval have the same cardinality than the reals? 
Even so an interval like $I=(-\pi/2,\pi/2)$ has finite length, 
one can bijectively map it to $\mathbb{R}$ with the $\tan$ function as $\tan: I \to \mathbb{R}$
is bijective. Similarly, one can see that any two intervals of positive length
have the same cardinality. 
\index{Ulam spiral}
It was a great moment of mathematics, when {\bf Georg Cantor} realized in 1874 that the interval 
$(0,1)$ does not have the same cardinality than the natural numbers. His argument is ingenious:
assume, we could count the points $a_1,a_2, \dots$. If $0.a_{i1} a_{i2} a_{i3} ... $ 
is the {\bf decimal expansion} of $a_i$, define the real number 
$b = 0.b_1 b_2 b_3 ...$, where $b_i = a_{ii}+1 \; {\rm mod} \; 10$. 
Because this number $b$ does not agree at the first decimal place with $a_1$, 
at the second place with $a_2$ and so on, the number $b$ does not appear in that 
enumeration of all reals. 
It has positive distance at least $10^{-i}$ from the $i$'th number 
(and any representation of the number by a decimal expansion which is equivalent). 
This is a contradiction.  The new cardinality, the {\bf continuum} 
is also denoted $\aleph_1$.  
The reals are {\bf uncountable}. This gives elegant proofs like the existence of 
{\bf transcendental number}, numbers which are not algebraic, meaning that they are
not the root of any polynomial with integer coefficients: algebraic numbers can be counted. 
Similarly as one can establish a bijection between the natural numbers $\mathbb{N}$ 
and the integers $\mathbb{Z}$, there is a bijection
$f$ between the interval $I$ and the unit square: 
if $x = 0.x_1 x_2 x_3 \dots $ is the decimal expansion of $x$ then 
$f(x) = (0.x_1 x_3 x_5 \dots,0.x_2 x_4 x_6 \dots)$ is the bijection. 
Are there cardinalities larger than $\aleph_1$? 
Cantor answered also this question. He showed that for an infinite set, the set of all subsets has a 
larger cardinality than the set itself. 
How does one see this? Assume there is a bijection $x \to A(x)$ which maps each 
point to a set $A(x)$. Now look at the set
 $B = \{ x \; | \; x \notin A(x) \; \}$ and let $b$ be the point in $X$ which corresponds
to $B$. If $y \in B$, then $y \notin B(x)$. 
On the other hand, if $y \notin B$, then $y \in B$. The set $B$ does appear in the
"enumeration" $x \to A(x)$ of all sets. The set of all subsets of $N$ has the same cardinality than
the continuum: $A \to \sum_{j \in A} 1/2^j$ provides a map from $P(N)$ to $[0,1]$. The set of all {\bf finite subsets} of $N$ 
however can be counted. The set of all subsets of the real numbers has cardinality $\aleph_2$, etc. 
Is there a cardinality between $\aleph_0$ and $\aleph_1$? In other words, is there a set which can not be 
counted and which is strictly smaller than the continuum in the sense that one can not find a bijection between it and $R$? 
This was the first of the 23 problems posed by Hilbert in 1900. The answer is surprising: one has a choice.
One can accept either the "yes" or the "no" as a new axiom. In both cases, Mathematics is still fine. The nonexistence of a 
cardinality between $\aleph_0$ and $\aleph_1$ is called the {\bf continuum hypothesis} and is usually abbreviated CH. 
It is independent of the other axioms making up mathematics. 
This was the work of {\bf Kurt G\"odel} in 1940 and {\bf Paul Cohen} in 1963.
The story of exploring the consistency and completeness of axiom systems of all of mathematics is exciting. 
Euclid axiomatized geometry, Hilbert's program was more ambitious. He aimed at 
a set of axiom systems for all of mathematics. The challenge to prove Euclid's 5'th postulate
is paralleled by the quest to prove the CH. But the later is much more fundamental
because it deals with {\bf all of mathematics} and not only with some geometric space.
Here are the {\bf Zermelo-Frenkel Axioms} (ZFC) including the Axiom of choice (C) as established by 
{\bf Ernst Zermelo} in 1908 and {\bf Adolf Fraenkel} and {\bf Thoral Skolem} in 1922.  

\begin{small}
\fcolorbox{yellow2}{yellow2}{
\begin{tabular}{ll}
{\bf Extension}    &  If two sets have the same elements, they are the same. \\
{\bf Image}        &  Given a function and a set, then the image of the function is a set too. \\
{\bf Pairing}      &  For any two sets, there exists a set which contains both sets.  \\
{\bf Property}     &  For any property, there exists a set for which each element has the property. \\
{\bf Union}        &  Given a set of sets, there exists a set which is the union of these sets.  \\
{\bf Power}        &  Given a set, there exists the set of all subsets of this set.  \\
{\bf Infinity}     &  There exists an infinite set. \\
{\bf Regularity}   &  Every nonempty set has an element which has no intersection with the set. \\
{\bf Choice}       &  Any set of nonempty sets leads to a set which contains an element from each.  
\end{tabular}
}
\end{small}

There are other systems like ETCS, which is the {\bf elementary theory of the category of sets}. In category
theory, not the sets but the categories are the building blocks. Categories do not form a set in general. 
It elegantly avoids the Russel paradox too. 
The {\bf axiom of choice (C)} has a nonconstructive
nature which can lead to seemingly paradoxical results like the {\bf Banach Tarski paradox}:
one can cut the unit ball into 5 pieces, rotate and translate the pieces to assemble two identical balls 
of the same size than the original ball.  G\"odel and Cohen showed that
the axiom of choice is logically independent of the other axioms ZF.
Other axioms in ZF have been shown to be independent, like the {\bf axiom of infinity}. A {\bf finitist} would
refute this axiom and work without it. It is surprising what one can do with finite sets. 
The {\bf axiom of regularity} excludes Russellian sets like the set $X$ of all sets which do not contain themselves.
The {\bf Russell paradox} is: Does $X$ contain $X$? It is popularized as the {\bf Barber riddle}: 
a barber in a town only shaves the people who do not shave themselves. Does the barber shave himself? 
{\bf G\"odels theorems} of 1931 deal
with {\bf mathematical theories} which are strong enough to do basic arithmetic in them.

\parbox{16.8cm}{
\fcolorbox{yellow1}{yellow1}{ \parbox{8cm}{
{\bf First incompleteness theorem:} \\
In any theory there are true statements which can not be proved within the theory.
}}
\hspace{4mm}
\fcolorbox{yellow1}{yellow1}{ \parbox{8cm}{
{\bf Second incompleteness theorem:} \\
In any theory, the consistency of the theory can not be proven within the theory.  
}}
}

The proof uses an encoding of mathematical sentences which allows to state liar paradoxical statement 
"this sentence can not be proved". While the later is an odd recreational entertainment gag, it is the 
core for a theorem which makes striking statements about mathematics. These theorems are not limitations of mathematics;
they illustrate its infiniteness. How awful if one could build axiom system and enumerate mechanically all 
possible truths from it.

 \pagebreak
\pagebreak {\bf E-320: Teaching Math with a Historical Perspective \hfill Oliver Knill, 2010-2018}

\chapter{Lecture 8: Probability theory}

{\bf Probability theory} is the science of chance.  It starts with {\bf combinatorics} and leads to a theory of 
{\bf stochastic processes}. Historically, probability theory initiated from gambling problems as in
{\bf Girolamo Cardano's} gamblers manual in the 16th century. A great moment of mathematics occurred, when 
{\bf Blaise Pascal} and {\bf Pierre Fermat} jointly laid a foundation of mathematical probability theory.  \\
It took a while to formalize ``randomness" precisely. Here is the setup as which it had been put forward
by {\bf Andrey Kolmogorov}: all possible experiments of a situation are modeled by a set $\Omega$, the "laboratory".
A measurable subset of experiments is called an ``event". Measurements are 
done by real-valued functions $X$. These functions are called {\bf random variables} and are used to 
{\bf observe the laboratory}. \\
As an example, let us model the process of throwing a coin 5 times. An experiment is a 
word like $httht$, where $h$ stands for ``head" and $t$ represents ``tail".  The laboratory consists of all
such 32 words. We could look for example at the event $A$ that the first two coin tosses are tail. 
It is the set $A=\{ttttt,tttth,tttht,ttthh,tthtt,tthth,tthht,tthhh\}$.
We could look at the random variable which assigns to a word the number of heads. For every experiment, 
we get a value, like for example, $X[tthht]=2$.  \\
In order to make statements about randomness, the concept of a {\bf probability measure} is needed. 
This is a function $P$ from the set of all events to the interval $[0,1]$. It should have the property that 
$P[\Omega]=1$ and $P[A_1 \cup A_2  \cup \cdots ] = P[A_1] + P[A_2]+ \cdots$, if $A_i$ is a sequence of
disjoint events.  \\
The most natural probability measure on a finite set $\Omega$ is $P[A] = \|A\|/\|\Omega\|$, where $\|A\|$ 
stands for the number of elements in $A$. It is the ``number of good cases" divided by the ``number of all cases". 
For example, to count the probability of the event $A$ that we throw $3$ heads during the 5 coin tosses,
we have $|A|=10$ possibilities. Since the entire laboratory has $|\Omega|=32$ possibilities, the probability 
of the event is $10/32$. In order to study these probabilities, one needs {\bf combinatorics}: 
\index{combinatorics}

\begin{center}
\fcolorbox{yellow1}{yellow1}{
\begin{tabular}{ll}
{\bf How many ways are there to:}       &  {\bf The answer is:}  \\ \hline
rearrange or permute $n$ elements       &  $n! = n (n-1) ... 2 \cdot 1$ \\
choose $k$ from $n$ with repetitions    &  $n^k$  \\
pick $k$ from $n$ if order matters      &  $\frac{n!}{(n-k)!}$ \\
pick $k$ from $n$ with order irrelevant &  $\left( \begin{array}{c}  n \\  k \end{array} \right) =  \frac{n!}{k! (n-k)!}$  \\ \hline
\end{tabular}
} \end{center}

The {\bf expectation} of a random variable $E[X]$ is defined as the sum 
$m=\sum_{\omega \in \Omega} X(\omega ) P[\{ \omega \}]$. In our coin toss experiment, this is
$5/2$. The {\bf variance} of $X$ is the expectation of $(X-m)^2$. In our coin experiments, it is $5/4$. 
The square root of the variance is the {\bf standard deviation}. This is the expected deviation 
from the mean. An event happens {\bf almost surely} if the event has probability $1$. \\
An important case of a random variable is $X(\omega) = \omega$ on $\Omega=R$ equipped with probability
$P[A] = \int_A \frac{1}{\sqrt{\pi}} e^{-x^2} \; dx$, the {\bf standard normal
distribution}. Analyzed first by {\bf Abraham de Moivre} in 1733, it was studied by {\bf Carl Friedrich Gauss} in 1807 
and therefore also called {\bf Gaussian distribution}. \\  
Two random variables $X,Y$ are called {\bf uncorrelated}, if $E[X Y] = E[X] \cdot E[Y]$. If for any 
functions $f,g$ also $f(X)$ and $g(Y)$ are uncorrelated, then $X,Y$ are called {\bf independent}. Two random 
variables are said to have the same distribution, if for any $a<b$, the events $\{ a \leq X \leq b \; \}$ and 
$\{ a \leq Y \leq b \; \}$ are independent. If $X,Y$ are uncorrelated, then the relation ${\rm Var}[X] + {\rm Var}[Y] = {\rm Var}[X+Y]$
holds which is just {\bf Pythagoras theorem}, because uncorrelated can be understood 
geometrically: $X-E[X]$ and $Y-E[Y]$ are orthogonal.
A common problem is to study the sum of independent random variables $X_n$ with identical
distribution. One abbreviates this IID. Here are the three most important theorems which we formulate in the case, where
all random variables are assumed to have expectatation $0$ and standard deviation $1$.
Let $S_n = X_1 + ... + X_n$ be the $n$'th sum of the IID 
random variables. It is also called a {\bf random walk}.  \\
\index{expectation}
\index{variance}
\index{standard deviation}

\fcolorbox{yellow1}{yellow1}{ \parbox{16.8cm}{ LLN {\bf Law of Large Numbers} assures that  $S_n/n$ converges to $0$. } } 

\fcolorbox{yellow1}{yellow1}{ \parbox{16.8cm}{ CLT {\bf Central Limit Theorem}:$S_n/\sqrt{n}$ approaches the Gaussian distribution. }} 

\fcolorbox{yellow1}{yellow1}{ \parbox{16.8cm}{ LIL {\bf Law of Iterated Logarithm:} $S_n/\sqrt{2n \log\log(n)}$ accumulates in $[-1,1]$.}}

\hspace{4mm} 
\index{law of large numbers}
\index{central limit theorem}
\index{law of iterated logarithm}

The LLN shows that one can find out about the expectation by averaging experiments. 
The CLT explains why one sees the standard normal distribution so often.
The LIL finally gives us a precise estimate how fast $S_n$ grows. 
Things become interesting if the random variables are no more independent. Generalizing LLN,CLT,LIL
to such situations is part of ongoing research. \\

Here are two open questions in probability theory: 

\begin{center}
\fcolorbox{blue1}{blue1}{ \parbox{16cm}{
Are numbers like $\pi,e,\sqrt{2}$ {\bf normal}: do all digits appear with the same frequency? \\
What growth rates $\Lambda_n$ can occur in $S_n/\Lambda_n$ having limsup $1$ and liminf $-1$? 
}} \end{center}

For the second question, there are examples for $\Lambda_n=1,\lambda_n=\log(n)$ and of course $\lambda_n=\sqrt{n \log \log(n)}$
from LIL if the random variables are independent. Examples of random variables which are not independent are $X_n = \cos(n \sqrt{2})$. \\

{\bf Statistics} is the science of modeling random events in a probabilistic setup. Given data points, we want
to find a {\bf model} which fits the data best. This allows to {\bf understand the past}, {\bf predict the future} 
or {\bf discover laws of nature}. The most common task is to find the {\bf mean} and the {\bf standard deviation} of some
data. The mean is also called the {\bf average} and given by $m=\frac{1}{n} \sum_{k=1}^n x_k$. The variance is
$\sigma^2= \frac{1}{n} \sum_{k=1}^{n} (x_k-m)^2$ with standard deviation $\sigma$. \\

A sequence of random variables $X_n$ define a so called {\bf stochastic process}. Continuous versions of such 
processes are where $X_t$ is a curve of random random variables. 
An important example is {\bf Brownian motion}, which is a model of a random particles. \\

Besides gambling and analyzing data, also {\bf physics} was an important motivator
to develop probability theory. An example is statistical mechanics, where the laws of nature are studied 
with probabilistic methods. A famous physical law is {\bf Ludwig Boltzmann's} relation 
$S=k \log(W)$ for entropy, a formula which decorates Boltzmann's tombstone. 
The {\bf entropy} of a probability measure $P[\{k\}] = p_k$ on a finite set $\{1,...,n\}$ 
is defined as $S = - \sum_{i=1}^n p_i \log(p_i)$. Today, we would reformulate Boltzmann's law and say
that it is the expectation $S=E[\log(W)]$ of the logarithm of the ``Wahrscheinlichkeit" random variable 
$W(i) = 1/p_i$ on $\Omega = \{ 1,...,n \; \}$. Entropy is important because nature tries to maximize it
\index{entropy}

 \pagebreak
\pagebreak {\bf E-320: Teaching Math with a Historical Perspective \hfill Oliver Knill, 2010-2018}

\chapter{Lecture 9: Topology}

{\bf Topology} studies properties of geometric objects
which do not change under continuous reversible deformations. 
In topology, a coffee cup with a single handle is the same as
a doughnut. One can deform one into the other without punching 
any holes in it or ripping it apart. Similarly, a plate and a 
croissant are the same. But a croissant is not equivalent to a doughnut.
On a doughnut, there are closed curves which can not be pulled together to 
a point. For a topologist the letters $O$ and $P$ are the equivalent but 
different from the letter $B$. The mathematical setup is beautiful: 
a {\bf topological space} is a set $X$ with a set $O$ of subsets of $X$ containing both 
$\emptyset$ and $X$ such that finite intersections and arbitrary 
unions in $O$ are in $O$. Sets in $O$ are called {\bf open sets} and $O$ is called a 
{\bf topology}. The complement of an open set is called 
{\bf closed}. Examples of topologies are the {\bf trivial topology }
$O=\{\emptyset,X \}$, where no open sets besides the empty set 
and $X$ exist or the {\bf discrete topology} $O = \{A \; | \;  A \subset X \}$,
where every subset is open. But these are in general not interesting.
An important example on the plane $X$ is the collection $O$ of sets 
$U$ in the plane $X$ for which every point is the center of a small 
disc still contained in $U$. A special class of topological 
spaces are {\bf metric spaces}, where a set $X$ is equipped with a 
{\bf distance function} $d(x,y)=d(y,x) \geq 0$ which satisfies 
the {\bf triangle inequality} $d(x,y)+d(y,z) \geq d(x,z)$ and for 
which $d(x,y)=0$ if and only if $x=y$.  A set $U$ in a metric 
space is open if to every $x$ in $U$, there is a 
{\bf ball} $B_{r}(x) = \{ y| d(x,y)<r \}$ of positive radius $r$ 
contained in $U$. Metric spaces are topological spaces but not vice versa:
the trivial topology for example is not in general. 
For doing {\bf calculus} on a topological space $X$, 
each point has a neighborhood called {\bf chart} which is topologically 
equivalent to a disc in Euclidean space. Finitely many
neighborhoods covering $X$ form an {\bf atlas} of $X$. If the charts 
are glued together with identification maps on the intersection
one obtains a {\bf manifold}. Two dimensional examples are the {\bf sphere}, the 
{\bf torus}, the projective plane or the {\bf Klein bottle}. Topological spaces 
$X,Y$ are called {\bf homeomorphic} meaning ``topologically equivalent" 
if there is an invertible map from $X$ to $Y$ such that this map
induces an invertible map on the corresponding topologies. 
How can one decide whether two spaces are equivalent in this sense?
The surface of the coffee cup for example is equivalent in this sense to the 
surface of a doughnut but it is not equivalent to the surface of a sphere. 
Many properties of geometric spaces can be understood by discretizing
it like with a graph. A graph is a finite collection of 
vertices $V$ together with a finite set of edges $E$, where each edge
connects two points in $V$. For example, the set $V$ of cities in the US where the edges 
are pairs of cities connected by a street is a graph. 
The {\bf K\"onigsberg bridge problem} was a trigger puzzle 
for the study of graph theory. {\bf Polyhedra} were an other start in graph theory.
It study is loosely related to the analysis of surfaces.
The reason is that one can see polyhedra as discrete versions of surfaces. 
In computer graphics for example, surfaces are rendered as finite graphs, 
using triangularizations. 
\index{Koenigsberg bridge problem}
The {\bf Euler characteristic} of a convex polyhedron is a remarkable topological 
invariant. It is $V-E+F   = 2$, 
where $V$ is the number of vertices, $E$ the number of edges and $F$ the number of
{\bf faces}. This number is equal to $2$ for connected polyhedra in which every 
closed loop can be pulled together to a point. This formula for the Euler 
characteristic is also called {\bf Euler's gem}. It comes with a rich history. 
{\bf Ren\'e Descartes} stumbled upon it and written it down in a 
secret notebook. It was Leonard Euler in 1752 was the first to proved the
formula for convex polyhedra.
A convex polyhedron is called a {\bf Platonic solid}, if all vertices are 
on the unit sphere, all edges have the same length and all faces are 
congruent polygons. A theorem of {\bf Theaetetus} states that there are only 
five Platonic solids: 
[Proof: Assume the faces are regular $n$-gons and $m$ of them meet at each vertex.
Beside the Euler relation $V+E+F=2$, a polyhedron also satisfies the relations 
$n F=2 E$ and $m V = 2 E$ which come from counting vertices
or edges in different ways. This gives $2E/m - E + 2E/n=2$ or $1/n + 1/m = 1/E + 1/2$.
From $n \geq 3$ and $m \geq 3$ we see that it is impossible that both $m$ and $n$ are larger than $3$.
There are now nly two possibilities: either $n=3$ or $m=3$. In the case $n=3$ we have $m=3,4,5$
in the case $m=3$ we have $n=3,4,5$. The five possibilities
$(3,3),(3,4),(3,5),(4,3),(5,3)$ represent the five Platonic solids.]
The pairs $(n,m)$ are called the {\bf Schl\"afly symbol} of the polyhedron:
\begin{small}
\parbox{16.8cm}{
\parbox{8.2cm}{
\begin{center} \fcolorbox{yellow2}{yellow2}{
\begin{tabular}{llllll} \hline
Name         &  V  &   E &  F  &  V-E+F  & Schl\"afli \\ \hline
tetrahedron  &  4  &   6 &  4  &     2   & $\{3,3 \}$ \\
hexahedron   &  8  &  12 &  6  &     2   & $\{4,3 \}$ \\
octahedron   &  6  &  12 &  8  &     2   & $\{3,4 \}$ \\ \hline
\end{tabular} } \end{center}
} \parbox{8.2cm}{
\begin{center} \fcolorbox{yellow2}{yellow2}{
\begin{tabular}{llllll} \hline
Name         &  V  &   E &  F  &  V-E+F  & Schl\"afli \\ \hline
\hspace{3mm} &     &     &     &         &            \\
dodecahedron & 20  &  30 & 12  &     2   & $\{5,3 \}$ \\
icosahedron  & 12  &  30 & 20  &     2   & $\{3,5 \}$ \\ \hline
\end{tabular} } \end{center} 
}}
\end{small}

The Greeks proceeded geometrically: Euclid showed
in the "Elements" that each vertex can have either 3,4 or 5 equilateral triangles attached,
3 squares or 3 regular pentagons. (6 triangles, 4 squares or 4 pentagons would lead to 
a total angle which is too large because each corner must have at least 3 different edges).
{\bf Simon Antoine-Jean L'Huilier} refined in 1813 Euler's formula to situations with holes:
\fcolorbox{yellow1}{yellow1}{ \parbox{5cm}{
$V-E+F = 2-2g \; ,$}} 
where $g$ is the number of holes. For a doughnut it is $V-E+F=0$.
Cauchy first proved that there are 4 non-convex regular {\bf Kepler-Poinsot} polyhedra.
\begin{small}
\begin{center} \fcolorbox{yellow2}{yellow2}{
\begin{tabular}{llllll} \hline
Name                          &  V  &   E &  F  &  V-E+F  & Schl\"afli\\ \hline
small stellated dodecahedron  & 12  &  30 & 12  &    -6   & $\{5/2,5 \}$ \\
great dodecahedron            & 12  &  30 & 12  &    -6   & $\{5,5/2 \}$ \\
great stellated dodecahedron  & 20  &  30 & 12  &     2   & $\{5/2,3 \}$ \\
great icosahedron             & 12  &  30 & 20  &     2   & $\{3,5/2 \}$ \\ \hline
\end{tabular} } \end{center}
\end{small}
If two different face types are allowed but each vertex still look the same, one obtains
13 {\bf semi-regular polyhedra.} They were first studied by {\bf Archimedes} in 287 BC.
Since his work is lost, {\bf Johannes Kepler} is considered the first
since antiquity to describe all of them them in his "Harmonices Mundi". 
The {\bf Euler characteristic} for surfaces is $\chi=2-2g$ where $g$ is the number of holes. 
The computation can be done by triangulating the surface. 
The Euler characteristic characterizes smooth compact surfaces if they are orientable. A non-orientable 
surface, the {\bf Klein bottle} can be obtained by gluing ends of the  M\"obius strip.
Classifying higher dimensional manifolds is more difficult and finding good invariants 
is part of modern research. Higher analogues of polyhedra are called {\bf polytopes} 
(Alicia Boole Stott). {\bf Regular polytopes} are the analogue 
of the Platonic solids in higher dimensions. Examples:
\begin{small}
\begin{center} \fcolorbox{yellow2}{yellow2}{
\begin{tabular}{lll}  \hline
dimension  & name                 & Schl\"afli symbols            \\    \hline
2:         & Regular polygons     & $\{3 \},\{4 \},\{5\}, ...$ \\
3:         & Platonic solids      & $\{3,3\},\{3,4\},\{3,5\},\{4,3\},\{5,3\}$  \\
4:         & Regular 4D polytopes & $\{3,3,3\},\{4,3,3\},\{3,3,4\},\{3,4,3\},\{5,3,3\},\{3,3,5\}$ \\
$\geq 5$:  & Regular polytopes    & $\{3,3,3,\dots,3\},\{4,3,3, \dots ,3\},\{3,3,3,\dots,3,4\}$ \\ \hline
\end{tabular} } \end{center}
\end{small}
{\bf Ludwig Schl\"lafly} saw in 1852 exactly six convex regular convex 4-polytopes or
{\bf polychora}, where "Choros" is Greek for "space". Schlaefli's polyhedral formula is
\fcolorbox{yellow1}{yellow1}{ \parbox{5cm}{
$V-E+F-C = 0$ }} 
holds, where $C$ is the number of 3-dimensional {\bf chambers}. 
In dimensions 5 and higher, there are only 3 types of polytopes: the higher dimensional analogues of 
the tetrahedron, octahedron and the cube. A general formula 
\fcolorbox{yellow1}{yellow1}{ \parbox{5cm}{ $\sum_{k=0}^{d-1} (-1)^k v_k = 1-(-1)^{d}$ }} 
gives the Euler characteristic of a convex polytop in $d$ dimensions with $k$-dimensional parts $v_k$. 
\index{chamber}
\index{polychora}
\index{polytopes}

 \pagebreak
\pagebreak {\bf E-320: Teaching Math with a Historical Perspective \hfill Oliver Knill, 2010-2018}

\chapter{Lecture 10: Analysis}

{\bf Analysis} is a science of measure and optimization. As a rather diverse 
collection of mathematical fields, it contains {\bf real and complex analysis}, {\bf functional analysis}, 
{\bf harmonic analysis} and {\bf calculus of variations}.
Analysis has relations to calculus, geometry, topology, probability theory and dynamical systems.
We focus here mostly on "the geometry of fractals" which can be seen as part of dimension theory.
Examples are Julia sets which belong to the subfield of "complex analysis" of "dynamical systems".
"Calculus of variations" is illustrated by the Kakeya needle set in "geometric measure theory",
"Fourier analysis" appears when looking at functions which have fractal graphs, 
"spectral theory" as part of functional analysis is represented by the "Hofstadter butterfly".
We somehow describe the topic using "pop icons". \\

A {\bf fractal} is a set with non-integer dimension. An example is
the {\bf Cantor set}, as discovered in 1875 by Henry Smith. 
Start with the unit interval. Cut the middle third, then cut the middle third 
from both parts then the middle parts of the four parts etc. The limiting set is the Cantor set.
The mathematical theory of fractals belongs to {\bf measure theory} and can also 
be thought of a playground for real analysis or topology. The term {\bf fractal} had been introduced 
by Benoit Mandelbrot in 1975. Dimension can be defined in different ways. The simplest is the
{\bf box counting definition} which works for most household fractals: if  we need $n$ squares of 
length $r$ to cover a set, then 
\fcolorbox{yellow1}{yellow1}{ \parbox{5cm}{
$d=-\log(n)/\log(r)$}} converges to the dimension of the set with 
$r \to 0$. A curve of length $L$ for example needs
$L/r$ squares of length $r$ so that its dimension is $1$. 
A region of area $A$ needs $A/r^2$ squares of length $r$ to be covered 
and its dimension is $2$. The Cantor set needs to be covered with 
$n=2^m$ squares of length $r=1/3^m$. Its dimension is 
$-\log(n)/\log(r) = - m \log(2)/(m \log(1/3)) = \log(2)/\log(3)$. 
Examples of fractals are the graph of the Weierstrass function 1872, 
the Koch snowflak (1904), the Sierpinski carpet (1915) or
the Menger sponge  (1926). \\
{\bf Complex analysis} extends calculus to the complex. It deals with functions $f(z)$ defined
in the complex plane. Integration is done along paths. Complex analysis completes the
understanding about functions. It also provides more examples of fractals by iterating functions 
like the {\bf quadratic map} $f(z) = z^2+c$:  \\
One has already iterated functions before like the Newton method (1879). The
Julia sets were introduced in 1918, the Mandelbrot set   in 1978 and the
Mandelbar set     in 1989.
Particularly famous are the {\bf Douady rabbit} and the {\bf dragon}, the {\bf dendrite}, the {\bf airplane}. 
{\bf Calculus of variations} is calculus in infinite dimensions. Taking derivatives is called taking "variations". 
Historically, it started with the problem to find the curve of fastest fall leading to the 
{\bf Brachistochrone} curve $\vec{r}(t) = (t-\sin(t),1-\cos(t))$. 
In calculus, we find maxima and minima of functions. In calculus of variations, we extremize on much
larger spaces. Here are examples of problems: 

\parbox{16.8cm}{
\parbox{6cm}{
\fcolorbox{yellow2}{yellow2}{
\begin{tabular}{ll}
Brachistochrone        &  1696 \\ 
Minimal surface        &  1760 \\
Geodesics              &  1830 \\
Isoperimetric problem  &  1838 \\
Kakeya Needle problem  &  1917 \\ 
\end{tabular}
}}}

{\bf Fourier theory} decomposes a function into basic components of various frequencies 
$f(x) = a_1 \sin(x) + a_2 \sin(2x) + a_3 \sin(3x) + \cdots$.  The numbers $a_i$ are called 
the {\bf Fourier coefficients}. Our ear does such a decomposition, when we listen to music. 
By distinguish different frequencies, our ear produces a Fourier analysis. 

\parbox{16.8cm}{
\parbox{6cm}{
\fcolorbox{yellow2}{yellow2}{
\begin{tabular}{ll}
Fourier series      &  1729   \\
Fourier transform (FT)  &  1811   \\
Discrete FT         &  Gauss?  \\
Wavelet transform   &  1930   \\   
\end{tabular}
}}}

The Weierstrass function mentioned above is given as a series
$\sum_n a^n \cos(\pi b^n x)$ with $0<a<1,  ab> 1+3\pi/2$. The dimension of its graph 
is believed to be $2+\log(a)/\log(b)$ but no rigorous computation of the dimension was done yet. 
{\bf Spectral theory} analyzes linear maps $L$. The {\bf spectrum}
are the real numbers $E$ such that $L-E$ is not invertible. A Hollywood celebrity among all linear
maps is the {\bf almost Matthieu operator} $L(x)_n = x_{n+1} + x_{n-1} + (2-2 \cos(c n)) x_n$: 
if we draw the spectrum for for each $c$, we see the 
{\bf Hofstadter butterfly}. For fixed $c$ the map describes the behavior of an electron 
in an almost periodic crystal. 
An other famous system is the {\bf quantum harmonic oscillator},
$L(f) = f ''(x) + f(x)$, the {\bf vibrating drum} $L(f) = f_{xx} + f_{yy}$, where $f$ is the 
amplitude of the drum and $f=0$ on the boundary of the drum.  \\
\index{Hofstadter butterly}
\index{Mathieu operator}

\parbox{16.8cm}{
\parbox{6cm}{
\fcolorbox{yellow2}{yellow2}{
\begin{tabular}{ll}
Hydrogen atom            &   1914      \\  
Hofstadter butterfly     &   1976      \\
Harmonic oscillator      &   1900      \\
Vibrating drum           &   1680      \\  
\end{tabular}}
}
}

All these examples in analysis look unrelated at first. 
Fractal geometry ties many of them together: spectra are often fractals, 
minimal configurations have fractal nature, like in solid state physics or in {\bf diffusion limited aggregation}
or in other critical phenomena like {\bf percolation} phenomena, {\bf cracks} in solids or the formation of 
{\bf lighting bolts}
In Hamiltonian mechanics, minimal energy configurations are often fractals like {\bf Mather theory}. And solutions
to minimizing problems lead to fractals in a natural way like when you have the task to turn around a needle 
on a table by 180 degrees and minimize the area swept out by the needle. The minimal turn leads to a Kakaya set, which is
a fractal. Finally, lets mention some unsolved problems in analysis: 
does the {\bf Riemann zeta function} $f(z) = \sum_{n=1}^\infty 1/n^z$ have all nontrivial roots on the axis $Re(z)=1/2$? 
This question is called the {\bf Riemann hypothesis} and is 
the most important open problem in mathematics.  It is an example of a question in {\bf analytic number theory} which
also illustrates how analysis has entered into number theory. Some mathematicians think that spectral theory might solve it.
Also the Mandelbrot set $M$ is not understood yet: the "holy grail" in the field of complex dynamics 
is the problem whether it $M$ is locally connected. From the Hofstadter butterfly one knows that it has measure zero. What is its dimension? 
An other open question in spectral theory is the "can one hear the sound of a drum" problem which asks 
whether there are two convex drums which are not congruent but which 
have the same spectrum. In the area of calculus of variations, just one problem: how long is the shortest curve in space
such that its convex hull (the union of all possible connections between two points on the curve)  contains the unit ball.

 \pagebreak
\pagebreak {\bf E-320: Teaching Math with a Historical Perspective \hfill Oliver Knill, 2010-2018}

\chapter{Lecture 11: Cryptography}

{\bf Cryptography} is the theory of {\bf codes}. Two important aspects of the field are 
the {\bf encryption} rsp. {\bf decryption} of information and {\bf error correction}. Both are crucial 
in daily life. When getting access to a computer, viewing a bank statement or when taking money from the 
ATM, encryption algorithms are used. When phoning, surfing the web, accessing data 
on a computer or listening to music, error correction algorithms are used. 
Since our lives have become more and more digital: music, movies, books, journals, 
finance, transportation, medicine, and communication have become digital, we rely on strong 
error correction to avoid errors and encryption to assure things can not be tempered
with. Without error correction, airplanes would crash: small errors in 
the memory of a computer would produce glitches in the navigation and control program. 
In a computer memory every hour a couple of bits are altered, for example by cosmic rays. Error 
correction assures that this gets fixed. Without error correction music would sound like a 1920 gramophone record.
Without encryption, everybody could intrude electronic banks and transfer money. 
Medical history shared with your doctor would all be public. 
Before the digital age, error correction was assured by extremely redundant information storage. 
Writing a letter on a piece of paper displaces billions of billions of molecules in ink.
Now, changing any single bit could give a letter a different meaning. 
Before the digital age, information was kept in well guarded safes which were 
physically difficult to penetrate. Now, information is locked up in computers which are 
connected to other computers. Vaults, money or voting ballots are secured by mathematical algorithms 
which assure that information can only be accessed by authorized users. 
Also life needs error correction: information in the genome is stored in a 
{\bf genetic code}, where a error correction makes sure that life can survive. 
A cosmic ray hitting the skin changes the DNA of a cell, but in general this is harmless.
Only a larger amount of radiation can render cells cancerous. \\
How can an encryption algorithm be safe? One possibility is to invent a new
method and keep it secret. An other is to use a well known encryption method 
and rely on the {\bf difficulty of mathematical computation tasks} to assure that the method is safe.
History has shown that the first method is unreliable. Systems which rely on "security through 
obfuscation" usually do not last. The reason is that it is tough to 
keep a method secret if the encryption tool is distributed.  
Reverse engineering of the method is often possible, for example using
plain text attacks. Given a map $T$, a third party can compute pairs $x,T(x)$ and by 
choosing specific texts figure out what happens.  \\
The {\bf Caesar cypher} permutes the letters of the alphabet. We can for example replace every 
letter $A$ with $B$, every letter $B$ with $C$ and so on until finally $Z$ is replaced with $A$. 
The word "Mathematics" becomes so encrypted as "Nbuifnbujdt". Caesar would shift the 
letters by $3$. The right shift just discussed was used by his Nephew Augustus. {\bf Rot13} shifts by 13, 
and {\bf Atbash cypher} reflects the alphabet, switch $A$ with $Z$, $B$ with $Y$ etc. The last two examples are
involutive: encryption is decryption. 
More general cyphers are obtained by permuting the alphabet. Because of $26!=403291461126605635584000000 \sim 10^{27}$
permutations, it appears first that a brute force attack is not possible. But Cesar cyphers can be cracked 
very quickly using  statistical analysis. If we know the frequency with which letters appear and match the frequency 
of a text we can figure out which letter was replaced with which.
The {\bf Trithemius cypher} prevents this simple analysis by changing the permutation in each step.
It is called a polyalphabetic substitution cypher. 
Instead of a simple permutation, there are many permutations. After transcoding a letter, we also change the key. 
Lets take a simple example. Rotate for the first letter the
alphabet by $1$, for the second letter, the alphabet by $2$, for the third letter, the alphabet by $3$
etc.  The word "Mathematics" becomes now "Ncwljshbrmd". Note that the second "a" has been translated 
to something different than $a$. A frequency analysis is now more difficult.
The {\bf Viginaire cypher} adds even more complexity: instead of shifting
the alphabet by 1, we can take a key like "BCNZ", then shift the first letter by 1, the second letter by 3
the third letter by 13, the fourth letter by 25 the shift the 5th letter by 1 again. While this cypher
remained unbroken for long, a more sophisticated frequency analysis which involves first finding
the length of the key makes the cypher breakable. With the emergence of computers, even more sophisticated
versions like the German {\bf enigma} had no chance. \\
{\bf Diffie-Hellman key exchange} allows Ana and Bob want to agree on a secret key
over a public channel. The two palindromic friends agree on a prime number $p$ and a base $a$. 
This information can be exchanged over an open channel.
Ana chooses now a secret number $x$ and sends $X=a^x$ modulo $p$ to Bob
over the channel. Bob chooses a secret number $y$ and sends $Y=a^y$ modulo $p$ 
to Ana. Ana can compute $Y^x$ and Bob can compute $X^y$ but both are equal to $a^{x y}$. 
This number is their common secret. The key point is that eves dropper Eve, can not compute this number. The
only information available to Eve are $X$ and $Y$, as well as the base $a$ and $p$. Eve knows that
$X=a^x$ but can not determine $x$. The key difficulty in this code is the  {\bf discrete log problem}: 
getting $x$ from $a^x$ modulo $p$ is believed to be difficult for large $p$. \\
The {\bf Rivest-Shamir-Adleman public key system} uses
a {\bf RSA public key} $(n,a)$ with an integer $n=pq$ and $a<(p-1) (q-1)$, where $p,q$ are prime. 
Also here, $n$ and $a$ are public. Only the factorization of $n$ is 
kept secret. Ana publishes this pair. Bob who wants to email Ana a message $x$, sends her 
$y = x^a \; {\rm mod} \; n$. Ana, who has computed $b$ with 
$ab = 1 \; {\rm mod} \; (p-1)(q-1)$ can read the secrete email $y$ because
$y^b = x^{a b} = x^{(p-1) (q-1)} = x \; {\rm mod} n$.
\index{RSA}
But Eve, has no chance because the only thing Eve knows is
$y$ and $(n,a)$. It is believed that without the {\bf factorization} of $n$, it is not 
possible to determine $x$. The message has been transmitted securely.
The core difficulty is that {\bf taking roots} in the ring 
$Z_n = \{ 0, \dots ,n-1 \; \}$ is difficult without knowing the factorization of $n$. 
With a factorization, we can quickly take arbitrary roots. 
If we can take square roots, then we can also factor: assume we have a 
product $n=pq$ and we know how to take square roots of $1$. 
If $x$ solves $x^2=1 \; {\rm mod} \; n$ and $x$ is different from $1$, then 
$x^2-1= (x-1) (x+1)$ is zero modulo $n$. This means that $p$ divides $(x-1)$ or $(x+1)$. 
To find a factor, we can take the greatest common divisor of $n,x-1$. 
Take $n=77$ for example. We are given the root $34$ of $1$. ( $34^2 = 1156$ has reminder $1$ when 
divided by $34$). The greatest common divisor of $34-1$ and $77$ is $11$ is a factor of $77$. 
Similarly, the greatest common divisor of $34+1$ and $77$ is $7$ divides $77$. Finding roots
modulo a composite number and factoring the number is equally difficult.  

\begin{tabular}{llll}
Cipher          & Used for                 &  Difficulty            & Attack         \\ \hline
Cesar           & transmitting messages    &  many permutations     & Statistics     \\ 
Viginere        & transmitting messages    &  many permutations     & Statistics     \\ 
Enigma          & transmitting messages    &  no frequency analysis & Plain text     \\
Diffie-Helleman & agreeing on secret key   &  discrete log mod p    & Unsafe primes  \\
RSA             & electronic commerce      &  factoring integers    & Factoring      \\
\end{tabular} 

The simplest {\bf error correcting code} uses 3 copies of the same 
information so single error can be corrected. With 3 watches for example, one watch can fail.
But this basic error correcting code is not efficient. 
It can correct single errors by tripling the size. Its efficiency is $33$ percent.
\index{error correction}

 \pagebreak
\pagebreak {\bf E-320: Teaching Math with a Historical Perspective \hfill Oliver Knill, 2010-2018}

\chapter{Lecture 12: Dynamical systems}

{\bf Dynamical systems theory} is the science of time evolution. If time is {\bf continuous} the evolution
is defined by a {\bf differential equation} $\dot{x}=f(x)$. If time is {\bf discrete} then we look at the
{\bf iteration of a map} $x \to T(x)$.  \\

The goal of the theory is to {\bf predict the future} of the system when the present
state is known. A {\bf differential equation} is an equation of the form $d/dt x(t) = f(x(t))$, where the unknown quantity 
is a path $x(t)$ in some ``phase space". We know the {\bf velocity} $d/dt x(t) = \dot{x}(t)$ 
at all times and the initial configuration $x(0))$, we can to compute the {\bf trajectory} $x(t)$. 
What happens at a future time? Does $x(t)$ stay in a bounded region or escape to infinity? Which 
areas of the phase space are visited and how often? Can we reach a certain part of the space when starting
at a given point and if yes, when. An example of such a question is to predict, whether an asteroid 
located at a specific location will hit the earth or not. An other example is to predict the weather of the 
next week. \\

An examples of a dynamical systems in one dimension is the differential equation
$$  x'(t)  = x(t) (2-x(t)), x(0)=1 $$
It is called the {\bf logistic system} and describes population growth. 
This system has the solution $x(t) = 2 e^t/(1+e^{2t})$ as you can see by computing the left and right hand side.  \\

A {\bf map} is a rule which assigns to a quantity $x(t)$ a new quantity $x(t+1) = T(x(t))$. The state $x(t)$
of the system determines the situation $x(t+1)$ at time $t+1$. An example is 
is the {\bf Ulam map } $T(x) = 4 x(1-x)$ on the interval $[0,1]$. This is an example,
where we have no idea what happens after a few hundred iterates even if we would know the initial position with 
the accuracy of the Planck scale. \\

Dynamical system theory has applications  all fields of mathematics. It can be used to find roots of equations
like for 
$$T(x) = x - f(x)/f'(x)  \; . $$
A system of number theoretical nature is the {\bf Collatz map} 
$$   T(x) = \frac{x}{2} \; {\rm (even \; x)}, 3x+1 \; {\rm else}  \; . $$
A system of geometric nature is the {\bf Pedal map} which assigns to a triangle
{\bf the pedal triangle}.  \\
\index{Pedal triangle map}
About 100 years ago, {\bf Henry Poincar\'e} was able to 
deal with {\bf chaos} of low dimensional systems. While {\bf statistical mechanics} 
had formalized the evolution of large systems with probabilistic methods already, the new
insight was that simple systems like a {\bf three body problem} or a {\bf billiard map} can produce very 
complicated motion. It was Poincar\'e who saw that even for such low dimensional and completely 
deterministic systems, random motion can emerge. While physisists have dealt with chaos 
earlier by assuming it or artificially feeding it into equations like the {\bf Boltzmann equation},
the occurrence of stochastic motion in geodesic flows or billiards or restricted three body problems 
was a surprise. These findings needed half a century to sink in and only with the emergence of computers 
in the 1960ies, the awakening happened. Icons like Lorentz helped to popularize the findings 
and we owe them the {\bf "butterfly effect"} picture: a wing of a butterfly can produce a tornado in Texas in a
few weeks. The reason for this statement is that the complicated equations to simulate the weather reduce
under extreme simplifications and truncations to a simple differential equation $\dot{x} =\sigma (y-x),
\dot{y}=rx - y -xz, \dot{z} = xy - bz$, the {\bf Lorenz system}. For $\sigma=10,r=28,b=8/3$, Ed Lorenz discovered 
in 1963 an interesting long time behavior and an aperiodic "attractor". Ruelle-Takens called it a 
{\bf strange attractor}. It is a {\bf great moment} in mathematics to realize that attractors of simple 
systems can become fractals on which the motion is chaotic. It suggests that such behavior
is abundant. What is chaos? If a dynamical system shows {\bf sensitive dependence 
on initial conditions}, we talk about {\bf chaos}. We will experiment with the 
two maps $T(x) = 4x (1-x)$ and $S(x) = 4x - 4x^2$ which starting with the same initial conditions will 
produce different outcomes after a couple of iterations. 

The sensitive dependence on initial conditions is measured by how fast the derivative $dT^n$
of the $n$'th iterate grows. The exponential growth rate $\gamma$ is called the {\bf Lyapunov exponent}. 
A small error of the size $h$ will be amplified to $h e^{\gamma n}$ after $n$ iterates. In 
the case of the Logistic map with $c=4$, the Lyapunov exponent is $\log(2)$ and an error of 
$10^{-16}$ is amplified to $2^n \cdot 10^{-16}$. For time $n=53$ already the error is of the order $1$. This explains the above
experiment with the different maps. The maps $T(x)$ and $S(x)$ round differently on the level $10^{-16}$. 
After 53 iterations, these initial fluctuation errors have grown to a macroscopic size. \\
Here is a famous open problem which has resisted many attempts to solve it: 
Show that the map $T(x,y) = (c \sin(2\pi x) + 2 x -y, x)$ with $T^n(x,y) = (f_n(x,y),g_n(x,y))$
has sensitive dependence on initial conditions on a set of positive area. 
More precisely, verify that for $c>2$ and all $n$
$\frac{1}{n} \int_0^1 \int_0^1 \log|\partial_x f_n(x,y)| \; dx dy \geq \log(\frac{c}{2})$. 
The left hand side converges to the average of the Lyapunov exponents which 
is in this case also the {\bf entropy} of the map. For some systems, one can compute the entropy. 
The logistic map with $c=4$ for example, which is also called the {\bf Ulam map}, has 
entropy $\log(2)$. The {\bf cat map} 
$$  T(x,y) = (2x+y,x+y)  \; {\rm mod 1} \;   $$
has positive entropy 
$\log|(\sqrt{5}+3)/2|$. This is the logarithm of the larger eigenvalue of the 
matrix implementing $T$. \\ 
While questions about simple maps look artificial at first, the mechanisms
prevail in other systems: in astronomy, when studying planetary motion or electrons in the 
van Allen belt, in mechanics when studying coupled pendulum or nonlinear oscillators, in fluid dynamics 
when studying vortex motion or turbulence, in geometry, when studying the evolution of light on a 
surface, the change of weather or tsunamis in the ocean.
Dynamical systems theory started historically with the problem to understand the {\bf motion of planets}. 
Newton realized that this is governed by a differential equation, the {\bf n-body problem}
$$ x_j''(t) = \sum_{i=1}^n \frac{c_{ij} (x_i-x_j)}{|x_i-x_j|^3}  \; , $$
where $c_{ij}$ depends on the masses and the gravitational constant. If one body is the sun and 
no interaction of the planets is assumed and using the common center of gravity as the origin, 
this reduces to the {\bf Kepler problem}
$x''(t) = -C x/|x|^3$, where planets move on {\bf ellipses}, the radius vector sweeps equal area in each time 
and the period squared is proportional to the semi-major axes cubed. A great moment in astronomy was when
Kepler derived these laws empirically. An other great moment in mathematics is Newton's
theoretically derivation from the differential equations.

 \pagebreak
\pagebreak {\bf E-320: Teaching Math with a Historical Perspective \hfill Oliver Knill, 2010-2018}

\chapter{Lecture 13: Computing}

{\bf Computing} deals with algorithms and the art of programming. 
While the subject intersects with
computer science, information technology, the theory is by nature very mathematical. But
there are new aspects: computers have opened the field of {\bf experimental mathematics} and
serve now as the {\bf laboratory} for new mathematics. Computers are not only able to {\bf simulate}
more and more of our physical world, they allow us to {\bf explore} new worlds. \\
A mathematician pioneering new grounds with computer experiments does similar work than an
experimental physicist. Computers have smeared 
the boundaries between physics and mathematics. 
According to Borwein and Bailey, experimental mathematics consists of:

\begin{center}
\parbox{16.8cm}{
\parbox{8.5cm}{
Gain insight and intuition. \\
Find patterns and relations \\
Display mathematical principles \\
Test and falsify conjectures
}
\parbox{8.5cm}{
Explore possible new results \\
Suggest approaches for proofs \\
Automate lengthy hand derivations \\
Confirm already existing proofs
}
}
\end{center}

When using computers to prove things, reading and verifying the computer program is part of the proof. 
If Goldbach's conjecture would be known to be true for all $n>10^{18}$, the conjecture should be accepted because
numerical verifications have been done until $2 \cdot 10^{18}$ until today.  
The first famous theorem proven with the help of a computer was the "4 color theorem" in 1976. 
Here are some pointers in the history of computing: 

\begin{small}
\begin{center}
\fcolorbox{yellow2}{yellow2}{
\parbox{16.8cm}{
\parbox{5.5cm}{
\begin{tabular}{ll}
2700BC  &   Sumerian Abacus \\
200BC   &   Chinese Abacus \\
150BC   &   Astrolabe \\
125BC   &   Antikythera \\
1300    &   Modern Abacus \\
1400    &   Yupana   \\
1600    &   Slide rule \\
1623    &   Schickard computer \\
1642    &   Pascal Calculator  \\
1672    &   Leibniz multiplier \\
1801    &   Punch cards \\
1822    &   Difference Engine \\
1876    &   Mechanical integrator \\
\end{tabular}
}
\parbox{5.5cm}{
\begin{tabular}{ll}
1935    &   Zuse 1 programmable \\
1941    &   Zuse 3     \\
1943    &   Harvard Mark I  \\
1944    &   Colossus \\
1946    &   ENIAC   \\
1947    &   Transistor \\
1948    &   Curta Gear Calculator \\
1952    &   IBM 701 \\
1958    &   Integrated circuit \\
1969    &   Arpanet \\
1971    &   Microchip \\
1972    &   Email \\
1972    &   HP-35 calculator \\
\end{tabular}
}
\parbox{5.5cm}{
\begin{tabular}{ll} 
1973    &   Windowed OS \\
1975    &   Altair 8800 \\
1976    &   Cray I  \\
1977    &   Apple II \\
1981    &   Windows I     \\
1983    &   IBM PC \\
1984    &   Macintosh  \\
1985    &   Atari      \\
1988    &   Next \\
1989    &   HTTP  \\
1993    &   Web browser, PDA \\
1998    &   Google  \\
2007    &   iPhone 
\end{tabular}
}
}
}
\end{center} 
\end{small}

We live in a time where technology explodes exponentially.{\bf Moore's law} from 1965 predicted that semiconductor
technology doubles in capacity and overall performance every 2 years. This has happened since. 
Futurologists like Ray Kurzweil conclude from this  technological singularity in which artificial intelligence 
might take over. An important question is how to decide whether a computation is "easy" or "hard".
In 1937, {\bf Alan Turing} introduced the idea of a {\bf Turing machine}, 
a theoretical model of a computer which allows to quantify complexity.
It has finitely many states $S=\{s_1,...,s_n,h \; \}$ and works on an tape of $0-1$ sequences. The state $h$ is the 
"halt" state. If it is reached, the machine stops. The machine has rules
which tells what it does if it is in state $s$ and reads a letter $a$. Depending on $s$ and $a$, it writes $1$ or $0$
or moves the tape to the left or right and moves into a new state. Turing showed
that anything we know to compute today can be computed with Turing machines. For any known machine, there is a polynomial
$p$ so that a computation done in $k$ steps with that computer can be done in $p(k)$ steps on a Turing machine. 
What can actually be computed? {\rm Church's thesis} of 1934 states that everything which can be computed can 
be computed with Turing machines. Similarly as in mathematics itself, there are limitations of computing. Turing's 
setup allowed him to enumerate all possible Turing machine and use them as input of an other machine. Denote by $TM$ the
set of all pairs $(T,x)$, where $T$ is a Turing machine and $x$ is a finite input. Let $H \subset TM$ denote
the set of Turing machines $(T,x)$ which halt with the tape $x$ as input. Turing looked at the decision problem: is there a machine
which decides whether a given machine $(T,x)$ is in $H$ or not. An ingenious Diagonal argument of Turing shows that the
answer is "no". [Proof: assume there is a machine 
$HALT$ which returns from the input $(T,x)$ the output ${\rm HALT(T,x)=true}$, if $T$ halts with the input $x$ and
otherwise returns ${\rm HALT(T,x)=false}$. Turing constructs a Turing machine ${\rm DIAGONAL}$, which 
does the following: 1) Read x. 2) Define Stop=HALT(x,x) 3) While Stop=True repeat Stop:=True; 4) Stop. \\
Now, DIAGONAL is either in H or not. If DIAGONAL is in H, then the variable Stop is true which means that the machine
DIAGONAL runs for ever and DIAGONAL is not in H. But if DIAGONAL is not in H, then the variable Stop is false which 
means that the loop 3) is never entered and the machine stops. The machine is in H.] \\
Lets go back to the problem of distinguishing "easy" and "hard" problems: 
One calls {\bf P} the class of decision problems that are solvable in polynomial time and
{\bf NP} the class of decision problems which can efficiently be tested if the solution is given.
These categories do not depend on the computing model used. The question 
\fcolorbox{yellow1}{yellow1}{"N=NP?"}
is the most important open problem in theoretical computer science.  
It is one of the seven {\bf millenium problems} and it is widely believed that $P \neq NP$. 
If a problem is such that every other NP problem can be reduced to it, it is called {\bf NP-complete}. 
Popular games like Minesweeper or Tetris are NP-complete. If $P \neq NP$, then there is no efficient 
algorithm to beat the game. The intersection of NP-hard and NP is the class of NP-complete problems.
An example of an NP-complete problem is the {\bf balanced number partitioning problem}:
given $n$ positive integers, divide them into two subsets $A,B$, so that the sum in $A$ and the sum in $B$ 
are as close as possible. A first shot: chose the largest remaining number and distribute it to alternatively 
to the two sets. \\
We all feel that it is harder to {\bf find a solution to a problem} rather than to {\bf verify a solution}. 
If $N \neq NP$ there are one way functions, functions which are easy to compute but hard to verify. 
For some important problems, we do not even know whether they are in NP. Examples are the
{\bf the integer factoring problem}.
An efficient algorithm for the first one would have enormous consequences.
Finally, lets look at some mathematical problems in artificial intelligence AI: \\

\fcolorbox{green1}{green1}{
\parbox{16.8cm}{
\begin{tabular}{ll}
problem solving        &  playing games like chess, performing algorithms, solving puzzles \\
pattern matching       &  speech, music, image, face, handwriting, plagiarism detection, spam \\
reconstruction         &  tomography, city reconstruction, body scanning \\
research               &  computer assisted proofs, discovering theorems, verifying proofs \\
data mining            &  knowledge acquisition, knowledge organization, learning \\
translation            &  language translation, porting applications to programming languages \\
creativity             &  writing poems, jokes, novels, music pieces, painting, sculpture \\
simulation             &  physics engines, evolution of bots, game development, aircraft design \\
inverse problems       &  earth quake location, oil depository, tomography \\
prediction             &  weather prediction, climate change, warming, epidemics, supplies \\
\end{tabular}
}}

 \pagebreak

\fontsize{12}{15} \selectfont

\section*{About this document}

It should have become obvious that I'm reporting on many of these theorems as
a {\bf tourist} and not as a {\bf local}. In some few areas I could qualify 
as a {\bf tour guide} but hardly as a local. The references contain only parts
which have been consulted but it does not imply that I know all of that 
source. My own background was in dynamical systems theory and mathematical physics. 
Both of these subjects by nature have many connections with other branches of mathematics.  \\

The motivation to try such a project came  through teaching
a course called {\bf Math E 320} at the Harvard extension school. This
math-multi-disciplinary course is part of the ``math for teaching program",
and tries to map out the major parts of mathematics and visit some selected 
placed on 12 continents. \\

It is wonderful to visit other places and see connections. One can learn 
new things, relearn old ones and marvel again about how large and diverse 
mathematics is but still to notice how many similarities there are between 
seemingly remote areas. 
A goal of this project is also to get back up to speed up to the level of 
a first year grad student (one forgets a lot of things over the years) and
maybe pass the qualifying exams (with some luck). \\

This summer 2018 project also illustrates the challenges
when trying to tour the most important mountain peaks in the 
mathematical landscape with limited time. Already the identification of
major peaks and attaching a ``height" can be challenging. Which theorems
are the most important? Which are the most fundamental? Which theorems
provide fertile seeds for new theorems? I recently got asked by some 
students what I consider the most important theorem in mathematics 
(my answer had been the ``Atiyah-Singer theorem"). \\
\index{Atiyah-Singer theorem}

Theorems are the entities which build up mathematics. 
Mathematical ideas show their merit only through theorems. 
Theorems not only help to bring ideas to live, they in turn 
allow to solve problems and justify the language or theory. 
But not only the results alone, also the history and the connections with the
mathematicians who created the results are fascinating. \\

The first version of this document got started in May 2018 
and was posted in July 2018. Comments, suggestions or corrections are welcome. 
I hope to be able to extend, update and clarify it and explore 
also still neglected continents in the future if time permits. \\

It should be pretty obvious that one can hardly do justice 
to all mathematical fields and that much more would be needed to 
cover the essentials. A more serious project would be to identify a dozen 
theorems in each of the major MSC classification fields. The current MSC2020 classification
system has now 64 major entries and thousands of sub-entries listed on 120 pages \cite{AMS2020}.
But even ``thousand and one theorem" list would only be the tip of the iceberg.
Such a list exists already: on Wikipedia, there are currently about 1000 theorems 
discussed. The one-document project getting closest to this project is 
maybe the beautiful book \cite{Neunhauserer}.

\pagebreak

\section{Document history}

The first draft was posted on July 22, 2018 \cite{FundamentalTheorems}.
On July 23, 2018, a short list of theorems was made available on 
\cite{Top10FundamentalTheorems}.
This document history section got started on July 25-27, 2018. 

\begin{tiny}
\begin{itemize}
\item July 28 2018: Entry 36 had been a repeated prime number theorem entry. Its alternative is now the Fredholm alternative.
      Also added are the Sturm theorem and Smith normal form.
\item July 29: The two entries about Lidskii theorem and Radon transform are added.
\item July 30: An entry about linear programming.
\item July 31: An entry about random matrices.
\item August 2: An entry about entropy of diffeomorphisms
\item August 4: 104-108 entries: linearization, law of small numbers, Ramsey, Fractals and Poincare duality. 
\item August 5: 109-111 entries: Rokhlin and Lax approximation, Sobolev embedding
\item August 6: 112: Whitney embedding.
\item August 8: 113-114: AI and Stokes entries
\item August 12: 115 and 116: Moment entry  and martingale theorem
\item August 13: 117 and 118: theorema egregium and Shannon theorem
\item August 14: 119 mountain pass
\item August 15: 120, 121,122,123 exponential sums, sphere theorem, word problem and finite simple groups
\item August 16: 124, 125, 126, Rubik, Sard and Elliptic curves, 
\item August 17: 127, 128, 129 billiards, uniformization, Kalman filter
\item August 18: 130,131 Zarisky and Poincare's last theorem
\item August 19: 132, 133 Geometrization, Steinitz
\item August 21: 134, 135 Hilbert-Einstein, Hall marriage
\item August 22: 136-130
\item August 24: 141-142
\item August 25: 143-144
\item August 27: 145-149
\item August 28: 150-151
\item August 31: 152
\item September 1: 153-155
\item September 2: 156
\item September 8: 157,158
\item September 14 2018: 159-161
\item September 25 2018: 162-164
\item March 17 2019: 165-169
\item March 20, 2019: section on paradigms
\item March 21, 2019: 170
\item March 27, 2019, 171
\item June 20, 2019, 172
\item August 6, 2020, 173-174, deepness section started
\item August 8, 2020, 175-177, more on deepness section
\item August 18, 2020, 178,179,
\item August 19, 2020, 180,181,182 
\item August 20, 2020, section on essential math, 183-185
\item August 24, 2020, 186,187
\item August 25, 2020, 188,189,190,191
\item August 26, 2020, 192, 193
\item August 27, 2020, 194 - 200
\item August 28, 2020, 201, 202
\item August 30, 2020, 203
\item August 31, 2020, 204,205
\item September 5, 2020, 206,207
\item September 6-8, 2020, 208-212
\item September 9, 2020, 213-214
\item September 10, 2020, 215,216,217,218
\item September 21, 2020, 219
\item October 2, 2020, 220-221
\item October 8, 2020, 222-223
\item October 12, 2020, 224-225
\item November 4, 2020, 226-227
\item November 5, 2020, 228-231
\item November 6, 2020, 232
\item November 16, 2020, 233-234
\item November 25, 2020, 235-236
\item December 3, 2020, 237-238
\item December 4, 2020, 239
\item January 20, 2021, 240-243
\item May 11, 2021, 244
\item February 2, 2022, 245-250
\end{itemize}
\end{tiny}

\section{Top choice}

The short list of 10 theorems mentioned in the youtube clip were:

\begin{itemize}
\item Fundamental theorem of arithmetic (prime factorization)
\item Fundamental theorem of geometry (Pythagoras theorem)
\item Fundamental theorem of logic (incompleteness theorem)
\item Fundamental theorem of topology (rule of product)
\item Fundamental theorem of computability (Turing computability)
\item Fundamental theorem of calculus (Stokes theorem)
\item Fundamental theorem of combinatorics, (pigeonhole principle)
\item Fundamental theorem of analysis (spectral theorem)
\item Fundamental theorem of algebra (polynomial factorization)
\item Fundamental theorem of probability (central limit theorem)
\end{itemize}

Let me try to justify this shortlist. It should go without saying that similar arguments could 
be stated for any other choice, except maybe for the
five classical fundamental theorems: Arithmetic, Geometry (which is undisputed Pythagoras), 
Calculus and Algebra, where one can hardly argue much: except for the Pythagorean theorm,
their given name already suggests that they are considered fundamental. Here is some reflection:

\begin{itemize}
\item {\bf Analysis}. Why chose the spectral theorem and not say the more general {\bf Jordan normal form theorem}? 
This is not an easy call but the {\bf Jordan normal form theorem} is less simple to state and 
furthermore, that it does not stress the importance of {\bf normality} giving the possibility 
for a {\bf functional calculus}. Also, the spectral theorem holds in infinite dimensions
for operators on Hilbert spaces. If one looks at mathematical physics for example,
then it is the {\bf functional calculus of operators} which is really made use of; 
the Jordan normal form theorem appears rarely in comparison. 
In infinite dimensions, a Jordan normal form theorem would be much more difficult
as the operator $Au(n) = u(n+1)$ on $l^2(\mathbb{Z})$ is both unitary as well as a ``Jordan form matrix". 
The spectral theorem however sails through smoothly to infinite dimensions and even applies with adaptations
to {\bf unbounded self-adjoint operators} which are important in physics. And as it is a core part of 
{\bf analysis}, it is also fine to see the theorem as part of analysis. The main reason of course is 
that the fundamental theorem of algebra is already occupied by a theorem. One could object that ``analysis" is 
already represented by the fundamental theorem of calculus but calculus is so important that it can represent
its own field. The idea of the fundamental theorem of calculus goes beyond calculus. It is essentially
a {\bf cancellation property}, a {\bf telescopic sum} or {\bf Pauli principle} ($d^2=0$ for exterior derivatives)
which makes the principle work. Calculus is the idea of an exterior derivative, the idea of cohomology, 
a link between algebra and geometry. One can see calculus also as a theory of ``time". In some sense, the fundamental 
theorem of calculus also represents the field of {\bf differential equations} and this is what 
``time is all about". 

\item {\bf Probability}. One can ask also why to pick the {\bf central limit theorem} and not say the 
{\bf Bayes formula} or then the deeper {\bf law of iterated logarithm}.
One objection against the Bayes formula is that it is essentially a definition, 
like the basic arithmetic properties ``commutativity, distributivity or associativity" in 
an algebraic structure like a ring. One does not present the identity $a+b=b+a$ for example as a fundamental theorem.
Yes, the Bayes theorem has an unusual high appeal to scientists as it appears like a {\bf magic
bullet}, but for a mathematician, the statement just does not have enough
beef: it is a definition, not a theorem. Not to belittle the Bayes theorem, like the notion
of {\bf entropy} or the notion of {\bf logarithm}, it is a {\bf genius concept}. But it is not an 
actual theorem, as the cleverness of the statement of Bayes lies in the {\bf definition} and so the
clarification of conditional probability theory.
For the central limit theorem, it is pretty clear that it should be high up on any list of theorems,
as the name suggests: it is central. But also, it actually is {\bf stronger} than some versions of the law of large 
numbers. The strong law is also super seeded by Birkhoff's ergodic theorem which is much more general. 
One could argue to pick the {\bf law of iterated logarithm} or some {\bf Martingale theorem} instead
but there is something appealing in the central limit theorem which goes over to other set-ups. One
can formulate the central limit theorem also for random variables taking values in a compact topological
group like when doing statistics with spherical data \cite{FLE}. An other pitch for the central 
limit theorem is that it is a {\bf fixed point of a renormalization map} $X \to \overline{X + X}$ (where the 
right hand side is the sum of two independent copies of $X$) in the space of random variables. 
This map {\bf increases entropy} and the fixed point is is a random variable whose distribution function 
$f$ has the {\bf maximal entropy} $- \int_{\RR} f(x) \log(f(x)) \; dx$ among all probability density functions. 
The entropy principle justifies essentially all known probability density functions. Nature just likes to 
maximize entropy and minimize energy or more generally - in the presence of energy - to minimize the free energy.

\item {\bf Topology}. Topology is about geometric properties which do not change under 
continuous deformation or more generally under homotopies.
Quantities which are invariant under homeomorphisms are interesting. 
Such quantities should add up under disjoint unions of 
geometries and multiply under products. The Euler characteristic is {\bf the} prototype.
Taking products is fundamental for building up Euclidean spaces
(also over other fields, not only the real numbers) which locally patch up more complicated spaces. 
It is the essence of vector spaces that after building a basis, one has a product of Euclidean spaces.
Field extensions can be seen therefore as product spaces. 
How does the counting principle come in? As stated, it actually is quite strong and calling it a
``fundamental principle of topology" can be justified if the product of topological spaces is defined
properly: if $1$ is the one-point space, one can see the statement $G \times 1 = G_1$ 
as the {\bf Barycentric refinement} of $G$, implying that the Euler characteristic is a 
Barycentric invariant and so that it is a ``counting tool" which can be pushed to the continuum,
to manifolds or varieties. And the compatibility with the product is the key to make it work. 
Counting in the form of Euler characteristic goes throughout mathematics, 
combinatorics, differential geometry or algebraic geometry. Riemann-Roch or Atiyah-Singer
and even dynamical versions like the Lefschetz fixed point theorem (which generalizes the Brouwer fixed 
point theorem) or the even more general Atiyah-Bott theorem can be seen as {\bf extending the basic counting principle}:
the {\bf Lefschetz number} $\chi(X,T)$ is a dynamical Euler characteristic which in the static case $T=Id$ 
reduces to the Euler characteristic $\chi(X)$. 
In ``school mathematics", one calls the principle the ``fundamental principle of counting"
or ``rule of product". It is put in the following way: ``If we have $k$ ways to do one thing  and $m$ ways to do an 
other thing, then we have $k*m$ ways to do both". It is so simple that one can argue that it is over represented
in teaching but it is indeed important.   \cite{BiggsRoots} makes the point that it should be considered a 
{\bf founding stone of combinatorics}.

Why is the multiplicative property more fundamental 
than the {\bf additive counting principle}. It is again that the additive property is essentially 
placed in as a definition of what a {\bf valuation} is. It is in the {\bf in-out-formula} 
$\chi(A \cup B) + \chi(A \cap B)=\chi(A)+\chi(B)$. Now, this inclusion-exclusion formula is also important in combinatorics
but it is already in the {\bf definition} of what we call counting or ``adding things up".
The multiplicative property on the other hand is not a definition; it actually
is quite non-trivial. It characterizes classical mathematics as {\bf quantum mechanics} or {\bf non-commutative 
flavors of mathematics} have shown that one can extend things.
So, if the ``rule of product" (which is taught in elementary school) is 
beefed up to be more geometric and interpreted to Euler characteristic, it becomes fundamental. 

\item {\bf Combinatorics}. The pigeonhole principle stresses the importance of 
{\bf order structure}, partially ordered sets (posets) and cardinality or comparisons of cardinality. 
The point for posets is made in \cite{JamesPropp} who writes 
{\it The biggest lesson I learned from Richard Stanley's work is, combinatorial objects want to be partially ordered!}
The use of injective functions to express cardinality is a key part of Cantor. 
Like some of the ideas of Grothendieck it is of ``infantile simplicity" (quote Grothendieck about schemes) 
but powerful.  It allowed for the stunning result that there are different infinities. One of the
reason for the success of Cantor's set theory is the immediate applicability. For any new theory,
one has to ask: ``does it tell me something I did not know?" In ``set theory" the larger cardinality
of the reals (uncountable) than the cardinality of the algebraic numbers (countable) gave immediately
the existence of {\bf transcendental numbers}. This is very elegant. The pigeonhole principle 
similarly gives combinatorial results which are non trivial and elegant. 
Currently, searching for ``the fundamental theorem of combinatorics" gives the {\bf ``rule of product"}.
As explained above, we gave it a geometric spin and placed it into topology. Now, combinatorics and 
topology have always been very hard to distinguish. Euler, who somehow booted up topology by reducing the 
{\bf K\"onigsberg problem} to a problem in graph theory did that already. Combinatorial topology is
essentially part of topology. Today, some very geometric topics like algebraic geometry
have been placed within pure {\bf commutative algebra} (this is how I myself was exposed to algebraic geometry)
On the other hand, some very hard core combinatorial problems like the upper bound conjecture have been 
proven with algebro-geometric methods like toric varieties which are geometric. In any case, 
order structures are important everywhere and the pigeonhole principle justifies the importance of order structures.

\item {\bf Computation}. There is no official ``fundamental theorem of computer science" but
the {\bf Turing completeness theorem} comes up as a top candidate when searching on engines. 
Turing formalized using Turing machines in a precise way, what computing is, and even what a 
proof is. It nails down {\bf mathematical activity} of running an algorithm or argument in a 
mathematical way. It is also pure as it is {\bf not hardware dependent}. One can also only appreciate
Turing's definition if one sees how different programming languages can look like and also in
logic, what type of different frame works have been invented. Turing breaks all this complexity with
a machine which can be itself part of mathematics leading to the {\bf Halte problem} illustrating the
basic limitations of computation. {\bf Quantum computing} would add a
hardware component and might break through the {\bf Turing-Church thesis} that everything we can 
compute can be computed with Turing machines in the same complexity class. 
G\"odel and Turing are related and the Turing incompleteness theorem has a 
similar flavor than the G\"odel incompleteness theorems. There is an other angle to it and that is 
the question of {\bf complexity}. I would predict that most mathematicians
would currently favor the Platonic view of the Church thesis and predict that 
also new paradigms like quantum computing will never go beyond {\bf Turing computability} or 
even not break through complexity barriers like {\bf P-NP thresholds}.
It is just that the Turing completeness theorem is too beautiful to be spoiled by a different
type of complexity tied to a physical world. The point of view is that anything we see in 
the  physical world can in principle be computed with a machine {\bf without changing the
complexity class}. But that picture could be as naive as Hilbert's dream one hundred years ago. 
Still, whatever happens in the future, the Turing completeness theorem remains a theorem. Theorems
stay true. 

\item {\bf Logic}. One can certainly argue whether it would be justified to have G\"odel's theorem replaced
by a theorem in category theory like the Yoneda lemma. The Yoneda result
is not easy to state and it does not produce yet an ``Aha moment" like G\"odel's theorem does 
(the liars paradox explains the core of G\"odel's theorem, and it was successfully popularized in
\cite{Hofstadter}.) Maybe The Yoneda theorem will hit the pop culture in the future, when all 
mathematics has been naturally and pedagogically well expressed in categorical language. 
I'm personally not sure whether this will ever happen: not everything which is nice also
had been penetrating large parts of mathematics: an example is given by {\bf non-standard analysis}, 
which makes calculus orders of magnitudes easier and which is related also to {\bf surreal numbers}, 
which are the most ``natural" numbers.  Both concepts have not entered calculus or 
algebra textbooks and there are reasons: the subjects need mathematical 
maturity and one can easily make mistakes. (I myself use non-standard analysis on an 
intuitive level as presented by Nelson \cite{Nelson77,Robert} and think of a compact 
set as a finite set for example which for example, where basic
theorems almost require no proof like the Bolzano theorem telling that a continuous 
function on a compact set takes
a maximum). But using non-standard analysis would be a ``no-no" both in 
teaching as well when formulating mathematical
thoughts for others who are not familar with the three additional axioms IST within ZFC of Nelson. 
It is non-standard and true to its name. 
An example where something was once pop-culture but then was sidelined are quaternions.
It might be a topic which has a comeback. Fashion is hard to predict. 
Also, much of category theory still feels just like a huge conglomerate of definitions. There
is lots of dough in the form of definitions and little raisins in the form of theorems. 
Historically also the language of set theory have been overkill especially in education, 
where it has lead to ``new math" controversies in the 1960ies. 
The work of Russel and Whitehead demonstrates, 
how clumsy things can become if boiled down to the small pieces. 
We humans like to think and programming in higher order structures,
rather than doing assembly coding, we like to work in object oriented languages which give more insight. 
But we like and make use of that higher order codes can be boiled down to assembly closer to what the 
basic instructions are. 
This is similar in mathematics and also in future, a topologist working in 4 manifold theory will hardly 
think about all the definitions in terms of sets for similar reasons than a modern computer algebra system
does not break down all the objects into lists and lists of lists (even so, that's what it often is).
Category theory has a chance to change the landscape because
it is close to computer science and to natural data structures. It is more pictorial and flexible
than set theory alone. It definitely has been very successful to find new structures and see connections within
different fields like computer science \cite{PierceCategory}. It also has lead to more flexible axiom systems.
\end{itemize}

\printindex

\chapter{Bibliography}

\bibliographystyle{plain}

\end{document}